\documentclass[11pt,twoside]{book} 

\usepackage{amssymb,graphicx,fancyhdr,amsmath,booktabs,pdflscape,longtable,multirow,listings,enumitem,amsthm,float,mathtools}
\usepackage{pgf,tikz}
\usepackage[NoDate]{currvita}
\usepackage[ansinew]{inputenc}
\usepackage[T1]{fontenc}
\usetikzlibrary{arrows,decorations.pathmorphing,backgrounds,positioning,fit}

\tikzset{place/.style={rectangle,draw=black!50,thick,inner sep=4pt,minimum size=6mm,rounded corners=2mm}}
\tikzset{kreis/.style={circle,draw=black!50,thick,inner sep=1pt,minimum size=11mm}}
\tikzset{kreis1/.style={circle,draw=black!50,thick,inner sep=1pt,minimum size=7mm}}
\tikzset{point/.style={circle,draw=black!100,thick,inner sep=0pt,minimum size=2pt}}

\definecolor{cTangoButter1}{RGB}{252,233,79}
\definecolor{cTangoButter2}{RGB}{237,212,0}
\definecolor{cTangoButter3}{RGB}{196,160,0}
\definecolor{cTangoChameleon1}{RGB}{138,226,52}
\definecolor{cTangoChameleon2}{RGB}{115,210,22}
\definecolor{cTangoChameleon3}{RGB}{78,154,6}
\definecolor{cTangoOrange1}{RGB}{252,175,62}
\definecolor{cTangoOrange2}{RGB}{245,121,0}
\definecolor{cTangoOrange3}{RGB}{206,92,0}
\definecolor{cTangoSkyBlue1}{RGB}{114,159,207}
\definecolor{cTangoSkyBlue2}{RGB}{52,101,164}
\definecolor{cTangoSkyBlue3}{RGB}{32,74,135}
\definecolor{cTangoPlum1}{RGB}{173,127,168}
\definecolor{cTangoPlum2}{RGB}{117,80,123}
\definecolor{cTangoPlum3}{RGB}{92,53,102}
\definecolor{cTangoChocolate1}{RGB}{233,185,110}
\definecolor{cTangoChocolate2}{RGB}{193,125,17}
\definecolor{cTangoChocolate3}{RGB}{143,89,2}
\definecolor{cTangoScarletRed1}{RGB}{239,41,41}
\definecolor{cTangoScarletRed2}{RGB}{204,0,0}
\definecolor{cTangoScarletRed3}{RGB}{164,0,0}
\definecolor{cTangoAluminum1}{RGB}{238,238,236}
\definecolor{cTangoAluminum2}{RGB}{211,215,207}
\definecolor{cTangoAluminum3}{RGB}{186,189,182}
\definecolor{cTangoAluminum4}{RGB}{136,138,133}
\definecolor{cTangoAluminum5}{RGB}{85,87,83}
\definecolor{cTangoAluminum6}{RGB}{46,52,54}

\renewcommand{\baselinestretch}{1.15}   % Zeilenabstand (Größenangabe nicht in cm!)
\renewcommand{\chaptermark}[1]{\markboth{ #1}{}}    % \markboth{\thechapter{} #1}
\renewcommand{\sectionmark}[1]{\markright{\thesection\ #1}}
\begin{document}

%--------------------------   Satzzählung, Stil der Theoreme und Neudefinitionen  ---------------------------------

\theoremstyle{plain}
\newtheorem{thm}{Theorem}[chapter]
\newtheorem{lemma}[thm]{Lemma}
\newtheorem{kor}[thm]{Corollary}
\newtheorem{prop}[thm]{Proposition}
\newtheorem{conj}[thm]{Conjecture}

\theoremstyle{definition}
\newtheorem{defi}[thm]{Definition}
\newtheorem{ex}[thm]{Example}
\newtheorem{alg}[thm]{Algorithm}

\theoremstyle{remark}
\newtheorem{rem}[thm]{Remark}
\newtheorem{notation}[thm]{Notation}
\newtheorem*{ack}{Acknowledgement}

\newcommand{\N}{\mathbb{N}}
\newcommand{\Z}{\mathbb{Z}}
\newcommand{\Q}{\mathbb{Q}}
\newcommand{\R}{\mathbb{R}}
\newcommand{\C}{\mathbb{C}}
\newcommand{\mc}{\mathcal}
\newcommand{\mf}{\mathfrak}
\newcommand{\mb}{\mathbb}
\newcommand{\x}{\underline{x}}
\newcommand{\GKZ}{GKZ }
\newcommand{\uu}{\underline{v}}

%--------------------------   Beginn Einleitungsteil  --------------------------------------------------------------

%--------------------------   Titelseite und Zusammenfassung in deutsch und englisch  ------------------------------

\pagestyle{empty}
\begin{titlepage}
\center
 \vspace*{6em}
 \Huge{Picard-Fuchs equations \\of special one-parameter families \\of invertible polynomials}
 \vspace{14em}
 \\
 \large{Der Fakult\"at für Mathematik und Physik\\
 der Gottfried Wilhelm Leibniz Universität Hannover\\
 zur Erlangung des Grades\\
 Doktorin der Naturwissenschaften\\
 Dr. rer. nat.\\
 vorgelegte Dissertation\\
 von\\
 Dipl.-Math. Swantje G\"ahrs\\
 geboren am 08. Mai 1983 in Hamburg
 }
% \newpage
% \vspace*{50em}
%\begin{flushleft}
% Referent: Prof. Dr. Wolfgang Ebeling, Hannover\\
% Korreferentin: Prof. Dr. Noriko Yui, Kingston\\
% Tag der Promotion: 21. Oktober 2011
%\end{flushleft}
\end{titlepage}

\newpage

\chapter*{Abstract}
The thesis deals with calculating the Picard-Fuchs equation of special one-parameter families of invertible polynomials. In particular, for an invertible polynomial $g(x_1,\dots,x_n)$ we consider the family $f(x_1,\dots,x_n)=g(x_1,\dots,x_n)+s\cdot\prod x_i$, where $s$ denotes the parameter. For the families of hypersurfaces defined by these polynomials, we compute the Picard-Fuchs equation, i.e. the ordinary differential equation which solutions are exactly the period integrals. For the proof of the exact appearance of the Picard-Fuchs equation we use a combinatorial version of the Griffiths-Dwork method and the theory of \GKZ systems. As consequences of our work and facts from the literature, we show the relation between the Picard-Fuchs equation, the Poincar\'{e} series and the monodromy in the space of period integrals.

\textbf{Keywords:} Picard-Fuchs equation, invertible polynomials, Griffiths-Dwork method, \GKZ systems, Poincar\'{e} series, Monodromy
\newpage
\thispagestyle{empty}

%--------------------------   Inhaltsverzeichnis  ------------------------------------------------------------------

\pagestyle{plain}
{\renewcommand{\baselinestretch}{0.7}
\tableofcontents
}
\newpage
\thispagestyle{empty}
\cleardoublepage

%\frontmatter  
%\mainmatter
\setcounter{page}{1}
%\pagenumbering{roman}
\pagenumbering{arabic}
\chapter*{Introduction}
\addcontentsline{toc}{chapter}{Introduction}
\pagestyle{fancyplain}
%\headwidth
\rhead[\fancyplain{}{\sffamily\footnotesize Introduction}]%
      {\fancyplain{}{\thepage}}
\chead[\fancyplain{}{}]{\fancyplain{}{}}
\lhead[\fancyplain{}{\thepage}]%
      {\fancyplain{}{\sffamily\footnotesize Introduction}}
\cfoot{}

%erst allgemeine einführung ins thema

In this thesis we investigate the Picard-Fuchs equation of special one-parameter families of Calabi-Yau varieties. Calabi-Yau varieties have been studied in much detail, especially in Mirror Symmetry. Much of the early interest in this field focused on toric varieties. This is mostly due to Batyrev \cite{MR1269718}, who showed that for hypersurfaces in toric varieties duality in the sense of Mirror Symmetry can be reduced to polar duality between polytopes of toric varieties. This was the starting point for many achievements in Mirror Symmetry of Calabi-Yau varieties. The work of Batyrev, however, does not cover the families in weighted projective space that we consider in this thesis. In particular, Batyrev requires an ambient space that is Gorenstein. This implies that every weight divides the sum of all weights. In the case of hypersurfaces in weighted projective spaces this restricts the class covered by Batyrev's approach to polynomials of Brieskorn-Pham type.

The hypersurfaces we investigate are defined by invertible polynomials. These are weighted homogeneous polynomials $g(x_1,\dots,x_n)\in\C[x_1,\dots,x_n]$, which are a sum of exactly $n$ monomials, such that the weights $q_1,\dots,q_n$ of the variables $x_1,\dots,x_n$ are unique up to a constant and the affine hypersurface defined by the polynomial has an isolated singularity at $0$. The class of invertible polynomials includes all polynomials of Brieskorn-Pham type, but is much bigger. These polynomials were already studied by Berglund and H\"ubsch \cite{MR1191434}, who showed that a mirror manifold is related to a dual polynomial. For an invertible polynomial $g(x_1,\dots,x_n)=\sum_{j=1}^n\prod_{i=1}^n x_i^{E_{ij}}$ the dual polynomial $g^t(x_1,\dots,x_n)$ is defined by transposing the exponent matrix $E=(E_{ij})_{i,j}$ of the original polynomial, so $g^t(x_1,\dots,x_n)=\sum_{j=1}^n\prod_{i=1}^n x_i^{E_{ji}}$. If the polynomial is of Brieskorn-Pham type then the polynomial is in the above sense self-dual. This work was made precise by Krawitz, et al. (cf. \cite{MR2606631}, \cite{Kr}), where an isomorphism is given between the FJRW-ring of the polynomial (cf. \cite{FJR}) and a quotient of the Milnor ring of the dual polynomial. In addition Chiodo and Ruan \cite{CR} have made progress by stating the so called Landau-Ginzburg/Calabi-Yau correspondence for invertible polynomials. Among other things this includes the statement that the Chen-Ruan orbifold cohomology of the Mirror Partners interchanges. Recently, Borisov \cite{Bor} developed a theory combining his work with Batyrev on toric varieties \cite{MR1416334} in Mirror Symmetry and the work of Krawitz on invertible polynomials in Mirror Symmetry \cite{Kr}.
\newpage
In this thesis we analyse the Picard-Fuchs equations of the one-parameter families of hypersurfaces. The Picard-Fuchs equation is a differential equation that is satisfied by the periods of the family, i.e. the integrals of a form over a basis of cycles. These differential equations have been studied by many people and this can lead to several aspects of Mirror Symmetry. For example, Morrison \cite{MR1191426} used the Picard-Fuchs equations of hypersurfaces to calculate the mirror map and Yukawa couplings for mirror manifolds. In \cite{MR2369490} Chen, Yang and Yui study the monodromy for Picard-Fuchs equations of Calabi-Yau threefolds in terms of monodromy groups. These give two potential applications of the results of this thesis to further research.

We consider a special one-parameter family over an invertible polynomial and calculate the Picard-Fuchs equation for this family. In detail we start with an invertible polynomial $g(x_1,\dots,x_n)$, and in addition we require that the weights $q_1,\dots,q_n$ of $g$ add up to the degree $d$ of $g$. This is called the Calabi-Yau condition, because in \cite{MR704986} Dolgachev showed that under this condition the canonical bundle of the hypersurface $\{g(x_1,\dots,x_n)=0\}\subset \mathbb{P}(q_1,\dots,q_n)$ is trivial. The special one-parameter family we are dealing with is given by
\begin{align*}
f(x_1,\dots,x_n):=g(x_1,\dots,x_n)+s\prod_{i=1}^n x_i,
\end{align*}
where $s$ is a parameter. We calculate the Picard-Fuchs equation for this one-parameter family by using the Griffiths-Dwork method, which provides an algorithm to calculate the Picard-Fuchs equation (cf. \cite{MR1677117}). Unfortunately this method of calculations can be quite computationally expensive. Therefore we develop a combinatorial approach for the Griffiths-Dwork method. This approach, among other things, allows us to prove the order of the Picard-Fuchs equation. With this statement and the computation of the \GKZ system satisfied by the same periods, we can prove a general formula for the Picard-Fuchs equation. For a one-parameter family $f$ defined above the Picard-Fuchs equation is given by
\begin{align*}
0=\prod_{i=1}^n \widehat{q}_i^{\widehat{q}_i}s^{\widehat{d}}\prod_{i=1}^{n}\prod_{j=0}^{\widehat{q}_i-1}(\delta+\frac{j\cdot \widehat{d}}{\widehat{q}_i})\prod_{\ell\in I}(\delta+\ell)^{-1}-(-\widehat{d})^{\widehat{d}}\prod_{j=0}^{\widehat{d}-1}(\delta-j)\prod_{\ell\in I}(\delta-\ell)^{-1},
\end{align*}
where $\widehat{q}_1,\dots,\widehat{q}_n$ are the weights of the dual polynomial $g^t$, $\widehat{d}$ is the degree of $g^t$, and $I=\{0,\dots,\widehat{d}-1\}\cap\bigcup_{i=1}^n \left\{0,\frac{\widehat{d}}{\widehat{q}_i},\frac{2\widehat{d}}{\widehat{q}_i},\dots,\frac{(\widehat{q}_i-1)\widehat{d}}{\widehat{q}_i}\right\}$.

One interesting observation is that the Picard-Fuchs equation consists only of the data given by the dual polynomial, namely the dual weights and the dual degree. As pointed out to us by Stienstra, this Picard-Fuchs equation was already obtained in a work by Corti and Golyshev \cite{CG} in the context of local systems and Landau-Ginzburg pencils. Our combinatorial approach however is very constructive and yields not only the Picard-Fuchs equation itself but computes a basis of the important part of the cohomology. These computations will again relate to the duality between the polynomials. In addition we are able to show for certain values of the parameter a 1-1 correspondence between the roots of the Picard-Fuchs equation of $f$, the Poincar\'{e} series of the dual polynomial $g^t$ and the monodromy in the solution space of the Picard-Fuchs equation.

One important class that will be studied in detail in this thesis is the case of the 14 exceptional unimodal hypersurface singularities that are part of Arnold's strange duality \cite{MR0420689}. The duality between these singularities was known before Mirror Symmetry, but was shown to fit into the language of Mirror Symmetry (cf. \cite{MR1420220}). We will not only calculate the Picard-Fuchs equation here, but also investigate the structure of the cohomology which is used in the calculations for the Picard-Fuchs equation.
\subsubsection*{Structure of the thesis}
This thesis starts with some preliminaries in Chapter \ref{eins}. We recall definitions and relevant statements needed for discussions and fix notation here. The first section of this chapter is devoted to invertible polynomials and the duality among them. After that, in the second section we concentrate on the Griffiths-Dwork method to calculate the Picard-Fuchs equations. In Subsection \ref{gd-klassisch} we present the method in general for hypersurfaces in weighted projective spaces and in Subsection \ref{gd-comb} we introduce a new combinatorial notation which will be used to reconstruct the Griffiths-Dwork method in detail for one-parameter families $f$ of the form mentioned above.

The main goal of Chapter \ref{haupt} is to prove the order of the Picard-Fuchs equation for a one-parameter family. To achieve this goal we will investigate the structure of the underlying combinatorics. This will shorten the calculations for the Griffiths-Dwork method, and also allows us to construct explicitly the forms involved in the calculations of the Picard-Fuchs equation. One important ingredient for the whole procedure to work nicely will be that the Calabi-Yau condition holds, so the weights of the polynomial add up to the degree. In the last section (Section \ref{detailed-example-gd}) of this chapter we calculate the complete Picard-Fuchs equation with the Griffiths-Dwork method in one example.

Chapter \ref{vier} combines several results achieved so far in the thesis with results that can be found in the literature. The main theorem (see Section \ref{pifu}) presents the Picard-Fuchs equation for a one-parameter family associated to an invertible polynomial in general. We already calculated the order of the Picard-Fuchs equation and together with the \GKZ system, computed in Section \ref{calculation-GKZ}, the theorem is proved. We take advantage of the constructive proof of the order of the Picard-Fuchs equation in Section \ref{cohom}, where we investigate the cohomology of the hypersurface defined by the one-parameter family. In Section \ref{ASD} we calculate the Picard-Fuchs equations explicitly for the famous class of examples of the 14 exceptional unimodal hypersurface singularities. In addition to this being an important class, this was the origin of our work and most of the interesting phenomena can already be seen. In the last section (Section \ref{poincare}) we investigate the 1-1 correspondence between the Picard-Fuchs equation of a one-parameter family, the Poincar\'{e} series of the transposed polynomial and the monodromy in the solution space of the Picard-Fuchs equation.

Finally, the Appendix is divided into two parts. The first part (Appendix \ref{k3}) shows another class of examples for which we calculated the Picard-Fuchs equation. This special class of examples is extracted from the list of 93 hypersurfaces stated in \cite{MR1066667}. In Appendix \ref{singular} one can find the code of the algorithm that is provided by the Griffiths-Dwork method. The algorithm is written for Singular, but it can easily be adapted for any computer algebra system.

\textbf{Acknowledgements:} First of all I would like to thank Prof. Wolfgang Ebeling for supervising me and for always giving me the opportunity to take advantage of his knowledge, his good intuition and his enthusiasm. I would also like to thank Prof. Noriko Yui and Prof. Ragnar-Olaf Buchweitz, who were great hosts during my stay at the Fields Institute in Toronto and gave me the possibility to speak in front of a wide audience, and in particular I want to thank Prof. Yui for discussing my work in a very encouraging way. I want to gratefully mention Prof. Jan Stienstra, who showed me the relation between my work and the work of Golyshev and Corti and gave me a final input to the proof of the main theorem, and Prof. Victor Golyshev, who took the time to answer my questions on his paper. I am obliged to Dr. Nathan Broomhead who helped me solving several problems in his patient manner and in addition all other members of the Insitute of Algebraic Geometry at the Leibniz Universit\"at Hannover.\\

\newpage
\thispagestyle{empty}

%--------------------------   Ende der Einleitung   -----------------------------------------------------------------

%--------------------------   Beginn des Hauptteils: Seitenzahlen arabisch bei 1 beginnend --------------------------

\mainmatter
\setcounter{page}{5}
%\pagenumbering{arabic}
\pagestyle{fancyplain}
%\headwidth
\renewcommand{\chaptermark}[1]{\markboth{ #1}{}}    % \markboth{\thechapter{} #1}
\renewcommand{\sectionmark}[1]{\markright{\thesection\ #1}}
\rhead[\fancyplain{}{\sffamily\footnotesize\leftmark}]%
      {\fancyplain{}{\thepage}}
\chead[\fancyplain{}{}]{\fancyplain{}{}}
\lhead[\fancyplain{}{\thepage}]%
      {\fancyplain{}{\sffamily\footnotesize\rightmark}}
\cfoot{}

\chapter{Background on invertible polynomials and the Griffiths-Dwork method}\label{eins}
\section{Invertible polynomials}

We start this chapter by defining invertible polynomials and proving some properties we need later.

\begin{defi}
Let $$g(\underline{x})=\sum_{j=1}^m c_j\prod_{i=1}^n x_i^{E_{ij}}\in\mathbb{C}[\underline{x}]$$ be a quasihomogeneous polynomial with weights $q_1,\dots,q_n\in\Z$, where $\underline{x}=(x_1,\dots,x_n)$ and $E_{ij}\in\mathbb{N}$. Then $g(\underline{x})$ is an invertible polynomial, if the following conditions hold:
\begin{enumerate}[label=(\roman*)]
\item $\# \text{ variables }=\# \text{ summands }$, i.e. $n=m$,
\item the weights $q_1,\dots,q_n$ are unique up to scaling by a multiple of $gcd(q_1,\dots,q_n)^{-1}$ and
\item the Milnor ring $\mathbb{C}[\underline{x}]/J(g)$ has a finite basis, where $J(g)=\langle \frac{\partial g}{\partial x_1},\dots,\frac{\partial g}{\partial x_n}\rangle$ is the Jacobian ideal. This is equivalent to $\underline{0}$ being an isolated singularity of $\{g(\x)=0\}\subseteq\C^n$.
\end{enumerate}
\end{defi}

\begin{rem}
We want to state some conventions we are using throughout the thesis:
\begin{itemize}
\item We require the weights to be reduced, i.e. $\gcd(q_1,\dots,q_n)=1$. This way the weights are unique.
\item Some authors call the polynomial $g(\underline{x})$ invertible, if the first two conditions are satisfied and a non-degenerate invertible polynomial, if $g(\underline{x})$ satisfies all three conditions.
\item From now on we assume that the coefficients $c_j$ are all equal to $1$. This can always be achieved by an easy coordinate transformation.
\item The weights are also defined by the smallest numbers $q_1,\dots,q_n\in\N$ and $d\in\N$ satisfying the equation
\begin{align*}
\begin{pmatrix} E_{11} & \cdots & E_{1n} \\ \vdots & & \vdots \\ E_{n1} & \cdots & E_{nn} \end{pmatrix} \begin{pmatrix} q_1\\ \vdots \\ q_n \end{pmatrix} = \begin{pmatrix} d \\ \vdots \\ d \end{pmatrix}
\end{align*}
or concisely $E\cdot\underline{q}=\underline{d}$. We call $E$ the \emph{exponent matrix}
\end{itemize}
\end{rem}

M.~Kreuzer and H.~Skarke showed that the polynomials which are invertible are a composition of only two types.

\begin{thm}\emph{(Kreuzer and Skarke \cite{MR1188500})}\label{kreuzer-skarke}
Every invertible polynomial is a sum of polynomials with distinct variables of the following two types
\begin{align*}
\text{loop:   }& x_1^{k_1}x_2+x_2^{k_2}x_3+\dots+x_{m-1}^{k_{m-1}}x_m+x_m^{k_m}x_1 &\text{ for }m\geq2\\
\text{chain:   }& x_1^{k_1}x_2+x_2^{k_2}x_3+\dots+x_{m-1}^{k_{m-1}}x_m+x_m^{k_m} &\text{ for }m\geq1
\end{align*}
\end{thm}

\begin{ex}
We want to list two very famous class of examples here.
\begin{enumerate}[label=(\roman*)]
\item A polynomial is of Brieskorn-Pham type if it is of the form $g(\x)=\sum_{i=1}^n x_i^{k_i}$ with $k_i\in\N$. In this case the polynomial is always invertible and the exponent matrix is a diagonal matrix with the exponents $k_i$ on the diagonal. It follows that $q_i=\frac{\textrm{lcm}(k_1,\dots,k_n)}{k_i}$ and $d=\textrm{lcm}(k_1,\dots,k_n)$.
\item For the 14 exceptional unimodal singularities, invertible polynomials can be chosen. Table \ref{ASD-normal} lists their name, invertible polynomial, reduced weights and degree in the first four columns. In the last columns the dual singularity due to Arnol'd \cite{MR0420689} is listed. In the next definition we will see how this duality fits into the context of invertible polynomials which also explains the rest of the table. The example of Arnold's strange duality will be studied in detail in Section \ref{ASD}.
\end{enumerate}
\end{ex}

In their paper \cite{MR1191434} P.~Berglund and T.~H\"ubsch proposed a way to define dual pairs of invertible polynomials by transposing the exponent matrix.

\begin{defi}
If $g(\underline{x})=\sum_{j=1}^n\prod_{i=1}^n x_i^{E_{ij}}$ is an invertible polynomial then the Berglund-H\"ubsch transpose is given by
$$g^t(\underline{x})=\sum_{j=1}^n\prod_{i=1}^n x_i^{E_{ji}}.$$
\end{defi}

\begin{ex}
As noticed before the dual singularities in the examples of Arnold's strange duality are given by transposed polynomials:
\begin{table}[H]
\begin{center}
\begin{tabular}{lcccccl}
\toprule
Name & $g(x,y,z)$ & Weights & Deg & Dual weights & $g^t(x,y,z)$ & Dual\\
\midrule
E$_{12}$ & $x^7+y^3+z^2$ & (6,14,21) & 42 & (6,14,21) & $x^7+y^3+z^2$ & E$_{12}$\\
&&&&&&\\
E$_{13}$ & $x^5y+y^3+z^2$ & (4,10,15) & 30 & (6,8,15) & $x^5+xy^3+z^2$ & Z$_{11}$\\
Z$_{11}$ & $x^5+xy^3+z^2$ & (6,8,15) & 30 & (4,10,15) & $x^5y+y^3+z^2$ & E$_{13}$\\
&&&&&&\\
E$_{14}$ & $x^4z+y^3+z^2$ & (3,8,12) & 24 & (6,8,9) & $x^4+y^3+xz^2$ & Q$_{10}$\\
Q$_{10}$ & $x^4+y^3+xz^2$ & (6,8,9) & 24 & (3,8,12) & $x^4z+y^3+z^2$ & E$_{14}$\\
&&&&&&\\
Z$_{12}$ & $x^4y+xy^3+z^2$ & (4,6,11) & 22 & (4,6,11) & $x^4y+xy^3+z^2$ & Z$_{12}$\\
&&&&&&\\
W$_{12}$ & $x^5+y^2z+z^2$ & (4,5,10) & 20 & (4,10,5) & $x^5+y^2+yz^2$ & W$_{12}$\\
&&&&&&\\
Z$_{13}$ & $x^3z+xy^3+z^2$ & (3,5,9) & 18 & (4,6,7) & $x^3y+y^3+xz^2$ & Q$_{11}$\\
Q$_{11}$ & $x^3y+y^3+xz^2$ & (4,6,7) & 18 & (3,5,9) & $x^3z+xy^3+z^2$ & Z$_{13}$\\
&&&&&&\\
W$_{13}$ & $x^4y+y^2z+z^2$ & (3,4,8) & 16 & (4,6,5) & $x^4+xy^2+yz^2$ & S$_{11}$\\
S$_{11}$ & $x^4+y^2z+xz^2$ & (4,5,6) & 16 & (3,8,4) & $x^4z+y^2+yz^2$ & W$_{13}$\\
&&&&&&\\
Q$_{12}$ & $x^3z+y^3+xz^2$ & (3,5,6) & 15 & (3,5,6) & $x^3z+y^3+xz^2$ & Q$_{12}$\\
&&&&&&\\
S$_{12}$ & $x^3y+y^2z+xz^2$ & (3,4,5) & 13 & (3,5,4) & $x^3z+xy^2+yz^2$ & S$_{12}$\\
&&&&&&\\
U$_{12}$ & $x^4+y^2z+yz^2$ & (3,4,4) & 12 & (3,4,4) & $x^4+y^2z+yz^2$ & U$_{12}$\\
\bottomrule
\end{tabular}
\end{center}
\caption{Arnold's strange duality}
\label{ASD-normal}
\end{table}
\end{ex}

\begin{rem}
Notice that taking the transpose does not change the type of the polynomial. The exponent matrix is a direct sum of matrices, where every summand belongs to a polynomial of chain or loop type. Therefore we can transpose every chain and loop separately:
\begin{itemize}
\item $g(\x)=x_1^{k_1}x_2+\dots+x_{m-1}^{k_{m-1}}x_m+x_m^{k_m}x_1 \Rightarrow g^t(\x)=x_m x_1^{k_1}+x_1 x_2^{k_2}\dots+x_{m-1}x_m^{k_m}$ and
\item $g(\x)=x_1^{k_1}x_2+\dots+x_{m-1}^{k_{m-1}}x_m+x_m^{k_m} \Rightarrow g^t(\x)=x_1^{k_1}+x_1 x_2^{k_2}\dots+x_{m-1}x_m^{k_m}$.
\end{itemize}
\end{rem}

\begin{defi}\label{1para}
Let $g(\x)$ be an invertible polynomial. We set $f(\x)$ to be the one-parameter family associated to $g(\x)$ via
$$f(\underline{x})=g(\underline{x})+s\prod_{i=1}^n x_i,$$ where $s$ denotes the parameter.
\end{defi}

This one-parameter family $f(\x)$ will be one of the main objects of interest in this thesis.
Because we still want this family to be quasihomogeneous, we require that the weights of $g(\underline{x})$ add up to the degree of $g(\underline{x})$.
In \cite{MR704986} I. Dolgachev showed that this is the condition for a quasihomogeneous polynomial to define a Calabi-Yau hypersurface.

\begin{prop}\emph{(\cite{MR704986})}\label{CYcond} Let $g(\underline{x})$ be a quasihomogeneous polynomial with weights $q_1,\dots,q_n$. Then $g(\underline{x})$ defines a hypersurface in $\mathbb{P}(q_1,\dots,q_n)$ that is Calabi-Yau, if $$\sum_{i=1}^nq_i=d=\deg g(\underline{x}).$$
\end{prop}

\begin{lemma}\label{dualCY-cond}
If the Calabi-Yau condition holds for the weights of an invertible polynomial then it also holds for the weights of the transposed polynomial.
\end{lemma}

\begin{notation}
For an invertible polynomial $g(\x)$ we denote the reduced weights with $q_1,\dots,q_n$ and $\deg g=d$. For the dual polynomial $g^t$ the weights are $\widehat{q}_1,\dots,\widehat{q}_n$ and $\deg g^t=\widehat{d}$. The diagonal entries of the exponent matrix $E$ are $k_1,\dots,k_n$. Notice that this are the same for $g$ and $g^t$.
\end{notation}

\begin{proof}
The Calabi-Yau condition is equivalent to the condition
\begin{align*}
\det\begin{pmatrix} \begin{array}{c}1\\\vdots\end{array}&E\\1&\begin{array}{cc}\cdots&1\end{array}  \end{pmatrix}=0
\end{align*}
This is due to the fact that the weights are unique up to scaling and therefore the linear relation between the rows of the above matrix have to be given by multiplying with the vector $(-d,q_1,\dots,q_n)^t$. Now it is obvious that if the above condition holds for $E$ it also holds for $E^t$.
\end{proof}

We want to investigate another relation of the dual weights here, which has to do with the partial derivatives of the one-parameter family $f(\x)$. This connection between the Jacobian ideal of $f(\x)$ and the exponent matrix of $g(\x)$ occurs again in the next chapter.

\begin{rem}\label{auftreten-duale-gewichte}
First assume $g(\x)=x_1^{k_1}x_2+\dots+x_{m-1}^{k_{m-1}}x_m+x_m^{k_m}x_1$ is a loop of length $m$, then we have $f(\x)=x_1^{k_1}x_2+\dots+x_{m-1}^{k_{m-1}}x_m+x_m^{k_m}x_1+sx_1\cdots x_m$. The dual weights of $g(\x)$ can be calculated via the equation
\begin{align*}
\begin{pmatrix} k_1&0&\cdots&~&0&1\\1&k_2&0&\cdots&~&0\\0&1&k_3&0&\cdots&0\\\vdots&~&\ddots&\ddots&~&\vdots\\0&\cdots&0&1&k_{m-1}&0\\0&\cdots&~&0&1&k_m \end{pmatrix}\cdot \begin{pmatrix} \widehat{q}_1\\ \vdots\\~\\~\\ \vdots\\\widehat{q}_n \end{pmatrix}= \begin{pmatrix} \widehat{d}\\ \vdots\\~\\~\\ \vdots\\\widehat{d} \end{pmatrix}
\end{align*}
Because the dual weights of $g(\x)$ also satisfy the Calabi-Yau equation due to Lemma \ref{dualCY-cond} we get the following equation:
\begin{align*}
0=&\left(\begin{pmatrix} k_1&0&\cdots&~&0&1\\1&k_2&0&\cdots&~&0\\0&1&k_3&0&\cdots&0\\\vdots&~&\ddots&\ddots&~&\vdots\\0&\cdots&0&1&k_{m-1}&0\\0&\cdots&~&0&1&k_m \end{pmatrix}-\begin{pmatrix} 1&\cdots&~&~&\cdots&1\\\vdots&\ddots&~&~&~&\vdots\\~&~&~&~&~&~\\~&~&~&~&~&~\\\vdots&~&~&~&\ddots&\vdots\\1&\cdots&~&~&\cdots&1\end{pmatrix}\right)\cdot \begin{pmatrix} \widehat{q}_1\\ \vdots\\~\\~\\ \vdots\\\widehat{q}_n \end{pmatrix} \\
=&\begin{pmatrix} k_1-1&-1&\cdots&~&-1&0\\0&k_2-1&-1&\cdots&~&-1\\-1&0&k_3-1&-1&\cdots&-1\\ \vdots&~&\ddots&\ddots&~&\vdots\\-1&\cdots&-1&0&k_{m-1}-1&-1\\-1&\cdots&~&-1&0&k_m-1 \end{pmatrix}\cdot \begin{pmatrix} \widehat{q}_1\\ \vdots\\~\\~\\ \vdots\\ 
\widehat{q}_n \end{pmatrix}
\end{align*}
This is interesting, because the $i$th column of this matrix is connected to $\frac{\partial f}{\partial x_i}=k_ix_i^{k_i-1}x_{i+1}+x_{i-1}+s\prod_{j\neq i}x_j$ in the sense that the $i$th column is given by subtracting the exponent vector of the summand $k_ix_i^{k_i-1}x_{i+1}$ by the exponent vector of the summand $s\prod_{j\neq i}x_j$. Of course the index is taken modulo $m$.\\
The same happens if $g(\x)=x_1^{k_1}x_2+\dots+x_{m-1}^{k_{m-1}}x_m+x_m^{k_m}$ is a polynomial of chain type. The equation from above is now given by
\begin{align*}
0=&\left(\begin{pmatrix} k_1&0&\cdots&~&0&0\\1&k_2&0&\cdots&~&0\\0&1&k_3&0&\cdots&0\\\vdots&~&\ddots&\ddots&~&\vdots\\0&\cdots&0&1&k_{m-1}&0\\0&\cdots&~&0&1&k_m \end{pmatrix}-\begin{pmatrix} 1&\cdots&~&~&\cdots&1\\\vdots&\ddots&~&~&~&\vdots\\~&~&~&~&~&~\\~&~&~&~&~&~\\\vdots&~&~&~&\ddots&\vdots\\1&\cdots&~&~&\cdots&1\end{pmatrix}\right)\cdot \begin{pmatrix} \widehat{q}_1\\ \vdots\\~\\~\\ \vdots\\\widehat{q}_n \end{pmatrix}\\
=&\begin{pmatrix} k_1-1&-1&\cdots&~&-1&-1\\0&k_2-1&-1&\cdots&~&-1\\-1&0&k_3-1&-1&\cdots&-1\\\vdots&~&\ddots&\ddots&~&\vdots\\-1&\cdots&-1&0&k_{m-1}-1&-1\\-1&\cdots&~&-1&0&k_m-1 \end{pmatrix}\cdot \begin{pmatrix} \widehat{q}_1\\ \vdots\\~\\~\\ \vdots\\\widehat{q}_n \end{pmatrix}
\end{align*}
and the partial derivatives of $f(\x)$ are $\frac{\partial f}{\partial x_1}=k_1x_1^{k_1}+s\prod_{j\neq 1}x_j$, $\frac{\partial f}{\partial x_i}=k_ix_i^{k_i-1}x_{i+1}+x_{i-1}+s\prod_{j\neq i}x_j$ for $i=2,\dots,m-1$ and $\frac{\partial f}{\partial x_m}=k_mx_m^{k_m-1}+x_{m-1}+s\prod_{j\neq m}x_j$. If we again subtract the exponent vector of the summand $s\prod_{j\neq i}x_j$, for $i=1,\dots,m$, in the partial derivative from the exponent vector of the summand $k_1x_1^{k_1}$, $k_ix_i^{k_i-1}x_{i+1}$ or $k_mx_m^{k_m-1}$ in the partial derivative then the result is exactly given by the columns of the matrix above.
\end{rem}
\section{Introduction to the Griffiths-Dwork method}

This section is divided into two parts. In the first part we explain the Griffiths-Dwork method, which is a well-known method to calculate the Picard-Fuchs equation of a one-parameter family of hypersurfaces. In the second part we will introduce our own combinatorial notation to describe the Griffiths-Dwork method. Using this we are able to describe the Griffiths-Dwork method in Chapter \ref{haupt} in a sufficiently abstract way. So we are able to present facts about the Picard-Fuchs equation in general.\\
Before we start describing the Griffiths-Dwork method, we want to recall the definition of a Picard-Fuchs equation.

\begin{defi}
The Picard-Fuchs equation of a one-parameter family $f(\x)$ of hypersurfaces is defined as the ordinary differential equation with differential operator $\delta=s\frac{\partial}{\partial s}$, where $s$ is the parameter, which has as solutions exactly the period integrals. So the solutions are given by $\int_{\gamma_i} \omega$ for a basis $\{\gamma_i\}$ of $H_{n-2}(V(f))$ and $\omega\in H^{n-2}(V(f))$.
\end{defi}

\subsection{The Griffiths-Dwork method}\label{gd-klassisch}

In this section we want to repeat how the Griffiths-Dwork method works. The method is due to Griffiths \cite{MR0260733}, Dwork \cite{MR0159823} and the generalization to weighted polynomials was done by Dolgachev \cite{MR704986}. A good reference for this method in general is Chapter 5.3 of the book by Cox and Katz \cite{MR1677117}. We will not do everything in general but restrict ourselves to the cases which are important to us. In particular we will only consider the one-parameter family $f(\x)=g(\x)+s\prod_{i=1}^n x_i\in\C(s)[\x]$ defined in \ref{1para} and explain the Griffiths-Dwork method for these families.\\
We will recall and fix the notation we already used in the last section. We use this notation throughout the rest of the thesis.

\begin{notation}
The polynomial $g(\x)=\sum_{j=1}^n\prod_{i=1}^n x_i^{E_{ij}}$ is an invertible polynomial. The diagonal entries of the exponent matrix $E$ are denoted by $k_1,\dots,k_n$. The polynomial $g(\x)$ is quasihomogeneous with weights $q_1,\dots,q_n$ and degree $d$. We define $f(\x)=g(\x)+s\prod_{i=1}^n x_i$ to be a one-parameter family with parameter $s$.
\end{notation}

We want to calculate the Picard-Fuchs equation for $f(\x)$. The first step is to describe the cohomology $H^{n-2}(V(f))$ in more detail. For this we use the residue map
\begin{align*}
\text{Res}:H^{n-1}(\mathbb{P}(q_1,\dots,q_n)\backslash V(f))\twoheadrightarrow PH^{n-2}(V(f))\subseteq H^{n-2}(V(f)),
\end{align*}
where $PH^{n-2}(V(f))=\{\eta\in H^{n-2}(V(f))|\,\eta\cdot H=0\text{ for the hyperplane class }H\}$ is the primitive cohomology. Note that $PH^{n-2}(V(f))= H^{n-2}(V(f))$ if $n$ is even.
The advantage of this map is that the cohomology $H^{n-1}(\mathbb{P}(q_1,\dots,q_n)\backslash V(f))$ was explicitly described by Griffiths in \cite{MR0260733} and we can use the Residue map to carry this description over to the cohomology of the hypersurface. In detail, the classes in $H^{n-1}(\mathbb{P}(q_1,\dots,q_n)\backslash V(f))$ can be represented by forms of the form $\frac{Q\Omega_0}{f^l}$, where $\Omega_0=\sum_{i=1}^n(-1)^jd_jx_jdx_1\wedge\dots\wedge\widehat{dx_j}\wedge\dots\wedge dx_n$, $l\in\mathbb{N}$ and $Q\in\C(s)[\x]$ is a polynomial with $\deg Q=(\deg f)(l-1)$.
We can now define the residue map as follows:
\begin{align}\label{res}
\int_{\gamma}\text{Res }\frac{Q\Omega_0}{f^l}=\int_{T_{\gamma}}\frac{Q\Omega_0}{f^l}
\end{align}
for an $(n-2)$-cycle $\gamma$ and $T_{\gamma}$ a tubular neighbourhood around $\gamma$.

Let us go back to our goal. The Picard-Fuchs equation is of the form
\begin{align*}
0=(p_r(s)\delta^r+\dots+p_1(s)\delta+p_0(s))\left(\int_{\gamma_i} \omega\right),
\end{align*}
where $p_i(s)\in\C(s)$. Suppose $\omega\in PH^{n-2}(V(f))$, so $\omega=\text{Res }\frac{Q\Omega_0}{f^l}$ for some $Q,f,l$. Using equation (\ref{res}) we get
\begin{align*}
0=&(p_r(s)\delta^r+\dots+p_1(s)\delta+p_0(s))\left(\int_{\gamma_i} \text{Res }\frac{Q\Omega_0}{f^l}\right)\\
=&\int_{T_{\gamma}}\left(p_r(s)\delta^r\frac{Q\Omega_0}{f^l}+\dots+p_1(s)\delta\frac{Q\Omega_0}{f^l}+p_0(s)\frac{Q\Omega_0}{f^l}\right)\\
=&\int_{\gamma_i}\left(p_r(s)\delta^r\omega+\dots+p_1(s)\delta\omega+p_0(s)\omega\right),
\end{align*}
because the integral commutes with the differential operator. This means if we find a differential equation satisfied by the $(n-2)$-form $\omega$, then this also holds for the period integrals of $\omega$. From now on we will write $\frac{Q\Omega_0}{f^l}$ instead of $\text{Res }\frac{Q\Omega_0}{f^l}$ for a form in $PH^{n-2}(V(f))$. The idea is now to calculate $\delta^i\frac{Q\Omega_0}{f^l}$ for $i=0,1,\dots$ until we find a linear relation between these forms. Unfortunately this is not so easy to do, because as $i$ increases, the pole order $l$ also increases. The Griffiths-Dwork method tells us how to solve this problem. The primitive cohomology can be compared with the Milnor ring by the following isomorphism
\begin{align*}
\left(\mathbb{C}(s)[\x]/J(f)\right)_{(\deg f)(l-1)} & \cong PH^{n-l-1,l-1}(V(f))\\
Q & \mapsto \frac{Q\Omega_o}{f^l},
\end{align*}
where the subscript $(\deg f)(l-1)$ denotes the graded piece of degree $(\deg f)(l-1)$ in the Milnor ring.
The fundamental ingredient to this isomorphism is the Griffiths formula (cf. Theorem 4.3 in \cite{MR0260733}) that tells us how and when to reduce the pole order of an $(n-2)$-form:
\begin{align}\label{griffiths-formula}
\frac{(l-1)\sum_{j=1}^nG_j\frac{\partial f}{\partial x_j}\Omega_0}{f^l}=\frac{\sum_{j=1}^n\frac{\partial G_j}{\partial x_j}\Omega_0}{f^{l-1}}\text{ (modulo exact forms)}.
\end{align}
This can be seen easily from the following calculation:
\begin{multline*}
\frac{(l-1)\sum_{j=1}^nG_j\frac{\partial f}{\partial x_j}\Omega_0}{f^l}-\frac{\sum_{j=1}^n\frac{\partial G_j}{\partial x_j}\Omega_0}{f^{l-1}}\\
=d\left(\frac{1}{f^{l-1}}\sum_{i<j}\left(x_iG_j-x_jG_i\right)dx_0\wedge\dots\wedge\widehat{dx_i}\wedge\dots\wedge\widehat{dx_j}\wedge\dots\wedge dx_n\right).
\end{multline*}

The big advantage is now that all computations can be done with a Gr\"obner basis in the Milnor ring and the Picard-Fuchs equation can be calculated with the following algorithm.

\begin{alg}\emph{(cf. Cox and Katz \cite{MR1677117}, Section 5.3)}\label{alg}
With the following steps one can determine the Picard-Fuchs equation for the one-parameter family $f(\x)$ with parameter $s$:
\begin{enumerate}[label=(\roman*)]
\item Find a basis $B$ of the Milnor Ring $\mathbb{C}(s)[\x]/J(f)$ in degree $d(l-1)$ for $1\leq l\leq n-1$ (this is equivalent to having a basis of the primitive cohomology).
\item Write $\delta^i\omega=\left(s\frac{\partial}{\partial s}\right)^i\omega$ in the basis $B$ for all $0\leq i\leq |B|$. This is done by writing $\delta^i\omega$ as a sum of a part that is in the basis and a part that is in the Jacobian ideal and can therefore be reduced with the Griffiths formula. After reducing this process can be repeated until pole order $0$ is reached.
\item Now there are $|B|$ basis elements and $|B|+1$ derivatives of $\omega$, so there is a linear relation between them. The linear relation between the $\delta^i\omega$ gives the Picard-Fuchs equation of $f$.
\end{enumerate}
\end{alg}

\begin{rem}
One could ask why it is still interesting to investigate the Griffiths-Dwork method in even more detail. The reason is that some of the calculations done in the above algorithm are very expensive. Furthermore, it very often happens that in the calculations not all elements of the basis of the Milnor ring are needed. The goal is therefore to find an abstract way to describe the steps in the Griffiths-Dwork method and try to restrict the calculations to a minimum.
\end{rem}

\begin{rem}\label{delta-omega}
If we have an invertible polynomial $g(\x)$ that satisfies the Calabi-Yau condition \ref{CYcond} and the one-parameter family we are looking at is $f(\x)=g(\x)+s\prod x_i$, then we can easily calculate $\delta^i\omega$ for $\omega=\frac{s\Omega_0}{f}$ and all $i\geq 0$:
\begin{align*}
\delta\omega=&\frac{s\Omega_0}{f}-\frac{s^2\prod x_i\Omega_0}{f^2}\\
\delta^2\omega=&\frac{s\Omega_0}{f}-3\frac{s^2\prod x_i\Omega_0}{f^2}+\frac{2s^3(\prod x_i)^2\Omega_0}{f^3}\\
\delta^3\omega=&\frac{s\Omega_0}{f}-7\frac{s^2\prod x_i\Omega_0}{f^2}+6\frac{2s^3(\prod x_i)^2\Omega_0}{f^3}-\frac{6s^4(\prod x_i)^3\Omega_0}{f^4}\\
\dots&\\
\delta^i\omega=&\sum_{m=0}^{i}(-1)^m r_m^i \frac{m!s^{m+1}(\prod x_j)^m\Omega_0}{f^{m+1}}, \\
\end{align*}
where $r_m^i=-r_{m-1}^{i-1}+(m+1)r_m^{i-1}$ for $i,m\geq 1,m<i$ with $r_0^m=r_m^m=1$ for all $m$ and $r_1^0=1$.\\
This means in the second step of the Algorithm \ref{alg} we have to write every $\frac{m!s^{m+1}(\prod x_j)^m\Omega_0}{f^{m+1}}$ in the basis $B$ of the Milnor ring.
\end{rem}

\begin{rem}
In practice we are going to interchange the first and the second step in Algorithm \ref{alg}. The basic idea is to first write the $\delta^i\omega$ with monomials of degree $\leq d(n-2)$ and then see which of them are linearly independent in the Milnor ring and choose a basis this way.
\end{rem}

\subsection{The Griffiths-Dwork method for invertible polynomials and its combinatorics}\label{gd-comb}

In this section we want to show a possibility how to give a diagrammatic version of the Griffiths-Dwork method. This means we will develop diagrams for all the steps in the Griffiths-Dwork method. This is helpful later to reduce the algorithm to the important parts and do the steps in a clear way to see what is happening there.

From now on we restrict ourselves completely to invertible polynomials. So $g(\x)=\sum_{i=1}^{n}\prod_{i=1}^n x_i^{E_{ij}}$ is an invertible polynomial with weights $q_1,\dots,q_n$ and $\deg g=d$. We denote by $k_1,\dots,k_n$ the diagonal entries of the matrix $E=(E_{ij})_{i,j}$. These are the only entries $\neq 0,1$ in this matrix. The one-parameter family $f$ is given by $f(\x)=g(\x)+s\prod_{i=1}^{n}x_i$. We will assume that the weights of $g$ satisfy the Calabi-Yau condition \ref{CYcond}, so that $f$ is still weighted homogeneous.

First we will have a closer look at the Jacobian ideal $J(f)$. We start with an invertible polynomial, so every variable can appear in at most two terms of $g$ or equivalently 3 terms of $f$. The possibilities for the terms that contain the variable $x_i$ in $g$ are the following:
\begin{enumerate}[label=(\roman*)]
\item $x_i^{k_i}$, which occurs if $x_i$ is in a chain of length $1$,
\item $x_i^{k_i}+x_ix_{i-1}^{k_{i-1}}$, which occurs if $x_i$ is the end of a chain of length $\geq 2$,
\item $x_i^{k_i}x_{i+1}$, which occurs if $x_i$ is the beginning of a chain of length $\geq 2$ or
\item $x_i^{k_i}x_{i+1}+x_ix_{i-1}^{k_{i-1}}$, which occurs if $x_i$ is in the middle of a chain of length $\geq 3$ or in a loop of arbitrary length.
\end{enumerate}
Therefore there are only 4 possibilities for the partial derivative of $f$ with respect to $x_i$:
\begin{enumerate}[label=(\roman*)]
\item $\frac{\partial}{\partial x_i}f=k_ix_i^{k_i-1}+s\prod_{j\neq i}x_j$,
\item $\frac{\partial}{\partial x_i}f=k_ix_i^{k_i-1}+x_{i-1}^{k_{i-1}}+s\prod_{j\neq i}x_j$,
\item $\frac{\partial}{\partial x_i}f=k_ix_i^{k_i-1}x_{i+1}+s\prod_{j\neq i}x_j$, or
\item $\frac{\partial}{\partial x_i}f=k_ix_i^{k_i-1}x_{i+1}+x_{i-1}^{k_{i-1}}+s\prod_{j\neq i}x_j$.
\end{enumerate}

\subsection*{Partial derivatives in a diagrammatic way}
Now we want to write these partial derivatives in a diagrammatic way. To do this we will not write down monomials, but restrict ourselves to the exponents. So instead of writing $\prod_{i=1}^nx_i^{a_i}$ we write the tuple of exponents $(a_1,a_2,\dots,a_n)$.

In the next step we want to write the sum of monomials in the diagrammatic notation. Because the Jacobian ideal of $f$ is only generated by sums of two or three monomials, we concentrate on how to write these sums of two or three monomials in the Jacobian ideal in a good way.\\
First let us assume our partial derivative with respect to $x_i$ is a sum of two monomials, e.g. $x_i$ occurs in $f$ as $k_ix_i^{k_i-1}+s\prod_{j\neq i}x_j$ or $k_ix_i^{k_i-1}x_{i+1}+s\prod_{j\neq i}x_j$. Then we can identify the two involved monomials by two points given by the exponents and indicate the fact that they are in a sum by an arrow that points to the monomial which has the parameter $s$ as coefficient. The coefficient $k_i$ of the other monomial will be written at the beginning of the arrow. Later on if it is not important we will omit all coefficients to reduce the notation to the essential information. So we would write the partial derivative $k_ix_i^{k_i-1}+s\prod_{j\neq i}x_j$ as
\begin{figure}[H]
\begin{center}
\begin{tikzpicture}
\node[place] (0) {$0,\dots,0,k_i-1,0,\dots,0$};
\node[place] (1) [right=70pt of 0]{$1,\dots,1,0,1,\dots,1$};
\draw [-latex,shorten >=2pt, shorten <=2pt,thick] (0) -- (1)
 node[very near start,above] {{\small $k_i$}};
\end{tikzpicture}
\end{center}
\caption{Diagram associated to $k_ix_i^{k_i-1}+s\prod_{j\neq i}x_j$}
\end{figure}
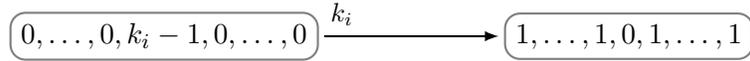
and similarly the sum $k_ix_i^{k_i-1}x_{i+1}+s\prod_{j\neq i}x_j$ would be represented by
\begin{figure}[H]
\begin{center}
\begin{tikzpicture}
\node[place] (0) {$0,\dots,0,k_i-1,1,0,\dots,0$};
\node[place] (1) [right=70pt of 0]{$1,\dots,1,0,1,1,\dots,1$};
\draw [-latex,shorten >=2pt, shorten <=2pt,thick] (0) -- (1)
 node[very near start,above] {{\small $k_i$}};
\end{tikzpicture}
\end{center}
\caption{Diagram associated to $k_ix_i^{k_i-1}x_{i+1}+s\prod_{j\neq i}x_j$}
\end{figure}
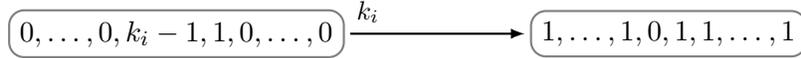
Now let us study what happens if we multiply a partial derivative by a monomial. Multiplication with a monomial does not change the number of summands. So we still end up with a sum of two or three monomials and the exponent of the monomial just gets added to the exponents of the partial derivative. For example if we multiply $k_ix_i^{k_i-1}+s\prod_{j\neq i}x_j$ by a monomial $m=c\prod x_i^{a_i}$ then the product is represented in our new notation, by
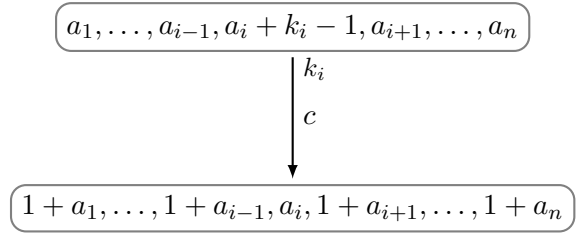
\begin{figure}[H]
\begin{center}
\begin{tikzpicture}
\node[place] (0) {$a_1,\dots,a_{i-1},a_i+k_i-1,a_{i+1},\dots,a_n$};
\node[place] (1) [below=50pt of 0]{$1+a_1,\dots,1+a_{i-1},a_i,1+a_{i+1},\dots,1+a_n$};
\draw [-latex,shorten >=2pt, shorten <=2pt,thick] (0) -- (1)
 node[very near start,right] {{\small $k_i$}}
 node[midway,right] {$c$};
\end{tikzpicture}
\end{center}
\caption{Multiplication of $k_ix_i^{k_i-1}+s\prod_{j\neq i}x_j$ by a monomial $m=c\prod x_i^{a_i}$}
\end{figure}
We keep track of the coefficient $c$ next to the middle of the arrow. Again if the information is not important, we will omit the coefficient. A sum of two monomials can be written as a monomial times a partial derivative, if the difference between the two monomials is $(1,\dots,1,-k_i+1,1,\dots,1)$ and there is a partial derivative of the form $k_ix_i^{k_i-1}+s\prod_{j\neq i}x_j$ or if the difference is $(1,\dots,1,k_i-1,0,1,\dots,1)$, the $(i+1)$th entry is $>0$ and there is a partial derivative of the form $k_ix_i^{k_i-1}x_{i+1}+s\prod_{j\neq i}x_j$. Of course the coefficients have to fit, but this will not need extra attention here.

Now let us discuss the case that the partial derivative is a sum of three monomials, so it is either $k_ix_i^{k_i-1}+x_{i-1}^{k_{i-1}}+s\prod_{j\neq i}x_j$ or $k_ix_i^{k_i-1}x_{i+1}+x_{i-1}^{k_{i-1}}+s\prod_{j\neq i}x_j$. As in the case of two monomials, we connect all monomials that form the partial derivative. This is best understood in an example.
First consider the case $k_ix_i^{k_i-1}+x_{i-1}^{k_{i-1}}+s\prod_{j\neq i}x_j$, this leads to the picture
\begin{figure}[H]
\begin{center}
\begin{tikzpicture}
\node[place] (0) {$0,\dots,0,0,k_i-1,0,\dots,0$};
\node[place] (1) [right=50pt of 0]{$1,\dots,1,1,0,1,\dots,1$};
\node[point] (p) [right=25pt of 0]{};
\node[place] (2) [above=20pt of p] {$0,\dots,0,k_{i-1},0,0\dots,0$};
\draw [-latex,shorten >=2pt, shorten <=2pt,thick] (0) -- (1)
 node[very near start,above] {{\small $k_i$}};
\draw [-,thick] (p) -- (2);
\end{tikzpicture}
\end{center}
\caption{Diagram associated to $k_ix_i^{k_i-1}+x_{i-1}^{k_{i-1}}+s\prod_{j\neq i}x_j$}
\end{figure}
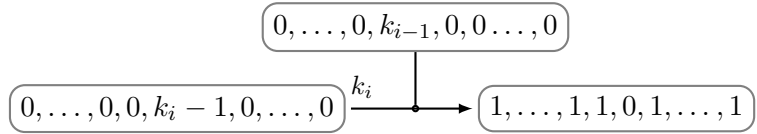
and for the case $k_ix_i^{k_i-1}x_{i+1}+x_{i-1}^{k_{i-1}}+s\prod_{j\neq i}x_j$ the picture is given by
\begin{figure}[H]
\begin{center}
\begin{tikzpicture}
\node[place] (0) {$0,\dots,0,0,k_i-1,1,0,\dots,0$};
\node[place] (1) [right=50pt of 0]{$1,\dots,1,1,0,1,1,\dots,1$};
\node[point] (p) [right=25pt of 0]{};
\node[place] (2) [above=20pt of p] {$0,\dots,0,k_{i-1},0,0,0\dots,0$};
\draw [-latex,shorten >=2pt, shorten <=2pt,thick] (0) -- (1)
 node[very near start,above] {{\small $k_i$}};
\draw [-,thick] (p) -- (2);
\end{tikzpicture}
\end{center}
\caption{Diagram associated to $k_ix_i^{k_i-1}x_{i+1}+x_{i-1}^{k_{i-1}}+s\prod_{j\neq i}x_j$}
\end{figure}
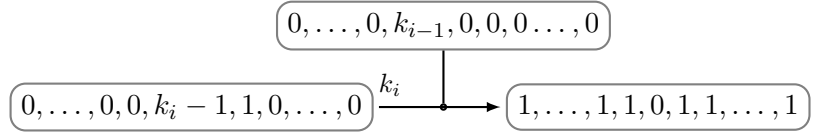
We put the arrow in the direction where the difference is given by $(1,\dots,1,-k_i+1,1,\dots,1)$ or $(1,\dots,1,-k_i+1,0,1,\dots,1)$ to indicate the sum between the monomial with the biggest $x_i$ exponent and the monomial with coefficient $s$. \\
In the same way as before we can describe what happens if we multiply such a derivative consisting of three monomials by another monomial $m=c\prod x_i^{a_i}$. The result of multiplying $k_ix_i^{k_i-1}+x_{i-1}^{k_{i-1}}+s\prod_{j\neq i}x_j$ by $m$ would be
\begin{figure}[H]
\begin{center}
\begin{tikzpicture}
\node[place] (0) {$a_1,\dots,a_{i-1},a_i+k_i-1,a_{i+1},\dots,a_n$};
\node[place] (1) [right=46pt of 0]{$1+a_1,\dots,1+a_{i-1},a_i,1+a_{i+1},\dots,1+a_n$};
\node[point] (p) [right=23pt of 0]{};
\node[place] (2) [above=20pt of p] {$a_1,\dots,a_{i-2},a_{i-1}+k_{i-1},a_i,\dots,a_n$};
\draw [-latex,shorten >=2pt, shorten <=2pt,thick] (0) -- (1)
 node[very near start,above] {{\small $k_i$}}
 node[midway, below] {$c$};
\draw [-,thick] (p) -- (2);
\end{tikzpicture}
\end{center}
\caption{Multiplication of $k_ix_i^{k_i-1}+x_{i-1}^{k_{i-1}}+s\prod_{j\neq i}x_j$ by a monomial $m=c\prod x_i^{a_i}$}
\end{figure}
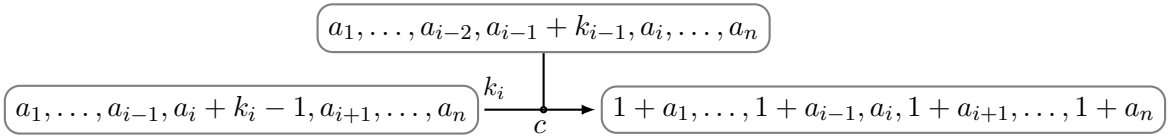
In the same way as before we can state the conditions that the sum of three monomials is given by multiplying a monomial to a partial derivative. In the case of $k_ix_i^{k_i-1}+x_{i-1}^{k_{i-1}}+s\prod_{j\neq i}x_j$ the pairwise differences between the three monomials has to be $(1,\dots,1,-k_i+1,1,\dots,1)$, $(1,\dots,1,-k_{i-1}+1,0,1,\dots,1)$ and $(0,\dots,0,k_{i-1},-k_i+1,0,\dots,0)$. In the case of $k_ix_i^{k_i-1}x_{i+1}+x_{i-1}^{k_{i-1}}+s\prod_{j\neq i}x_j$ the differences of the summands have to be $(1,\dots,1,-k_i+1,0,1,\dots,1)$,$(1,\dots,1,-k_{i-1}+1,0,1,\dots,1)$ and $(0,\dots,0,k_{i-1},-k_i+1,-1,\dots,0)$. Again we ignore the fact that the coefficients have to match.
\begin{rem}
We will sometimes say that an arrow of one of the two above types which has three adjacent vertices creates an extra vertex. This is meant in the sense that if we want to connect two vertices with distance $(1,\dots,1,-k_i+1,1,\dots,1)$ or $(1,\dots,1,-k_i+1,0,1,\dots,1)$, then an extra vertex has to be created in order to get all the differences correct.
\end{rem}
\begin{rem}\label{degree-partial-derivative}
One should notice that the vertices adjacent to the same arrow all have the same weighted degree, because all summands of the partial derivative have the same weighted degree. This means in particular that if you know the degree of one vertex, you know the degree of all others.
\end{rem}

\subsection*{The role of chains and loops}
\begin{rem}\label{chain-loop}
The property of a chain and a loop is also represented in the partial derivatives. \\
In a loop $x_1^{k_1}x_2+\dots+x_m^{k_m}x_1$ all partial derivatives are of the form $\frac{\partial f}{\partial x_i}=k_ix_i^{k_i-1}x_{i+1}+x_{i-1}^{k_{i-1}}+s\prod_{j\neq i}x_j$. Notice that the difference between the exponents of $x_{i-1}^{k_{i-1}}$ and $s\prod_{j\neq i}x_j$ is $(1,\dots,1,-k_{i-1}+1,0,1,\dots,1)$ which is exactly the difference of the arrow which belongs to the partial derivative $\frac{\partial f}{\partial x_{i-1}}=k_{i-1}x_{i-1}^{k_{i-1}-1}x_i+x_{i-2}^{k_{i-2}}+s\prod_{j\neq i-1}x_j$ with respect to $x_{i-1}$. We have to be a little bit careful that all exponents involved are positive, which means that $\frac{\partial f}{\partial x_i}$ has to be at least multiplied by $x_i$. In our notation this means that if the numbers are big enough, i.e. all entries of the vertex at the arrow tip are $\geq 1$, the partial derivatives of a polynomial of loop type form a loop. We show here the smallest possible example. In general this works if all entries are at least as big as shown here:
\begin{figure}[H]
\begin{center}
\begin{tikzpicture}
\node[place] (1) {$1,\dots,1$};
\node[place] (0) [left=70pt of 1]{$k_1,1,0,\dots,0$};
\node[place] (2) [above left=70pt of 1,xshift=40pt] {$1,0,\dots,0,k_m$};
\draw [-latex,shorten >=2pt, shorten <=2pt,thick] (0) -- (1)
 node[point] (p1) [pos=0.5,label=70:{\small $\frac{\partial f}{\partial x_1}$}] {}
 node[very near start,above] {{\small $k_1$}};
\draw[-,thick,shorten >=2pt] (p1)--(2);
\draw [-latex,shorten >=2pt, shorten <=2pt,thick] (2) -- (1)
 node[point] (p2) [pos=0.5] {}
 node[midway,above right,yshift=5pt] {{\small $\frac{\partial f}{\partial x_m}$}}
 node[very near start,right] {{\small $k_m$}};
\node[place] (3) [above right=70pt of 1,xshift=-40pt] {$0,\dots,0,k_{m-1},1$};
\draw[-,thick,shorten >=2pt] (p2)--(3);
\draw [-latex,shorten >=2pt, shorten <=2pt,thick] (3) -- (1)
 node[very near start,right] {{\small $k_{m-1}$}}
 node[point] (p3) [pos=0.5] {}
 node[midway,left,xshift=-2pt] {{\small $\frac{\partial f}{\partial x_{m-1}}$}};
\node[place] (4) [below left=70pt of 1,xshift=40pt] {$0,k_2,1,0,\dots,0$};
\draw [-latex,shorten >=2pt, shorten <=2pt,thick] (4) -- (1)
 node[point] (p4) [pos=0.5] {}
 node[midway,above,yshift=2pt] {{\small $\frac{\partial f}{\partial x_2}$}}
 node[very near start,left] {{\small $k_2$}};
\draw[-,thick,shorten >=2pt] (p4)--(0);
\node[place] (5) [right=70pt of 1] {$0,\dots,0,k_{m-2},1,0$};
\draw [-latex,shorten >=2pt, shorten <=2pt,thick] (5) -- (1)
 node[point] (p5) [pos=0.5] {}
 node[very near start,below] {{\small $k_{m-2}$}};
\draw[-,thick,shorten >=2pt] (p3)--(5);
\draw[-,out=-90,in=-25,loosely dotted,thick] (p5) to (p4);
\end{tikzpicture}
\end{center}
\caption{Partial derivatives of a loop}
\label{loop}
\end{figure}
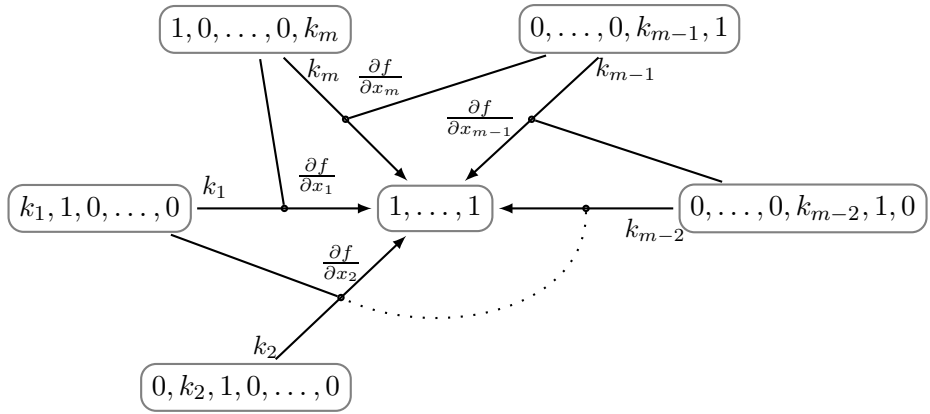
In the case of a chain $x_1^{k_1}x_2+\dots+x_{m-1}^{k_{m-1}}x_m+x_m^{k_m}$ of length $\geq 3$, there are three types of partial derivatives involved. The partial derivatives of the variables in the middle of the chain are of the form $k_ix_i^{k_i-1}x_{i+1}+x_{i-1}^{k_{i-1}}+s\prod_{j\neq i}x_j$ and the beginning and end are of the form $k_1x_1^{k_1-1}x_{2}+s\prod_{j\neq 1}x_j$ and $k_mx_m^{k_m-1}+x_{m-1}^{k_{m-1}}+s\prod_{j\neq m}x_j$ respectively. In the same way as for polynomials of loop type these partial derivatives match in our notation to give a chain. Again we show an example with the smallest non-negative entries.
\begin{figure}[H]
\begin{center}
\begin{tikzpicture}
\node[place] (1) {$1,\dots,1$};
\node[place] (0) [left=70pt of 1]{$0,\dots,0,k_m$};
\node[place] (2) [above left=70pt of 1,xshift=40pt] {$0,\dots,0,k_{m-1},1$};
\draw [-latex,shorten >=2pt, shorten <=2pt,thick] (0) -- (1)
 node[point] (p1) [pos=0.5,label=below:{\small $\frac{\partial f}{\partial x_m}$}] {}
 node[very near start,above] {{\small $k_m$}};
\draw[-,thick,shorten >=2pt] (p1)--(2);
\draw [-latex,shorten >=2pt, shorten <=2pt,thick] (2) -- (1)
 node[point] (p2) [pos=0.5,label=-90:{\small $\frac{\partial f}{\partial x_{m-1}}$}] {}
 node[very near start,right] {{\small $k_{m-1}$}};
\node[place] (3) [above right=70pt of 1,xshift=-40pt] {$0,\dots,0,k_{m-2},1,0$};
\draw[-,thick,shorten >=2pt] (p2)--(3);
\draw [-latex,shorten >=2pt, shorten <=2pt,thick] (3) -- (1)
 node[point] (p3) [pos=0.5,label=right:{\small $\frac{\partial f}{\partial x_{m-2}}$}] {}
 node[very near start,right] {{\small $k_{m-2}$}};
\node[place] (4) [below=60pt of 1] {$k_1,1,0,\dots,0$};
\draw [-latex,shorten >=2pt, shorten <=2pt,thick] (4) -- (1)
 node[circle] (p4) [pos=0.5,label=left:{\small $\frac{\partial f}{\partial x_1}$}] {}
 node[very near start,left] {{\small $k_1$}};
\draw[-,out=300,in=0,loosely dotted,thick] (p3) to (p4);
\end{tikzpicture}
\end{center}
\caption{Partial derivatives of a chain of length $m>2$}
\label{chain}
\end{figure}
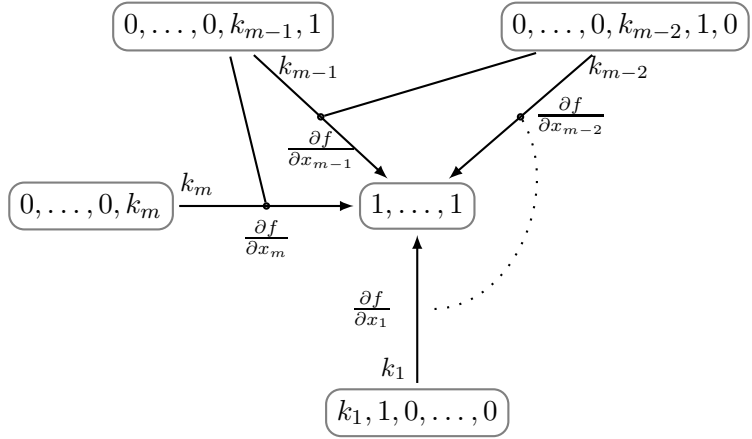
\end{rem}
Of course this also works for chains of length $2$, where only two partial derivatives are involved. One is of the type $k_{i-1}x_{i-1}^{k_{i-1}-1}x_{i}+s\prod_{j\neq i-1}x_j$ and the other $k_ix_i^{k_i-1}+x_{i-1}^{k_{i-1}}+s\prod_{j\neq i}x_j$. How this matches is also shown in the following picture:
\begin{figure}[H]
\begin{center}
\begin{tikzpicture}
\node[place] (0) {$0,\dots,0,0,k_i,0,\dots,0$};
\node[place] (1) [right=100pt of 0]{$1,\dots,1$};
\node[point] (p) [right=50pt of 0]{};
\node[place] (2) [above=60pt of p] {$0,\dots,0,k_{i-1},1,0\dots,0$};
\draw [-latex,shorten >=2pt, shorten <=2pt,thick] (0) -- (1)
 node[above left,midway] {$\frac{\partial f}{\partial x_i}$}
 node[very near start,above] {{\small $k_i$}};
\draw [-,thick] (p) -- (2);
\draw [-latex,shorten >=2pt, shorten <=2pt,thick] (2) -- (1)
 node[right,midway] {$\frac{\partial f}{\partial x_{i-1}}$}
 node[very near start,right] {{\small $k_{i-1}$}};
\end{tikzpicture}
\end{center}
\caption{Partial derivatives of a chain of length $2$}
\end{figure}
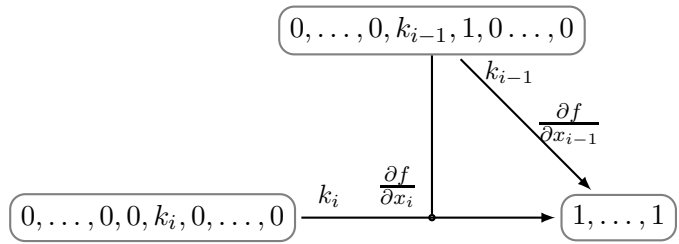

Later one of the important parts is to know when an arrow is generated by a partial derivative. We have seen before that for this to happen the difference between the monomials has to be appropriate. To shorten the notation we define the following.
\begin{notation}\label{partial_i}
$\partial_i$ is an abbreviation of the partial derivative in the new notation. In detail $\partial_i$ is a short notation for the arrow connecting all vertices of the partial derivative of $f$ with respect to $x_i$.
\end{notation}

\subsection*{Writing a monomial with the generators of the Jacobian ideal}
Now let us see how our new way of writing the partial derivatives fits into the Griffiths-Dwork method:
We concentrate on the second part of the Algorithm \ref{alg} and repeat what there is to do. We assume we have a basis of the Milnor ring in the appropriate degrees. From Remark \ref{delta-omega}, we know that all $\delta^i\omega$ are linear combinations of $\frac{m!s^{m+1}(\prod x_j)^m\Omega_0}{f^{m+1}}$ with $m\leq i$. So all we have to do is write every $m!s^{m+1}(\prod x_j)^m$ in the basis of the Milnor ring using the Griffiths formula. We should notice that for $m<n-1$ the monomial $m!s^{m+1}(\prod x_j)^m$ is not in the Jacobian ideal. To make things easier we can without loss of generality assume that they are basis elements. But for $m\geq n-1$ the monomial $m!s^{m+1}(\prod x_j)^m$ is in the Jacobian ideal and can be written as a linear combination of the partial derivatives. So for $m\geq n-1$ one can find polynomials $p^m_i(\x)$ for $i=1,\dots,n$ such that $m!s^{m+1}(\prod x_j)^m=\sum_{i=1}^n p^m_i(\x)\frac{\partial f}{\partial x_i}$.\\
Let us see what this means in our notation. First we represent the monomial $m!s^{m+1}(\prod x_j)^m$ by $(m,\dots,m)$. On the other side we have a sum of $p^m_i(\x)\frac{\partial f}{\partial x_i}$. Here we get an arrow (maybe with an additional vertex) for every monomial in $p^m_i(\x)$, but these arrows are not completely independent in the following sense: If we expand $p^m_i(\x)\frac{\partial f}{\partial x_i}$ than every monomial apart from $(\prod x_j)^m$ has to appear at least twice, because all monomials apart from $(\prod x_j)^m$ have to disappear after adding up everything. In the new notation this means that there are always at least two arrows meeting at a vertex. Putting this information for every vertex together, we can say that because there are only finitely many monomials involved, all arrows that represent the sum $\sum_{i=1}^n p^m_i(\x)\frac{\partial f}{\partial x_i}$ form not necessarily oriented cycles with the exception of a line meeting one of the cycles at the one end and the point $(m,\dots,m)$, which corresponds to the monomial $(\prod x_j)^m$, at the other end. Notice that all cycles are connected otherwise they can be omitted.\\
Summarizing the above we know that the monomial $(\prod x_j)^m$ written as a linear combination of partial derivatives must consist of connected cycles and maybe an additional line from one of the cycles ending at $(m,\dots,m)$. However, a representation of $(\prod x_j)^m$ by the partial derivatives is not given a priori and it is not necessarily unique. Therefore the goal is to find such a linear combination of partial derivatives. Is it possible to arrange the arrows so that we end up with a linear combination of partial derivatives giving $(\prod x_j)^m$? The answer is yes, because otherwise the monomial is not in the Jacobian ideal. A first solution on how this arrangement of arrows looks like will be given in the next chapter in Proposition \ref{shortestway}. After that we develop an explicit method of finding such a representation with the arrows of partial derivatives.

\subsection*{Using the Griffiths formula}
Now let us return to Algorithm \ref{alg} and assume we have found a way to write $(m,\dots,m)$ with the arrows representing the partial derivatives. The Griffiths formula (\ref{griffiths-formula}) tells us that we can reduce the monomial $(\prod x_j)^m$ to a sum of monomials of degree $d(m-1)$ by differentiating the coefficient polynomials in the representation by the partial derivatives. In other words, if we can write $m!s^{m+1}(\prod x_j)^m=m\sum_{i=1}^n p^m_i(\x)\frac{\partial f}{\partial x_i}$, then in the primitive cohomology the monomial $m!s^{m+1}(\prod x_j)^m$ can be identified with the polynomial $\sum_{i=0}^n \frac{\partial}{\partial x_i}p^m_i(\x)$. Now we want to translate this behaviour to our notation.\\
For every arrow belonging to the partial derivative we have to extract the information what the coefficient monomial is, i.e. the monomial the partial derivative was multiplied with to give this arrow, and then take the appropriate partial derivative. This means that using the Griffiths formula every arrow gets contracted to a point in the following way:
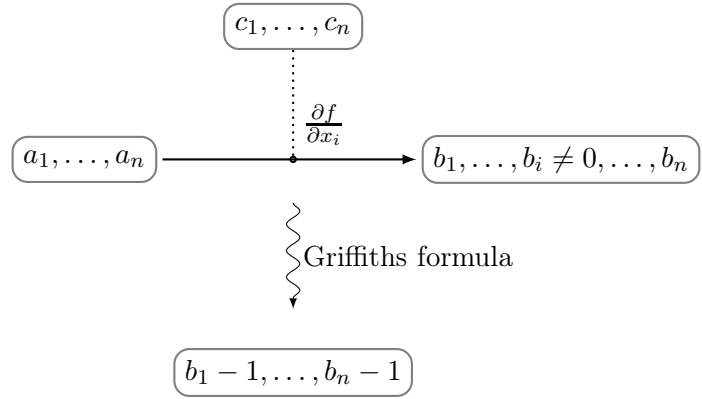
\begin{figure}[H]
\begin{center}
\begin{tikzpicture}
\node[place] (0) {$a_1,\dots,a_n$};
\node[place] (1) [right=100pt of 0]{$b_1,\dots,b_i\neq 0,\dots,b_n$};
\node[point] (p) [right=50pt of 0]{};
\node[place] (2) [above=40pt of p] {$c_1,\dots,c_n$};
\draw [-latex,shorten >=2pt, shorten <=2pt,thick] (0) -- (1);
\draw [-,thick,dotted] (p) -- (2)
 node [above right,at start] {$\frac{\partial f}{\partial x_i}$};
\node[place] (g) [below=70pt of p]{$b_1-1,\dots,b_n-1$};
\draw [-latex reversed,shorten >=18pt,decorate,decoration={snake,pre=moveto,pre length=15pt,post=moveto,post length=15pt}] (p) -- (g) 
 node [right,midway] {Griffiths formula};
\end{tikzpicture}
\end{center}
\caption{Griffiths formula}
\end{figure}
This states that no matter which of the 4 types of partial derivatives we have, as long as $b_i\neq 0$ (all other $b_j\neq 0$ anyway) the Griffiths formula maps the whole arrow to the point at the arrow tip subtracted by $(1,\dots,1)$. Let us see why this is correct: The point at the arrow tip represents the exponents of the coefficient monomial times $\prod_{j\neq i} x_i$ therefore the coefficient monomial is given by $x_i^{b_i}\prod_{j\neq i} x_j^{b_j-1}$. If we differentiate this monomial with respect to $x_i$ and $b_i\neq 0$ we end up with the monomial $\prod_{j=1}^n x_j^{b_j-1}$. If the entry $b_i=0$ then the conclusion that the coefficient monomial is $\prod_{j\neq i} x_j^{b_j-1}$ is still true, but if we differentiate this with respect to $x_i$ it simply vanishes.

\subsection*{Summary of the diagrammatic Griffiths-Dwork method}
We want to do a final summary of how the Griffiths-Dwork method and in particular Algorithm \ref{alg} work in our diagrammatic interpretation. We will skip the step of choosing a basis. This issue will be addressed later in Chapter \ref{haupt} and also we will not deal with finding the linear relation. But these two steps are not the hard part of Algorithm \ref{alg}. The most difficult part is to use the Griffiths formula \ref{griffiths-formula}. So we want to write the derivatives $\delta^i\omega$ in the basis for a given $i$. We have seen in Remark \ref{delta-omega} that they consist only of the monomials $(\prod x_j)^m$ for $m\leq i$. So we can restrict ourselves to write $(\prod x_j)^m$ in a basis. For $m\geq n-2$ these monomials are in the Jacobian ideal $(\prod x_j)^m\in J(f)$ and for $m\leq n-1$ we can assume they are basis elements, because they are definitely linear independent. So we concentrate on writing $(\prod x_j)^m$ for $m\geq n-2$ in the basis. First thing we have to do to achieve this goal is writing $(\prod x_j)^m$ in terms of partial derivatives or, in our new diagrammatic notation, find a way starting at $(m,\dots,m)$ and using the arrows belonging to the partial derivatives in such a way that there are always at least two adjacent arrows at each vertex. As mentioned before this is not a sufficient condition, so if we have found such an arrangement of arrows we have to check that it really works and the way does not end up to be trivial in the way that the coefficients we ignored are trivial. If we have found a valid arrangement giving $(\prod x_j)^m$ in terms of partial derivatives we can use the Griffiths formula, which diagrammatically involves replacing every arrow by a single point, namely the vertex at the arrow tip subtracted by $(1,\dots,1)$, if the appropriate entry is $\geq 1$, or if this entry is $0$ the arrow just vanishes when we use the Griffiths formula. After using the Griffiths formula once we end up with vertices corresponding to monomials which are either in the basis or can be written as something in the basis plus something in the Jacobian ideal. We then have to repeat the same procedure for everything in the Jacobian ideal until we end up with monomials in the basis. This is the idea, however, in practice it becomes slightly more complicated. We will come back to this in Chapter \ref{haupt} in the cases important to us. In particular we will see that we can restrict ourselves to understanding the whole procedure for the monomial $(\prod x_j)^m$ for just one $m$ and from this we can deduce what happens in all the other cases.

\newpage
\thispagestyle{empty}
\chapter{Calculations for the Picard-Fuchs equation with the Griffiths-Dwork method}\label{haupt}

In this chapter we will analyse the Picard-Fuchs equation for our one-parameter family in a lot more detail. We will simulate all steps of the Griffiths-Dwork method in our new notation. The goal of this chapter is to calculate the order of the Picard-Fuchs equation. The proof of the order of the Picard-Fuchs equation will also play a role in Chapter \ref{vier}, where we will use this result together with the calculation of the \GKZ System to show exactly what the Picard-Fuchs equation looks like.\\
We will use the same notation as before, but we want to recall it again here and use it throughout this chapter without further notice.

\begin{notation}\label{notation}
Let $g(\x)=g(x_1,x_2,\dots,x_n):=\sum_{i=1}^{n}\prod_{i=1}^n x_i^{E_{ij}}$ be an invertible polynomial with reduced weights $q_1,q_2,\dots,q_n$ and $\deg g=d$ for which the Calabi-Yau condition, $d=\sum_{i=1}^n q_i$, holds. The diagonal entries of the exponent matrix $E=(E_{ij})_{i,j}$ are defined as $k_1,\dots,k_n$. We denote by $g^t(\x)$ the transposed polynomial of $g$, the dual reduced weights belonging to $g^t$ are denoted by $\widehat{q}_1,\widehat{q}_2,\dots,\widehat{q}_n$ and the degree by $\deg g^t=\widehat{d}$.\\
The invertible polynomial consists of loops and chains of arbitrary length. For a variable $x_i$ we always take $x_{i-1}$ and $x_{i+1}$ to be the neighbouring variables in the loop or chain. The indices are without further notice taken modulo the length of the loop or chain.\\
We always denote by $f(x_1,\dots,x_n)$ the one-parameter family with parameter $s$ defined by $f(\x)=f(x_1,\dots,x_n):=g(x_1,\dots,x_n)+s\prod_{i=1}^n x_i$.
\end{notation}

We will prove in Chapter \ref{vier} that the Picard-Fuchs equation has a very special form, which is only dependent on the dual weights and the dual degree of the invertible polynomial. This special form can be seen in the following theorem:

\pagebreak

\textbf{Theorem \ref{conj}.}
\textit{The Picard-Fuchs equation for the one-parameter family $f(x_1,\dots,x_n)=g(x_1,\dots,x_n)+s\prod x_i$ is given by}
\begin{align*}
0=\prod_{i=1}^n \widehat{q}_i^{\widehat{q}_i}s^{\widehat{d}}\prod_{i=1}^{n}\prod_{j=0}^{\widehat{q}_i-1}(\delta+\frac{j\cdot \widehat{d}}{\widehat{q}_i})\prod_{\ell\in I}(\delta+\ell)^{-1}-(-\widehat{d})^{\widehat{d}}\prod_{j=0}^{\widehat{d}-1}(\delta-j)\prod_{\ell\in I}(\delta-\ell)^{-1},
\end{align*}
\textit{where $I=\{0,\dots,\widehat{d}-1\}\cap\bigcup_{i=1}^n \left\{0,\frac{\widehat{d}}{\widehat{q}_i},\frac{2\widehat{d}}{\widehat{q}_i},\dots,\frac{(\widehat{q}_i-1)\widehat{d}}{\widehat{q}_i}\right\}$.
}

So in this chapter we will prove that the order of the Picard-Fuchs equation in the theorem above is correct. 

\section{Combinatorial ideas for the order of the Picard-Fuchs equation}

Before we state the main theorem (Theorem \ref{thm}) of this chapter which presents what the order of the Picard-Fuchs equation is, we will have a closer look at the Griffiths-Dwork method. In Section \ref{gd-comb} we already discussed how one can see the Griffiths-Dwork method in an diagrammatic way. But as promised we will be more concrete in this Chapter. The first thing we want to make concrete is how to write $(\prod x_i)^{n-1}$ with the generators of the Jacobian ideal. We will see in the proof of Theorem \ref{thm} that this information is everything one needs to determine $(\prod x_i)^m$ for arbitrary $m$. We already know from Section \ref{gd-comb} that we need a path using the arrows representing the partial derivatives and we know from Remark \ref{degree-partial-derivative} that all vertices, or correspondingly monomials, on this path have the same degree. So in this case all have degree $d(n-1)$. In the first lemma we will study all monomials of this degree.

\begin{lemma}\label{3d}
Suppose we have a monomial of weighted degree $d(n-1)$, the weights satisfy the Calabi-Yau condition and the monomial is not $\prod_{i=1}^n x_i^{k_i-1}$. Then there is an $i\in\{1,\dots,n\}$ such that $x_i$ has an exponent $\geq k_i$.
\end{lemma}

\begin{proof}
Assume $m(\x)$ is a monomial, where for all $i$ the exponent of $x_i$ is $\leq k_i-1$, then the weighted degree of this monomial satisfies
\begin{align*}
\deg m\leq \sum_{i=1}^n q_i(k_i-1)=\sum_{i=1}^n q_i k_i -d\leq nd-d=d(n-1).
\end{align*}
This means the degree of $m(\x)$ is smaller than $d(n-1)$ except if $q_ik_i=d$ for all $i$ (polynomial of Brieskorn-Pham type) and $m(\x)=\prod x_i^{k_i-1}$.
\end{proof}

Using Lemma \ref{3d} we can give a quite concrete construction of the path using the arrows representing the partial derivatives. 

\begin{prop}\label{shortestway}
The shortest non-trivial way from $(n-1,\dots,n-1)$ to itself in any diagram that can be constructed using the arrows corresponding to the partial derivatives has length $\widehat{d}$. Furthermore there is a shortest way that has at most one vertex with a zero entry. We call such a way {\em Jacobi path}.
\end{prop}

\begin{proof}
The idea is to first restrict ourselves to a main path, which will be a cycle starting and ending at $(n-1,\dots,n-1)$ and then we will take care of the extra vertices we created. So first of all we forget about the extra vertices and treat every partial derivative as if it would consist of one arrow and two vertices.
If we use every step $\partial_i$ (cf. Notation \ref{partial_i}) of the Jacobian ideal exactly $\widehat{q}_i$ times, then in total we end up at the starting point. This can be seen in the following calculation. We will look at the entries separately and show that after adding all the steps all entries are zero. In the case when $\widehat{q}_ik_i=\widehat{d}$, the $i$th entry of adding up all steps is
\begin{align*}
-\widehat{q}_i(k_i-1)+\sum_{j\neq i}\widehat{q}_j=-\widehat{q}_ik_i+\sum_{j=1}^n\widehat{q}_j=-\widehat{d}+\widehat{d}=0.
\end{align*}
Notice that if $\widehat{q}_ik_i=\widehat{d}$ this means that in our polynomial $g(\x)$ the variable $x_i$ is only appearing as the term $x_i^{k_i}$ or $x_i^{k_i}x_{i-1}$.\\
Now let $\widehat{q}_ik_i\neq\widehat{d}$. This means that in $g(\x)$ the variable $x_i$ is either appearing as the term $x_i^{k_i}+x_{i-1}^{k_{i-1}}x_i$ or $x_i^{k_i}x_{i+1}+x_{i-1}^{k_{i-1}}x_i$. So $x_i$ is the end of the chain or in the middle of a chain or in a loop. In these cases the $i$th entry of adding up all steps is
\begin{align*}
-\widehat{q}_i(k_i-1)+\sum_{j\neq i,i-1}\widehat{q}_j=-\widehat{q}_ik_i+\sum_{j\neq i-1}\widehat{q}_j=-\widehat{d}+\sum_{j=1}^n\widehat{q}_j=0,
\end{align*}
because $\partial_{i-1}$ will have $0$ as $i$th entry.\\
According to Lemma \ref{3d} we know that we can arrange the arrows in a way that at each vertex all entries are $>0$. The lemma states that we can always use some arrow, because at least for one $x_i$ the exponent is bigger than $k_i$, i.e. one entry in the vertex is bigger than $k_i$, and we can exclude that we have to use $\partial_i$ more than $\widehat{q}_i$ times, because this would produce a negative entry somewhere which is not possible.\\
Now we have to take care of the extra vertices but this is very easy because according to Lemma \ref{3d} all entries are $\geq 1$ in the case where an extra vertex appears. So we can use Remark \ref{chain-loop}, which shows that we can use the rest of the chain or the full loop here. This means that all arrows of the rest of the chain or loop fit here with the arrow tip pointing at the same vertex as can be seen in Figure \ref{loop} and \ref{chain}. After doing this, every vertex has at least two adjacent arrows.\\
There is obviously no non-trivial shorter way, because this means that there is a linear relation between the rows of $E$ and therefore the weights would not be reduced.\\
The last thing we have to exclude is that this path is trivial. This is not possible, because we produced a path where every arrow on the main path points in the same direction. At every vertex except $(n-1,\dots,n-1)$ the coefficients are chosen in a way that they add up to zero, but the arrow tip always carries the $s$ as a coefficient so in order to add up to zero the coefficient of the next arrow has to have a higher $s$ exponent. This means that at the point $(n-1,\dots,n-1)$ there is one arrow with a coefficient that contains $s^0$ and one arrow with a coefficient containing $s^{\widehat{d}}$. So they can never add up to zero, but we can easily normalise all coefficients to get the coefficient $1$ at this vertex.
\end{proof}

\begin{rem}
In regular notation the output of Proposition \ref{shortestway} is that the monomial $\prod_{i=1}^n x_i^{n-1}$ can be written exactly as a sum $\sum_{i=1}^n \left(p_i(\x)+h_i(\x)\right)\frac{\partial f}{\partial x_i}$, where $p_i(\x)$ is a sum of exactly $\widehat{q}_i$ monomials and with only one exception all of these monomials include all variables with a positive exponent and $h_i(\x)$ includes all monomials that come from the arrows of the chains and loops to adopt the extra vertices we created on the Jacobi path.
\end{rem}

\begin{rem}
We want to draw some attention on the fact that the dual weights and the dual degree come into play in Proposition \ref{shortestway}. The reason for that was already established in Remark \ref{auftreten-duale-gewichte}. There we found out that the differences in the exponent vectors are given by the matrix $E^t-\left(\begin{smallmatrix}1&\cdots& 1\\\vdots&~&\vdots\\1&\cdots&1 \end{smallmatrix}\right)$. Also in the equations in Remark \ref{auftreten-duale-gewichte} we can see that the relations between these columns are given by the dual weights:
\begin{align*}
\left(E^t-\begin{pmatrix}1&\cdots& 1\\\vdots&~&\vdots\\1&\cdots&1 \end{pmatrix}\right)\begin{pmatrix}\widehat{q}_1\\\vdots\\\widehat{q}_n\end{pmatrix}=\begin{pmatrix}0\\\vdots\\0\end{pmatrix}
\end{align*}
This explains immediately why the dual weights give the relation between the steps done by the partial derivatives.
\end{rem}

Now we know exactly how many times we need each arrow $\partial_i$, corresponding to a partial derivative, to have the shortest possible way of writing the monomial $\prod_{i=1}^n x_i$ with the partial derivatives. However, we still don't know in which order to use them. To achieve this, we will define a few sets that will tell us exactly where to use each partial derivative (see in Lemma \ref{position}) and will also be used in the formulation of the main theorem, telling us the order of the Picard-Fuchs equation.

\begin{defi}\label{vundu}
\begin{align*}
D:&=\{1,2,\dots,\widehat{d}\} & &\\
Q_i:&=\{\frac{\widehat{d}}{\widehat{q}_i},2\frac{\widehat{d}}{\widehat{q}_i},\dots,(\widehat{q}_i-1)\frac{\widehat{d}}{\widehat{q}_i},\widehat{d}\} & &\\
Q_i^{\Z}:&= Q_i\cap \Z & Q_i^{\Q}:&=Q_i\cap(\Q\setminus\Z)\\
V:&=D\setminus (\bigcup_{i=1}^nQ_i^{\Z}) & v:&=|V|\\
u:&=\sum_{i=1}^n|Q_i|-|\bigcup_{i=1}^nQ_i^{\Z}| &&
\end{align*}
\end{defi}

\begin{rem}\label{v=u+3}
Since $|D|=\widehat{d}$ and $|Q_i|=\widehat{q}_i$, we have
\begin{align*}
v&=|V|=\widehat{d}-|\bigcup_{i=1}^nQ_i^{\Z}|\\
&=\sum_{i=1}^n(\widehat{q}_i)-|\bigcup_{i=1}^nQ_i^{\Z}|=\sum_{i=1}^n|Q_i|-|\bigcup_{i=1}^nQ_i^{\Z}|\\
&=u
\end{align*}
and $u=v\geq \varphi(\widehat{d})$, where $\varphi$ is Euler's phi function. In addition we have $u\geq n-1$.
\end{rem}

Now we have everything to state the main theorem of the chapter:

\begin{thm}\label{thm}
Let $g(\x)$ be an invertible polynomial and $f(\x)=g(\x)+s\prod_i x_i$ a one-parameter family with parameter $s$. Then the Picard-Fuchs equation of $f(\x)$ has order $u$.
\end{thm}

Before we prove this theorem we will work out some details about the polynomial $g(\x)$ and the one-parameter family $f(\x)$. We will need this information to prove Theorem \ref{thm}. Especially we need to know in more detail how our Jacobi path looks like. We know which steps have to be done, but we do not know in which order they are used. But this is important to keep track of which monomials vanish when we use the Griffiths formula.
As before we will concentrate on the monomial $(\prod_{i=1}^n x_i)^{n-1}$ and try to figure out everything in this case first.

We already know that the Jacobi path has length $\widehat{d}$. So there are $\widehat{d}$ positions on our path that have to be filled. We want to figure out now at which position a partial derivative produces a vertex where at least one entry is $1$. This is important, because we know with only one exception, that all entries on our Jacobi path are $\geq 1$. So if a partial derivative produces a vertex where at least one entry is $1$ then this is the earliest position where this partial derivative will be used. The following lemma tells us in detail where these earliest positions are.

\begin{lemma}\label{position}
The smallest position a partial derivative $\partial_i$ can be used is where it produces an entry $1$ in the vertex at the arrow tip. To state the smallest positions we distinguish between two cases. The smallest positions for the partial derivative $\partial_i$ are
\begin{enumerate}[label=(\roman*)]
\item $q-n+2$ for $q\in Q_i$ and $\widehat{q}_i k_i=\widehat{d}$ or
\item $\lfloor q \rfloor-n+2$ for $q\in Q_i^{\Q}$ and $q-n+1$ for $q\in Q_i^{\Z}$, if $\widehat{q}_i k_i\neq \widehat{d}$.
\end{enumerate}
\end{lemma}

\begin{rem}
The numbers $q-n+2$, $q-n+1$ and $\lfloor q \rfloor-n+2$ can be $\leq 0$. If this happens, this obviously means that we cannot use this partial derivative at a position where we produce an entry $1$. We have to move these partial derivatives at least to position $1$. But this will become clear later.
\end{rem}

\begin{proof}
Let us assume we are in case (i), so $\widehat{q}_i k_i=\widehat{d}$. This means that in $g(\x)$ the variable $x_i$ appears only as $x_i^{k_i}$ or $x_i^{k_i}x_{i+1}$. It follows now that all $\partial_l$ for $l\neq i$ have $1$ as $i$th entry. So if we assume that all other positions are taken, then the first time we can use $\partial_i$ is when the $i$th entry is $k_i=\frac{\widehat{d}}{\widehat{q}_i}$, but because we started with the monomial $(n-1,\dots,n-1)$ and at every step $\neq \partial_i$ we add $1$ in the $i$th entry this happens after $k_i-(n-2)$ steps. Now the next time we can use $\partial_i$ is after $k_i$ steps of adding $1$ to the $i$th entry. So in total we can use $\partial_i$ at the positions $q-n+2$ for $q\in Q_i$. This proves case (i).\\
For case (ii) we assume $\widehat{q}_i k_i+\widehat{q}_{i-1}=\widehat{d}$. The numbers in $Q_i$ are evenly spread between $0$ and $\widehat{d}$. These numbers minus $n-2$ nearly give the smallest positions of $\partial_i$. But we have to investigate this a little bit more to see what is happening. Similar to the other case it follows that the only terms that consist of $x_i$ in $g(\x)$ are $x_i^{k_i}+x_ix_{i-1}^{k_{i-1}}$ or $x_i^{k_i}x_{i+1}+x_ix_{i-1}^{k_{i-1}}$. But this time not all of the other partial derivatives add $1$ to the $i$th entry. The partial derivative $\partial_{i-1}$ adds $0$ to the $i$th entry. We can use $\partial_i$ for the first time if the $i$th entry is $k_i$, so we have to use $k_i-n+2=\frac{\widehat{d}-\widehat{q}_{i-1}}{\widehat{q}_i}-n+2$ of the partial derivatives with respect to $x_j$ where $j\neq i,i-1$. In addition we have to count how often $\partial_{i-1}$ got used before we use $\partial_i$. Since for both partial derivatives the numbers in $Q_i$ and $Q_{i-1}$ are evenly spread, the relation between $\widehat{q}_{i}$ and $\widehat{q}_{i-1}$ tells us exactly the relative position of both numbers on the Jacobi path. The term $\frac{\widehat{q}_{i-1}}{\widehat{q}_i}$ tells us exactly how often $\partial_{i-1}$ is used before we used $\partial_i$ for the first time. If this is not a natural number we have to round down and get $\lfloor q \rfloor-n+2$ as position for $\partial_i$ and as before the other positions are the multiples of these. If $\frac{\widehat{q}_{i-1}}{\widehat{q}_i}$ is a natural number, then this means that $\partial_{i-1}$ can be used at the same position. But because this does not contribute anything in the entry $i$, we can use $\partial_i$ one position earlier. This proves part (ii).
\end{proof}

\begin{rem}\label{letzte-positionen}
Notice that all partial derivatives $\partial_i$ with $\widehat{q}_i k_i=\widehat{d}$ are at position $\widehat{d}-n+2$ and all partial derivatives $\partial_i$ with $\widehat{q}_i k_i\neq\widehat{d}$ are at position $\widehat{d}-n+1$, which agrees with the fact that after using every partial derivative $\partial_i$ exactly $\widehat{q}-1$ times we always end with the monomial $\prod_{\widehat{q}_ik_i\neq\widehat{d}}x_i^{k_i-1}\prod_{\widehat{q}_ik_i=\widehat{d}}x_i^{k_i}$.
\end{rem}

We want to draw a picture illustrating the Jacobi path and indicating where to use the partial derivatives. We want to explain this using an example.

\begin{ex}\label{exIV}
Consider the one-parameter family $f(x_1,x_2,x_3,x_4)=x_1^{18}+x_2^2x_3+x_3^3x_4+x_4^3+sx_1x_2x_3x_4$ with weights $(q_1,q_2,q_3,q_4)=(1,7,4,6)$ and weighted degree $\deg f=d=18$. The dual weights and the dual weighted degree are given by $(\widehat{q}_1,\widehat{q}_2,\widehat{q}_3,\widehat{q}_4)=(1,9,3,5)$ and $\widehat{d}=18$.\\
The sets needed for calculating the path in the Jacobian ideal are given as follows:
\begin{align}\label{qi-ex}
Q_1=&\{18\},& Q_2=&\{2,4,6,8,10,12,14,16,18\},\\ \notag
Q_3=&\{6,12,18\},& Q_4=&\{\frac{18}{5},\frac{36}{5},\frac{54}{5},\frac{72}{5},18\}
\end{align}
We know that the Jacobi path has length $18=\widehat{d}$, so we will make a table where every position can be entered. The table stops at position $16=\widehat{d}-n+2$ because that is the biggest number that can occur as a smallest position.
\begin{figure}[H]
\begin{center}
\begin{tikzpicture}
\node[circle] (0) {\textcolor{cTangoAluminum3}{$0$}};
\node[circle] (1) [right=3pt of 0] {$1$};
\node[circle] (2) [right=3pt of 1] {$2$};
\node[circle] (3) [right=3pt of 2] {$3$};
\node[circle] (4) [right=3pt of 3] {$4$};
\node[circle] (5) [right=3pt of 4] {$5$};
\node[circle] (6) [right=3pt of 5] {$6$};
\node[circle] (7) [right=3pt of 6] {$7$};
\node[circle] (8) [right=3pt of 7] {$8$};
\node[circle] (9) [right=3pt of 8] {$9$};
\node[circle] (10) [right=3pt of 9] {$10$};
\node[circle] (11) [right=3pt of 10] {$11$};
\node[circle] (12) [right=3pt of 11] {$12$};
\node[circle] (13) [right=3pt of 12] {$13$};
\node[circle] (14) [right=3pt of 13] {$14$};
\node[circle] (15) [right=3pt of 14] {$15$};
\node[circle] (16) [right=3pt of 15] {$16$};
\node[circle] (hilf1) [below left=12pt of 0] {}; 
\node[circle] (hilf2) [below right=12pt of 16] {};
\draw[-,thick] (hilf1)--(hilf2);
\node[circle] (j1) [below=20pt of 0] {$\partial_2$};
\node[circle] (j2) [below=20pt of 2] {$\partial_2$};
\node[circle] (j3) [below=20pt of 4] {$\partial_2$};
\node[circle] (j4) [below=20pt of 6] {$\partial_2$};
\node[circle] (j5) [below=20pt of 8] {$\partial_2$};
\node[circle] (j6) [below=20pt of 10] {$\partial_2$};
\node[circle] (j7) [below=20pt of 12] {$\partial_2$};
\node[circle] (j18) [below=20pt of 14] {$\partial_2$};
\node[circle] (j10) [below=20pt of 3] {$\partial_3$};
\node[circle] (j11) [below=20pt of 9] {$\partial_3$};
\node[circle] (j12) [below=20pt of 1] {$\partial_4$};
\node[circle] (j13) [below=20pt of 5] {$\partial_4$};
\node[circle] (j14) [below=10pt of j5] {$\partial_4$};
\node[circle] (j15) [below=10pt of j7] {$\partial_4$};
\node[circle] (j9) [below=20pt of 16] {$\partial_1$};
\node[circle] (j16) [below=20pt of 15] {$\partial_3$};
\node[circle] (j8) [below=10pt of j16] {$\partial_4$};
\node[circle] (j17) [below=10pt of j9] {$\partial_2$};
\end{tikzpicture}
\end{center}
\caption{Smallest positions for the partial derivatives}
\end{figure}
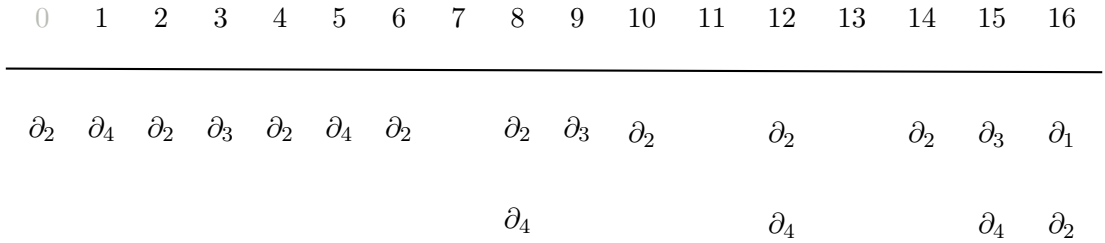
\end{ex}

The idea is that because we know where the smallest position is where we can use $\partial_i$, we get a Jacobi path by just shifting the positions until we have one partial derivative at every position. Proposition \ref{shortestway} tells us that this is always possible. Basically we should only shift the partial derivatives to the right, because shifting to the left is only possible by one position and this might disconnect the path. The only position where we allow to shift to the left will be from position $\widehat{d}-n+2$ to position $\widehat{d}-n+1$. The reason for this becomes clear later. We will have a look what to do in the example and indicate the shifts by arrows in the table.

\begin{ex}\emph{(Continuation Example \ref{exIV})}
The shifting we have to do can be seen in the following picture:
\begin{figure}[H]
\begin{center}
\begin{tikzpicture}
\node[circle] (0) {\textcolor{cTangoAluminum3}{$0$}};
\node[circle] (1) [right=20pt of 0] {$1$};
\node[circle] (2) [right=20pt of 1] {$2$};
\node[circle] (3) [right=20pt of 2] {$3$};
\node[circle] (4) [right=20pt of 3] {$4$};
\node[circle] (5) [right=20pt of 4] {$5$};
\node[circle] (6) [right=20pt of 5] {$6$};
\node[circle] (7) [right=20pt of 6] {$7$};
\node[circle] (8) [right=20pt of 7] {$8$};
\node[circle] (9) [right=20pt of 8] {$9$};
\node[circle] (10) [below=80pt of 0] {$10$};
\node[circle] (11) [right=20pt of 10] {$11$};
\node[circle] (12) [right=20pt of 11] {$12$};
\node[circle] (13) [right=20pt of 12] {$13$};
\node[circle] (14) [right=20pt of 13] {$14$};
\node[circle] (15) [right=20pt of 14] {$15$};
\node[circle] (16) [right=20pt of 15] {$16$};
\node[circle] (17) [right=20pt of 16] {$17$};
\node[circle] (18) [right=20pt of 17] {$18$};
\node[circle] (hilf1) [below left=12pt of 0] {}; 
\node[circle] (hilf2) [below right=12pt of 9] {};
\draw[-,thick] (hilf1)--(hilf2);
\node[circle] (hilf3) [below right=12pt of 18] {}; 
\node[circle] (hilf4) [below left=12pt of 10] {};
\draw[-,thick] (hilf3)--(hilf4);
\node[circle] (j1) [below=20pt of 0] {$\partial_2$};
\node[circle] (j2) [below=20pt of 2] {$\partial_2$};
\node[circle] (j3) [below=20pt of 4] {$\partial_2$};
\node[circle] (j4) [below=20pt of 6] {$\partial_2$};
\node[circle] (j5) [below=20pt of 8] {$\partial_2$};
\node[circle] (j6) [below=20pt of 10] {$\partial_2$};
\node[circle] (j7) [below=20pt of 12] {$\partial_2$};
\node[circle] (j18) [below=20pt of 14] {$\partial_2$};
\node[circle] (j10) [below=20pt of 3] {$\partial_3$};
\node[circle] (j11) [below=20pt of 9] {$\partial_3$};
\node[circle] (j12) [below=20pt of 1] {$\partial_4$};
\node[circle] (j13) [below=20pt of 5] {$\partial_4$};
\node[circle] (j14) [below=10pt of j5] {$\partial_4$};
\node[circle] (j15) [below=10pt of j7] {$\partial_4$};
\node[circle] (j9) [below=20pt of 16] {$\partial_1$};
\node[circle] (j16) [below=20pt of 15] {$\partial_3$};
\node[circle] (j8) [below=10pt of j16] {$\partial_4$};
\node[circle] (j17) [below=10pt of j9] {$\partial_2$};
\node[circle] (h1) [right=20pt of j4] {};
\node[circle] (h2) [right=20pt of j6] {};
\node[circle] (h3) [right=20pt of j7] {};
\node[circle] (h4) [right=20pt of j11] {};
\node[circle] (h5) [left=20pt of j6] {};
\node[circle] (h6) [right=20pt of j9] {};
\node[circle] (h7) [right=20pt of h6] {};
\draw[-latex,color=cTangoPlum3] (j1)--(j12);
\draw[-latex,color=cTangoPlum3] (j12)--(j2);
\draw[-latex,color=cTangoPlum3] (j2)--(j10);
\draw[-latex,color=cTangoPlum3] (j10) to [out=-45,in=-135](j13);
\draw[-latex,color=cTangoPlum3] (j13)--(j4);
\draw[-latex,color=cTangoPlum3] (j4)--(h1);
\draw[-latex,color=cTangoPlum3] (j14)--(j11);
\draw[-latex,color=cTangoPlum3] (j11)--(h4);
\draw[-latex,color=cTangoPlum3] (h5)to [out=-45,in=-135](h2);
\draw[-latex,color=cTangoPlum3] (j15)--(h3);
\draw[-latex,color=cTangoPlum3] (j16)--(j9);
\draw[-latex,color=cTangoPlum3] (j8)to [in=-100,out=-45](h6);
\draw[-latex,color=cTangoPlum3] (j9)--(h7);
\draw[-latex,color=cTangoPlum3] (j17)--(j16);
\end{tikzpicture}
\end{center}
\caption{Shifting of positions on the Jacobi path}
\label{ex-shift}
\end{figure}
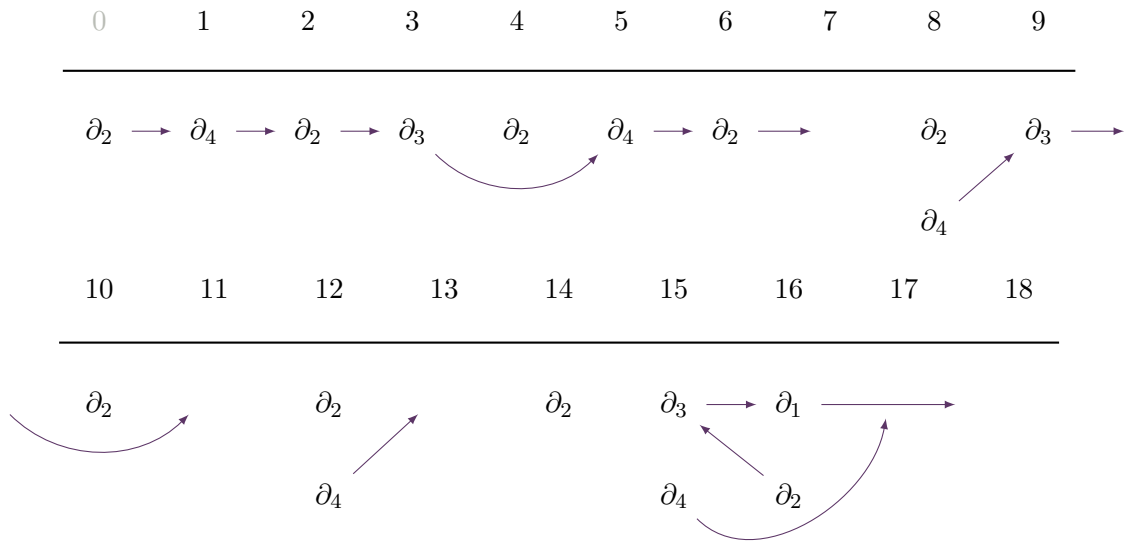
\end{ex}

At this point it is not entirely clear why we choose to shift exactly like this. But the important point here is that we have exactly one partial derivative at each position and apart from the partial derivative $\partial_2$ at position $16$ we shifted all arrows to the right.\\
The above picture tells us which partial derivative we have to use at every position. So we can use the above picture to write down the Jacobi path in the diagrammatic notation introduced in Section \ref{gd-comb}. The only thing we have to take care of in addition is to complete the loops and chains so that the coefficients can be chosen such that after adding up everything else but $(\prod_{i=1}^n x_i)^{n-1}$ vanishes. From Remark \ref{chain-loop} we know that this is nearly always possible without any difficulties. We will go back to the Example now to see how this works.

\begin{ex}\emph{(Continuation Example \ref{exIV})}\label{jpexIV}
We will write down how the Jacobi path looks like explicitly. We start with the monomial $(3,3,3,3)$ and use the partial derivatives as shown in Figure \ref{ex-shift}, which gives the order
\begin{align*}
\partial_2,\partial_4,\partial_2,\partial_2,\partial_3,\partial_4,\partial_2,\partial_2,\partial_4,\partial_2,\partial_3,\partial_2,\partial_4,\partial_2,\partial_2,\partial_3,\partial_4,\partial_1.
\end{align*}
This leads to the following picture, where we wrote down every monomial on the Jacobi path. Notice that we neglected all the coefficients here.
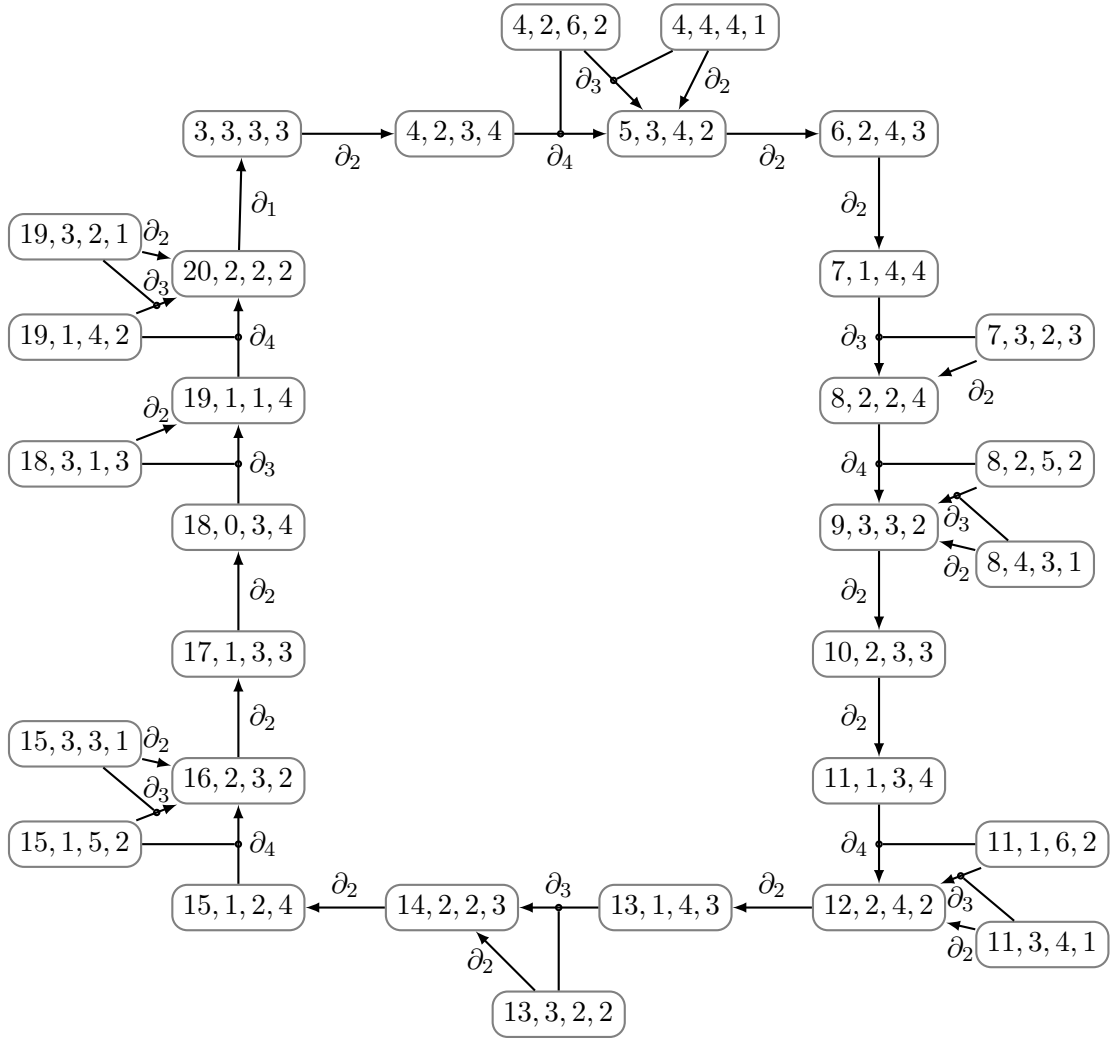
\begin{figure}[H]
\begin{center}
\begin{tikzpicture}
\node[place] (0) {$3,3,3,3$};
\node[place] (1) [right=35pt of 0] {$4,2,3,4$}; 
\node[place] (2) [right=35pt of 1] {$5,3,4,2$}; 
\node[place] (3) [right=35pt of 2] {$6,2,4,3$}; 
\node[place] (4) [below=35pt of 3] {$7,1,4,4$}; 
\node[place] (5) [below=30pt of 4] {$8,2,2,4$}; 
\node[place] (6) [below=30pt of 5] {$9,3,3,2$}; 
\node[place] (7) [below=30pt of 6] {$10,2,3,3$}; 
\node[place] (8) [below=30pt of 7] {$11,1,3,4$}; 
\node[place] (9) [below=30pt of 8] {$12,2,4,2$}; 
\node[place] (10) [left=30pt of 9] {$13,1,4,3$}; 
\node[place] (11) [left=30pt of 10] {$14,2,2,3$}; 
\node[place] (12) [left=30pt of 11] {$15,1,2,4$}; 
\node[place] (13) [above=30pt of 12] {$16,2,3,2$}; 
\node[place] (14) [above=30pt of 13] {$17,1,3,3$}; 
\node[place] (15) [above=30pt of 14] {$18,0,3,4$}; 
\node[place] (16) [above=30pt of 15] {$19,1,1,4$}; 
\node[place] (17) [above=30pt of 16] {$20,2,2,2$}; 
\draw[-latex,thick] (0)--(1)
 node[midway,below] {$\partial_2$};
\draw[-latex,thick] (1)--(2)
 node[midway,below] {$\partial_4$}
 node[point,midway] (h1) {};
\draw[-latex,thick] (2)--(3)
 node[midway,below] {$\partial_2$};
\draw[-latex,thick] (4)--(5)
 node[midway,left] {$\partial_3$}
 node[point,midway] (h2) {};
\draw[-latex,thick] (3)--(4)
 node[midway,left] {$\partial_2$};
\draw[-latex,thick] (5)--(6)
 node[midway,left] {$\partial_4$}
 node[point,midway] (h3) {};
\draw[-latex,thick] (6)--(7)
 node[midway,left] {$\partial_2$};
\draw[-latex,thick] (7)--(8)
 node[midway,left] {$\partial_2$};
\draw[-latex,thick] (8)--(9)
 node[midway,left] {$\partial_4$}
 node[point,midway] (h4) {};
\draw[-latex,thick] (10)--(11)
 node[midway,above] {$\partial_3$}
 node[point,midway] (h5) {};
\draw[-latex,thick] (9)--(10)
 node[midway,above] {$\partial_2$};
\draw[-latex,thick] (11)--(12)
 node[midway,above] {$\partial_2$};
\draw[-latex,thick] (12)--(13)
 node[midway,right] {$\partial_4$}
 node[point,midway] (h6) {};
\draw[-latex,thick] (13)--(14)
 node[midway,right] {$\partial_2$};
\draw[-latex,thick] (15)--(16)
 node[midway,right] {$\partial_3$}
 node[point,midway] (h7) {};
\draw[-latex,thick] (16)--(17)
 node[midway,right] {$\partial_4$}
 node[point,midway] (h8) {};
\draw[-latex,thick] (14)--(15)
 node[midway,right] {$\partial_2$};
\draw[-latex,thick] (17)--(0)
 node[midway,right] {$\partial_1$};
\node[place] (p1) [above=30pt of h1] {$4,2,6,2$};
\draw[-,thick] (h1)--(p1);
\draw[-latex,thick] (p1)--(2)
 node[point,midway] (h11) {}
 node[midway,left] {$\partial_3$};
\node[place] (p11) [right=15pt of p1] {$4,4,4,1$};
\draw[-,thick] (h11)--(p11);
\draw[-latex,thick] (p11)--(2)
 node[midway,right] {$\partial_2$};
\node[place] (p2) [right=35pt of h2] {$7,3,2,3$};
\draw[-,thick] (h2)--(p2);
\draw[-latex,thick] (p2)--(5)
 node[midway,below right] {$\partial_2$};
\node[place] (p3) [right=35pt of h3] {$8,2,5,2$};
\draw[-,thick] (h3)--(p3);
\draw[-latex,thick] (p3)--(6)
 node[point,midway] (h33) {}
 node[midway,below] {$\partial_3$};
\node[place] (p33) [below=20pt of p3] {$8,4,3,1$};
\draw[-,thick] (h33)--(p33);
\draw[-latex,thick] (p33)--(6)
 node[midway,below] {$\partial_2$};
\node[place] (p4) [right=35pt of h4] {$11,1,6,2$};
\draw[-,thick] (h4)--(p4);
\draw[-latex,thick] (p4)--(9)
 node[point,midway] (h44) {}
 node[midway,below] {$\partial_3$};
\node[place] (p44) [below=20pt of p4] {$11,3,4,1$};
\draw[-,thick] (h44)--(p44);
\draw[-latex,thick] (p44)--(9)
 node[midway,below] {$\partial_2$};
\node[place] (p5) [below=30pt of h5] {$13,3,2,2$};
\draw[-,thick] (h5)--(p5);
\draw[-latex,thick] (p5)--(11)
 node[midway,left] {$\partial_2$};
\node[place] (p6) [left=35pt of h6] {$15,1,5,2$};
\draw[-,thick] (h6)--(p6);
\draw[-latex,thick] (p6)--(13)
 node[point,midway] (h66) {}
 node[midway,above] {$\partial_3$};
\node[place] (p66) [above=20pt of p6] {$15,3,3,1$};
\draw[-,thick] (h66)--(p66);
\draw[-latex,thick] (p66)--(13)
 node[midway,above] {$\partial_2$};
\node[place] (p7) [left=35pt of h7] {$18,3,1,3$};
\draw[-,thick] (h7)--(p7);
\draw[-latex,thick] (p7)--(16)
 node[midway,above] {$\partial_2$};
\node[place] (p8) [left=35pt of h8] {$19,1,4,2$};
\draw[-,thick] (h8)--(p8);
\draw[-latex,thick] (p8)--(17)
 node[point,midway] (h88) {}
 node[midway,above] {$\partial_3$};
\node[place] (p88) [above=20pt of p8] {$19,3,2,1$};
\draw[-,thick] (h88)--(p88);
\draw[-latex,thick] (p88)--(17)
 node[midway,above] {$\partial_2$};
\end{tikzpicture}
\end{center}
\caption{The Jacobi path for Example \ref{exIV}}
\label{jpexIVpath}
\end{figure}
\end{ex}

At this point we can build a Jacobi path for every invertible polynomial. So we will go one step further. We will make use of the Griffiths formula. From Section \ref{gd-comb} we know that using the Griffiths formula means contracting every arrow to the point at the arrow tip minus $(1,\dots,1)$ or if the point has $0$ as an entry the arrow vanishes completely. But from Proposition \ref{shortestway} we know that the Jacobi path starting at $(n-1,\dots,n-1)$ has at most one vertex with $0$ as an entry. This means that after using the Griffiths formula at most one vertex will vanish. The vertices that are still there after the use of the Griffiths formula have the same differences as before. So we can basically put the arrows in again, but we have to be a little bit careful. If the partial derivative belongs to a loop or a chain, then all entries belonging to another variable of the loop or to the rest of the chain have to be $>0$. If an arrow fits in between two vertices together with the rest of the chain, or the loop as can be seen in Figure \ref{loop} and \ref{chain}, then we can adjust all coefficients. This means we only need one basis element for every part of a path. But we have to be careful that our shifting did not disconnect two vertices which are not linearly dependent. Therefore we want to investigate what a good and what a bad way of shifting is. This means we want to find out how to shift the partial derivatives such that if a path gets disconnected there is no other connection between two vertices. The last $n$ positions play a special role here and we will take care of them in the end. We distinguish between the following cases:

\begin{enumerate}[label=(\roman*)]
\item Two partial derivatives $\partial_{i_1}$ and $\partial_{i_2}$ are at the same position $p$, where
\begin{enumerate}[label=(\alph*)]
\item $\partial_{i_1}$ and $\partial_{i_2}$ are not neighbouring elements in a chain or a loop or
\item $\partial_{i_1}$ and $\partial_{i_2}$ are neighbouring elements in a chain or a loop.
\end{enumerate}
\item Two partial derivatives $\partial_{i_1}$ and $\partial_{i_2}$ are at two succeeding positions $p$ and $p+1$ and the first one $\partial_{i_1}$ gets shifted, where
\begin{enumerate}[label=(\alph*)]
\item $\partial_{i_1}$ and $\partial_{i_2}$ are not neighbouring elements in a chain or a loop or
\item $\partial_{i_1}$ and $\partial_{i_2}$ are neighbouring elements in a chain or a loop.
\end{enumerate}
\end{enumerate}

Step by step we will show how to shift in all these cases. Before we take care of all special cases, we will state a lemma that shows us that for some partial derivatives it is always possible to shift them.

\begin{lemma}\label{dreieck-shift}
Let $M\subseteq\{1,\dots,n\}$ be the set of all indices with $\widehat{q}_m k_m\neq\widehat{d}$ for $m\in M$. Then, without one exception, every monomial $\prod_{i=1}^nx_i^{\alpha_i}$, with $\alpha_m=k_m$ for $m\in M$ and $\alpha_i<k_i$ for $i\in \{1,\dots,n\}\setminus M$, has degree $<d(n-1)$. The only exception is the monomial $\prod_{m\in M}x_m^{k_m}\prod_{i\notin M}x_i^{k_i-1}$.
\end{lemma}

\begin{proof}
Assume the statement is false. This means that there is a monomial $\prod_{i=1}^nx_i^{\alpha_i}$ with $\alpha_m=k_m$ for $m\in M$ and $\alpha_i<k_i$ for $i\in \{1,\dots,n\}\setminus M$ that has weighted degree $d(n-1)$. Notice that $\widehat{d}\neq\widehat{q}_m k_m$ means that either $x_m^{k_m}+x_{m-1}^{k_{m-1}}x_m$ or $x_m^{k_m}x_{m+1}+x_{m-1}^{k_{m-1}}x_m$ is in the polynomial $g(\x)$, where the indices are taken modulo the length of the appropriate chain or loop. It follows that $q_m+q_{m-1}k_{m-1}=d$. If we calculate the degree of the monomial $\prod_{i=1}^nx_i^{\alpha_i}$, we get:
\begin{align*}
\deg\left( \prod_{i=1}^nx_i^{\alpha_i}\right)&=\sum_{i=1}^n q_i\alpha_i=\sum_{m\in M}q_mk_m+\sum_{i\notin M} q_i\alpha_i\\
&\leq \sum_{m\in M}q_mk_m+\sum_{i\notin M} q_i(k_i-1)=-d+\sum_{m\in M}q_m+\sum_{i=1}^n q_ik_i\\
&=-d+\sum_{m\in M}(\underbrace{q_m+q_{m-1}k_{m-1}}_{=d})+\sum_{i;i+1\notin M} \underbrace{q_ik_i}_{=d}\\
&=-d+d|M|+d(n-|M|)= d(n-1).
\end{align*}
It follows that the degree of the monomial $\prod_{i=1}^nx_i^{\alpha_i}$ with $\alpha_m=k_m$ for $m\in M$ and $\alpha_i<k_i$ for $i\in \{1,\dots,n\}\setminus M$ is $<d(n-1)$ unless $\alpha_i=k_i-1$ for all $i\notin M$. 
%In addition we have that the degree of a monomial $\prod_{i=1}^nx_i^{\alpha_i}$ is always $<d(n-1)$ if $\alpha_i\leq k_i$ for all $i$ and the set of the indices with $\alpha_i=k_i$ is smaller than $M$.
\end{proof}

We want to relate Lemma \ref{dreieck-shift} to what we know. The lemma states in particular that if $\partial_i$ creates an extra vertex, then there is, with one exception, no monomial where $x_i$ has exponent $k_i$ and for all other $x_j$ it is smaller than $k_j$. This means that on the Jacobi path there is no position, apart from $\widehat{d}-n+1$ (cf. Remark \ref{letzte-positionen}), where we can use $\partial_i$ exclusively. In other words: There is always the possibility to shift $\partial_i$ if it produces an extra vertex.\\
To make the notation a little bit easier, we state two extra definitions.

\begin{defi}
For a fixed Jacobi path, we denote by 
\begin{align*}
\kappa(p):=\min_{1\leq i\leq n}\{a_i\,|\,(a_1,\dots,a_n)\text{ is the }p\text{th vertex on the Jacobi path}\}.
\end{align*}
So $\kappa(p)$ is the smallest entry of the 
vertex at position p.\\
The second number we define is $\partial(p)$. For every position on a fixed Jacobi path $\partial(p):=\partial_i$, if $\partial_i$ is the arrow connecting the vertices at position $p$ and $p+1$ and $\partial(p):=0$, if there is no arrow connecting the vertices at position $p$ and $p+1$.
\end{defi}

Now we want to investigate how to shift in case (i). So we have two partial derivatives $\partial_{i_1}$ and $\partial_{i_2}$ that are possible at the same position $p$. If the variables $x_{i_1}$ and $x_{i_2}$ are not neighbours in a loop or a chain, then they are independent of each other. This is subcase (a). Assume we shifted $\partial_{i_2}$, so we assume $\partial(p)=\partial_{i_1}$ and $\partial(p+1)=\partial_{i_2}$. This means that $\kappa(p+1)=1$ and $\kappa(p+2)=2$. If we used the Griffiths formula once, we have $\kappa(p+2)=1$ and therefore we still have $\partial(p+1)=\partial_{i_2}$. After using the Griffiths formula twice the vertex at position $p+1$ will vanish, because after the first use of the Griffiths formula, we had $\kappa(p+1)=0$. For the arrow at position $p$ it can make a difference what partial derivative we use here. If $\kappa(p)=\kappa(p+1)$, then after the first use of the Griffiths formula there is an arrow between the two vertices if and only if the partial derivative belongs to a chain of length $1$ or is the beginning of a chain. So if only one of the two partial derivatives belongs to the middle or end of a chain or to a loop than this should be shifted to position $p+1$ otherwise it does not matter. Now we get to subcase (b), which means that the two partial derivatives $\partial_{i_1}$ and $\partial_{i_2}$ at position $p$ belong to neighbouring variables in a loop or chain. This is the case that needs most work. First of all we will prove that if this is the case, then the rest of the loop or the beginning of the chain is also at this position, this will be done in the following lemma. After that we will state directly what the best way of shifting for a chain or a loop is.

\begin{lemma}\label{ganzzahlig-loop-chain}
Let $x_1^{k_1}x_2+\dots+x_m^{k_m}x_1$ be a loop of length $m$ in $g(\x)$.
\begin{enumerate}[label=(\roman*)]
\item $q\in Q_i^{\Z}$ for some $i\in\{1,\dots,m\}$ $\Rightarrow$ $q\in Q_i^{\Z}$ for all $i\in\{1,\dots,m\}$.
\item $q=\lfloor \tilde{q}_i\rfloor=\lfloor \tilde{q}_{i+1}\rfloor$ with $\tilde{q}_i \in Q_i,\tilde{q}_{i+1}\in Q_{i+1}$ $\Rightarrow$ $q\in Q_i^{\Z}$ for all $i\in\{1,\dots,m\}$.
\end{enumerate}
Let $x_1^{k_1}x_2+\dots+x_{m-1}^{k_{m-1}}x_m+x_m^{k_m}$ be a chain of length $m$ in $g(\x)$.
\begin{enumerate}[label=(\roman*),resume]
\item $q\in Q_i^{\Z}$ for $i\in\{1,\dots,m\}$ $\Rightarrow$ $q\in Q_j^{\Z}$ for all $j\in\{1,\dots,i\}$.
\item $q=\lfloor \tilde{q}_i\rfloor=\lfloor \tilde{q}_{i+1}\rfloor$ with $\tilde{q}_i \in Q_i,\tilde{q}_{i+1}\in Q_{i+1}$ $\Rightarrow$ $q\in Q_j^{\Z}$ for all $j\in\{1,\dots,i\}$.
\end{enumerate}
\end{lemma}

\begin{proof}
{\em (i):} If $x_1^{k_1}x_2+\dots+x_m^{k_m}x_1$ is in $g(\x)$ then this means that 
\begin{align*}
\widehat{q}_1 k_1+\widehat{q}_m=\widehat{d},\widehat{q}_m k_m+\widehat{q}_{m-1}=\widehat{d},\dots,\widehat{q}_2 k_2+\widehat{q}_{1}=\widehat{d}.
\end{align*}
If $q=\frac{c\widehat{d}}{\widehat{q}_i}\in Q_i^{\Z}$ then $\widehat{q}_i\mid c\widehat{d}$ and it follows immediately that $\widehat{q}_i\mid c\widehat{q}_{i-1}, \dots,c\widehat{q}_1,c\widehat{q}_m,\dots, c\widehat{q}_{i+1}, c\widehat{d}$. Define $\beta_j:=\frac{c\widehat{q}_j}{\widehat{q}_i}$ then we have $\beta_j\frac{\widehat{d}}{\widehat{q}_j}=\beta_j\frac{c\widehat{d}}{\beta_j\widehat{q}_i}=q$ and therefore $q\in Q_j^{\Z}$ for all $j\in\{1,\dots,m\}$.\\
{\em (ii):} Let $\tilde{q}_i=\frac{c_i\widehat{d}}{\widehat{q}_i}$ and $\tilde{q}_{i+1}=\frac{c_{i+1}\widehat{d}}{\widehat{q}_{i+1}}$, then $q=\frac{c_i\widehat{d}-a_i}{\widehat{q}_i}=\frac{c_{i+1}\widehat{d}-a_{i+1}}{\widehat{q}_{i+1}}$ with $a_i,a_{i+1}\in\Z$, $a_i<\widehat{q}_i$ and $a_{i+1}<\widehat{q}_{i+1}$. Because $\widehat{d}=\widehat{q}_{i+1}k_{i+1}+\widehat{q}_i$ the following calculation holds
\begin{align*}
\widehat{q}_{i+1}(c_i\widehat{d}-a_i)&=\widehat{q}_i(c_{i+1}\widehat{d}-a_{i+1})\\
\widehat{q}_{i+1}(c_i\widehat{d}-a_i)&=(\widehat{d}-\widehat{q}_{i+1}k_{i+1})(c_{i+1}\widehat{d}-a_{i+1})\\
\widehat{q}_{i+1}(c_i\widehat{d}-a_i)&=\widehat{d}(c_{i+1}\widehat{d}-a_{i+1})-\widehat{q}_{i+1}k_{i+1}(c_{i+1}\widehat{d}-a_{i+1})\\
c_i\widehat{d}-a_i&=\widehat{d}q-k_{i+1}(c_{i+1}\widehat{d}-a_{i+1}).
\end{align*}
It follows that $\widehat{d}\mid (a_i-a_{i+1}k_{i+1})<\widehat{q}_{i+1}k_{i+1}+\widehat{q}_i=\widehat{d}$. Therefore $a_i=a_{i+1}=0$ and with (i) the result follows.\\
{\em (iii):}
If $x_1^{k_1}x_2+\dots+x_{m-1}^{k_{m-1}}x_m+x_m^{k_m}$ is in $g(\x)$ then it follows that 
\begin{align*}
\widehat{q}_1 k_1=\widehat{d},\widehat{q}_2 k_2+\widehat{q}_1=\widehat{d},\dots,\widehat{q}_i k_i+\widehat{q}_{i-1}=\widehat{d}.
\end{align*} 
If $q=\frac{c\widehat{d}}{\widehat{q}_i}\in Q_i^{\Z}$ then
$\widehat{q}_i\mid c\widehat{d}$ and it follows that $\widehat{q}_i\mid c\widehat{q}_{i-1},\dots, c\widehat{q}_1$. Define $\beta_j:=\frac{c\widehat{q}_j}{\widehat{q}_i}$ then we have $\beta_j\frac{\widehat{d}}{\widehat{q}_j}=\beta_j\frac{c\widehat{d}}{\beta_j\widehat{q}_i}=q$ and therefore $q\in Q_j^{\Z}$ for all $j\in\{1,\dots,i\}$.\\
{\em (iv):} The proof is essentially the same as in {\em (ii)}. The only extra case is if $i=1$, but $q_1\in\Z$ and this just means $a_1=0$ from the beginning.
\end{proof}

We proved that if we are in case (ib), where two partial derivatives are at the same position and they belong to neighbouring variables in a chain or loop, then the loop is completely at this position or the chain until ending at one of the two variables or later is at this position. To be more precise, if two neighbouring variables $\partial_i$ and $\partial_{i-1}$ of a chain have the same number $q\in Q_i\cap Q_{i-1}$, then $q\in Q_j$ for $j\leq i$ as long as $x_j$ is part of the chain. According to Lemma \ref{position} this means that the beginning of the chain is at position $q-n+2$ and at least everything in the chain between the beginning and $\partial_i$ is at position $q-n+1$. We will show now in detail what to do if a loop or a chain of arbitrary length is at one position. In the first remark we will see how to shift a complete loop and what the linear dependencies between the monomials are. In Remark \ref{chain-gleicher-platz} we will do the same for a chain of arbitrary length. 

\begin{rem}\label{loop-gleicher-exponent}
Assume the loop $x_1^{k_1}x_2+\dots+x_m^{k_m}x_1$ is in $g(\x)$ and there is an element $q\in\bigcap_{i=1}^m Q_i$. This means all partial derivatives with respect to a variable in the loop have the same smallest position $q-n+1$. In the pictures below we want to see what happens if we use the partial derivatives in order, i.e. starting with $\partial_1$, then $\partial_2$ until in the end we use $\partial_m$. the polynomial $g(\x)$ might have variables $x_{m+1},\dots$ which are not in the loop, but we will omit all entries in the vertices that do not belong to the loop, i.e. all other variables in the monomials, because they will only increase by $1$ in every step and do not have any effect on the partial derivatives we use here. All partial derivatives have the same smallest position, so the starting monomial has to be $(k_1+c,\dots,k_m+c)$, where $c$ depends on how often the whole loop got shifted. We will start at $c=-1$ because this is the first time $0$ appears as entry, so if $c$ is bigger nothing interesting is happening until we have used the Griffiths formula several times. So we start with the vertex $(k_1-1,\dots,k_m-1)$ and use every partial derivative of the variables in the loop exactly once. The picture we get is now the following.

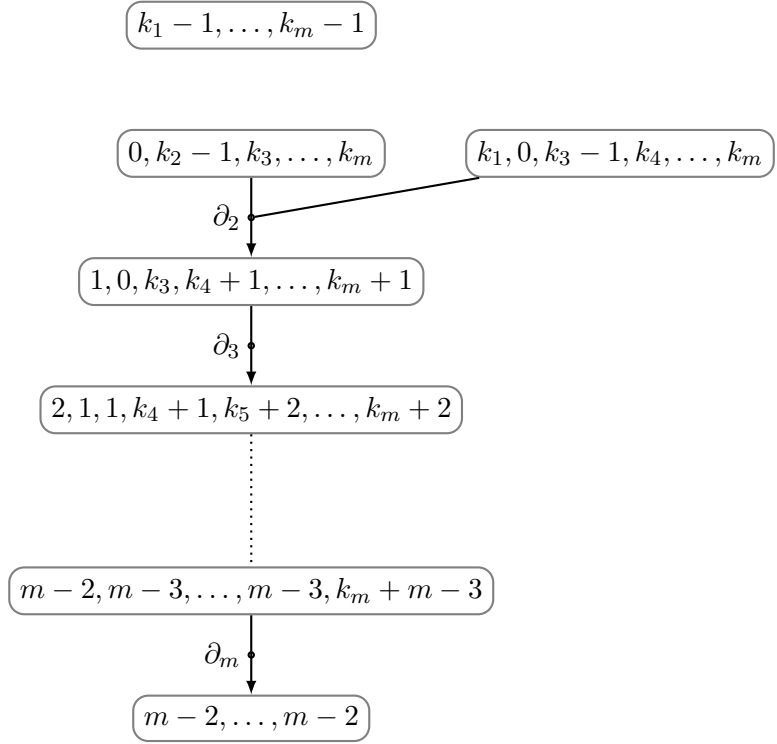
\begin{figure}[H]
\begin{center}
\begin{tikzpicture}
\node[place] (0) {$k_1-1,\dots,k_m-1$};
\node[place] (11) [below=30pt of 0] {$0,k_2-1,k_3,\dots,k_m$};
\node[place] (22) [below=30pt of 11] {$1,0,k_3,k_4+1,\dots,k_m+1$};
\node[place] (33) [below=30pt of 22] {$2,1,1,k_4+1,k_5+2,\dots,k_m+2$};
\node[place] (m-1m-1) [below=50pt of 33] {$m-2,m-3,\dots,m-3,k_m+m-3$};
\node[place] (m) [below=30pt of m-1m-1] {$m-2,\dots,m-2$};
\node[place] (12) [right=30pt of 11] {$k_1,0,k_3-1,k_4,\dots,k_m$};
\draw[-latex,thick] (11)--(22)
 node[midway,left] {$\partial_2$}
 node[point,midway] (p11) {};
\draw[-,thick] (p11)--(12);
\draw[-latex,thick] (22)--(33)
 node[midway,left] {$\partial_3$}
 node[point,midway] (p22) {};
\draw[dotted,thick] (33)--(m-1m-1);
\draw[-latex,thick] (m-1m-1)--(m)
 node[midway,left] {$\partial_m$}
 node[point,midway] (pm-1m-1) {};
\end{tikzpicture}
\end{center}
\caption{The case of a complete loop with the same smallest position}
\label{loop-kette}
\end{figure}

We want to point out that an arrow with a dot in the middle indicates that an extra vertex is created and all entries of the vertex at the arrow tip are $>0$. Therefore we can take care of the extra vertices with the rest of the partial derivatives as we have seen in Figure \ref{loop}. So there are a lot of extra vertices that we do not write down because this is the normal case as in Remark \ref{chain-loop}. In the above picture all arrows starting from the third vertex until the end point to a vertex with positive entries and therefore we marked the arrows with a dot in the middle. Also in this part of the picture the smallest number increases by $1$ in every step. So there is nothing to worry here. The part that needs more attention is between the first and the second vertex and the second and the third vertex. In both cases the vertex that the arrow would point to has a zero entry. So the situation looks the same in both cases. We can use the arrow $\partial_1$ between the first two vertices and we can use $\partial_2$ between the second and the third vertex. However, we are in both cases not able to use all other partial derivatives such that they point to the second or third vertex, which means that there would be a vertex with just one adjacent arrow. What we will see in the following part is that the problem can be fixed for the arrow between the second and the third vertex, but not for the arrow between the first two vertices. Therefore in Figure \ref{loop-kette} we draw the arrow $\partial_2$ but not the arrow $\partial_1$. In the next picture we see that there is a path which is on the one end connected with the extra vertex $(k_1,0,k_3-1,k_4,\dots,k_m)$ from Figure \ref{loop-kette} and on the other end with the vertex $(m-2,\dots,m-2)$, where the part of the path we are looking at ends. Now we draw the complete picture of the part of the path. After that we will show the interesting part of this picture again in more detail. The original path from Figure \ref{loop-kette} is shown in the first row of the following picture.

\begin{figure}[H]
\begin{center}
\begin{tikzpicture}
\node[place] (0) {$v_0$};
\node[place] (11) [right=30pt of 0] {$v_{1,1}$};
\node[place] (22) [right=30pt of 11] {$v_{2,2}$};
\node[place] (33) [right=30pt of 22] {$v_{3,3}$};
\node[place] (m-1m-1) [right=50pt of 33] {$v_{m-1,m-1}$};
\node[place] (m) [right=30pt of m-1m-1] {$v_{m}$};
\node[place] (12) [below=20pt of 11] {$\textcolor{cTangoAluminum4}{v_{1,2}}$};
\node[place] (13) [below=20pt of 12] {$\textcolor{cTangoAluminum4}{v_{1,3}}$};
\node[place] (1m-1) [below=40pt of 13] {$\textcolor{cTangoAluminum4}{v_{1,m-1}}$};
\node[place] (1m) [below=20pt of 1m-1] {$\textcolor{cTangoAluminum4}{v_{1,m}}$};
\node[place] (23) [below=20pt of 22] {$\textcolor{cTangoAluminum4}{v_{2,3}}$};
\node[place] (24) [below=20pt of 23] {$\textcolor{cTangoAluminum4}{v_{2,4}}$};
\node[place] (2m) [below=40pt of 24] {$\textcolor{cTangoAluminum4}{v_{2,m}}$};
\node[place] (21) [below=20pt of 2m] {$\textcolor{cTangoAluminum4}{v_{2,1}}$};
\node[place] (34) [below=20pt of 33] {$\textcolor{cTangoAluminum4}{v_{3,4}}$};
\node[place] (35) [below=20pt of 34] {$\textcolor{cTangoAluminum4}{v_{3,5}}$};
\node[place] (31) [below=40pt of 35] {$\textcolor{cTangoAluminum4}{v_{3,1}}$};
\node[place] (32) [below=20pt of 31] {$\textcolor{cTangoAluminum4}{v_{3,2}}$};
\node[place] (m-1m) [below=20pt of m-1m-1] {$v_{m-1,m}$};
\node[place] (m-11) [below=20pt of m-1m] {$v_{m-1,1}$};
\node[place] (m-1m-3) [below=40pt of m-11] {$v_{m-1,m-3}$};
\node[place] (m-1m-2) [below=20pt of m-1m-3] {$v_{m-1,m-2}$};
\draw[-latex,thick] (11)--(22)
 node[midway,above] {$\partial_2$}
 node[point,midway] (p11) {};
\draw[-,thick] (p11)--(12);
\draw[-latex,thick] (22)--(33)
 node[midway,above] {$\partial_3$}
 node[point,midway] (p22) {};
\draw[dotted,thick] (33)--(m-1m-1);
\draw[-latex,thick] (m-1m-1)--(m)
 node[midway,above] {$\partial_m$}
 node[point,midway] (pm-1m-1) {};
\draw[-,thick] (pm-1m-1) to [out=-90,in=25] (m-1m-2);
\draw[-latex,thick,color=cTangoAluminum4] (12)--(23)
 node[midway,above] {$\textcolor{cTangoAluminum4}{\partial_3}$}
 node[point,midway,color=cTangoAluminum4] (p12) {};
\draw[-,thick,color=cTangoAluminum4] (p12)--(13);
\draw[-latex,thick,color=cTangoAluminum4] (13)--(24)
 node[midway,above] {$\textcolor{cTangoAluminum4}{\partial_4}$}
 node[point,midway,color=cTangoAluminum4] (p13) {};
\draw[dotted,color=cTangoAluminum4] (13)--(1m-1);
\draw[dotted,color=cTangoAluminum4] (24)--(2m);
\draw[dotted,color=cTangoAluminum4] (35)--(31);
\draw[dotted] (m-11)--(m-1m-3);
\draw[dotted,thick,color=cTangoAluminum4] (34)--(m-1m);
\draw[dotted,thick,color=cTangoAluminum4] (35)--(m-11);
\draw[dotted,thick,color=cTangoAluminum4] (31)--(m-1m-3);
\draw[dotted,thick,color=cTangoAluminum4] (32)--(m-1m-2);
\draw[-latex,thick,color=cTangoAluminum4] (1m-1)--(2m)
 node[midway,above,xshift=-3pt] {$\textcolor{cTangoAluminum4}{\partial_m}$}
 node[point,midway,color=cTangoAluminum4] (p1m-1) {};
\draw[-,thick,color=cTangoAluminum4] (p1m-1)--(1m);
\draw[-latex,thick,color=cTangoAluminum4] (1m)--(21)
 node[midway,below] {$\textcolor{cTangoAluminum4}{\partial_1}$}
 node[point,midway,color=cTangoAluminum4] (p1m) {};
\draw[-,thick,color=cTangoAluminum4] (p1m) to [out=85,in=-50] (11);
\draw[-latex,thick,color=cTangoAluminum4] (23)--(34)
 node[midway,above] {$\textcolor{cTangoAluminum4}{\partial_4}$}
 node[point,midway,color=cTangoAluminum4] (p23) {};
\draw[-latex,thick,color=cTangoAluminum4] (24)--(35)
 node[midway,above] {$\textcolor{cTangoAluminum4}{\partial_5}$}
 node[point,midway,color=cTangoAluminum4] (p24) {};
\draw[-latex,thick,color=cTangoAluminum4] (2m)--(31)
 node[midway,above] {$\textcolor{cTangoAluminum4}{\partial_1}$}
 node[point,midway,color=cTangoAluminum4] (p2m) {};
\draw[-latex,thick,color=cTangoAluminum4] (21)--(32)
 node[midway,above] {$\textcolor{cTangoAluminum4}{\partial_2}$}
 node[point,midway,color=cTangoAluminum4] (p21) {};
\draw[-latex,thick] (m-1m)--(m)
 node[midway,above] {$\partial_1$}
 node[point,midway] (pm-1m) {};
\draw[-,thick] (pm-1m)--(m-1m-1);
\draw[-latex,thick] (m-11) to [out=30,in=-110] (m);
\node[point] (pm-11) [right=10pt of m-11,yshift=19pt,label=above:$\partial_2$] {};
\draw[-,thick] (pm-11)--(m-1m);
\draw[-latex,thick] (m-1m-3) to [out=30,in=-95] (m);
\node[point] (pm-1m-3) [right=10pt of m-1m-3,yshift=27pt,label=left:$\partial_{m-2}$] {};
\draw[-latex,thick] (m-1m-2) to [out=0,in=-90] (m);
\node[point] (pm-1m-2) [right=43pt of m-1m-2,yshift=48pt,label=right:$\partial_{m-1}$] {};
\draw[-,thick] (pm-1m-2)--(m-1m-2);
\end{tikzpicture}
\end{center}
\caption{The case of a complete loop at one position}
\label{kompletter-loop-bild}
\end{figure}

To show the structure of what is happening, we used the following abbreviations for the vertices:

\begin{align*}
v_0&=(k_1-1,\dots,k_m-1)\\
v_m&=(m-2,\dots,m-2)\\
v_{i,j}&=(v_{i,j}^1,\dots,v_{i,j}^m) \text{ with}\\
v_{i,j}^l&=\left\{\begin{array}{cl}k_l+i-2 &\text{ if } l\equiv j+1 \pmod{m}\\k_l+i-1 & \text{ if } l\equiv j+2,\dots,j-i \pmod{m}\\i-1& \text{ if }l\equiv j-i+1 \pmod{m}\\i-2 & \text{ if } l\equiv j-i+2,\dots,j \pmod{m}\end{array}\right.
\end{align*}

In the picture we marked in light grey the vertices which are not produced by the Jacobi path, so which are not included in Figure \ref{loop-kette}, but which we can put in extra in order to have two adjacent arrows at every vertex. Again the dot in the middle of an arrow indicates that there is actually an extra vertex and this extra vertex can be adopted with a normal loop as in Remark \ref{chain-loop}. To have a better view what is happening here, we draw a detailed picture of the second and the third vertex on the original path, which includes the extra vertices we have to put in additionally to take care of the extra vertex created by $\partial_2$.
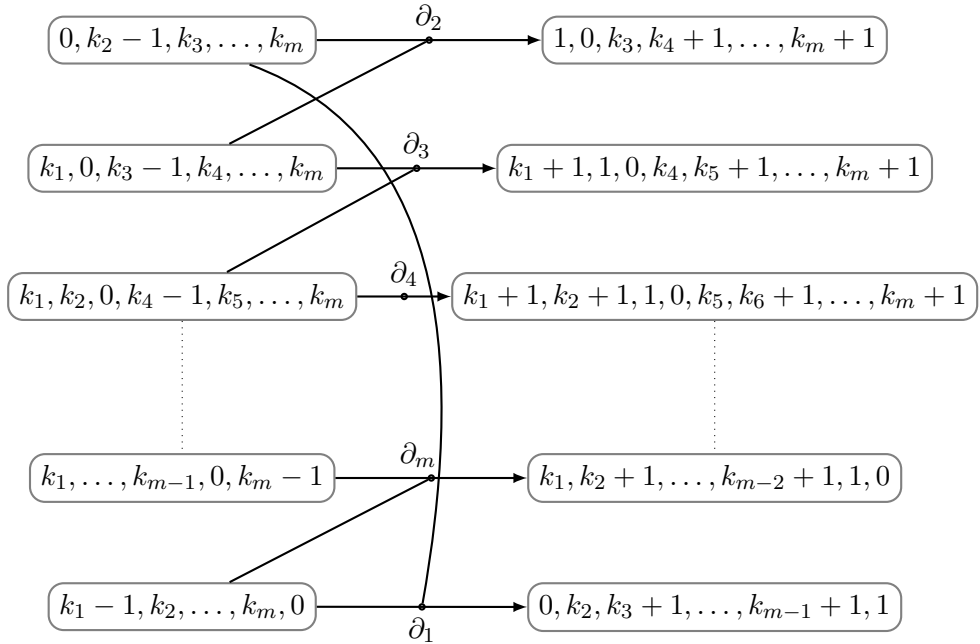
\begin{figure}[H]
\begin{center}
\begin{tikzpicture}
\node[place] (11) {$0,k_2-1,k_3,\dots,k_m$};
\node[place] (22) [right=85pt of 11] {$1,0,k_3,k_4+1,\dots,k_m+1$};
\node[place] (12) [below=30pt of 11] {$k_1,0,k_3-1,k_4,\dots,k_m$};
\node[place] (13) [below=30pt of 12] {$k_1,k_2,0,k_4-1,k_5,\dots,k_m$};
\node[place] (1m-1) [below=50pt of 13] {$k_1,\dots,k_{m-1},0,k_m-1$};
\node[place] (1m) [below=30pt of 1m-1] {$k_1-1,k_2,\dots,k_m,0$};
\node[place] (23) [below=30pt of 22] {$k_1+1,1,0,k_4,k_5+1,\dots,k_m+1$};
\node[place] (24) [below=30pt of 23] {$k_1+1,k_2+1,1,0,k_5,k_6+1,\dots,k_m+1$};
\node[place] (2m) [below=50pt of 24] {$k_1,k_2+1,\dots,k_{m-2}+1,1,0$};
\node[place] (21) [below=30pt of 2m] {$0,k_2,k_3+1,\dots,k_{m-1}+1,1$};
\draw[-latex,thick] (11)--(22)
 node[midway,above] {$\partial_2$}
 node[point,midway] (p11) {};
\draw[-,thick] (p11)--(12);
\draw[-latex,thick] (12)--(23)
 node[midway,above] {$\partial_3$}
 node[point,midway] (p12) {};
\draw[-,thick] (p12)--(13);
\draw[-latex,thick] (13)--(24)
 node[midway,above] {$\partial_4$}
 node[point,midway] (p13) {};
\draw[dotted] (13)--(1m-1);
\draw[dotted] (24)--(2m);
\draw[-latex,thick] (1m-1)--(2m)
 node[midway,above,xshift=-5pt] {$\partial_m$}
 node[point,midway] (p1m-1) {};
\draw[-,thick] (p1m-1)--(1m);
\draw[-latex,thick] (1m)--(21)
 node[midway,below] {$\partial_1$}
 node[point,midway] (p1m) {};
\draw[-,thick] (p1m) to [out=80,in=-20] (11);
\end{tikzpicture}
\end{center}
\caption{The second and third vertex of a complete loop at one position}
\end{figure}
Here we can see in detail that the extra vertex from Figure \ref{loop-kette} and all other vertices on the left side have two adjacent arrows. This means that it is always possible to adjust the coefficients and therefore the extra vertex is linear dependent to the vertices that are already on the path. The vertices on the right side have also another adjacent arrow: For every vertex on the right side the next arrow gives an independent path to the vertex $(m-2,\dots,m-2)$. This was shown in Figure \ref{kompletter-loop-bild} above.\\
We will recall now the important facts of this remark which we will also need later. The vertices on the Jacobi path can be seen in Figure \ref{loop-kette} and we want to summarize the relations of the vertices on this path. First of all starting from the third vertex until the end we have $\kappa(p)+1=\kappa(p+1)$. This is all information we need for these vertices. We also see that between the first and the second vertex there is a gap as soon as a $0$ appears in the second vertex. The special behaviour is between the second and the third vertex. Here we found out that even if a $0$ appears in the third vertex and although the arrow connecting the two vertices creates an extra vertex, the second vertex is still linear dependent to the rest of the path. The extra vertices we put in to make the second vertex linear dependent to the already existing path do not really play an extra role here. The reason is that we were able to find a path starting at every of these extra vertices ending at $(m-2,\dots,m-2)$. It follows that after using the Griffiths formula we will loose the beginning of these paths and they will always be connected to $(m-2,\dots,m-2)$. Therefore we can ignore all the extra vertices we put in and just keep in mind that the second vertex is linear dependent despite the fact that the smallest number of the third vertex is $0$.\\
Another important fact to notice is that nothing changes if we change the order in the positions of the partial derivatives, because we will always have an arrow that connects two vertices with the same smallest numbers. So the basic idea is the same.\\
The above picture does not work if $m=2$, but there something similar happens. One of the important facts in the above picture was that the second vertex did not stand alone. The same thing will happen in the special case of a loop of length $2$. This time we will not omit the coefficients, because they are the key ingredient in this case. We will mark the coefficients from inside the partial derivatives in purple and the ones we can choose in blue. So assume $x_1^{k_1}x_2+x_2^{k_2}x_1$ is in $g(\x)$. Then we get the following picture:
\begin{figure}[H]
\begin{center}
\begin{tikzpicture}
\node[place] (0) {$k_1,k_2$};
\node[place] (1) [right=150pt of 0,label=above:{\small \textcolor{cTangoSkyBlue1}{$0$}}] {$1,k_2$};
\node[place] (2) [right=150pt of 1] {$1,1$};
\draw[-latex,thick] (0)--(1)
 node[midway,below left] {$\partial_1$}
 node[midway,point,label=above:{\footnotesize \textcolor{cTangoSkyBlue1}{$\frac{k_2\alpha}{1-k_1k_2}$}}] (p1) {}
 node[very near start,above] {{\footnotesize $\textcolor{cTangoPlum1}{k_1}$}};
\draw[-latex,thick] (1)--(2)
 node[midway,below left] {$\partial_2$}
 node[midway,point,label=above:{\footnotesize \textcolor{cTangoSkyBlue1}{$\frac{sk_1(k_2-1)\alpha}{(1-k_1k_2)^2}$}}] (p2) {}
 node[very near start,above] {{\footnotesize $\textcolor{cTangoPlum1}{k_2}$}};
\node[place] (3) [below=60pt of p1,label=below:{\small \textcolor{cTangoSkyBlue1}{$0$}}] {$1,2k_2-1$};
\node[place] (4) [below=60pt of p2,label=below:{\small \textcolor{cTangoSkyBlue1}{$0$}}] {$k_1,1$};
\draw[-,thick] (p1)--(3);
\draw[-,thick] (p2)--(4);
\draw[-latex,thick] (3)--(1)
 node[midway,above] {$\partial_2$}
 node[midway,point,label=right:{\footnotesize \textcolor{cTangoSkyBlue1}{$\frac{\alpha}{1-k_1k_2}$}}] (p3) {}
 node[very near start,right] {{\footnotesize $\textcolor{cTangoPlum1}{k_2}$}};
\draw[-,thick] (p3)--(0);
\draw[-latex,thick] (4)--(2)
 node[midway,above] {$\partial_1$}
 node[midway,point,label=right:{\footnotesize \textcolor{cTangoSkyBlue1}{$\frac{s(k_2-1)\alpha}{(1-k_1k_2)^2}$}}] (p4) {}
 node[very near start,right] {{\footnotesize $\textcolor{cTangoPlum1}{k_1}$}};
\draw[-,thick] (p4)--(1);
\node[place] (1a) [below=130pt of 1,label=below:{\small \textcolor{cTangoSkyBlue1}{$0$}}] {$0,k_2-1$};
\draw [-latex reversed,shorten >=8pt,decorate,decoration={snake,pre=moveto,pre length=80pt,post=moveto,post length=5pt}] (1) -- (1a) 
 node [right,near end] {Griffiths formula};
\node[place] (2a) [below=130pt of 2,label=above:{\small \textcolor{cTangoSkyBlue1}{$\frac{s^2(k_1-1)(k_2-1)\alpha}{(1-k_1k_2)^2}$}}] {$0,0$};
\end{tikzpicture}
\end{center}
\caption{The case of a loop of length $2$ at one position}
\end{figure}
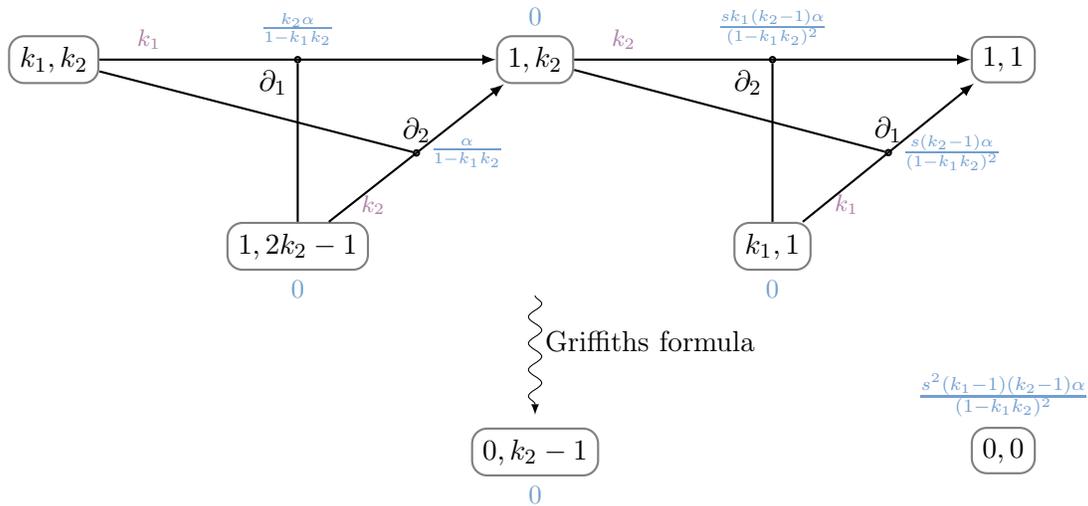
Here $\alpha$ is the coefficient that comes from the arrow before the loop starts. The following coefficients are chosen such that we get $0$ at this vertex. If this is not the first time we used the Griffiths formula, which means that the coefficients should not add up to $0$ at a vertex but to a certain constant, it is very easy to adjust the coefficients. 
We can see in the picture that the second vertex vanishes after the use of the Griffiths formula, because of the fact that the coefficient is $0$. This means again that the second vertex does not stand alone. In other words, if we have a loop of length $2$ everything is linear dependent as long as all entries are bigger than $0$. If we used the Griffiths formula and the entries become $0$, we get only one gap and not $2$ gaps as one might expect.
\end{rem}

With this remark we know what to do if there is a complete loop at one position. Now we will show what happens if a chain (or a part of a chain) is at one position. This partially involves part (ii), because the beginning of the loop and the partial derivative before are never at the same position. In the case that occurs in Lemma \ref{ganzzahlig-loop-chain} the partial derivative belonging to the beginning of the chain is at position $p+1$ and the rest of the part of the chain is at position $p$. 

\begin{rem}\label{chain-gleicher-platz}
From Lemma \ref{ganzzahlig-loop-chain}, we know that if two partial derivatives of neighbouring variables $x_i$ and $x_{i-1}$ in a chain have the same element $p$ in $Q_i$ and $Q_{i-1}$, then this element is also in all $Q_j$ for $j\in\{1,\dots,i\}$. This means that we can use $\partial_1$ at position $p$ and $\partial_j$ at the position $p-1$ for all $2\leq j\leq i$. So let $x_1^{k_1}x_2+\dots+x_{m-1}^{k_{m-1}}x_m+x_m^{k_m}$ be a chain in $g(\x)$ and assume $\partial_1,\dots,\partial_i$ are at the same position. Let $a_{i+1}>k_{i+1}$. Then the picture is the following
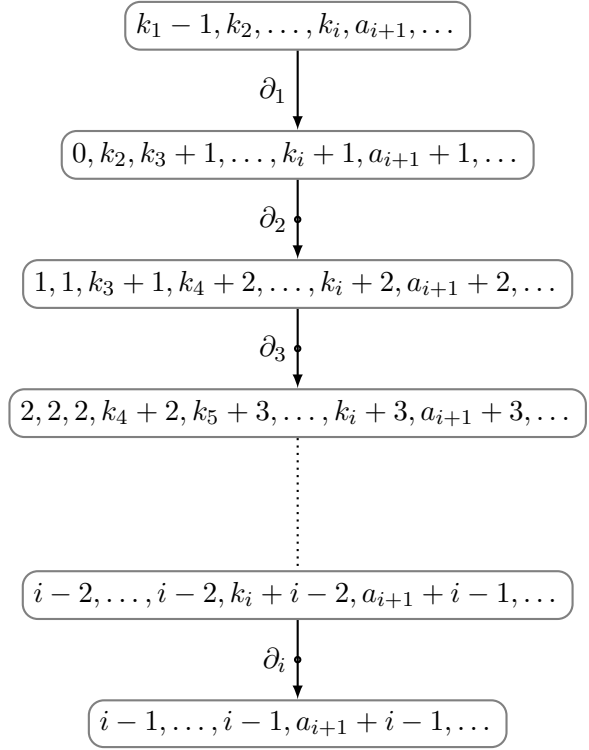
\begin{figure}[H]
\begin{center}
\begin{tikzpicture}
\node[place] (0) {$k_1-1,k_2,\dots,k_i,a_{i+1},\dots$};
\node[place] (1) [below=30pt of 0] {$0,k_2,k_3+1,\dots,k_i+1,a_{i+1}+1,\dots$};
\node[place] (2) [below=30pt of 1] {$1,1,k_3+1,k_4+2,\dots,k_i+2,a_{i+1}+2,\dots$};
\node[place] (3) [below=30pt of 2] {$2,2,2,k_4+2,k_5+3,\dots,k_i+3,a_{i+1}+3,\dots$};
\node[place] (4) [below=50pt of 3] {$i-2,\dots,i-2,k_i+i-2,a_{i+1}+i-1,\dots$};
\node[place] (5) [below=30pt of 4] {$i-1,\dots,i-1,a_{i+1}+i-1,\dots$};
\draw[-latex,thick] (0)--(1)
 node[left,midway] {$\partial_1$};
\draw[-latex,thick] (1)--(2)
 node[left,midway] {$\partial_2$}
 node[point,midway] (p1) {};
\draw[-latex,thick] (2)--(3)
 node[left,midway] {$\partial_3$}
 node[point,midway] (p2) {};
\draw[dotted,thick] (3)--(4);
\draw[-latex,thick] (4)--(5)
 node[left,midway] {$\partial_i$}
 node[point,midway] (p3) {};
\end{tikzpicture}
\end{center}
\caption{The case of a chain at one position}
\end{figure}
Again it makes no sense to change the order of the partial derivatives, because we know from Lemma \ref{dreieck-shift} that everything apart from the beginning of the chain gets shifted. As mentioned earlier it is important to use $\partial_1$ first, because otherwise one of the arrows creating an extra vertex would connect two vertices with the same smallest numbers and therefore disconnect the path earlier. The question in which order to use the rest of the partial derivatives is similar. The answer is that doing if we use them in a different order, then there is an arrow $\partial_j$ where for the vertex at the arrow tip  $(a_1,\dots,a_m)$ we have $a_j=a_i$ for an $i<j$ and $a_j$ is also the smallest number. If we use the Griffiths formula several times such that $a_j=a_i=0$, then the partial derivative $\partial_i$ does not fit here and therefore we are not able to take care of the extra vertex produced by $\partial_j$. For later purposes the important part of this remark is that the smallest number increases at every position.
\end{rem}

Now let us investigate the last two cases, listed under (ii). So the smallest possible positions for $\partial_{i_1}$ and $\partial_{i_2}$ are $p$ and $p+1$ respectively and the first partial derivative $\partial_{i_1}$ gets shifted. First assume that we are in subcase (a) and the two partial derivatives do not belong to neighbouring variables of a chain or loop. Then we can choose $\partial(p+1)=\partial_{i_1}$ and $\partial(p+2)=\partial_{i_2}$. This way $\kappa(p+2)=\kappa(p+3)$ and the vertices only get disconnected if $\partial_{i_2}$ is in a loop, a middle or end of a chain and the smallest number is $0$ and after the next Griffiths step everything vanishes. If we choose $\partial(p+1)=\partial_{i_2}$ and $\partial(p+2)=\partial_{i_1}$, then $\kappa(p+2)\leq\kappa(p+3)=2$, because $\partial_{i_1}$ was shifted further than $\partial_{i_2}$. So the connection gets cut earlier or at the same time. So we should use the first way of shifting to be on the safe side.\\
The last case to consider is (iib), so as before the two partial derivatives $\partial_{i_1}$ and $\partial_{i_2}$ can be used first at position $p$ and $p+1$ and the first partial derivative $\partial_{i_1}$ gets shifted. This time, however, they correspond to neighbouring variables in a loop or a chain. If $\partial_{i_2}$ is the beginning of a chain then this is the case of Remark \ref{chain-gleicher-platz} and we should choose $\partial(p+1)=\partial_{i_2}$ and $\partial(p+2)=\partial_{i_1}$. Otherwise, both partial derivatives are either from a loop, or from the middle or end of a chain, then we should choose $\partial(p+1)=\partial_{i_1}$ and $\partial(p+2)=\partial_{i_2}$. The argument is the same as in the case (iia). If we use the partial derivatives as indicated then $\kappa(p+2)=\kappa(p+3)$, so the diagram disconnects first when $\kappa(p+2)=\kappa(p+3)=0$. If we change the positions, then $\kappa(p+2)=1$ and $\kappa(p+3)=2$ and the connection is cut earlier, so we stick to the first way of ordering the partial derivatives.\\
Now we know where to shift everything and this explains most of the shifting we did in Example \ref{exIV} and the shifting we have to do in general. But as mentioned before we still have to look closer at the last $n$ positions. In Remark \ref{letzte-positionen} we saw that all partial derivatives are at the positions $\widehat{d}-n+2$ and $\widehat{d}-n+1$. At all the positions $p>\widehat{d}-n+2$ there is no partial derivative that produces $1$ as an entry here. So we have to spread all partial derivatives over the positions $\widehat{d}-n+1,\dots,\widehat{d}$. There we have to use every partial derivative exactly once. In order to see all linear dependencies between the monomials this should be done in the same way as before. So every chain and every loop itself should be used in the order suggested in Remark \ref{loop-gleicher-exponent} and \ref{chain-gleicher-platz}. Because separate chains and loops do not interact, the order between the loops and chains does not matter. Notice that at position $\widehat{d}-n+1$ there is either an arrow with an extra vertex and $\kappa(\widehat{d}-n+1)=1$ or an arrow without extra vertex and $\kappa(\widehat{d}-n+1)=0$. In the second case we have to shift this partial derivative to the left, which also explains the shift to the left in Example \ref{exIV}. This means, no matter what the rest of the path looks like, after the first use of the Griffiths formula there will be a gap behind the vertex at position $\widehat{d}-n$, so $\partial(\widehat{d}-n)=0$. In the first case the vertex at position $\widehat{d}-n+1$ will not vanish but the arrow will not fit with the complete loop anymore and in the second case the vertex at position $\widehat{d}-n+1$ simply vanishes. Because of Remark \ref{loop-gleicher-exponent} and \ref{chain-gleicher-platz} the smallest numbers never decrease between the position $\widehat{d}-n+1$ and $\widehat{d}$. Therefore the vertex at position $p+1$ always vanishes after the vertex at position $p$ and it follows that only the beginning of the path from $\widehat{d}-n+1$ to $\widehat{d}$ vanishes. So the last $n$ steps and what is left after all the Griffiths steps is always linear dependent to the vertex $(a,\dots,a),a\leq n-1$, where we started.

\section{Proof of the order of the Picard-Fuchs equation}

First we will prove a weak form of the main theorem of this chapter. We will show that one needs $u$ basis elements to write one special form of the forms appearing in the Picard-Fuchs equation of $f(\x)$. In the proof of Theorem \ref{thm} we will see that these are all basis elements we need for all forms appearing in the Picard-Fuchs equation, which means that $u$ is also the order of the Picard-Fuchs equation.

\begin{prop}\label{propo}
The form $\frac{s^n(\prod_{i=1}^n x_i)^{n-1}\Omega_0}{f^n}$ is a linear combination of $u$ basis elements of the primitive cohomology with coefficients in $\C(s)$. In other words, using the Griffiths formula, $(\prod_{i=1}^n x_i)^{n-1}$ can be written as a combination of $u$ basis elements of the Milnor ring $\C(s)/J(f)$ with coefficients in $\C(s)$.
\end{prop}

Notice that the $(\prod_{i=1}^n x_i)^{n-1}$ can never be a basis element in the Milnor ring itself, because we have $(\prod_{i=1}^n x_i)^{j}\in J(f)$ for $j\geq n-1$. Before we prove this proposition we want to calculate an example to have a better understanding what we have to do in the proof. We will continue the example we already used throughout the whole chapter. The notation for this can be found in Example \ref{exIV}. In order to show the proposition in this example, we will start with $(\prod_{i=1}^n x_i)^{n-1}$ and count how many basis elements we need to write $(\prod_{i=1}^n x_i)^{n-1}$ as a linear combination of them. We will do this by following the steps in the Griffiths-Dwork method.

\begin{ex}\emph{(Continuation Example \ref{exIV})}
We want to calculate the number of basis elements we need to write $(wxyz)^3$ in this example. We know that $(wxyz)^3$ is an element of the Jacobian ideal. Following the Griffiths-Dwork method we have to write down the Jacobi path, which we already did in Example \ref{jpexIV}. We start with this Jacobi path and use the Griffiths formula once. This means we subtract $(1,1,1,1)$ from every vertex. Again we do not mention the various coefficients one needs to do the actual calculations because they do not give any interesting input for the calculation of the number of basis elements. The important part for us is to count the disconnected parts of the Jacobi path in every degree. We want to note that for counting the basis elements we would not need to write down the extra vertices. It would be enough to focus on the $18$ vertices on the Jacobi path and remember which of the arrows produces an extra vertex. So after using the Griffiths formula once Figure \ref{jpexIVpath} becomes:
\begin{figure}[H]
\begin{center}
\begin{tikzpicture}
\node[place] (0) {\textcolor{cTangoPlum1}{$2,2,2,2$}};
\node[place] (1) [right=35pt of 0] {$3,1,2,3$}; 
\node[place] (2) [right=35pt of 1] {$4,2,3,1$}; 
\node[place] (3) [right=35pt of 2] {$5,1,3,2$}; 
\node[place] (4) [below=35pt of 3] {$6,0,3,3$}; 
\node[place] (5) [below=30pt of 4] {$7,1,1,3$}; 
\node[place] (6) [below=30pt of 5] {$8,2,2,1$}; 
\node[place] (7) [below=30pt of 6] {$9,1,2,2$}; 
\node[place] (8) [below=30pt of 7] {$10,0,2,3$}; 
\node[place] (9) [below=30pt of 8] {$11,1,3,1$}; 
\node[place] (10) [left=30pt of 9] {$12,0,3,2$}; 
\node[place] (11) [left=30pt of 10] {$13,1,1,2$}; 
\node[place] (12) [left=30pt of 11] {$14,0,1,3$}; 
\node[place] (13) [above=30pt of 12] {$15,1,2,1$}; 
\node[place] (14) [above=30pt of 13] {$16,0,2,2$}; 
\node[place] (17) [below=30pt of 0] {$19,1,1,1$}; 
\node[place] (16) [below=30pt of 17] {$18,0,0,3$}; 
\draw[-latex,thick] (0)--(1)
 node[midway,below] {$\partial_2$};
\draw[-latex,thick] (1)--(2)
 node[midway,below] {$\partial_4$}
 node[point,midway] (h1) {};
\draw[-latex,thick] (2)--(3)
 node[midway,below] {$\partial_2$};
\draw[-latex,thick] (4)--(5)
 node[midway,left] {$\partial_3$}
 node[point,midway] (h2) {};
\draw[-latex,thick] (3)--(4)
 node[midway,left] {$\partial_2$};
\draw[-latex,thick] (5)--(6)
 node[midway,left] {$\partial_4$}
 node[point,midway] (h3) {};
\draw[-latex,thick] (6)--(7)
 node[midway,left] {$\partial_2$};
\draw[-latex,thick] (7)--(8)
 node[midway,left] {$\partial_2$};
\draw[-latex,thick] (8)--(9)
 node[midway,left] {$\partial_4$}
 node[point,midway] (h4) {};
\draw[-latex,thick] (10)--(11)
 node[midway,above] {$\partial_3$}
 node[point,midway] (h5) {};
\draw[-latex,thick] (9)--(10)
 node[midway,above] {$\partial_2$};
\draw[-latex,thick] (11)--(12)
 node[midway,above] {$\partial_2$};
\draw[-latex,thick] (12)--(13)
 node[midway,right] {$\partial_4$}
 node[point,midway] (h6) {};
\draw[-latex,thick] (13)--(14)
 node[midway,right] {$\partial_2$};
\draw[-latex,thick] (16)--(17)
 node[midway,right] {$\partial_4$}
 node[point,midway] (h8) {};
\draw[-latex,thick] (17)--(0)
 node[midway,right] {$\partial_1$};
\node[place] (p1) [above=30pt of h1] {$3,1,5,1$};
\draw[-,thick] (h1)--(p1);
\draw[-latex,thick] (p1)--(2)
 node[point,midway] (h11) {}
 node[midway,left] {$\partial_3$};
\node[place] (p11) [right=15pt of p1] {$3,3,3,0$};
\draw[-,thick] (h11)--(p11);
\draw[-latex,thick] (p11)--(2)
 node[midway,right] {$\partial_2$};
\node[place] (p2) [right=35pt of h2] {$6,2,1,2$};
\draw[-,thick] (h2)--(p2);
\draw[-latex,thick] (p2)--(5)
 node[midway,below right] {$\partial_2$};
\node[place] (p3) [right=35pt of h3] {$7,1,4,1$};
\draw[-,thick] (h3)--(p3);
\draw[-latex,thick] (p3)--(6)
 node[point,midway] (h33) {}
 node[midway,below] {$\partial_3$};
\node[place] (p33) [below=20pt of p3] {$7,3,2,0$};
\draw[-,thick] (h33)--(p33);
\draw[-latex,thick] (p33)--(6)
 node[midway,below] {$\partial_2$};
\node[place] (p4) [right=35pt of h4] {$10,0,5,1$};
\draw[-,thick] (h4)--(p4);
\draw[-latex,thick] (p4)--(9)
 node[point,midway] (h44) {}
 node[midway,below] {$\partial_3$};
\node[place] (p44) [below=20pt of p4] {$10,2,3,0$};
\draw[-,thick] (h44)--(p44);
\draw[-latex,thick] (p44)--(9)
 node[midway,below] {$\partial_2$};
\node[place] (p5) [below=30pt of h5] {$12,2,1,1$};
\draw[-,thick] (h5)--(p5);
\draw[-latex,thick] (p5)--(11)
 node[midway,left] {$\partial_2$};
\node[place] (p6) [left=35pt of h6] {$14,0,4,1$};
\draw[-,thick] (h6)--(p6);
\draw[-latex,thick] (p6)--(13)
 node[point,midway] (h66) {}
 node[midway,above] {$\partial_3$};
\node[place] (p66) [above=20pt of p6] {$14,2,2,0$};
\draw[-,thick] (h66)--(p66);
\draw[-latex,thick] (p66)--(13)
 node[midway,above] {$\partial_2$};
\node[place] (p8) [left=35pt of h8] {$18,0,3,1$};
\draw[-,thick] (h8)--(p8);
\draw[-latex,thick] (p8)--(17)
 node[point,midway] (h88) {}
 node[midway,above] {$\partial_3$};
\node[place] (p88) [above=20pt of p8] {$18,2,1,0$};
\draw[-,thick] (h88)--(p88);
\draw[-latex,thick] (p88)--(17)
 node[midway,above] {$\partial_2$};
\end{tikzpicture}
\end{center}
\caption{The Jacobi path for Example \ref{exIV} after the first use of the Griffiths formula}
\end{figure}
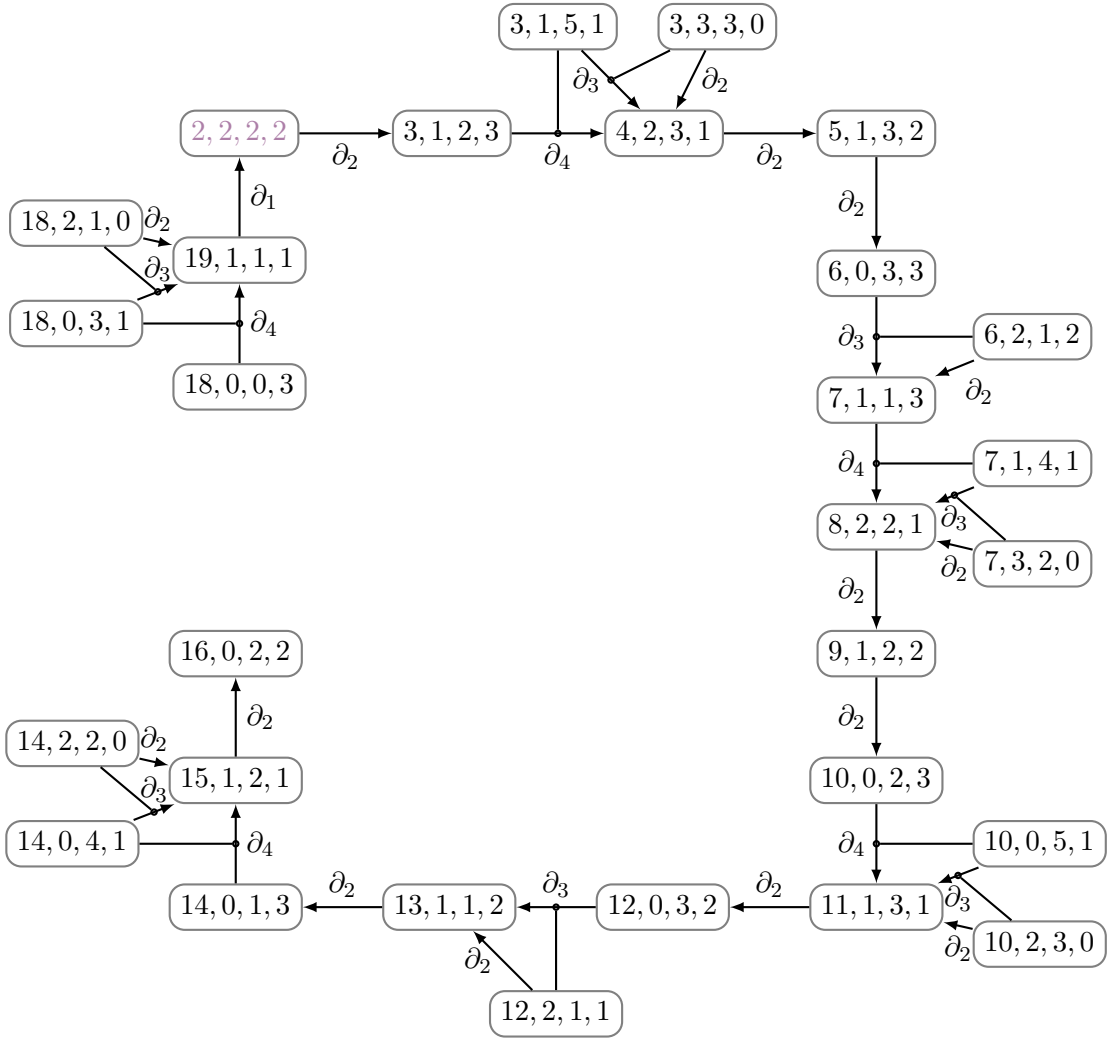
As we can see from the picture, all monomials are still connected. We only have one gap between the $14$th and the $16$th vertex. But we knew before that we will create a gap after the vertex at position $\widehat{d}-n=18-4=14$. This means we need one basis element in this degree (indicated by a purple colour) in order to be able to choose all coefficients appropriately. Let us assume we add this basis element with the appropriate coefficient, such that the resulting path is in the Jacobian ideal. Now we can use the Griffiths formula again and see how much gaps we get in the next degree. We end up with the following picture:
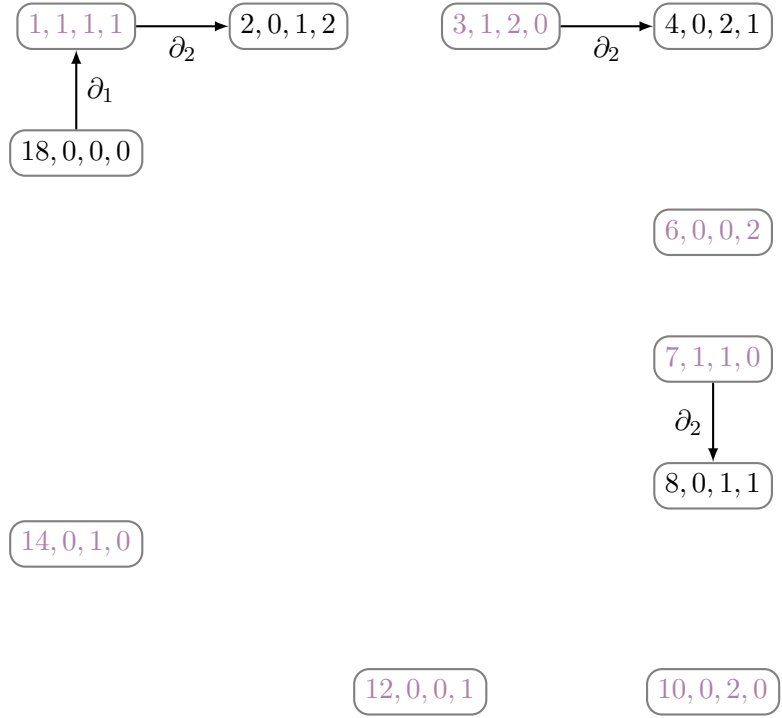
\begin{figure}[H]
\begin{center}
\begin{tikzpicture}
\node[place] (0) {\textcolor{cTangoPlum1}{$1,1,1,1$}};
\node[place] (1) [right=35pt of 0] {$2,0,1,2$}; 
\node[place] (2) [right=35pt of 1] {\textcolor{cTangoPlum1}{$3,1,2,0$}}; 
\node[place] (3) [right=35pt of 2] {$4,0,2,1$}; 
\node[place] (5) [below=60pt of 3] {\textcolor{cTangoPlum1}{$6,0,0,2$}}; 
\node[place] (6) [below=30pt of 5] {\textcolor{cTangoPlum1}{$7,1,1,0$}}; 
\node[place] (7) [below=30pt of 6] {$8,0,1,1$}; 
\node[place] (9) [below=60pt of 7] {\textcolor{cTangoPlum1}{$10,0,2,0$}}; 
\node[place] (11) [left=60pt of 9] {\textcolor{cTangoPlum1}{$12,0,0,1$}}; 
\node[place] (17) [below=30pt of 0] {$18,0,0,0$}; 
\node[place] (13) [below=130pt of 17] {\textcolor{cTangoPlum1}{$14,0,1,0$}}; 
\draw[-latex,thick] (0)--(1)
 node[midway,below] {$\partial_2$};
\draw[-latex,thick] (2)--(3)
 node[midway,below] {$\partial_2$};
\draw[-latex,thick] (6)--(7)
 node[midway,left] {$\partial_2$};
\draw[-latex,thick] (17)--(0)
 node[midway,right] {$\partial_1$};
\end{tikzpicture}
\end{center}
\caption{The Jacobi path for Example \ref{exIV} after the second use of the Griffiths formula}
\end{figure}
Again we coloured a choice for basis elements we need in purple. This means we need $7$ basis elements in this degree. If we use the Griffiths formula again, we are left with $(0,0,0,0)$, which therefore has to be a basis element as well. So in total we counted $1$ basis element in degree $0$ and degree $2\cdot 18$ and $7$ basis elements in degree $18$, which adds up to $9$ basis elements overall. If the theorem is true, then it should hold that $u=9$. So we will calculate $u$ with the $Q_i$ calculated in (\ref{qi-ex}):
\begin{align*}
u&=\widehat{d}-|\bigcup_{i=1}^nQ_i^{\Z}|=18-|\{18\}\cup\{2,4,6,8,10,12,14,16,18\}\cup\{6,12,18\}\cup\{18\}|\\
&=18-9=9.
\end{align*}
\end{ex}

\begin{lemma}\label{u=basis}
Let
\begin{align*}
\eta_1:&=\{p\,|\,\kappa(p)+1=\kappa(p+1),1\leq p\leq \widehat{d}\},\\
\eta_2:&=\{p\,|\,\kappa(p)+2=\kappa(p+1),1\leq p\leq \widehat{d}\},\\
\eta_3:&=\{p\,|\,\kappa(p)=\kappa(p+1),1\leq p\leq \widehat{d},\partial(p)=\partial_i \text{ satisfies Condition ($\ast$)}\}.
\end{align*}
Then
\begin{align*}
u=|\eta_1|+2|\eta_2|+|\eta_3|:=\eta.
\end{align*}

\pagebreak
\textbf{Condition ($\ast$)}
\begin{itemize}
\item $\partial_i$ creates an extra vertex, i.e. $\widehat{q}_ik_i\neq \widehat{d}$.
\item If $x_i$ is part of the loop $x_{i_1}^{k_{i_1}}x_{i_2}+\dots+x_{i_m}^{k_{i_m}}x_{i_1}$ and the partial derivatives are used in the order $\partial_{i_1},\partial_{i_2},\dots,\partial_{i_m}$ whenever they have the same smallest position, then $i$ should not equal $i_2$.
\end{itemize}
\end{lemma}

\begin{proof}
First recall that $u=\sum_{i=1}^n|Q_i|-|\bigcup_{i=1}^nQ_i^{\Z}|=\sum_{i=1}^{n}|Q_i^{\Q}|+\sum_{i=1}^n|Q_i^{\Z}|-|\bigcup_{i=1}^nQ_i^{\Z}|$.
{\em Part 1: $u\leq\eta$}\\
First assume $q\in Q_{i_j}^{\Q}$ for $j=1\dots,\ell$, this means we have the summand $\ell$ in $u$. But this also means that position $p=\lfloor q \rfloor -n+2$  is the smallest possible position for $\partial_{i_1},\dots,\partial_{i_\ell}$ and due to Lemma \ref{ganzzahlig-loop-chain} none of the $x_{i_j}$ are neighbouring variables in a loop or a chain. Therefore they are all independent and $\partial_{i_j}$ adds $1$ to all entries $i_1,\dots,i_\ell$ except $i_j$. Let us assume that $\partial_{i_1}$ got shifted to position $\widetilde{p}$ (according to Lemma \ref{dreieck-shift} it has to be shifted) and all others where shifted correspondingly, i.e. $\partial(\widetilde{p}+j-1)=\partial_{i_{j}}$ for $1\leq j\leq \ell$. This means that after completely shifting we get $\kappa(\widetilde{p}+\ell)=\kappa(\widetilde{p}+\ell-1)+1=\dots=\kappa(\widetilde{p}+2)+(\ell-2)=\kappa(\widetilde{p}+1)+(\ell-1)$ and therefore $\widetilde{p}+1,\dots,\widetilde{p}+\ell-1\in\eta_1$. In addition we get that $\kappa(\widetilde{p})\leq\kappa(\widetilde{p}+1)$, which means $\widetilde{p}\in\eta_3$, because each $\partial_{i_j}$ creates an extra vertex, or $\kappa(\widetilde{p})+1=\kappa(\widetilde{p}+1)$, which means $\widetilde{p}\in\eta_1$. In total this means that we also sum up $\ell$ in $\eta$.\\
Now assume $q\in Q_{i_j}^{\Z}$ for $j=1\dots,\ell+1$, this means that we sum up $\ell$ in $u$. But this also means that whole loops are at position $q-n+1$ or the beginning of a chain is at position $q-n+2$ and the rest of the chain stopping somewhere is at position $q-n+1$ as shown in Lemma \ref{ganzzahlig-loop-chain}. Remember that we saw in Lemma \ref{dreieck-shift} that all partial derivatives at position $q-n+1$ get shifted anyway. So assume that everything gets shifted to position $\widetilde{p}$. Now we additionally shift every loop and every chain according to Remark \ref{loop-gleicher-exponent} and \ref{chain-gleicher-platz} respectively and we want to look at $\kappa(\widetilde{p}),\dots,\kappa(\widetilde{p}+\ell)$. Therefore suppose a loop of length $m$ got shifted to $\widetilde{p}+a$. Then Remark \ref{loop-gleicher-exponent} tells us that $\kappa(\widetilde{p}+a+m)=\kappa(\widetilde{p}+a+m-1)+1=\dots=\kappa(\widetilde{p}+a+2)+m-2$ and $\kappa(\widetilde{p}+a+1)=\kappa(\widetilde{p}+a+2)$. In addition we get that $\kappa(\widetilde{p}+a)+2=\kappa(\widetilde{p}+a+1)$, because either there is another loop at the positions before $\widetilde{p}$ and as we can see in this loop $\kappa(\widetilde{p}+a+1)+m-2=\kappa(\widetilde{p}+a+m)$ so the smallest numbers increase by $m-2$ in $m-1$ steps and in addition the partial derivative at position $\widetilde{p}-1$ got shifted one less which leads to $\kappa(\widetilde{p}+a)+2=\kappa(\widetilde{p}+a+1)$. If there is a chain at the positions before $\widetilde{p}$, then the beginning of the chain was originally possible at position $q-n+2$ and the loop at position $q-n+1$ and now the loop is at one position later, therefore the loop got shifted two more and $\kappa(\widetilde{p}+a)+2=\kappa(\widetilde{p}+a+1)$. Now we can calculate $\eta$: We have $\widetilde{p}+a+m-1,\dots,\widetilde{p}+a+2\in\eta_1$ and $\widetilde{p}+a\in\eta_2$, which means adding up $m-2+2\cdot 1=m$ in $\eta$. Now let us assume that there is the beginning of a chain which now has length $m$ at position $\widetilde{p}+a$. Then Remark \ref{chain-gleicher-platz} tells us that $\kappa(\widetilde{p}+a+m)=\kappa(\widetilde{p}+a+m-1)+1=\dots=\kappa(\widetilde{p}+a+1)+m-1$. In addition $\kappa(\widetilde{p}+a)+1=\kappa(\widetilde{p}+a+1)$, because if there is a chain at the positions before, we just shifted one more and if there is a loop at the positions before, this got shifted the same amount, but as mentioned before the smallest numbers increased one less. In total this gives $\widetilde{p}+a+m-1,\dots,\widetilde{p}+a\in\eta_1$ and therefore we added $m$ in $\eta$. The last thing to do is to look at the chain or loop at position $\widetilde{p}$. The smallest numbers at position $\widetilde{p}+1$ and the later ones are as above. But the arrow at position $\widetilde{p}-1$ got shifted at least as much as the one at position $\widetilde{p}$, so $\kappa(\widetilde{p})\geq\kappa(\widetilde{p}+1)$. If there is a chain at position $\widetilde{p}$, we definitely get no contribution to $\eta$ from this position. If there is a loop at position $\widetilde{p}$, then $\widetilde{p}\in\eta_3$. So in general this means the chain or loop at position $\widetilde{p}$ adds up one less then its length and all others add up exactly their lengths in $\eta$. Since the lengths add up to $\ell+1$, we add up $\ell$ in $\eta$.\\
{\em Part 2: $u\geq\eta$}\\
We start with a position $p$ in $\eta_1$. This means that $\kappa(p)+1=\kappa(p+1)$ and therefore $\partial(p)=\partial_{i_1}$ and $\partial(p+1)=\partial_{i_2}$ have the same smallest possible position $\widetilde{p}$. This means either $\widetilde{p}+n-2\in Q_{i_1}\cap Q_{i_2}$, $\lfloor q \rfloor=\widetilde{p}+n-2$ with $q\in Q_{i_1}\cap Q_{i_2}$ or $\widetilde{p}+n-2\in Q_{i_1}$ and $\lfloor q \rfloor=\widetilde{p}+n-2$ with $q\in Q_{i_2}$. But in all cases we add up $1$ in $u$.\\
Now let $p\in\eta_2$. So $\kappa(p)+2=\kappa(p+1)$, this only happens if we shifted $\partial(p+1)$ over $\partial(p)$ which we only do if a full loop is at the same position and as we saw in part 1, this leads to adding up $2$ for this position in $u$.\\
The last possibility is $p\in\eta_3$. Here we have $\kappa(p)\geq\kappa(p+1)$ and $\partial(p)=\partial_i$ creates an extra vertex. It follows that either $\partial_i$ comes from a position $\widetilde{p}$ with $\lfloor q \rfloor=\widetilde{p}+n-2$ and $q\in Q_i^{\Q}$ or $\widetilde{p}+n-2\in Q_i^{\Z}$ and this is the first position of a complete loop. But in both cases we have added up $1$ in $u$.
\end{proof}

\begin{proof}[Proof of Proposition \ref{propo}]
We will now start proving Proposition \ref{propo} using all the results we achieved so far. The rough idea is the following: We count the holes that occur after using the Griffiths formula and relate them to the sets $\eta_1,\eta_2$ and $\eta_3$ and therefore with Lemma \ref{u=basis} to the number $u$, because for each hole on the Jacobi path we need an extra basis element.\\
So we investigate all cases when a path becomes disconnected. This depends on the smallest numbers occurring in a monomial, or correspondingly vertex. Given that a vertex vanishes if the smallest number was $0$ before using the Griffiths formula. If two vertices are neighbours in the Jacobi path we distinguish between $3$ cases of relations between their smallest numbers. Let $p$ and $p+1$ be two positions on the Jacobi path, then the following situations can occur:
\begin{enumerate}[label=(\roman*)]
\item $\kappa(p)=\kappa(p+1)$,
\item $\kappa(p)<\kappa(p+1)$ or
\item $\kappa(p)>\kappa(p+1)$.
\end{enumerate}
We should notice that the maximal gap between the smallest numbers in (ii) is $2$ because in the way we shift, the arrow used before is at most shifted two less than the arrow we use between the two vertices. We already investigated this in the proof of Lemma \ref{u=basis}.\\
We will count one basis element for every start of a disconnected part of the path. First consider case (iii), so $\kappa(p)>\kappa(p+1)$. If $\partial(p)$ does not create an extra vertex, then this will only shorten a path after using the Griffiths formula the appropriate number of times, but this will never be the beginning of a path. So we can neglect this case. But if $\partial(p)$ does create an extra vertex and we get to the point that $\kappa(p+1)=0$ the arrow $\partial(p)$ does not fit here together with the whole loop or the beginning of the chain. So the path gets disconnected and after the next use of the Griffiths formula the vertex $p+1$ vanishes. So we have to count one basis element for every time that $\kappa(p)>\kappa(p+1)$ and $\partial(p)$ creates an extra vertex. But this is done in $|\eta_3|$.\\
Now consider case (ii). As long as $\kappa(p)\geq 1$ the positions $p$ and $p+1$ are always connected by an arrow. So assume that we used the Griffiths formula several times until $\kappa(p)=0$ and $\kappa(p+1)>0$. Then there is still an arrow connecting the two vertices no matter which kind of partial derivative it is, but after the next use of the Griffiths formula the vertex at position $p$ will vanish and the one at position $p+1$ will still be there. This means that at position $p+1$ a disconnected part of the path starts, which shows that we need an extra basis element here. The vertex at position $p+1$ will stay until $\kappa(p+1)=0$, so we need an extra basis element in every degree until the vertex vanishes. So we need $1$ basis element if $\kappa(p)+1=\kappa(p+1)$ and $2$ basis elements if $\kappa(p)+2=\kappa(p+1)$. The set $\eta_1$ counts exactly the first case and $2\eta_2$ the second case.\\
The last case to consider is (i). Here we have to distinguish between several cases: First assume that the two vertices are connected by an arrow $\partial_i$ with $\widehat{q}_i k_i=\widehat{d}$, i.e. the arrow has no additional vertex. In this case we can put the arrow in as long as $\kappa(p+1)\geq0$ and then the vertices at position $p$ and $p+1$ vanish at the same time, so we don't need an extra basis element here. Remember that if $\partial_i$ belongs to the beginning of a chain of length $\geq 2$ there is an arrow as long as the $(i+1)$th entry of the vertex at position $p+1$ is $>0$. This can only occur if $\partial_{i+1}$ was used at the position before. But we already discussed in Remark \ref{chain-gleicher-platz} that if $\partial_{i+1}$ is at the position before $\partial_i$ we shift $\partial_{i+1}$ over $\partial_i$. So this case never occurs and we can assure that we never need an extra basis element if $\partial_i$ is between two vertices with the same smallest numbers and $\widehat{q}_i k_i=\widehat{d}$. Assume now that we are still in case (i), so the smallest numbers are the same, but the arrow $\partial_i$ between the two vertices produces an extra vertex. As long as $\kappa(p+1)\geq 1$ the arrow and all the other arrows from the chain or loop can be used here, which means that we can always choose the coefficients in a way that the additional vertices vanish and the vertex on the path has the appropriate coefficient. If we use the Griffiths formula until $\kappa(p+1)=0$ we might still be able to put in the arrow but the next arrow in the loop or chain does not fit anymore. So we need an extra basis element except if we are in the situation of Remark \ref{loop-gleicher-exponent}, where all partial derivatives of the loop or chain are at the same position. This is exactly what is counted in $|\eta_3|$ in addition to case (iii).\\
In total we see that counting basis elements is the same as $|\eta_1|+|2\eta_2|+|\eta_3$| and therefore the number of basis elements is $u$.
\end{proof}

\begin{rem}
We have $n-1$ basis elements for sure, because we need at least $1$ basis element in every degree. We choose this to be $\frac{s^k(\prod x_i)^{k-1}\Omega_0}{f^k}$ for $1\leq k\leq n-1$, because we need this anyway to write the first $(n-1)$ derivatives of $\omega$.
\end{rem}

Now we are able to put everything together and prove Theorem \ref{thm}.

\begin{proof}[Proof of Theorem \ref{thm}]
To prove Theorem \ref{thm} we will show that all the powers of $\prod_{i=1}^n x_i$ can be written as a combination of the same $u$ basis elements we needed for $(\prod_{i=1}^n x_i)^{n-1}$ seen in Proposition \ref{propo}. If we have done that, it is clear that $\delta^i\omega$ can be written as a combination of $u$ basis elements for all $i$. So if we take all $\delta^i\omega$ up to $i=u$, then we get a linear relation between them. So the Picard-Fuchs equation has order $u$.\\
We will show by induction that $(\prod_{i=1}^n x_i)^j$ can for all $j$ be written as a combination of the same basis elements. Therefore we first look at $(\prod_{i=1}^n x_i)^n$: We know how to write $(\prod_{i=1}^n x_i)^{n-1}$ in terms of the partial derivatives, so if we multiply this by $\prod_{i=1}^n x_i$ we get an expression for $(\prod_{i=1}^n x_i)^n$. In our notation this means adding $(1,\dots,1)$ to every monomial that appears on the Jacobi path. Now we can use the Griffiths formula and because all entries were bigger than $1$ all monomials are still there. Now we can use the Griffiths formula again, the only thing we have to add is a multiple of $(\prod_{i=1}^n x_i)^{n-1}$. From now on everything works as in the case for $(\prod_{i=1}^n x_i)^{n-1}$. This means we can write $(\prod_{i=1}^n x_i)^n$ as a linear combination of the $u$ basis elements and $(\prod_{i=1}^n x_i)^{n-1}$, but this monomial is itself a linear combination of the $u$ basis elements. So $(\prod_{i=1}^n x_i)^n$ can be written with the same $u$ basis elements as $(\prod_{i=1}^n x_i)^{n-1}$. Now look at $(\prod_{i=1}^n x_i)^j$ for a $j>n-1$ then again we get the expression of $(\prod_{i=1}^n x_i)^j$ in the partial derivatives by multiplying the expression of $(\prod_{i=1}^n x_i)^{n-1}$ by $(\prod_{i=1}^n x_i)^{j-n+1}$. And again we can use the Griffiths formula once and no monomial will vanish until we used the Griffiths formula $j-n+2$ times. This means that we can write $(\prod_{i=1}^n x_i)^j$ as a combination of the $u$ basis elements and all $(\prod_{i=1}^n x_i)^l$ with $l<j$, but by induction all these powers can themselves be written as a linear combination of the same $u$ basis elements. So in total we get that we can write all powers of $\prod_{i=1}^n x_i$ as a linear combination of the same $u$ basis elements. This leads to our statement that the Picard-Fuchs equation of $f(\x)$ has order $u$.
\end{proof}

\section{Detailed example for the Picard-Fuchs equation}\label{detailed-example-gd}
In this section we want to calculate in all details an example of computing the Picard-Fuchs equation with the Griffiths-Dwork method. We will choose a slightly smaller example as in the section before. Theorem \ref{thm} tells us immediately how many calculations we have to do, because we know how many basis elements we need and therefore how many of the $\delta^i\omega$ we have to calculate. We will prove the actual appearance of the Picard-Fuchs equation of Theorem \ref{conj} not by using these calculations, but we want to show how this can be done using the Griffiths-Dwork method and especially that with the help of Theorem \ref{thm} and our new diagrammatic notation it can be done relatively quick.

\begin{ex}
Let $g(x_1,x_2,x_3,x_4)=x_1^5x_2+x_2^4x_3+x_3^8+x_4^2$. This is the polynomial we are looking at. The reduced weights for this polynomial are given by $(q_1,q_2,q_3,q_4)=(5,7,4,16)$ and the degree is $d=32$. The transposed polynomial is given by $g^t(w,x,y,z)=w^5+wx^4+xy^8+z^2$ and has weights $(\widehat{q}_1,\widehat{q}_2,\widehat{q}_3,\widehat{q}_4)=(2,2,1,5)$ and the degree is $\widehat{d}=10$. So according to Theorem \ref{thm} the Picard-Fuchs equation of $f(x_1,x_2,x_3,x_4)=x_1^5x_2+x_2^4x_3+x_3^8+x_4^2+sx_1x_2x_3x_4$ has order 

\begin{align*}
u=\widehat{d}-|\bigcup_{i=1}^nQ_i^{\Z}|=10-|\{5,10\}\cup\{5,10\}\cup\{10\}\cup\{2,4,6,8,10\}|=10-6=4.
\end{align*}

This means we have to calculate $\delta^i\omega$ for $i=0,\dots,4$ and $\omega=\frac{s\Omega_0}{f}$. We know from Remark \ref{delta-omega} that the derivatives of $\omega$ can be written as a sum of $\frac{\ell! s^{\ell+1}(x_1x_2x_3x_4)^\ell\Omega_0}{f^{\ell+1}}$ for $0\leq\ell$. In detail we get:

\begin{align*}
\delta\omega&=\frac{s\Omega_0}{f}-\frac{s^2x_1x_2x_3x_4\Omega_0}{f^2}\\
\delta^2\omega&=\frac{s\Omega_0}{f}-3\frac{s^2x_1x_2x_3x_4\Omega_0}{f^2}+\frac{2s^3(x_1x_2x_3x_4)^2\Omega_0}{f^3}\\
\delta^3\omega&=\frac{s\Omega_0}{f}-7\frac{s^2x_1x_2x_3x_4\Omega_0}{f^2}+6\frac{2s^3(x_1x_2x_3x_4)^2\Omega_0}{f^3}-\frac{6s^4(x_1x_2x_3x_4)^3\Omega_0}{f^4}\\
\delta^4\omega&=\frac{s\Omega_0}{f}-15\frac{s^2x_1x_2x_3x_4\Omega_0}{f^2}+25\frac{2s^3(x_1x_2x_3x_4)^2\Omega_0}{f^3}-10\frac{6s^4(x_1x_2x_3x_4)^3\Omega_0}{f^4}\\
&\quad +\frac{24s^5(x_1x_2x_3x_4)^4\Omega_0}{f^5}.
\end{align*}

We are looking at the Milnor ring in degree $0,10$ and $20$. Therefore we can choose $\frac{s\Omega_0}{f}, \frac{s^2x_1x_2x_3x_4\Omega_0}{f^2}$ and $\frac{2s^3(x_1x_2x_3x_4)^2\Omega_0}{f^3}$ to be basis elements. We define them as $b_0,b_1$ and $b_2$ respectively. Then the above expression reduces to

\begin{align}\label{formel-delta}
\begin{split}
\omega&=b_0\\
\delta\omega&=b_0-b_1 \\
\delta^2\omega&=b_0-3b_1+b_2\\
\delta^3\omega&=b_0-7b_1+6b_2-\frac{6s^4(x_1x_2x_3x_4)^3\Omega_0}{f^4} \\
\delta^4\omega&=b_0-15b_1+25b_2-10\frac{6s^4(x_1x_2x_3x_4)^3\Omega_0}{f^4}+\frac{24s^5(x_1x_2x_3x_4)^4\Omega_0}{f^5}.
\end{split}
\end{align}

Now we have to figure out how to write $\frac{6s^4(x_1x_2x_3x_4)^3\Omega_0}{f^4}$ and $\frac{24s^5(x_1x_2x_3x_4)^4\Omega_0}{f^5}$ in a basis of the Milnor ring. From Proposition \ref{propo} we already know that we need $4$ basis elements. So the three we had so far are not enough. We will find out what the extra basis element should be in the process of calculating. We want to shorten the notation for numbers that occur throughout the whole calculation.

\begin{notation}
The number $\Delta$ is defined as $\Delta:=\prod \widehat{q}_i^{\widehat{q}_i}s^{\widehat{d}}-(-\widehat{d})^{\widehat{d}}=5^52^4s^{10}-10^{10}$ and will occur as a normalization factor in the calculations. We also want to define $c_q:=\prod \widehat{q}_i^{\widehat{q}_i}=5^52^4$ and $c_d:=(-\widehat{d})^{\widehat{d}}=10^{10}$ separately. So $\Delta=c_qs^{\widehat{d}}-c_d$.
\end{notation}

\subsection*{Calculating $(x_1x_2x_3x_4)^3$}

We start by calculating $\frac{6s^4(x_1x_2x_3x_4)^3\Omega_0}{f^4}$. The first step for this is to write down the Jacobi path with all coefficients. This is done in the following picture. In the last two sections we mostly ignored all coefficients, but now we have to calculate all of them. In the diagram the coefficients are marked in three ways, which is distinguished by three colours. The blue number near each vertex is the coefficient that the corresponding monomial should have after adding everything up. The green number next to each arrow is the coefficient the corresponding partial derivative needs such that after adding up we get the blue numbers as results. In addition in purple we marked the exponents $k_i$ which appear as coefficients inside the partial derivative.  

\begin{figure}[H]
\begin{center}
\begin{tikzpicture}
\node[place] (0) [label=above:\small{\textcolor{cTangoSkyBlue1}{$1$}}] {$3,3,3,3$};
\node[place] (1) [right=60pt of 0,label=above:\small{\textcolor{cTangoSkyBlue1}{$0$}}] {$4,4,4,2$};
\node[place] (2) [right=60pt of 1,label=above:\small{\textcolor{cTangoSkyBlue1}{$0$}}] {$5,5,5,1$};
\node[place] (3) [below=60pt of 2,label=120:\small{\textcolor{cTangoSkyBlue1}{$0$}}] {$1,5,6,2$};
\node[place] (4) [below=60pt of 3,label=120:\small{\textcolor{cTangoSkyBlue1}{$0$}}] {$2,2,6,3$};
\node[place] (5) [below=60pt of 4,label=below:\small{\textcolor{cTangoSkyBlue1}{$0$}}] {$3,3,7,2$};
\node[place] (6) [left=60pt of 5,label=below:\small{\textcolor{cTangoSkyBlue1}{$0$}}] {$4,4,8,1$};
\node[place] (7) [left=60pt of 6,label=below:\small{\textcolor{cTangoSkyBlue1}{$0$}}] {$0,4,9,2$};
\node[place] (8) [above=60pt of 7,label=-60:\small{\textcolor{cTangoSkyBlue1}{$0$}}] {$1,1,9,3$};
\node[place] (9) [above=60pt of 8,label=-60:\small{\textcolor{cTangoSkyBlue1}{$0$}}] {$2,2,2,4$};
\draw[-latex,thick] (0)--(1)
 node[very near start,above] {\small{\textcolor{cTangoPlum1}{$2$}}}
 node[midway,above] {\small{\textcolor{cTangoChameleon1}{$-\frac{10^9 5}{\Delta}$}}};
\draw[-latex,thick] (1)--(2)
 node[very near start,above] {\small{\textcolor{cTangoPlum1}{$2$}}}
 node[midway,above] {\small{\textcolor{cTangoChameleon1}{$\frac{10^8 5^2s}{\Delta}$}}};
\draw[-latex,thick] (2)--(3)
 node[very near start,right] {\small{\textcolor{cTangoPlum1}{$5$}}}
 node[midway,left] {\small{\textcolor{cTangoChameleon1}{$-\frac{10^7 5^2 2s^2}{\Delta}$}}};
\draw[-latex,thick] (3)--(4)
 node[very near start,right] {\small{\textcolor{cTangoPlum1}{$4$}}}
 node[midway,left] {\small{\textcolor{cTangoChameleon1}{$\frac{10^7 5^2 2s^3}{4\Delta}$}}}
 node[point,thick,midway] (h1) {};
\node[place] (p1) [right=40pt of h1] {$6,2,5,2$};
\draw[-,thick] (p1)--(h1);
\draw[-latex, thick] (p1)--(4)
 node[very near start,below] {\small{\textcolor{cTangoPlum1}{$5$}}}
 node[midway,below right] {\small{\textcolor{cTangoChameleon1}{$-\frac{10^7 5^2 2s^3}{4\cdot 5\Delta}$}}};
\draw[-latex,thick] (4)--(5)
 node[very near start,right] {\small{\textcolor{cTangoPlum1}{$2$}}}
 node[midway,left] {\small{\textcolor{cTangoChameleon1}{$-\frac{10^5 5^3 2^2s^4}{\Delta}$}}};
\draw[-latex,thick] (5)--(6)
 node[very near start,below] {\small{\textcolor{cTangoPlum1}{$2$}}}
 node[midway,below] {\small{\textcolor{cTangoChameleon1}{$\frac{10^4 5^4 2^2s^5}{\Delta}$}}};
\draw[-latex,thick] (6)--(7)
 node[very near start,below,xshift=2pt] {\small{\textcolor{cTangoPlum1}{$5$}}}
 node[midway,below] {\small{\textcolor{cTangoChameleon1}{$-\frac{10^3 5^4 2^3s^6}{\Delta}$}}};
\draw[-latex,thick] (7)--(8)
 node[very near start,left] {\small{\textcolor{cTangoPlum1}{$4$}}}
 node[midway,right] {\small{\textcolor{cTangoChameleon1}{$\frac{10^3 5^4 2^3s^7}{4\Delta}$}}}
 node[point,midway,thick] (h2) {};
\node[place] (p2) [left=60pt of h2] {$5,1,8,2$};
\draw[-,thick] (p2)--(h2);
\draw[-latex,thick] (p2)--(8)
 node[very near start,above] {\small{\textcolor{cTangoPlum1}{$5$}}}
 node[midway,above,xshift=-12pt,yshift=2pt] {\small{\textcolor{cTangoChameleon1}{$-\frac{10^3 5^4 2^3s^7}{4\cdot 5\Delta}$}}};
\draw[-latex,thick] (8)--(9)
 node[very near start,left] {\small{\textcolor{cTangoPlum1}{$8$}}}
 node[midway,right] {\small{\textcolor{cTangoChameleon1}{$-\frac{10^2 5^4 2^4s^8}{8\Delta}$}}}
 node[point,midway,thick] (h3) {};
\node[place] (p3) [left=60pt of h3] {$1,5,2,3$};
\draw[-,thick] (p3)--(h3);
\draw[-latex,thick] (p3)--(9)
 node[very near start,above] {\small{\textcolor{cTangoPlum1}{$4$}}}
 node[midway,above left,xshift=-11pt] {\small{\textcolor{cTangoChameleon1}{$\frac{10^2 5^4 2^4s^8}{4\cdot8\Delta}$}}}
 node[point,midway,thick] (h4) {};
\node[place] (p4) [above=50pt of p3] {$6,2,1,3$};
\draw[-,thick] (h4)--(p4);
\draw[-latex,thick] (p4)--(9)
 node[very near start,below] {\small{\textcolor{cTangoPlum1}{$5$}}}
 node[midway,above,xshift=4pt] {\small{\textcolor{cTangoChameleon1}{$-\frac{10^2 5^4 2^4s^8}{5\cdot4\cdot8\Delta}$}}};
\draw[-latex,thick] (9)--(0)
 node[very near start,left] {\small{\textcolor{cTangoPlum1}{$2$}}}
 node[midway,right] {\small{\textcolor{cTangoChameleon1}{$\frac{5^5 2^4s^9}{\Delta}$}}};
\end{tikzpicture}
\end{center}
\caption{The Jacobi path for $(x_1x_2x_3x_4)^3$}
\label{jp-eins}
\end{figure}
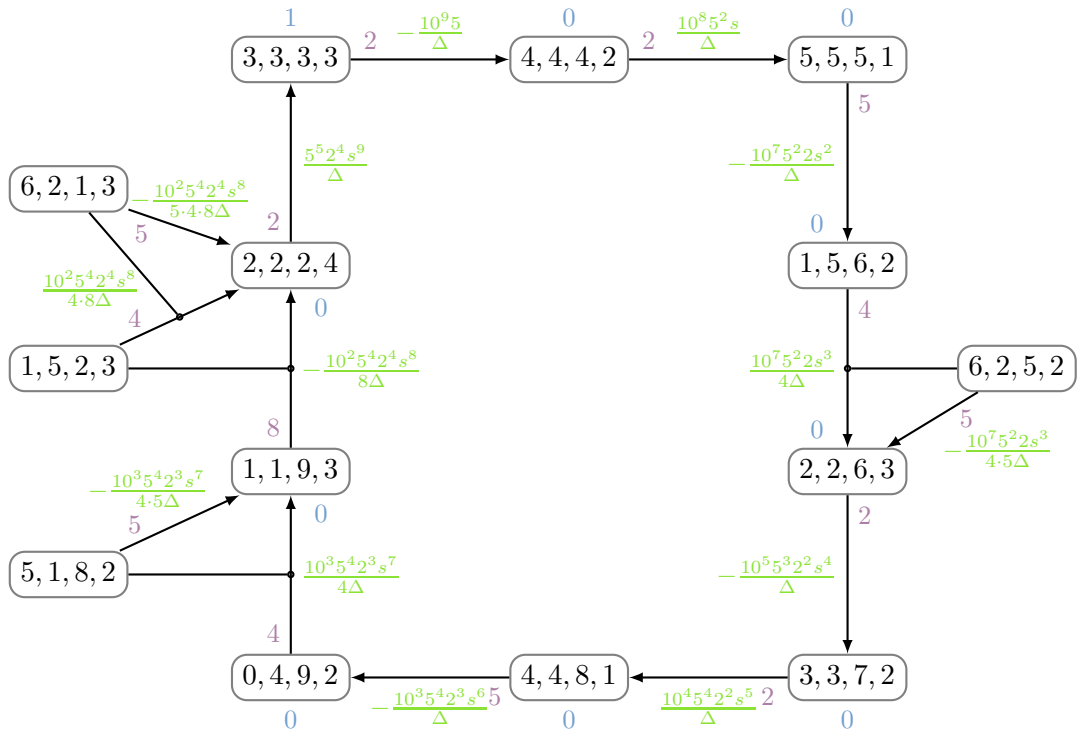

In this first step the goal is to get a description of $(wxyz)^3$. So after adding up all other monomials should vanish. Therefore the blue numbers , i.e. coefficients of the monomials, are all $0$ except the first one. In Figure \ref{jp-eins} the green numbers are relatively easy to find. he last one is always $\frac{c_qs^{\widehat{d}-1}}{\Delta}$ and the others can be calculated inductively. The basic idea is that from one arrow to the next one has to divide by $\widehat{d}$ and multiply by the appropriate $\widehat{q}_i$. One needs to put more effort in doing this in general and we will not do this here. Nevertheless, it is easy to find these coefficients, because they only consist of powers of $\widehat{d}$, $s$ and $\widehat{q}_i$ for $i=1,\dots,n$. 
If one translates the above picture in normal notation it tells us that the form $\frac{6s^4(x_1x_2x_3x_4)^3\Omega_0}{f^4}$ can be written in the following way as a linear combination of the partial derivatives:

\begin{align*}
\frac{6s^4(x_1x_2x_3x_4)^3\Omega_0}{f^4}&=\frac{6s^4\Omega_0}{f^4}\left(-\frac{10^75^22s^2}{\Delta}x_1x_2^4x_3^5x_4-\frac{10^75^22s^3}{4\cdot 5\Delta}x_1^2x_2x_3^5x_4^2-\frac{10^35^42^3s^6}{\Delta}x_2^3x_3^8x_4\right.\\
&\quad -\left.\frac{10^3 5^4 2^3s^7}{4\cdot 5\Delta}x_1x_3^8x_4^2-\frac{10^2 5^4 2^4s^8}{5\cdot4\cdot8\Delta}x_1^2x_2x_3x_4^3\right)\frac{\partial f}{\partial x_1}\\
&\quad +\frac{6s^4\Omega_0}{f^4}\left(\frac{10^7 5^2 2s^3}{4\Delta}x_1x_2^2x_3^5x_4^2+\frac{10^3 5^4 2^3s^7}{4\Delta}x_2x_3^8x_4^2\right.\\
&\quad +\left.\frac{10^2 5^4 2^4s^8}{4\cdot8\Delta}x_1x_2^2x_3x_4^3\right)\frac{\partial f}{\partial x_2}\\
&\quad +\frac{6s^4\Omega_0}{f^4}\left(-\frac{10^2 5^4 2^4s^8}{8\Delta}x_1x_2x_3^2x_4^3\right)\frac{\partial f}{\partial x_3}\\
&\quad +\frac{6s^4\Omega_0}{f^4}\left(-\frac{10^9 5}{\Delta}x_1^3x_2^3x_3^3x_4^2+\frac{10^8 5^2s}{\Delta}x_1^4x_2^4x_3^4x_4^1-\frac{10^5 5^3 2^2s^4}{\Delta}x_1^2x_2^2x_3^6x_4^2\right.\\
&\quad +\left.\frac{10^4 5^4 2^2s^5}{\Delta}x_1^3x_2^3x_3^7x_4^1+\frac{5^52^4s^9}{\Delta}x_1^2x_2^2x_3^2x_4^3\right)\frac{\partial f}{\partial x_4}.
\end{align*}
Now we use the Griffiths formula (\ref{griffiths-formula}). Because we have just written everything in terms of the partial derivatives we can do this directly and get the following result:
\begin{align*}
\frac{6s^4(x_1x_2x_3x_4)^3\Omega_0}{f^4}&=\frac{2s^4\Omega_0}{f^3}\left(-\frac{10^75^22s^2}{\Delta}x_2^4x_3^5x_4-2\frac{10^75^22s^3}{4\cdot 5\Delta}x_1x_2x_3^5x_4^2\right.\\
&\quad -\left.\frac{10^3 5^4 2^3s^7}{4\cdot 5\Delta}x_3^8x_4^2-2\frac{10^2 5^4 2^4s^8}{5\cdot4\cdot8\Delta}x_1x_2x_3x_4^3\right)\\
&\quad +\frac{2s^4\Omega_0}{f^3}\left(2\frac{10^7 5^2 2s^3}{4\Delta}x_1x_2x_3^5x_4^2+\frac{10^3 5^4 2^3s^7}{4\Delta}x_3^8x_4^2\right.\\
&\quad +\left.2\frac{10^2 5^4 2^4s^8}{4\cdot8\Delta}x_1x_2x_3x_4^3\right)\\
&\quad +\frac{2s^4\Omega_0}{f^3}\left(-2\frac{10^2 5^4 2^4s^8}{8\Delta}x_1x_2x_3x_4^3\right)\\
&\quad +\frac{2s^4\Omega_0}{f^3}\left(-2\frac{10^9 5}{\Delta}x_1^3x_2^3x_3^3x_4+\frac{10^8 5^2s}{\Delta}x_1^4x_2^4x_3^4-2\frac{10^5 5^3 2^2s^4}{\Delta}x_1^2x_2^2x_3^6x_4\right.\\
&\quad +\left.\frac{10^4 5^4 2^2s^5}{\Delta}x_1^3x_2^3x_3^7+3\frac{5^52^4s^9}{\Delta}x_1^2x_2^2x_3^2x_4^2\right).
\end{align*}

We can also very easily do this step in our picture and the advantage is that we directly get a decomposition of the result as a linear combination of partial derivatives. We want to stress that in this new picture the coefficients are closely related to the coefficients from before. The blue number, i.e. coefficient of the monomial, in the second picture is the green number next to the arrow pointing to the vertex in the last picture multiplied by a constant, because the Griffiths formula contracts an arrow to the vertex at the arrow tip. The factor one has to multiply is given by the $i$th entry if the arrow pointing to the vertex is the partial derivative with respect to $x_i$. The green numbers in the second picture can now be calculated such that all blue numbers are correct after adding up. This is possible everywhere but at a vertex corresponding to a basis element. At these vertices we will have to add a multiple of the basis element. The factor of the basis element is written down next to the vertex with the comment ``extra''.

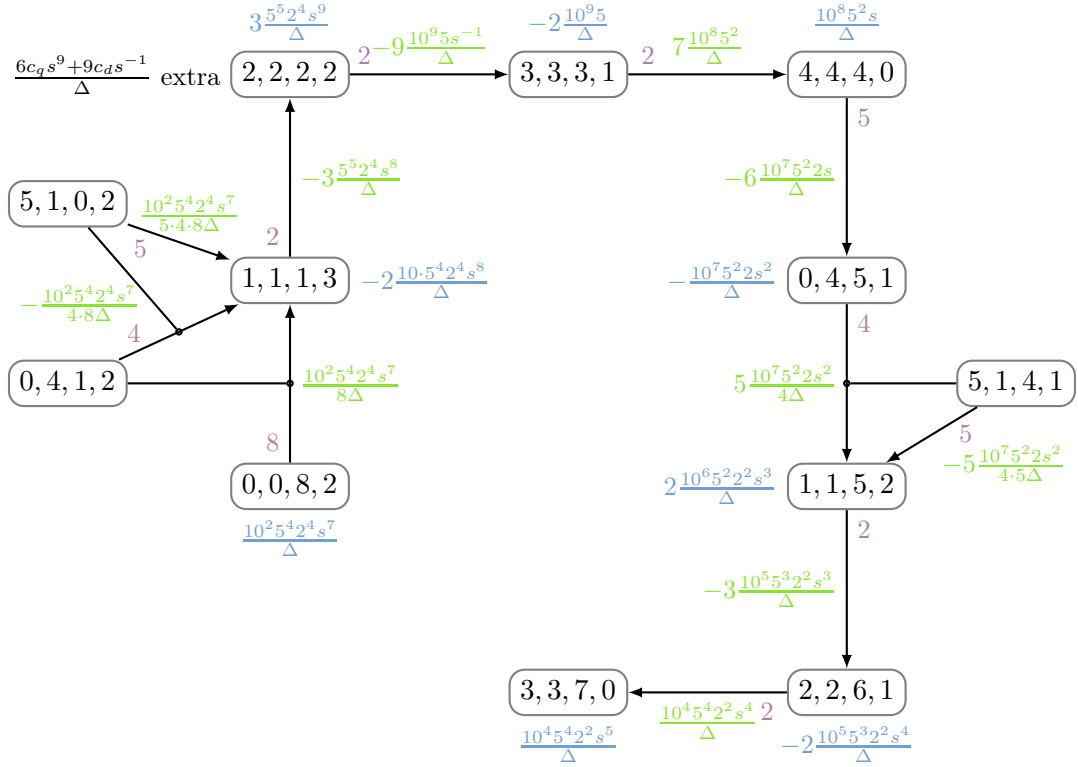
\begin{figure}[H]
\begin{center}
\begin{tikzpicture}
\node[place] (0) [label=above:\small{\textcolor{cTangoSkyBlue1}{$3\frac{5^52^4s^9}{\Delta}$}},label=left:\small{$\frac{6c_qs^9+9c_ds^{-1}}{\Delta}$ extra}] {$2,2,2,2$};
\node[place] (1) [right=60pt of 0,label=above:\small{\textcolor{cTangoSkyBlue1}{$-2\frac{10^9 5}{\Delta}$}}] {$3,3,3,1$};
\node[place] (2) [right=60pt of 1,label=above:\small{\textcolor{cTangoSkyBlue1}{$\frac{10^8 5^2s}{\Delta}$}}] {$4,4,4,0$};
\node[place] (3) [below=60pt of 2,label=left:\small{\textcolor{cTangoSkyBlue1}{$-\frac{10^7 5^2 2s^2}{\Delta}$}}] {$0,4,5,1$};
\node[place] (4) [below=60pt of 3,label=left:\small{\textcolor{cTangoSkyBlue1}{$2\frac{10^6 5^2 2^2s^3}{\Delta}$}}] {$1,1,5,2$};
\node[place] (5) [below=60pt of 4,label=below:\small{\textcolor{cTangoSkyBlue1}{$-2\frac{10^5 5^3 2^2s^4}{\Delta}$}}] {$2,2,6,1$};
\node[place] (6) [left=60pt of 5,label=below:\small{\textcolor{cTangoSkyBlue1}{$\frac{10^4 5^4 2^2s^5}{\Delta}$}}] {$3,3,7,0$};
\node[place] (9) [below=60pt of 0,label=right:\small{\textcolor{cTangoSkyBlue1}{$-2\frac{10\cdot5^4 2^4s^8}{\Delta}$}}] {$1,1,1,3$};
\node[place] (8) [below=60pt of 9,label=below:\small{\textcolor{cTangoSkyBlue1}{$\frac{10^2 5^4 2^4s^7}{\Delta}$}}] {$0,0,8,2$};
\draw[-latex,thick] (0)--(1)
 node[very near start,above,xshift=-2pt] {\small{\textcolor{cTangoPlum1}{$2$}}}
 node[midway,above] {\small{\textcolor{cTangoChameleon1}{$-9\frac{10^95s^{-1}}{\Delta}$}}};
\draw[-latex,thick] (1)--(2)
 node[very near start,above] {\small{\textcolor{cTangoPlum1}{$2$}}}
 node[midway,above] {\small{\textcolor{cTangoChameleon1}{$7\frac{10^85^2}{\Delta}$}}};
\draw[-latex,thick] (2)--(3)
 node[very near start,right] {\small{\textcolor{cTangoPlum1}{$5$}}}
 node[midway,left] {\small{\textcolor{cTangoChameleon1}{$-6\frac{10^75^22s}{\Delta}$}}};
\draw[-latex,thick] (3)--(4)
 node[very near start,right] {\small{\textcolor{cTangoPlum1}{$4$}}}
 node[midway,left] {\small{\textcolor{cTangoChameleon1}{$5\frac{10^75^22s^2}{4\Delta}$}}}
 node[point,thick,midway] (h1) {};
\node[place] (p1) [right=40pt of h1] {$5,1,4,1$};
\draw[-,thick] (p1)--(h1);
\draw[-latex, thick] (p1)--(4)
 node[very near start,below] {\small{\textcolor{cTangoPlum1}{$5$}}}
 node[midway,below right] {\small{\textcolor{cTangoChameleon1}{$-5\frac{10^75^22s^2}{4\cdot 5\Delta}$}}};
\draw[-latex,thick] (4)--(5)
 node[very near start,right] {\small{\textcolor{cTangoPlum1}{$2$}}}
 node[midway,left] {\small{\textcolor{cTangoChameleon1}{$-3\frac{10^55^32^2s^3}{\Delta}$}}};
\draw[-latex,thick] (5)--(6)
 node[very near start,below] {\small{\textcolor{cTangoPlum1}{$2$}}}
 node[midway,below] {\small{\textcolor{cTangoChameleon1}{$\frac{10^45^42^2s^4}{\Delta}$}}};
\draw[-latex,thick] (8)--(9)
 node[very near start,left] {\small{\textcolor{cTangoPlum1}{$8$}}}
 node[midway,right] {\small{\textcolor{cTangoChameleon1}{$\frac{10^2 5^4 2^4s^7}{8\Delta}$}}}
 node[point,midway,thick] (h3) {};
\node[place] (p3) [left=60pt of h3] {$0,4,1,2$};
\draw[-,thick] (p3)--(h3);
\draw[-latex,thick] (p3)--(9)
 node[very near start,above] {\small{\textcolor{cTangoPlum1}{$4$}}}
 node[midway,above left,xshift=-11pt] {\small{\textcolor{cTangoChameleon1}{$-\frac{10^2 5^4 2^4s^7}{4\cdot 8\Delta}$}}}
 node[point,midway,thick] (h4) {};
\node[place] (p4) [above=50pt of p3] {$5,1,0,2$};
\draw[-,thick] (h4)--(p4);
\draw[-latex,thick] (p4)--(9)
 node[very near start,below] {\small{\textcolor{cTangoPlum1}{$5$}}}
 node[midway,above,xshift=4pt] {\small{\textcolor{cTangoChameleon1}{$\frac{10^2 5^4 2^4s^7}{5\cdot 4\cdot 8\Delta}$}}};
\draw[-latex,thick] (9)--(0)
 node[very near start,left] {\small{\textcolor{cTangoPlum1}{$2$}}}
 node[midway,right] {\small{\textcolor{cTangoChameleon1}{$-3\frac{5^52^4s^8}{\Delta}$}}};
\end{tikzpicture}
\end{center}
\caption{The Jacobi path for $(x_1x_2x_3x_4)^3$ after the first use of the Griffiths formula}
\end{figure}

Interpreting the picture we get a description of $\frac{6s^4(x_1x_2x_3x_4)^3\Omega_0}{f^4}$ in terms of the basis in degree $20$, which is only given by $b_2=\frac{2s^3(x_1x_2x_3x_4)^2\Omega_0}{f^3}$, and a linear combination of the partial derivatives:

\begin{align*}
\frac{6s^4(x_1x_2x_3x_4)^3\Omega_0}{f^4}&=\frac{2s^3(x_1x_2x_3x_4)^2\Omega_0}{f^3}\left(\frac{6c_qs^{10}+9c_d}{\Delta}\right)\\
&\quad +\frac{2s^4\Omega_0}{f^3}\left(-6\frac{10^75^22s}{\Delta}x_2^3x_3^4-5\frac{10^75^22s^2}{4\cdot 5\Delta}x_1x_3^4x_4+\frac{10^2 5^4 2^4s^7}{5\cdot 4\cdot 8\Delta}x_1x_4^2\right)\frac{\partial f}{\partial x_1}\\
&\quad +\frac{2s^4\Omega_0}{f^3}\left(5\frac{10^75^22s^2}{4\Delta}x_2x_3^4x_4-\frac{10^2 5^4 2^4s^7}{4\cdot 8\Delta}x_2x_4^2\right)\frac{\partial f}{\partial x_2}\\
&\quad +\frac{2s^4\Omega_0}{f^3}\left(\frac{10^2 5^4 2^4s^7}{8\Delta}x_3x_4^2\right)\frac{\partial f}{\partial x_3}\\
&\quad +\frac{2s^4\Omega_0}{f^3}\left(-9\frac{10^95s^{-1}}{\Delta}x_1^2x_2^2x_3^2x_4+7\frac{10^85^2}{\Delta}x_1^3x_2^3x_3^3\right.\\
&\left.-3\frac{10^55^32^2s^3}{\Delta}x_1x_2x_3^5x_4+\frac{10^45^42^2s^4}{\Delta}x_1^2x_2^2x_3^6-3\frac{5^52^4s^8}{\Delta}x_1x_2x_3x_4^2\right)\frac{\partial f}{\partial x_4}.
\end{align*}

Again with the Griffiths formula we can lower the degree for everything in the Jacobian ideal, so for the linear combination of the partial derivatives. Doing this we end up with the following expression:

\begin{align*}
\frac{6s^4(x_1x_2x_3x_4)^3\Omega_0}{f^4}&=\left(\frac{6c_qs^{10}+9c_d}{\Delta}\right)b_2\\
&\quad +\frac{s^4\Omega_0}{f^2}\left(-5\frac{10^75^22s^2}{4\cdot 5\Delta}x_3^4x_4+\frac{10^2 5^4 2^4s^7}{5\cdot 4\cdot 8\Delta}x_4^2\right)\\
&\quad +\frac{s^4\Omega_0}{f^2}\left(5\frac{10^75^22s^2}{4\Delta}x_3^4x_4-\frac{10^2 5^4 2^4s^7}{4\cdot 8\Delta}x_4^2\right)\\
&\quad +\frac{s^4\Omega_0}{f^2}\left(\frac{10^2 5^4 2^4s^7}{8\Delta}x_4^2\right)\\
&\quad +\frac{s^4\Omega_0}{f^2}\left(-9\frac{10^95s^{-1}}{\Delta}x_1^2x_2^2x_3^2-3\frac{10^55^32^2s^3}{\Delta}x_1x_2x_3^5-6\frac{5^52^4s^8}{\Delta}x_1x_2x_3x_4\right).
\end{align*}

In total there are still $5$ monomials apart from $b_2$ in the formula, but this is not the end, because some of them are linear dependent to each other. This is not always easy to see in the normal notation, but it is very easy to see in our new notation. We can see this if we go back to our picture. After subtracting $\left(\frac{6c_qs^{10}+9c_d}{\Delta}\right)b_2$ we get a polynomial in the Jacobian ideal and we can use the Griffiths formula which leads to the following picture:
\begin{figure}[H]
\begin{center}
\begin{tikzpicture}
\node[place] (0) [label=above:\small{\textcolor{cTangoSkyBlue1}{$-6\frac{5^52^4s^8}{\Delta}$}},label=left:\small{$\frac{-7c_qs^8+9c_ds^{-2}}{\Delta}$ extra}] {$1,1,1,1$};
\node[place] (1) [right=60pt of 0,label=above:\small{\textcolor{cTangoSkyBlue1}{$-9\frac{10^9 5s^{-1}}{\Delta}$}}] {$2,2,2,0$};
\node[place] (9) [below=60pt of 0,label=right:\small{\textcolor{cTangoSkyBlue1}{$\frac{10\cdot5^4 2^4s^7}{\Delta}$}}] {$0,0,0,2$};
\node[place] (4) [below=60pt of 9,xshift=180pt,label=left:\small{\textcolor{cTangoSkyBlue1}{$5\frac{10^65^22^2s^2}{\Delta}$}},label=right:{\small{$8\frac{10^65^22^2s^2}{\Delta}$ extra}}] {$0,0,4,1$};
\node[place] (5) [below=60pt of 4,label=below:\small{\textcolor{cTangoSkyBlue1}{$-3\frac{10^5 5^3 2^2s^3}{\Delta}$}}] {$1,1,5,0$};
\draw[-latex,thick] (0)--(1)
 node[very near start,above,xshift=-2pt] {\small{\textcolor{cTangoPlum1}{$2$}}}
 node[midway,above] {\small{\textcolor{cTangoChameleon1}{$-9\frac{10^95s^{-2}}{\Delta}$}}};
\draw[-latex,thick] (4)--(5)
 node[very near start,right] {\small{\textcolor{cTangoPlum1}{$2$}}}
 node[midway,left] {\small{\textcolor{cTangoChameleon1}{$-3\frac{10^55^32^2s^2}{\Delta}$}}};
\draw[-latex,thick] (9)--(0)
 node[very near start,left] {\small{\textcolor{cTangoPlum1}{$2$}}}
 node[midway,right] {\small{\textcolor{cTangoChameleon1}{$\frac{5^52^4s^7}{\Delta}$}}};
\end{tikzpicture}
\end{center}
\caption{The Jacobi path for $(x_1x_2x_3x_4)^3$ after using the Griffiths formula twice}
\label{jptwice}
\end{figure}
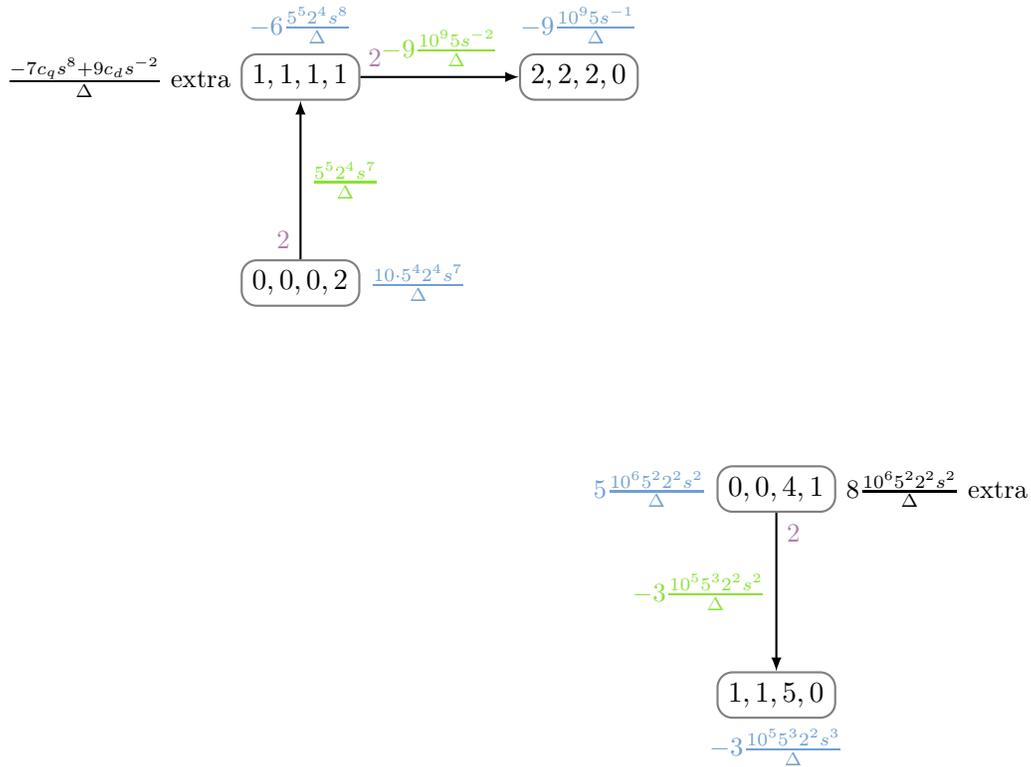

Before we do further calculations, we define a fourth basis element $b_3=\frac{10^65^22^2s^6x_3^4x_4\Omega_0}{f^2}$. We know from Theorem \ref{thm} that we need a fourth basis element and from the picture above, we can see this is a good way of choosing $b_3$, because we can immediately see the coefficients of this element in the picture. Choosing a good basis is another reason why our construction helps making the computations faster. If one has bigger examples and uses a computer algebra system to compute the coefficients, knowing a good basis makes it more efficient. So in total we have one basis element in degree $20$, which is $b_2$, two basis elements in degree $10$ which are $b_1$ and $b_3$ and we will get one basis element in degree $0$ denoted by $b_0$. Now translating back gives

\begin{align*}
\frac{6s^4(x_1x_2x_3x_4)^3\Omega_0}{f^4}&=\left(\frac{6c_qs^{10}+9c_d}{\Delta}\right)b_2\\
&\quad +\left(\frac{-7c_qs^{10}+9c_d}{\Delta}\right)b_1+\frac{8}{\Delta}b_3\\
&\quad +\frac{s^4\Omega_0}{f}\left(-9\frac{10^95s^{-2}}{\Delta}x_1x_2x_3-3\frac{10^55^32^2s^2}{\Delta}x_3^4+\frac{5^52^4s^7}{\Delta}x_4\right)\frac{\partial f}{\partial x_4}.
\end{align*}

It is very easy to use the Griffiths formula here, because there is only a multiple of the partial derivative with respect to $x_4$ left and because two of the terms vanish when one takes the partial derivative, we end up with a rather short expression for $\frac{6s^4(x_1x_2x_3x_4)^3\Omega_0}{f^4}$:

\begin{align}\label{formel-3,3,3,3}
\begin{split}
\frac{6s^4(x_1x_2x_3x_4)^3\Omega_0}{f^4}&=\left(\frac{6c_qs^{10}+9c_d}{\Delta}\right)b_2+\left(\frac{-7c_qs^{10}+9c_d}{\Delta}\right)b_1\\
&\quad +\frac{8}{\Delta}b_3+\left(\frac{c_qs^{10}}{\Delta}\right)b_0.
\end{split}
\end{align}

Of course the same happens in the picture. Every vertex except $(1,1,1,1)$ in Figure \ref{jptwice} has a zero entry and therefore vanishes after the use of the Griffiths formula. The only remaining vertex is $(0,0,0,0)$ and the coefficient is just the same as the coefficient next to the arrow pointing at $(1,1,1,1)$ in Figure \ref{jptwice}. So the picture is given by

\begin{figure}[H]
\begin{center}
\begin{tikzpicture}
\node[place] (0) [label=above:\small{\textcolor{cTangoSkyBlue1}{$\frac{c_qs^{7}}{\Delta}$}}] {$0,0,0,0$};
\end{tikzpicture}
\end{center}
\caption{The Jacobi path for $(x_1x_2x_3x_4)^3$ after using the Griffiths formula three times}
\end{figure}
It is obviously not necessary to translate back in between the pictures and formulas the whole time. We can do the whole calculations in the pictures and it is less work to draw all pictures first and after that translate back to the formulas. We will do this in the next step, where \pagebreak we calculate the expression in the basis elements for $\frac{24s^5(x_1x_2x_3x_4)^4\Omega_0}{f^5}$. We will also see that we can use some of the calculations we have already done for computing a linear combination of $\frac{6s^5(x_1x_2x_3x_4)^3\Omega_0}{f^4}$ in a basis of the Milnor ring for the case of $\frac{24s^5(x_1x_2x_3x_4)^4\Omega_0}{f^5}$. Especially we already chose all basis elements we need and now the goal is to find the correct coefficients here.

\subsection*{Calculating $(x_1x_2x_3x_4)^4$}

Again we start with writing down the Jacobi path, because $(x_1x_2x_3x_4)^4\in J(f)$ and therefore this is possible. This can be done very easy, because we just have to add $(1,1,1,1)$ to the Jacobi path in Figure \ref{jp-eins}. In addition the coefficients are all the same as in Figure \ref{jp-eins}. This does not hold for the rest of the pictures, because the factors we have to multiply to the coefficients after using the Griffiths formula are bigger as in the earlier pictures, because the entries in the vertices and therefore the exponents of the monomials are bigger. However, the coefficients will not be to far from each other. From Figure \ref{jp-eins} we get that the Jacobi path for $(4,4,4,4)$ is given by

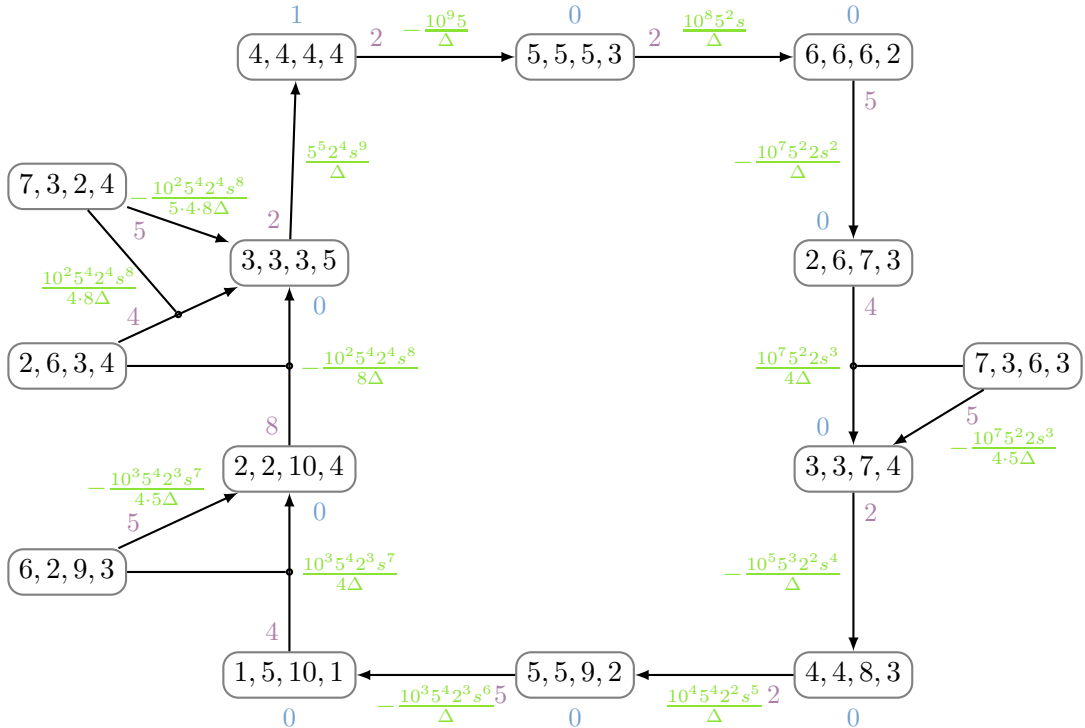
\begin{figure}[H]
\begin{center}
\begin{tikzpicture}
\node[place] (0) [label=above:\small{\textcolor{cTangoSkyBlue1}{$1$}}] {$4,4,4,4$};
\node[place] (1) [right=60pt of 0,label=above:\small{\textcolor{cTangoSkyBlue1}{$0$}}] {$5,5,5,3$};
\node[place] (2) [right=60pt of 1,label=above:\small{\textcolor{cTangoSkyBlue1}{$0$}}] {$6,6,6,2$};
\node[place] (3) [below=60pt of 2,label=120:\small{\textcolor{cTangoSkyBlue1}{$0$}}] {$2,6,7,3$};
\node[place] (4) [below=60pt of 3,label=120:\small{\textcolor{cTangoSkyBlue1}{$0$}}] {$3,3,7,4$};
\node[place] (5) [below=60pt of 4,label=below:\small{\textcolor{cTangoSkyBlue1}{$0$}}] {$4,4,8,3$};
\node[place] (6) [left=60pt of 5,label=below:\small{\textcolor{cTangoSkyBlue1}{$0$}}] {$5,5,9,2$};
\node[place] (7) [left=60pt of 6,label=below:\small{\textcolor{cTangoSkyBlue1}{$0$}}] {$1,5,10,1$};
\node[place] (8) [above=60pt of 7,label=-60:\small{\textcolor{cTangoSkyBlue1}{$0$}}] {$2,2,10,4$};
\node[place] (9) [above=60pt of 8,label=-60:\small{\textcolor{cTangoSkyBlue1}{$0$}}] {$3,3,3,5$};
\draw[-latex,thick] (0)--(1)
 node[very near start,above] {\small{\textcolor{cTangoPlum1}{$2$}}}
 node[midway,above] {\small{\textcolor{cTangoChameleon1}{$-\frac{10^9 5}{\Delta}$}}};
\draw[-latex,thick] (1)--(2)
 node[very near start,above] {\small{\textcolor{cTangoPlum1}{$2$}}}
 node[midway,above] {\small{\textcolor{cTangoChameleon1}{$\frac{10^8 5^2s}{\Delta}$}}};
\draw[-latex,thick] (2)--(3)
 node[very near start,right] {\small{\textcolor{cTangoPlum1}{$5$}}}
 node[midway,left] {\small{\textcolor{cTangoChameleon1}{$-\frac{10^7 5^2 2s^2}{\Delta}$}}};
\draw[-latex,thick] (3)--(4)
 node[very near start,right] {\small{\textcolor{cTangoPlum1}{$4$}}}
 node[midway,left] {\small{\textcolor{cTangoChameleon1}{$\frac{10^7 5^2 2s^3}{4\Delta}$}}}
 node[point,thick,midway] (h1) {};
\node[place] (p1) [right=40pt of h1] {$7,3,6,3$};
\draw[-,thick] (p1)--(h1);
\draw[-latex, thick] (p1)--(4)
 node[very near start,below] {\small{\textcolor{cTangoPlum1}{$5$}}}
 node[midway,below right] {\small{\textcolor{cTangoChameleon1}{$-\frac{10^7 5^2 2s^3}{4\cdot 5\Delta}$}}};
\draw[-latex,thick] (4)--(5)
 node[very near start,right] {\small{\textcolor{cTangoPlum1}{$2$}}}
 node[midway,left] {\small{\textcolor{cTangoChameleon1}{$-\frac{10^5 5^3 2^2s^4}{\Delta}$}}};
\draw[-latex,thick] (5)--(6)
 node[very near start,below] {\small{\textcolor{cTangoPlum1}{$2$}}}
 node[midway,below] {\small{\textcolor{cTangoChameleon1}{$\frac{10^4 5^4 2^2s^5}{\Delta}$}}};
\draw[-latex,thick] (6)--(7)
 node[very near start,below,xshift=2pt] {\small{\textcolor{cTangoPlum1}{$5$}}}
 node[midway,below] {\small{\textcolor{cTangoChameleon1}{$-\frac{10^3 5^4 2^3s^6}{\Delta}$}}};
\draw[-latex,thick] (7)--(8)
 node[very near start,left] {\small{\textcolor{cTangoPlum1}{$4$}}}
 node[midway,right] {\small{\textcolor{cTangoChameleon1}{$\frac{10^3 5^4 2^3s^7}{4\Delta}$}}}
 node[point,midway,thick] (h2) {};
\node[place] (p2) [left=60pt of h2] {$6,2,9,3$};
\draw[-,thick] (p2)--(h2);
\draw[-latex,thick] (p2)--(8)
 node[very near start,above] {\small{\textcolor{cTangoPlum1}{$5$}}}
 node[midway,above,xshift=-12pt,yshift=2pt] {\small{\textcolor{cTangoChameleon1}{$-\frac{10^3 5^4 2^3s^7}{4\cdot 5\Delta}$}}};
\draw[-latex,thick] (8)--(9)
 node[very near start,left] {\small{\textcolor{cTangoPlum1}{$8$}}}
 node[midway,right] {\small{\textcolor{cTangoChameleon1}{$-\frac{10^2 5^4 2^4s^8}{8\Delta}$}}}
 node[point,midway,thick] (h3) {};
\node[place] (p3) [left=60pt of h3] {$2,6,3,4$};
\draw[-,thick] (p3)--(h3);
\draw[-latex,thick] (p3)--(9)
 node[very near start,above] {\small{\textcolor{cTangoPlum1}{$4$}}}
 node[midway,above left,xshift=-11pt] {\small{\textcolor{cTangoChameleon1}{$\frac{10^2 5^4 2^4s^8}{4\cdot8\Delta}$}}}
 node[point,midway,thick] (h4) {};
\node[place] (p4) [above=50pt of p3] {$7,3,2,4$};
\draw[-,thick] (h4)--(p4);
\draw[-latex,thick] (p4)--(9)
 node[very near start,below] {\small{\textcolor{cTangoPlum1}{$5$}}}
 node[midway,above,xshift=4pt] {\small{\textcolor{cTangoChameleon1}{$-\frac{10^2 5^4 2^4s^8}{5\cdot4\cdot8\Delta}$}}};
\draw[-latex,thick] (9)--(0)
 node[very near start,left] {\small{\textcolor{cTangoPlum1}{$2$}}}
 node[midway,right] {\small{\textcolor{cTangoChameleon1}{$\frac{5^5 2^4s^9}{\Delta}$}}};
\end{tikzpicture}
\end{center}
\caption{The Jacobi path for $(x_1x_2x_3x_4)^4$}
\end{figure}

Now we use the Griffiths formula for the first time. The picture looks very similar to the one in the first case before we used the Griffiths formula. The vertices are exactly like they were before, but of course the coefficients are different. But one should notice that throughout all calculations the denominator of the coefficients is always the same, i.e. $\Delta$. The reason is that the denominator $\Delta$ always appears in the first step of writing $(wxyz)^i$ as a normalising factor and gets carried on afterwards. We will see below that this gives rise to the fact that $\Delta$ is also the leading coefficient of the Picard-Fuchs equation.

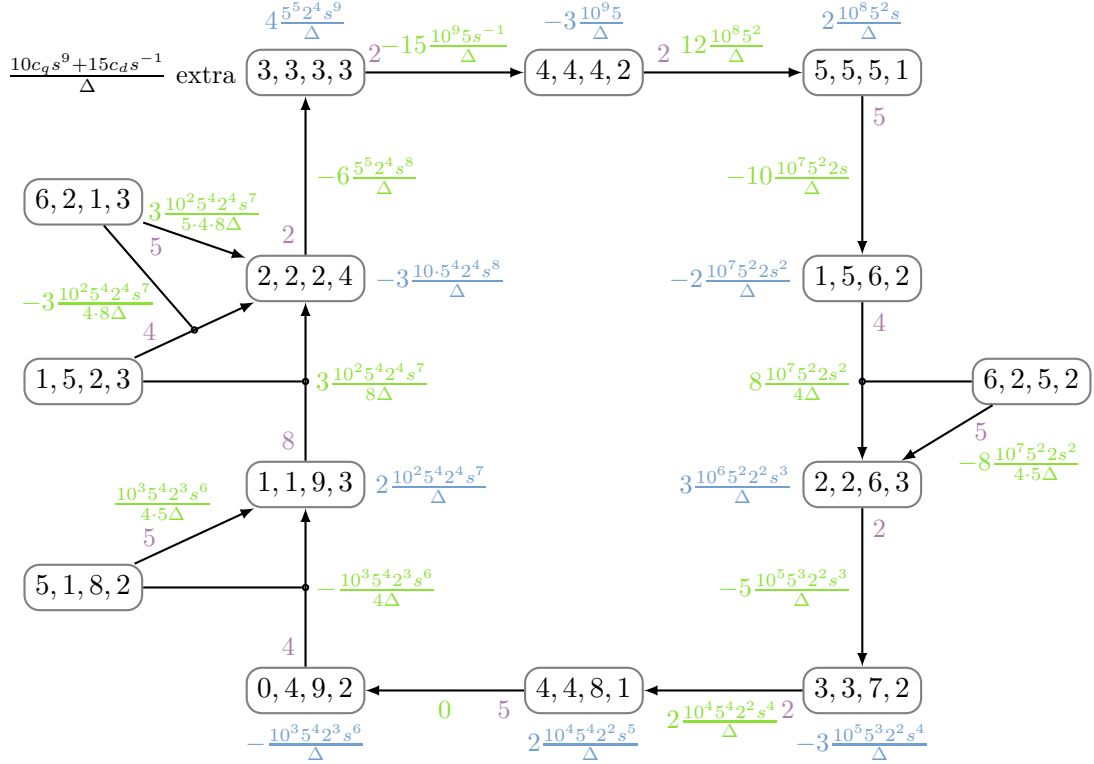
\begin{figure}[H]
\begin{center}
\begin{tikzpicture}
\node[place] (0) [label=above:\small{\textcolor{cTangoSkyBlue1}{$4\frac{5^5 2^4s^9}{\Delta}$}},label=left:\small{$\frac{10c_qs^9+15c_ds^{-1}}{\Delta}$ extra}] {$3,3,3,3$};
\node[place] (1) [right=60pt of 0,label=above:\small{\textcolor{cTangoSkyBlue1}{$-3\frac{10^9 5}{\Delta}$}}] {$4,4,4,2$};
\node[place] (2) [right=60pt of 1,label=above:\small{\textcolor{cTangoSkyBlue1}{$2\frac{10^8 5^2s}{\Delta}$}}] {$5,5,5,1$};
\node[place] (3) [below=60pt of 2,label=left:\small{\textcolor{cTangoSkyBlue1}{$-2\frac{10^7 5^2 2s^2}{\Delta}$}}] {$1,5,6,2$};
\node[place] (4) [below=60pt of 3,label=left:\small{\textcolor{cTangoSkyBlue1}{$3\frac{10^6 5^2 2^2s^3}{\Delta}$}}] {$2,2,6,3$};
\node[place] (5) [below=60pt of 4,label=below:\small{\textcolor{cTangoSkyBlue1}{$-3\frac{10^5 5^3 2^2s^4}{\Delta}$}}] {$3,3,7,2$};
\node[place] (6) [left=60pt of 5,label=below:\small{\textcolor{cTangoSkyBlue1}{$2\frac{10^4 5^4 2^2s^5}{\Delta}$}}] {$4,4,8,1$};
\node[place] (7) [left=60pt of 6,label=below:\small{\textcolor{cTangoSkyBlue1}{$-\frac{10^3 5^4 2^3s^6}{\Delta}$}}] {$0,4,9,2$};
\node[place] (8) [above=60pt of 7,label=right:\small{\textcolor{cTangoSkyBlue1}{$2\frac{10^2 5^4 2^4s^7}{\Delta}$}}] {$1,1,9,3$};
\node[place] (9) [above=60pt of 8,label=right:\small{\textcolor{cTangoSkyBlue1}{$-3\frac{10\cdot 5^4 2^4s^8}{\Delta}$}}] {$2,2,2,4$};
\draw[-latex,thick] (0)--(1)
 node[very near start,above,xshift=-4pt] {\small{\textcolor{cTangoPlum1}{$2$}}}
 node[midway,above] {\small{\textcolor{cTangoChameleon1}{$-15\frac{10^9 5s^{-1}}{\Delta}$}}};
\draw[-latex,thick] (1)--(2)
 node[very near start,above] {\small{\textcolor{cTangoPlum1}{$2$}}}
 node[midway,above] {\small{\textcolor{cTangoChameleon1}{$12\frac{10^8 5^2}{\Delta}$}}};
\draw[-latex,thick] (2)--(3)
 node[very near start,right] {\small{\textcolor{cTangoPlum1}{$5$}}}
 node[midway,left] {\small{\textcolor{cTangoChameleon1}{$-10\frac{10^7 5^2 2s}{\Delta}$}}};
\draw[-latex,thick] (3)--(4)
 node[very near start,right] {\small{\textcolor{cTangoPlum1}{$4$}}}
 node[midway,left] {\small{\textcolor{cTangoChameleon1}{$8\frac{10^7 5^2 2s^2}{4\Delta}$}}}
 node[point,thick,midway] (h1) {};
\node[place] (p1) [right=40pt of h1] {$6,2,5,2$};
\draw[-,thick] (p1)--(h1);
\draw[-latex, thick] (p1)--(4)
 node[very near start,below] {\small{\textcolor{cTangoPlum1}{$5$}}}
 node[midway,below right] {\small{\textcolor{cTangoChameleon1}{$-8\frac{10^7 5^2 2s^2}{4\cdot 5\Delta}$}}};
\draw[-latex,thick] (4)--(5)
 node[very near start,right] {\small{\textcolor{cTangoPlum1}{$2$}}}
 node[midway,left] {\small{\textcolor{cTangoChameleon1}{$-5\frac{10^5 5^3 2^2s^3}{\Delta}$}}};
\draw[-latex,thick] (5)--(6)
 node[very near start,below,xshift=2pt] {\small{\textcolor{cTangoPlum1}{$2$}}}
 node[midway,below] {\small{\textcolor{cTangoChameleon1}{$2\frac{10^4 5^4 2^2s^4}{\Delta}$}}};
\draw[-latex,thick] (6)--(7)
 node[very near start,below] {\small{\textcolor{cTangoPlum1}{$5$}}}
 node[midway,below] {\small{\textcolor{cTangoChameleon1}{$0$}}};
\draw[-latex,thick] (7)--(8)
 node[very near start,left] {\small{\textcolor{cTangoPlum1}{$4$}}}
 node[midway,right] {\small{\textcolor{cTangoChameleon1}{$-\frac{10^3 5^4 2^3s^6}{4\Delta}$}}}
 node[point,midway,thick] (h2) {};
\node[place] (p2) [left=60pt of h2] {$5,1,8,2$};
\draw[-,thick] (p2)--(h2);
\draw[-latex,thick] (p2)--(8)
 node[very near start,above] {\small{\textcolor{cTangoPlum1}{$5$}}}
 node[midway,above,xshift=-12pt,yshift=2pt] {\small{\textcolor{cTangoChameleon1}{$\frac{10^3 5^4 2^3s^6}{4\cdot 5\Delta}$}}};
\draw[-latex,thick] (8)--(9)
 node[very near start,left] {\small{\textcolor{cTangoPlum1}{$8$}}}
 node[midway,right] {\small{\textcolor{cTangoChameleon1}{$3\frac{10^2 5^4 2^4s^7}{8\Delta}$}}}
 node[point,midway,thick] (h3) {};
\node[place] (p3) [left=60pt of h3] {$1,5,2,3$};
\draw[-,thick] (p3)--(h3);
\draw[-latex,thick] (p3)--(9)
 node[very near start,above] {\small{\textcolor{cTangoPlum1}{$4$}}}
 node[midway,above left,xshift=-11pt] {\small{\textcolor{cTangoChameleon1}{$-3\frac{10^2 5^4 2^4s^7}{4\cdot8\Delta}$}}}
 node[point,midway,thick] (h4) {};
\node[place] (p4) [above=50pt of p3] {$6,2,1,3$};
\draw[-,thick] (h4)--(p4);
\draw[-latex,thick] (p4)--(9)
 node[very near start,below] {\small{\textcolor{cTangoPlum1}{$5$}}}
 node[midway,above,xshift=4pt] {\small{\textcolor{cTangoChameleon1}{$3\frac{10^2 5^4 2^4s^7}{5\cdot4\cdot8\Delta}$}}};
\draw[-latex,thick] (9)--(0)
 node[very near start,left] {\small{\textcolor{cTangoPlum1}{$2$}}}
 node[midway,right] {\small{\textcolor{cTangoChameleon1}{$-6\frac{5^5 2^4s^8}{\Delta}$}}};
\end{tikzpicture}
\end{center}
\caption{The Jacobi path for $(x_1x_2x_3x_4)^4$ after the first use of the Griffiths formula}
\end{figure}

Here we have to take $\frac{10c_qs^{10}+15 c_d}{\Delta}$ of the form $\frac{6s^4(x_1x_2x_3x_4)^3\Omega_0}{f^4}$ extra in order to have an expression in the Jacobian ideal. The form $\frac{6s^4(x_1x_2x_3x_4)^3\Omega_0}{f^4}$ itself can be written with the monomials appearing on the path as we already calculated, so we could adjust the coefficients in the picture by adding $\frac{10c_qs^{10}+15 c_d}{\Delta}$ times the coefficients from Figure \ref{jp-eins}, but then the picture gets much bigger and this is not necessary because we already know how to write $\frac{6s^4(x_1x_2x_3x_4)^3\Omega_0}{f^4}$ in the basis. So there is no need to put this information in the picture and do the calculations again. It is enough to remember the coefficient $\frac{10c_qs^{10}+15 c_d}{\Delta}$ and add $\frac{10c_qs^{10}+15 c_d}{\Delta}$ times formula (\ref{formel-3,3,3,3}) to the description of $\frac{24s^4(x_1x_2x_3x_4)^4\Omega_0}{f^5}$ in the end. We want to mention here that in the description of $\frac{6s^4(x_1x_2x_3x_4)^3\Omega_0}{f^4}$ in formula (\ref{formel-3,3,3,3}) all coefficients already have $\Delta$ as denominator, so after multiplying with $\frac{10c_qs^{10}+15 c_d}{\Delta}$ we have $\Delta^2$ as denominator. This is true in general: In the description of $\frac{\ell!s^{\ell+1}(x_1x_2x_3x_4)^\ell\Omega_0}{f^{\ell+1}}$ the denominator $\Delta^{\ell-1}$ will appear and $\ell-1$ is the biggest exponent that appears. This means that we have at least to multiply everything by $\Delta$ to get a relation between the partial derivatives. This explains why $\Delta$ is the leading coefficient of the Picard-Fuchs equation.\\
Now we can use the Griffiths formula for the above picture and the appropriate coefficient at the vertex $(3,3,3,3)$ for the second time. Again the monomials that occur are the same we had for the calculations of $(wxyz)^3$, but the coefficients are bigger. The Jacobi path we get after using the Griffiths formula two times is the following:

\begin{figure}[H]
\begin{center}
\begin{tikzpicture}
\node[place] (0) [label=above:\small{\textcolor{cTangoSkyBlue1}{$-18\frac{5^52^4s^8}{\Delta}$}},label=left:\small{$\frac{-25c_qs^8+80c_ds^{-2}}{\Delta}$ extra}] {$2,2,2,2$};
\node[place] (1) [right=60pt of 0,label=above:\small{\textcolor{cTangoSkyBlue1}{$-30\frac{10^9 5s^{-1}}{\Delta}$}}] {$3,3,3,1$};
\node[place] (2) [right=60pt of 1,label=above:\small{\textcolor{cTangoSkyBlue1}{$12\frac{10^8 5^2}{\Delta}$}}] {$4,4,4,0$};
\node[place] (3) [below=60pt of 2,label=left:\small{\textcolor{cTangoSkyBlue1}{$-10\frac{10^7 5^2 2s}{\Delta}$}}] {$0,4,5,1$};
\node[place] (4) [below=60pt of 3,label=left:\small{\textcolor{cTangoSkyBlue1}{$16\frac{10^6 5^2 2^2s^2}{\Delta}$}}] {$1,1,5,2$};
\node[place] (5) [below=60pt of 4,label=below:\small{\textcolor{cTangoSkyBlue1}{$-10\frac{10^5 5^3 2^2s^3}{\Delta}$}}] {$2,2,6,1$};
\node[place] (6) [left=60pt of 5,label=below:\small{\textcolor{cTangoSkyBlue1}{$2\frac{10^4 5^4 2^2s^4}{\Delta}$}}] {$3,3,7,0$};
\node[place] (9) [below=60pt of 0,label=right:\small{\textcolor{cTangoSkyBlue1}{$6\frac{10\cdot5^4 2^4s^7}{\Delta}$}}] {$1,1,1,3$};
\node[place] (8) [below=60pt of 9,label=below:\small{\textcolor{cTangoSkyBlue1}{$-\frac{10^2 5^4 2^4s^6}{\Delta}$}}] {$0,0,8,2$};
\draw[-latex,thick] (0)--(1)
 node[very near start,above,xshift=-4pt] {\small{\textcolor{cTangoPlum1}{$2$}}}
 node[midway,above] {\small{\textcolor{cTangoChameleon1}{$-80\frac{10^95s^{-2}}{\Delta}$}}};
\draw[-latex,thick] (1)--(2)
 node[very near start,above,xshift=-3pt] {\small{\textcolor{cTangoPlum1}{$2$}}}
 node[midway,above] {\small{\textcolor{cTangoChameleon1}{$50\frac{10^85^2s^{-1}}{\Delta}$}}};
\draw[-latex,thick] (2)--(3)
 node[very near start,right] {\small{\textcolor{cTangoPlum1}{$5$}}}
 node[midway,left] {\small{\textcolor{cTangoChameleon1}{$-38\frac{10^75^22}{\Delta}$}}};
\draw[-latex,thick] (3)--(4)
 node[very near start,right] {\small{\textcolor{cTangoPlum1}{$4$}}}
 node[midway,left] {\small{\textcolor{cTangoChameleon1}{$28\frac{10^75^22s}{4\Delta}$}}}
 node[point,thick,midway] (h1) {};
\node[place] (p1) [right=40pt of h1] {$5,1,4,1$};
\draw[-,thick] (p1)--(h1);
\draw[-latex, thick] (p1)--(4)
 node[very near start,below] {\small{\textcolor{cTangoPlum1}{$5$}}}
 node[midway,below right] {\small{\textcolor{cTangoChameleon1}{$-28\frac{10^75^22s}{4\cdot 5\Delta}$}}};
\draw[-latex,thick] (4)--(5)
 node[very near start,right] {\small{\textcolor{cTangoPlum1}{$2$}}}
 node[midway,left] {\small{\textcolor{cTangoChameleon1}{$-12\frac{10^55^32^2s^2}{\Delta}$}}};
\draw[-latex,thick] (5)--(6)
 node[very near start,below,xshift=2pt] {\small{\textcolor{cTangoPlum1}{$2$}}}
 node[midway,below] {\small{\textcolor{cTangoChameleon1}{$2\frac{10^45^42^2s^3}{\Delta}$}}};
\draw[-latex,thick] (8)--(9)
 node[very near start,left] {\small{\textcolor{cTangoPlum1}{$8$}}}
 node[midway,right] {\small{\textcolor{cTangoChameleon1}{$-\frac{10^2 5^4 2^4s^6}{8\Delta}$}}}
 node[point,midway,thick] (h3) {};
\node[place] (p3) [left=60pt of h3] {$0,4,1,2$};
\draw[-,thick] (p3)--(h3);
\draw[-latex,thick] (p3)--(9)
 node[very near start,above] {\small{\textcolor{cTangoPlum1}{$4$}}}
 node[midway,above left,xshift=-11pt] {\small{\textcolor{cTangoChameleon1}{$\frac{10^2 5^4 2^4s^6}{4\cdot 8\Delta}$}}}
 node[point,midway,thick] (h4) {};
\node[place] (p4) [above=50pt of p3] {$5,1,0,2$};
\draw[-,thick] (h4)--(p4);
\draw[-latex,thick] (p4)--(9)
 node[very near start,below] {\small{\textcolor{cTangoPlum1}{$5$}}}
 node[midway,above,xshift=4pt] {\small{\textcolor{cTangoChameleon1}{$-\frac{10^2 5^4 2^4s^6}{5\cdot 4\cdot 8\Delta}$}}};
\draw[-latex,thick] (9)--(0)
 node[very near start,left] {\small{\textcolor{cTangoPlum1}{$2$}}}
 node[midway,right] {\small{\textcolor{cTangoChameleon1}{$7\frac{5^52^4s^7}{\Delta}$}}};
\end{tikzpicture}
\end{center}
\caption{The Jacobi path for $(x_1x_2x_3x_4)^4$ after the second use of the Griffiths formula}
\end{figure}
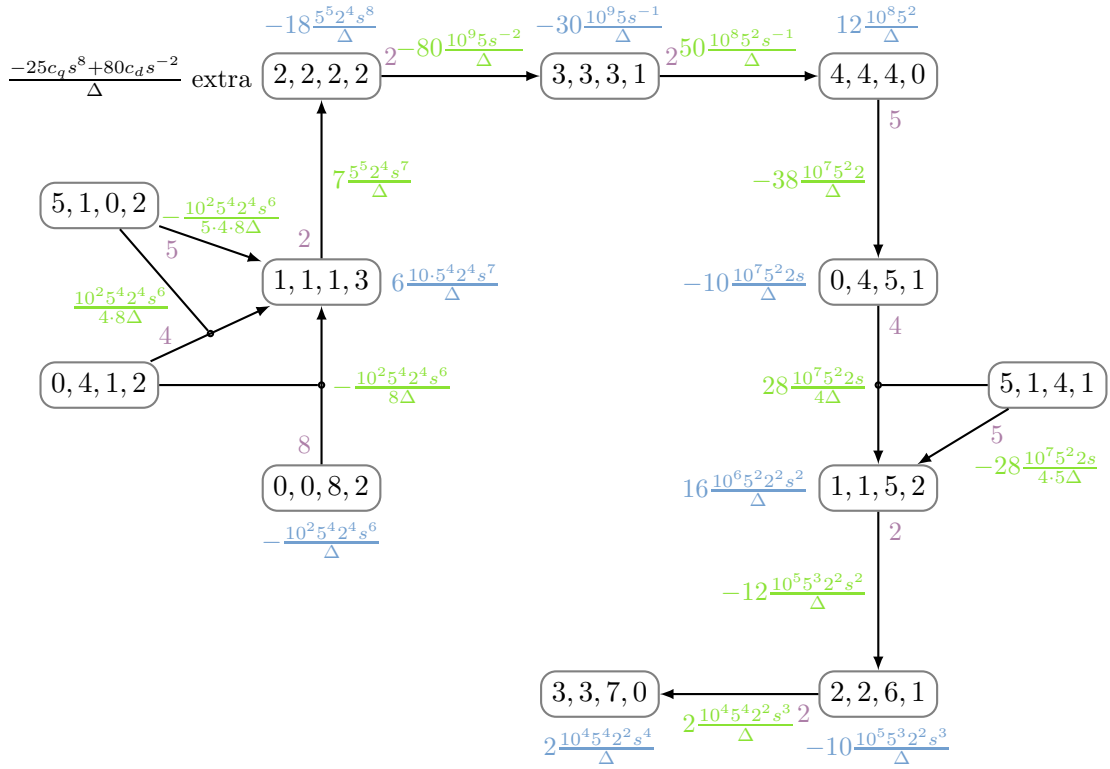

As expected we have one gap in this picture after the $6$th vertex. So we have to use the basis element $b_2=\frac{2s^3(x_1x_2x_3x_4)^2\Omega_0}{f^3}$ here. If we add $\frac{25c_qs^8-80c_ds^{-2}}{\Delta}$ of this basis element to the above picture, we are ending up with an expression in the Jacobian ideal and we are able to use the Griffiths formula again. There are only two steps until we have everything to write down the linear combination of $\frac{24s^5(x_1x_2x_3x_4)^4\Omega_0}{f^5}$ in the basis and because we already know what the basis is, we can concentrate on the coefficients at the corresponding vertices. Especially in the last step we only have to figure out the coefficient at the vertex $(0,0,0,0)$. Now we will show the pictures after the third and fourth use of the Griffiths formula:

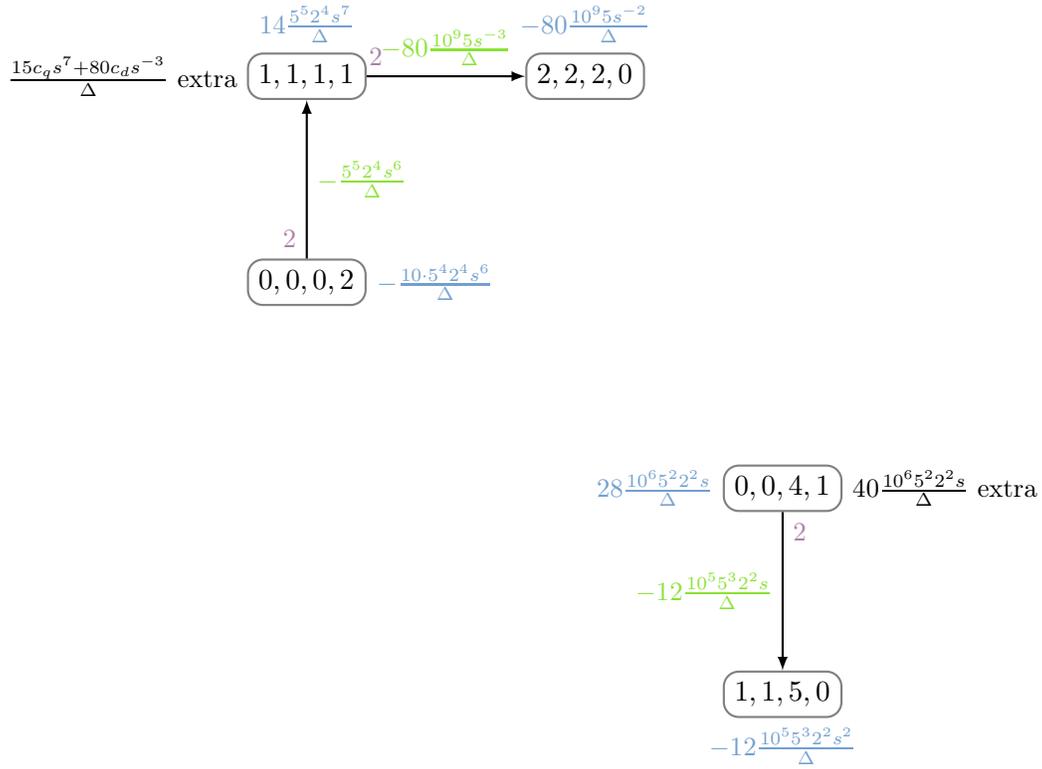
\begin{figure}[H]
\begin{center}
\begin{tikzpicture}
\node[place] (0) [label=above:\small{\textcolor{cTangoSkyBlue1}{$14\frac{5^52^4s^7}{\Delta}$}},label=left:\small{$\frac{15c_qs^7+80c_ds^{-3}}{\Delta}$ extra}] {$1,1,1,1$};
\node[place] (1) [right=60pt of 0,label=above:\small{\textcolor{cTangoSkyBlue1}{$-80\frac{10^9 5s^{-2}}{\Delta}$}}] {$2,2,2,0$};
\node[place] (9) [below=60pt of 0,label=right:\small{\textcolor{cTangoSkyBlue1}{$-\frac{10\cdot5^4 2^4s^6}{\Delta}$}}] {$0,0,0,2$};
\node[place] (4) [below=60pt of 9,xshift=180pt,label=left:\small{\textcolor{cTangoSkyBlue1}{$28\frac{10^65^22^2s}{\Delta}$}},label=right:{\small{$40\frac{10^65^22^2s}{\Delta}$ extra}}] {$0,0,4,1$};
\node[place] (5) [below=60pt of 4,label=below:\small{\textcolor{cTangoSkyBlue1}{$-12\frac{10^5 5^3 2^2s^2}{\Delta}$}}] {$1,1,5,0$};
\draw[-latex,thick] (0)--(1)
 node[very near start,above,xshift=-4pt] {\small{\textcolor{cTangoPlum1}{$2$}}}
 node[midway,above] {\small{\textcolor{cTangoChameleon1}{$-80\frac{10^95s^{-3}}{\Delta}$}}};
\draw[-latex,thick] (4)--(5)
 node[very near start,right] {\small{\textcolor{cTangoPlum1}{$2$}}}
 node[midway,left] {\small{\textcolor{cTangoChameleon1}{$-12\frac{10^55^32^2s}{\Delta}$}}};
\draw[-latex,thick] (9)--(0)
 node[very near start,left] {\small{\textcolor{cTangoPlum1}{$2$}}}
 node[midway,right] {\small{\textcolor{cTangoChameleon1}{$-\frac{5^52^4s^6}{\Delta}$}}};
\end{tikzpicture}
\end{center}
\caption{The Jacobi path for $(x_1x_2x_3x_4)^4$ after using the Griffiths formula three times}
\end{figure}

\begin{figure}[H]
\begin{center}
\begin{tikzpicture}
\node[place] (0) [label=above:\small{\textcolor{cTangoSkyBlue1}{$-\frac{c_qs^{6}}{\Delta}$}}] {$0,0,0,0$};
\end{tikzpicture}
\end{center}
\caption{The Jacobi path for $(x_1x_2x_3x_4)^4$ after using the Griffiths formula four times}
\end{figure}

Now we can put everything together and write down the expression of $\frac{24s^5(x_1x_2x_3x_4)^4\Omega_0}{f^5}$ in the basis $\{b_0,b_1,b_2,b_3\}$. First the formula we got from the Jacobi path starting and ending at $(4,4,4,4)$:

\begin{align*}
\frac{24s^5(x_1x_2x_3x_4)^4\Omega_0}{f^5}&= \left(\frac{10c_qs^{10}+15c_d}{\Delta}\right)\frac{6s^4(x_1x_2x_3x_4)^3\Omega_0}{f^4}\\
&\quad +\left(\frac{-10c_qs^{10}+80c_d}{\Delta}\right)b_2+\left(\frac{5c_qs^{10}+80c_d}{\Delta}\right)b_1\\
&\quad +\frac{40}{\Delta}b_3+\left(\frac{-c_qs^{10}}{\Delta}\right)b_0.\\
\end{align*}

And now we put in formula (\ref{formel-3,3,3,3}) and get the final expression in the basis:

\begin{align}\label{formel-4,4,4,4}
\begin{split}
\frac{24s^5(x_1x_2x_3x_4)^4\Omega_0}{f^5}&=\left(\frac{10c_qs^{10}+15c_d}{\Delta}\right)\left(\left(\frac{6c_qs^{10}+9c_d}{\Delta}\right)b_2\right.\\
&\quad +\left.\left(\frac{-7c_qs^{10}+9c_d}{\Delta}\right)b_1+\frac{8}{\Delta}b_3+\left(\frac{c_qs^{10}}{\Delta}\right)b_0\right)\\
&\quad +\left(\frac{-25c_qs^{10}+80c_d}{\Delta}\right)b_2+\left(\frac{15c_qs^{10}+80c_d}{\Delta}\right)b_1\\
&\quad +\frac{40}{\Delta}b_3+\left(\frac{-c_qs^{10}}{\Delta}\right)b_0\\
&=b_2\left(\frac{60c_qs^{20}+180c_dc_qs^{10}+135c_d^2}{\Delta^2}+\frac{-25c_qs^{10}+80c_d}{\Delta}\right)\\
&\quad +b_1\left(\frac{-70c_qs^{20}-15c_dc_qs^{10}+135c_d^2}{\Delta^2}+\frac{15c_qs^{10}+80c_d}{\Delta}\right)\\
&\quad +b_3\left(\frac{80c_qs^{10}+120c_d}{\Delta^2}+\frac{40}{\Delta}\right)\\
&\quad +b_0\left(\frac{10c_qs^{20}+15c_dc_qs^{10}}{\Delta^2}+\frac{-c_qs^{10}}{\Delta}\right).
\end{split}
\end{align}

From Theorem \ref{conj} we know that the Picard-Fuchs equation should be

\begin{align*}
0&=(5^52^4s^{10}\delta^3(\delta+5)-10^{10}(\delta-1)(\delta-3)(\delta-7)(\delta-9))(\omega)\\
&=(c_qs^{10}-c_d)\delta^4\omega+(5\cdot c_qs^{10}+20\cdot c_d)\delta^3\omega-130\cdot c_d\delta^2\omega+300\cdot c_d\delta\omega-189\cdot c_d\omega.
\end{align*}

Now we can put equation (\ref{formel-delta}) in this Picard-Fuchs equation and end up with the following formula to check

\begin{align*}
0&=(c_qs^{10}-c_d)\left(b_0-15b_1+25b_2-10\frac{6s^4(x_1x_2x_3x_4)^3\Omega_0}{f^4}+\frac{24s^5(x_1x_2x_3x_4)^4\Omega_0}{f^5}\right)\\
&\quad +(5\cdot c_qs^{10}+20\cdot c_d)\left(b_0-7b_1+6b_2-\frac{6s^4(x_1x_2x_3x_4)^3\Omega_0}{f^4}\right)\\
&\quad -130\cdot c_d(b_0-3b_1+b_2)+300\cdot c_d(b_0-b_1)-189\cdot c_db_0\\
&\stackrel{\ref{formel-delta}}{=}(c_qs^{10}-c_d)\frac{24s^5(x_1x_2x_3x_4)^4\Omega_0}{f^5}\\
&\quad +(-15\cdot c_qs^{10}-10\cdot c_d)\frac{6s^4(x_1x_2x_3x_4)^3\Omega_0}{f^4}\\
&\quad +(55\cdot c_qs^{10}-35\cdot c_d)b_2+(-50\cdot c_qs^{10}-35\cdot c_d)b_1+6\cdot c_qs^{10}b_0.
\end{align*}

Now we can put in the expression of $\frac{6s^4(x_1x_2x_3x_4)^3\Omega_0}{f^4}$ in the basis (\ref{formel-3,3,3,3}):

\begin{align*}
0&=\Delta\frac{24s^5(x_1x_2x_3x_4)^4\Omega_0}{f^5}\\
&\quad +\left(\frac{-90c_qs^{20}-195\cdot c_dc_qs^{10}-90c_d^2}{\Delta}+55\cdot c_qs^{10}-35\cdot c_d\right)b_2\\
&\quad +\left(\frac{105c_qs^{20}-65\cdot c_dc_qs^{10}-90c_d^2}{\Delta}-50\cdot c_qs^{10}-35\cdot c_d\right)b_1\\
&\quad +\left(\frac{-120\cdot c_qs^{10}-80\cdot c_d}{\Delta}\right)b_3+\left(\frac{-15c_qs^{20}-10\cdot c_dc_qs^{10}}{\Delta}+6\cdot c_qs^{10}\right)b_0.
\end{align*}

Finally we put in formula (\ref{formel-4,4,4,4}) which describes $\frac{24s^5(x_1x_2x_3x_4)^4\Omega_0}{f^5}$ in the basis and we check that the Picard-Fuchs equation holds:

\begin{align*}
0&=\left(\frac{-30c_qs^{20}-15\cdot c_dc_qs^{10}+45c_d^2}{\Delta}+30\cdot c_qs^{10}+45\cdot c_d\right)b_2\\
&\quad +\left(\frac{35c_qs^{20}-80\cdot c_dc_qs^{10}+45c_d^2}{\Delta}-35\cdot c_qs^{10}+45\cdot c_d\right)b_1\\
&\quad +\left(\frac{-40\cdot c_qs^{10}-40\cdot c_d}{\Delta}+40\right)b_3+\left(\frac{-5c_qs^{20}+5\cdot c_dc_qs^{10}}{\Delta}+5\cdot c_qs^{10}\right)b_0\\
&=\frac{1}{\Delta}(0\cdot b_2+0\cdot b_1+0\cdot b_0+0\cdot b_3).
\end{align*}
\end{ex}

\newpage
\thispagestyle{empty}
\chapter{The Picard-Fuchs equation for invertible polynomials and consequences}\label{vier}

In this Chapter we focus on the Picard-Fuchs equation of the one-parameter family $f(\x)$ and discuss some consequences of the results achieved so far. From the last chapter we already know the order of the Picard-Fuchs equation. In the first section of this chapter we calculate the \GKZ system and in the second section we see how this fits together with the result on the order of the Picard-Fuchs equation to prove Theorem \ref{conj}, which states the Picard-Fuchs equation for the one-parameter family $f(\x)$. In the Section \ref{pifu} we will also see how this relates to a paper by Corti and Golyshev \cite{CG}, where the same differential equation appears. This is also the starting point for Section \ref{cohom}, where we concentrate on relations between the cohomology of the hypersurface defined by the one-parameter family $f(\x)$ and the cohomology of the solution space of the Picard-Fuchs equation. In Section \ref{ASD} we will discuss the results in an important class of examples given by Arnold's strange duality. This was also the starting point of the research done in this thesis. Finally, in the last section we cover the relation between the zero sets of the Picard-Fuchs equation of $f$ for special choices of the parameter, the Poincar\'{e} series of the dual polynomial $g^t$ and the monodromy in the solution space of the Picard-Fuchs equation.

\section{The \GKZ system for invertible polynomials}\label{calculation-GKZ}

This section is devoted to \GKZ systems. We will give a short introduction to \GKZ systems and do the calculations for invertible polynomials afterwards.

\subsection*{Introduction to \GKZ systems}

In this first part we want to give a short introduction to \GKZ systems as far as we need it. The theory on \GKZ systems is much larger than the part we present here. Good references for an introduction as well as an overview on several aspects of \GKZ systems are the article by Stienstra \cite{MR2306158}, which has a large part on solutions of \GKZ systems, the book by Katz and Cox \cite{MR1677117}, which among other things embeds \GKZ systems in a bigger context, and the article of Hosono \cite{MR1653027}, which focuses on the case of toric varieties. The theory of \GKZ systems was originally established by a series of articles of Gelfand, Kapranov and Zelevinsky \cite{MR1011353,MR1264328,MR1080980,MR1166813} as a generalisation of hypergeometric differential equations. This also explains the name \GKZ systems.

\begin{notation}\label{notation-GKZ}
Let $\mc{A}\subset \Z^{n}$ be a finite subset which generates $\Z^{n}$ as an abelian group and for which there exists a group homomorphism $h:\Z^{n}\rightarrow \Z$ such that $h(\mc{A})={1}$, i.e. $\mc{A}$ lies in a $(n-1)$-dim. hypersurface. Let $\gamma\in\C^{n}$ be an arbitrary vector.\\
Let $|\mc{A}|=N$, then $\mb{L}:=\{(l_1,\dots,l_N)\in\Z^N: l_1a_1+\dots+l_Na_N=0, a_i\in\mc{A}\}$ denotes the lattice of linear relations among $\mc{A}$. Because of $\mc{A}$ lying in a hypersurface, $\sum l_i=0$ holds for $(l_1,\dots,l_N)\in\mb{L}$.
\end{notation}
 
\begin{rem}
We will calculate the \GKZ system for the one-parameter family $f(\x)$ later. Keep in mind that for these calculations $\mc{A}$ will be the set of all exponent vectors of our one-parameter family. The reasons for this will also become clear later.
\end{rem}

\begin{defi}
The \GKZ system (sometimes also called $\mc{A}$ system) for $\mc{A}$ and $\gamma$ is a system of differential equations for functions $\Phi$ of $N$ variables $v_1,\dots,v_N$ given by
\begin{align}
\prod_{l_i>0}\big(\frac{\partial}{\partial v_i}\big)^{l_i}\Phi&=\prod_{l_i<0}\big(\frac{\partial}{\partial v_i}\big)^{-l_i}\Phi \text{ for every }l\in\mb{L}\text{ and}\label{GKZ1}\\
\sum_{i=1}^{N} a_{ij} v_i\frac{\partial\Phi}{\partial v_i}&=\gamma_j\Phi \text{ for all }j=1,\dots,k+1\text{ and }(a_{i1},\dots,a_{i\,k+1})\in\mc{A}\label{GKZ2}.
\end{align}
\end{defi}

The above definition gives a system of partial differential equations. We stated the basic definition of a \GKZ system. From this point on we will continue in the special case of invertible polynomials.

\subsection*{Calculation of the \GKZ system for invertible polynomials}

We will now start calculating the \GKZ system for the one-parameter family $f(\x)=g(\x)+s\prod_i x_i$, where $g(\x)$ is an invertible polynomial. The notation in this section is the same as  before and can be found in \ref{notation} and \ref{notation-GKZ}. In addition we will define some extra notation:

\begin{notation}
We define the rows of the exponent matrix $E$ to be $\underline{e}_i=(e_{i\,1},\dots,e_{i\,n})$ for $i=1,\dots,n$. Then we can write $g(\x)$ as $g(\x)=\sum_{i=1}^n \x^{\underline{e}_i}$, where $\x^{\underline{e}_i}=\prod_{j=1}^n x_j^{e_{i\,j}}$. Now we define a general $n+1$-parameter family $f_{\uu}(\x)=f_{v_1,\dots,v_n}(\x)=\sum_{i=1}^n v_i\x^{\underline{e}_i}+s\x^{(1,\dots,1)}$ with parameters $v_1\dots,v_n$ and $s$. So in the previously used notation we have $N=n+1$ and we set $v_{n+1}:=s$. In this way the notation is consistent with the previous chapters, because we have that
\begin{align*}
f_{1,\dots,1}(\x)=\sum_{i=1}^n \x^{\underline{e}_i}+s\x^{(1,\dots,1)}=g(\x)+s\prod_{i=1}^n x_i=f(\x).
\end{align*}
\end{notation}

We will now start calculating the \GKZ system for $\mc{A}=\{\underline{e}_1^t,\dots,\underline{e}_n^t,(1,\dots,1)^t\}$ and $\gamma=(-1,\dots,-1)^t$. The reason for the choice of $\gamma$ will become clear when we look at the solutions of the \GKZ system. For the first equation (\ref{GKZ1}) we need to calculate the lattice of linear relations $\mathbb{L}$ among the vectors in $\mc{A}$. If we define $A$ to be the matrix with columns $\underline{e}_1^t,\dots,\underline{e}_n^t,(1,\dots,1)^t$, then $A$ is an $n \times (n+1)$-matrix and $\mathbb{L}$ is 1-dimensional. We know that
\begin{align*}
A\cdot \begin{pmatrix} \widehat{q}_1\\\vdots\\\widehat{q}_n\\-\widehat{d}\end{pmatrix}=E^t\cdot\begin{pmatrix} \widehat{q}_1\\\vdots\\\widehat{q}_n\end{pmatrix}-\begin{pmatrix} \widehat{d}\\\vdots\\\widehat{d}\end{pmatrix}=\begin{pmatrix} 0\\\vdots\\0\end{pmatrix}
\end{align*}
and therefore $\mathbb{L}=\langle(\widehat{q}_1,\dots,\widehat{q}_n,-\widehat{d})^t\rangle$. Now we are able to write down equation (\ref{GKZ1}) for this lattice $\mathbb{L}$:
\begin{align}\label{GKZeins}
\left(\frac{\partial}{\partial s}\right)^{\widehat{d}}\Phi=\left(\frac{\partial}{\partial v_1}\right)^{\widehat{q}_1}\cdot\dots\cdot\left(\frac{\partial}{\partial v_n}\right)^{\widehat{q}_n}\Phi.
\end{align}

In the end we want to compare the \GKZ system to the Picard-Fuchs equation from Theorem \ref{conj}. To do this we will write the \GKZ system with the differential operators $\delta=s\frac{\partial}{\partial s}$ and $\delta_i=v_i\frac{\partial}{\partial v_i}$ for $i=1,\dots,n$ by inserting $s^{-1}\delta=\frac{\partial}{\partial s}$ and $v_i^{-1}\delta_i=\frac{\partial}{\partial v_i}$.
\begin{align*}
\left(s^{-1}\delta\right)^{\widehat{d}}\Phi=\left(v_1^{-1}\delta_1\right)^{\widehat{q}_1}\cdot\dots\cdot\left(v_n^{-1}\delta_n\right)^{\widehat{q}_n}\Phi.
\end{align*}

Now we move $s^{-1}$ and $v_i^{-1}$ to the front and the product rule gives us an easy way to interchange the differential operators $\delta, \delta_i$ with the variables $s,v_i$:

\begin{align}
\begin{split}\label{delta-vertauschen}
\delta s^p&=s^p(\delta+p) \quad\text{ for } p\in\Z  \text{ and}\\
\delta_iv_i^p&=v_i^p(\delta_i+p) \quad\text{ for } i=1,\dots,n \text{ and } p\in\Z.
\end{split}
\end{align}
\pagebreak
Using these equations we can move every $s$ and every $v_i$ very quickly to the front of the equation:
\begin{align*}
(s^{-1}\delta)^{\widehat{d}}&=s^{-1}\delta s^{-1}\delta\dots s^{-1}\underbracket[1pt]{\delta s^{-1}}_{s^{-1}(\delta-1)}\delta\\
&=s^{-1}\delta s^{-1}\delta\dots s^{-1}\underbracket[1pt]{\delta s^{-2}}_{s^{-2}(\delta-2)}(\delta-1)\delta\\
&=\dots\\
&=s^{-\widehat{d}}(\delta-(\widehat{d}-1))\cdot\dots\cdot(\delta-1)\delta
\end{align*}
and in the same way we get
\begin{align*}
(v_i^{-1}\delta_i)^{\widehat{q}_i}&=v_i^{-1}\delta_i v_i^{-1}\delta_i\dots v_i^{-1}\underbracket[1pt]{\delta v_i^{-1}}_{v_i^{-1}(\delta_i-1)}\delta_i\\
&=v_i^{-1}\delta_i v_i^{-1}\delta_i\dots v_i^{-1}\underbracket[1pt]{\delta_i v_i^{-2}}_{v_i^{-2}(\delta_i-2)}(\delta_i-1)\delta_i\\
&=\dots\\
&=v_i^{-\widehat{q}_i}(\delta_i-(\widehat{q_i}-1))\cdot\dots\cdot(\delta_i-1)\delta_i.
\end{align*}

Putting this all together the first equation of the \GKZ system is given by

\begin{align}\label{gkz-gleichung1}
s^{-\widehat{d}}(\delta-(\widehat{d}-1))\cdot\dots\cdot(\delta-1)\delta\Phi=\prod_{i=1}^{n}v_i^{-\widehat{q}_i}(\delta_i-(\widehat{q}_i-1))\cdot\dots\cdot(\delta_i-1)\delta_i\Phi.
\end{align}

We will work with this equation later on and calculate the second part (\ref{GKZ2}) of the \GKZ system next.
The second system of equations of the \GKZ system is given by putting $\gamma=(-1,\dots,-1)^t$ in (\ref{GKZ2}):

\begin{align}\label{GKZzwei}
A\cdot\begin{pmatrix}v_1\frac{\partial}{\partial v_1}\\\vdots\\v_n\frac{\partial}{\partial v_n}\\s\frac{\partial}{\partial s}\end{pmatrix}\Phi=A\cdot\begin{pmatrix}\delta_1\\\vdots\\\delta_n\\\delta\end{pmatrix}\Phi=E^t\begin{pmatrix}\delta_1\\\vdots\\\delta_n\end{pmatrix}\Phi+\begin{pmatrix}1\\\vdots\\1\end{pmatrix}\delta\Phi=\begin{pmatrix}-1\\\vdots\\-1\end{pmatrix}\Phi
\end{align}

Before we do any further calculations, we focus on solutions of the \GKZ system. There is a whole theory on solutions of GKZ-System which, for example, is explained in \cite{MR2306158}. We however, do not need the full strength of this, because to compare the \GKZ system to the Picard-Fuchs equation in Theorem \ref{conj}, it is enough to know that the form $\omega=\frac{s\Omega_0}{f(\x)}$ is a solution of the \GKZ system shown in equation (\ref{gkz-gleichung1}) and (\ref{GKZzwei}). This is the goal, but we will start with a slightly different solution in the next lemma.

\begin{lemma}\label{solution-gkz}
The form $\Phi=\frac{\Omega_0}{f_{\uu}(\x)}$ is a solution for the above \GKZ system.
\end{lemma}

\begin{proof}
We will calculate the differentials for $\Phi=\frac{\Omega_0}{f_{\uu}(\x)}$ and see that the equations (\ref{GKZeins}) and (\ref{GKZzwei}) hold. For equation (\ref{GKZeins}) we need the partial derivatives with respect to $s$ and $v_i$ for $i=1,\dots,n$. They are easy to calculate:

\begin{align*}
\left(\frac{\partial}{\partial s}\right)^{\widehat{d}}\Phi&=(-1)^{\widehat{d}}\widehat{d}!\left(\prod x_i\right)^{\widehat{d}}\frac{\Omega_0}{\left(f_{\uu}(\x)\right)^{\widehat{d}+1}}\\
\frac{\partial}{\partial v_i}\Phi&=-\x^{\underline{e}_i}\frac{\Omega_0}{\left(f_{\uu}(\x)\right)^2}\\
\prod_{i=1}^{n}\left(\frac{\partial}{\partial v_i}\right)^{\widehat{q}_i}\Phi&=(-1)^{\sum\widehat{q}_i}\left(\sum \widehat{q}_i\right)!\x^{\sum\widehat{q}_i\underline{e}_i}\frac{\Omega_0}{\left(f_{\uu}(\x)\right)^{1+\sum \widehat{q}_i}}.
\end{align*}

Because of the Calabi-Yau condition we have $\sum \widehat{q}_i=\widehat{d}$ and from the definition of the dual weights and degree we get $\sum\widehat{q}_i\underline{e}_i=E^t\cdot (\widehat{q}_1,\dots,\widehat{q}_n)^t=(\widehat{d},\dots,\widehat{d})^t$. Therefore we have

\begin{align*}
\left(\frac{\partial}{\partial s}\right)^{\widehat{d}}\Phi&=(-1)^{\widehat{d}}\widehat{d}!\left(\prod x_i\right)^{\widehat{d}}\frac{\Omega_0}{\left(f_{\uu}(\x)\right)^{1+\widehat{d}}}\\
&=(-1)^{\sum\widehat{q}_i}\left(\sum \widehat{q}_i\right)!\x^{\sum\widehat{q}_i\underline{e}_i}\frac{\Omega_0}{\left(f_{\uu}(\x)\right)^{1+\sum \widehat{q}_i}}\\
&=\prod_{i=1}^{n}\left(\frac{\partial}{\partial v_i}\right)^{\widehat{q}_i}\Phi.
\end{align*}

This proves that $\Phi$ is a solution for equation (\ref{GKZeins}).
Now we check the second equation, where we need $\delta\Phi$ and $\delta_i\Phi$ for $i=1,\dots,n$, because the system of equations is given by:

\begin{align*}
\begin{pmatrix}1\\ \vdots \\ 1\end{pmatrix}\delta\Phi+E^t\begin{pmatrix}\delta_1\\\vdots\\\delta_n\end{pmatrix}\Phi+\begin{pmatrix}1\\\vdots\\1\end{pmatrix}\Phi=\begin{pmatrix}0\\\vdots\\0\end{pmatrix}\Phi.
\end{align*}

So for every $j=1,\dots,n$ we have the following equation:

\begin{align*}
\delta\Phi+\sum_{i=1}^n e_{ij}\delta_i\Phi+\Phi&=-s\x^{(1,\dots,1)}\frac{\Omega_0}{\left(f_{\uu}(\x)\right)^2}+\sum_{i=1}^ne_{ij}\left(-v_i\x^{\underline{e}_i}\right)\frac{\Omega_0}{\left(f_{\uu}(\x)\right)^2}+\frac{\Omega_0}{f_{\uu}(\x)}\\
&=-\frac{\left(s\x^{(1,\dots,1)}+\sum_{i=1}^ne_{ij}v_i\x^{\underline{e}_i}\right)\Omega_0}{\left(f_{\uu}\right)^2(\x)}+\frac{\Omega_0}{f_{\uu}(\x)}\\
&=-\frac{x_j\frac{\partial}{\partial x_j}f_{\uu}(\x)\Omega_0}{\left(f_{\uu}(\x)\right)^2}+\frac{\Omega_0}{f_{\uu}(\x)}\\
&=0,
\end{align*}

where the last expression is an exact form due to the Griffiths formula and is therefore zero.
\end{proof}

As mentioned before $\Phi$ is not the solution we want to have. A solution that would fit our purposes would be $\omega_{\uu}=\frac{s\Omega_0}{f_{\uu}(\x)}$, because $\omega_{1,\dots,1}=\frac{s\Omega_0}{f(\x)}=\omega$. So we insert $\Phi=s^{-1}\omega_{\uu}$ in the equations \ref{gkz-gleichung1} and \ref{GKZzwei}. So equation (\ref{gkz-gleichung1}) leads to:

\begin{align*}
s^{-\widehat{d}}(\delta-(\widehat{d}-1))\cdot\dots\cdot(\delta-1)\delta s^{-1}\omega_{\uu}&=\prod_{i=1}^{n}v_i^{-\widehat{q}_i}(\delta_i-(\widehat{q}_i-1))\cdot\dots\cdot(\delta_i-1)\delta_is^{-1}\omega_{\uu}.
\end{align*}

We can use equation (\ref{delta-vertauschen}) as earlier to move the variable $s$ to the front and get the following equation.

\begin{align}
s^{-\widehat{d}}(\delta-\widehat{d})\cdot\dots\cdot(\delta-1)\omega_{\uu}&=\prod_{i=1}^{n}v_i^{-\widehat{q}_i}(\delta_i-(\widehat{q}_i-1))\cdot\dots\cdot(\delta_i-1)\delta_i\omega_{\uu}.\label{finales-gkz-mit-u}
\end{align}

By putting $\Phi=s^{-1}\omega_{\uu}$ in equation (\ref{GKZzwei}) and using (\ref{delta-vertauschen}) again we get:

\begin{align*}
\begin{pmatrix}1\\\vdots\\1\end{pmatrix}\delta\Phi+E^t\begin{pmatrix}\delta_1\\\vdots\\\delta_n\end{pmatrix}\Phi&=\begin{pmatrix}-1\\\vdots\\-1\end{pmatrix}\Phi\\
\begin{pmatrix}1\\\vdots\\1\end{pmatrix}\delta s^{-1}\omega_{\uu}+E^t\begin{pmatrix}\delta_1\\\vdots\\\delta_n\end{pmatrix}s^{-1}\omega_{\uu}&=\begin{pmatrix}-1\\\vdots\\-1\end{pmatrix}s^{-1}\omega_{\uu}\\
s^{-1}\begin{pmatrix}1\\\vdots\\1\end{pmatrix}(\delta-1) \omega_{\uu}+s^{-1}E^t\begin{pmatrix}\delta_1\\\vdots\\\delta_n\end{pmatrix}\omega_{\uu}&=s^{-1}\begin{pmatrix}-1\\\vdots\\-1\end{pmatrix}\omega_{\uu}\\
\begin{pmatrix}1\\\vdots\\1\end{pmatrix}\delta \omega_{\uu}+E^t\begin{pmatrix}\delta_1\\\vdots\\\delta_n\end{pmatrix}\omega_{\uu}&=\begin{pmatrix}0\\\vdots\\0\end{pmatrix}\omega_{\uu}.
\end{align*}

Solving this equation for $(\delta_1,\dots,\delta_n)^t$ gives

\begin{align*}
\begin{pmatrix} \delta_1\\\vdots\\\delta_n \end{pmatrix}=-(E^t)^{-1}\begin{pmatrix}1\\\vdots\\1\end{pmatrix}\delta=\begin{pmatrix}\frac{\widehat{q}_1}{\widehat{d}}\\\vdots\\\frac{\widehat{q}_n}{\widehat{d}}\end{pmatrix}\delta.
\end{align*}

In other words we can write each of the differential operators $\delta_1,\dots,\delta_n$ in terms of $\delta$. For all $i=1,\dots,n$ we have

\begin{align*}
\delta_i=-\frac{\widehat{q}_i}{\widehat{d}}\delta.
\end{align*}

We can use this equation to write equation (\ref{finales-gkz-mit-u}) as an ordinary differential equation with differential operator $\delta$:

\begin{align*}
s^{-\widehat{d}}(\delta-\widehat{d})\cdot\dots\cdot(\delta-1)\omega_{\uu}&=\prod_{i=1}^{n}v_i^{-\widehat{q}_i}(\delta_i-(\widehat{q}_i-1))\cdot\dots\cdot(\delta_i-1)\delta_i\omega_{\uu}\\
&=\prod_{i=1}^{n}v_i^{-\widehat{q}_i}\left(-\frac{\widehat{q}_i}{\widehat{d}}\delta-(\widehat{q}_i-1)\right)\cdot\dots\cdot\left(-\frac{\widehat{q}_i}{\widehat{d}}\delta-1\right)\left(-\frac{\widehat{q}_i}{\widehat{d}}\delta\right)\omega_{\uu}\\
&=\prod_{i=1}^{n}\left(-\frac{\widehat{q}_i}{\widehat{d}}v_i^{-1}\right)^{\widehat{q}_i}\left(\delta+\frac{(\widehat{q}_i-1)\widehat{d}}{\widehat{q}_i}\right)\cdot\dots\cdot\left(\delta+\frac{\widehat{d}}{\widehat{q}_i}\right)\delta\omega_{\uu}.
\end{align*}

Now we set $v_i=1$, which brings us back to our one-parameter family $f(\x)$. Because the solutions of the differential equation before are given by $\omega_{\uu}$, we get a differential equation for $\omega=\frac{s\Omega_0}{f(\x)}$. So our final expression is given by

\begin{align*}
s^{-\widehat{d}}(\delta-\widehat{d})\cdot\dots\cdot(\delta-1)\omega=\prod_{i=1}^{n}\left(-\frac{\widehat{q}_i}{\widehat{d}}\right)^{\widehat{q}_i}\left(\delta+\frac{(\widehat{q}_i-1)\widehat{d}}{\widehat{q}_i}\right)\cdot\dots\cdot\left(\delta+\frac{\widehat{d}}{\widehat{q}_i}\right)\delta\omega.
\end{align*}

or to have the same appearance as in Theorem \ref{conj}:

\begin{align}\label{finales-gkz}
0=s^{\widehat{d}}\prod_{i=1}^{n}\left(\widehat{q}_i\right)^{\widehat{q}_i}\delta\left(\delta+\frac{\widehat{d}}{\widehat{q}_i}\right)\cdot\dots\cdot\left(\delta+\frac{(\widehat{q}_i-1)\widehat{d}}{\widehat{q}_i}\right)\omega-(-\widehat{d})^{-\widehat{d}}(\delta-1)\cdot\dots\cdot(\delta-\widehat{d})\omega.
\end{align}

\section{The Picard-Fuchs equation}\label{pifu}

In the last chapter we already proved the order of the Picard-Fuchs equation. If we look at examples such as those in Section \ref{detailed-example-gd} and \ref{ASD} and in Appendix \ref{k3}, we can also conjecture exactly what the Picard-Fuchs equation looks. We can use the \GKZ system that we calculated in the last section to confirm that this is true.

\pagebreak
\begin{thm}\label{conj}
Let $g(x_1,\dots,x_n)$ be an invertible polynomial with weighted degree $\deg g=d$ and reduced weights $q_1,\dots,q_n$ for which the Calabi-Yau condition,  $d=\sum q_i$, holds. Let $g^t(x_1,\dots,x_n)$ be the transposed polynomial with reduced weights $\widehat{q}_1,\dots,\widehat{q}_n$ and degree $\deg g^t=\widehat{d}$. Then the Picard-Fuchs equation for the one-parameter family $f(x_1,\dots,x_n)=g(x_1,\dots,x_n)+s\prod x_i$ is given by
\begin{align*}
0=\prod_{i=1}^n \widehat{q}_i^{\widehat{q}_i}s^{\widehat{d}}\prod_{i=1}^{n}\prod_{j=0}^{\widehat{q}_i-1}(\delta+\frac{j\cdot \widehat{d}}{\widehat{q}_i})\prod_{\ell\in I}(\delta+\ell)^{-1}-(-\widehat{d})^{\widehat{d}}\prod_{j=0}^{\widehat{d}-1}(\delta-j)\prod_{\ell\in I}(\delta-\ell)^{-1},
\end{align*}
where $I=\{0,\dots,\widehat{d}-1\}\cap\bigcup_{i=1}^n \left\{0,\frac{\widehat{d}}{\widehat{q}_i},\frac{2\widehat{d}}{\widehat{q}_i},\dots,\frac{(\widehat{q}_i-1)\widehat{d}}{\widehat{q}_i}\right\}$.
\end{thm}

\begin{proof}
From Lemma \ref{solution-gkz} we know that $\omega=\frac{s\Omega}{f(\x)}$ is a solution for the equation (\ref{finales-gkz}). It follows that all period integrals are solutions of (\ref{finales-gkz}) and therefore the Picard-Fuchs equation divides
\begin{align*}
0=s^{\widehat{d}}\prod_{i=1}^{n}\left(\widehat{q}_i\right)^{\widehat{q}_i}\delta\left(\delta+\frac{\widehat{d}}{\widehat{q}_i}\right)\cdot\dots\cdot\left(\delta+\frac{(\widehat{q}_i-1)\widehat{d}}{\widehat{q}_i}\right)\omega-(-\widehat{d})^{-\widehat{d}}(\delta-1)\cdot\dots\cdot(\delta-\widehat{d})\omega. 
\end{align*}
We also know from Theorem \ref{thm} that the order of the Picard-Fuchs equation is given by 
\begin{align*}
u=\widehat{d}-\left(\{0,1,\dots,\widehat{d}-1\}\cap\bigcup_{i=1}^n\{0,\frac{\widehat{d}}{\widehat{q}_i},\dots,\frac{(\widehat{q}_i-1)\widehat{d}}{\widehat{q}_i}\}\right).
\end{align*}
So we try to find common factors in the summands of (\ref{finales-gkz}) until the order of the equation is $u$. If we multiply equation (\ref{finales-gkz}) by $s^{-\widehat{d}}$ and use the commutation relations (\ref{delta-vertauschen}) to pass it through the differential operators we get
\begin{align*}
0=\prod_{i=1}^{n}\left(\widehat{q}_i\right)^{\widehat{q}_i}\delta\left(\delta+\frac{\widehat{d}}{\widehat{q}_i}\right)\cdot\dots\cdot\left(\delta+\frac{(\widehat{q}_i-1)\widehat{d}}{\widehat{q}_i}\right)\omega-(-\widehat{d})^{-\widehat{d}}(\delta+(\widehat{d}-1))\cdot\dots\cdot(\delta+1)\delta s^{-\widehat{d}}\omega. 
\end{align*}
Now it is easy to see that every linear factor $\delta+j$ with $j\in\{0,1,\dots,\widehat{d}-1\}\cap\bigcup_{i=1}^n\{0,\frac{\widehat{d}}{\widehat{q}_i},\dots,\frac{(\widehat{q}_i-1)\widehat{d}}{\widehat{q}_i}\}$ is in both summands and can therefore be deleted.
This leads us to the equation
\begin{align*}
0=\prod_{i=1}^n \widehat{q}_i^{\widehat{q}_i}s^{\widehat{d}}\prod_{i=1}^{n}\prod_{j=0}^{\widehat{q}_i-1}(\delta+\frac{j\cdot \widehat{d}}{\widehat{q}_i})\prod_{\ell\in I}(\delta+\ell)^{-1}-(-\widehat{d})^{\widehat{d}}\prod_{j=0}^{\widehat{d}-1}(\delta-j)\prod_{\ell\in I}(\delta-\ell)^{-1}
\end{align*}
where $I=\{0,\dots,\widehat{d}-1\}\cap\bigcup_{i=1}^n \left\{0,\frac{\widehat{d}}{\widehat{q}_i},\frac{2\widehat{d}}{\widehat{q}_i},\dots,\frac{(\widehat{q}_i-1)\widehat{d}}{\widehat{q}_i}\right\}$.\\
Finally, this equation is divisible by the Picard-Fuchs equation and has the degree of the Picard-Fuchs equation (cf. Theorem \ref{thm}).
\end{proof}

We give another class of examples here, which are the simple elliptic singularities. There are only 3 examples and their Picard-Fuchs equation is known.

\begin{ex}
In the following table we can see 3 polynomials that define the simple elliptic singularities, their weights, the degree and the Picard-Fuchs equation, which can easily be calculated with Theorem \ref{conj} or with the Griffiths-Dwork method directly.

\begin{table}[H]
\begin{center}
\begin{tabular}{lcccl}
\toprule
Name & Invertible polynomial & Degree & Weights & Picard-Fuchs equation\\
\midrule
$\tilde{\text{E}}_{6}$ & $x^3+y^3+z^3$ & $3$ & $(1,1,1)$ & $s^3\delta^2+3^3(\delta-1)(\delta-2)$ \\
\midrule
$\tilde{\text{E}}_{7}$ & $x^4+y^4+z^2$ & $4$ & $(1,1,2)$ & $s^4\delta^2-4^3(\delta-1)(\delta-3)$ \\
\midrule
$\tilde{\text{E}}_{8}$ & $x^6+y^3+z^2$ & $6$ & $(1,2,3)$ & $s^6\delta^2-2\cdot6^3(\delta-1)(\delta-5)$ \\
\bottomrule
\end{tabular}
\end{center}
\caption{Simple elliptic singularities and their Picard-Fuchs equations}
\label{elliptic curves}
\end{table}
\end{ex}

A similar result, but approached from a different point of view, can be found in a paper by Corti and Golyshev \cite{CG}. In this paper the differential equation that they look at is the same as our Picard-Fuchs equation, but they start with a local system, which is given in the following way:

\begin{align}\label{local-system}
Y=\left\{\begin{array}{cl} \prod_{i=1}^n y_i^{w_i}&=\lambda\\\sum_{i=1}^n y_i&=1 \end{array} \right. \quad\subset(\C^{*})^{n}\times\C^{*}
\end{align}

If we insert $y_i=-s^{-1}\x^{\underline{e}_i-(1,\dots,1)}$ and $w_i=\widehat{q}_i$, then we get that $Y$ consists of the following two equations:

\begin{align*}
\lambda=\prod_{i=1}^n y_i^{w_i}&=\prod_{i=1}^n \left(-s^{-1}\x^{\underline{e}_i-(1,\dots,1)}\right)^{\widehat{q}_i}=\left((-s)^{-\sum \widehat{q}_i}\right)\x^{\sum \widehat{q}_i\underline{e}_i-(\sum \widehat{q}_i,\dots,\sum \widehat{q}_i)}\\
&=(-s)^{-\widehat{d}}\x^{(\widehat{d},\dots,\widehat{d})-(\widehat{d},\dots,\widehat{d})}=(-s)^{-\widehat{d}}\\
1=\sum_{i=1}^n y_i&=\sum_{i=1}^n \left(-s^{-1}\x^{\underline{e}_i-(1,\dots,1)}\right)=-s^{-1}\x^{-(1,\dots,1)}\sum_{i=1}^n \x^{\underline{e}_i}.
\end{align*}

So, from the first equation we get $(-s)^{-\widehat{d}}=\lambda$ and the second equation can easily be rewritten as

\begin{align*}
0=\sum_{i=1}^n \x^{\underline{e}_i}+s\x^{(1,\dots,1)}=f(\x).
\end{align*}

This shows the direct connection to our hypersurface $V(f)$. It is very easy to write the Picard-Fuchs equation with differential operator $\mathcal{D}=\lambda\frac{\partial}{\partial \lambda}$, because the relation between $\mathcal{D}$ and $\delta$ is just given by 
\begin{align*}
\delta=s\frac{\partial}{\partial s}=-\widehat{d}(-s)^{-\widehat{d}}\frac{\partial}{\partial (-s)^{-\widehat{d}}}=-\widehat{d}\lambda\frac{\partial}{\partial \lambda}=-\widehat{d}\mathcal{D}. 
\end{align*}
So in terms of $\mathcal{D}$ the Picard-Fuchs equation is given by
\begin{align}
0&=\prod_{i=1}^n \widehat{q}_i^{\widehat{q}_i}s^{\widehat{d}}\prod_{i=1}^{n}\prod_{j=0}^{\widehat{q}_i-1}(\delta+\frac{j\cdot \widehat{d}}{\widehat{q}_i})\prod_{\ell\in I}(\delta+\ell)^{-1}-(-\widehat{d})^{\widehat{d}}\prod_{j=0}^{\widehat{d}-1}(\delta-j)\prod_{\ell\in I}(\delta-\ell)^{-1}\notag\\
&=\prod_{i=1}^n \widehat{q}_i^{\widehat{q}_i}\lambda^{-1}\prod_{i=1}^{n}\prod_{j=0}^{\widehat{q}_i-1}(-\widehat{d}\mathcal{D}+\frac{j\cdot \widehat{d}}{\widehat{q}_i})\prod_{\ell\in I}(-\widehat{d}\mathcal{D}+\ell)^{-1}-\widehat{d}^{\widehat{d}}\prod_{j=0}^{\widehat{d}-1}(-\widehat{d}\mathcal{D}-j)\prod_{\ell\in I}(-\widehat{d}\mathcal{D}-\ell)^{-1}\notag\\
0&=\prod_{i=1}^n \widehat{q}_i^{\widehat{q}_i}\prod_{i=1}^{n}\prod_{j=0}^{\widehat{q}_i-1}(\mathcal{D}-\frac{j}{\widehat{q}_i})\prod_{\ell\in I}(\mathcal{D}-\frac{\ell}{\widehat{d}})^{-1}-\widehat{d}^{\widehat{d}}\lambda\prod_{j=0}^{\widehat{d}-1}(\mathcal{D}+\frac{j}{\widehat{d}})\prod_{\ell\in I}(\mathcal{D}+\frac{\ell}{\widehat{d}})^{-1}\label{pf-in-t}
\end{align}

which agrees with formula (1) in \cite{CG}. 

In Theorem 1.1 of the article \cite{CG} it is stated that the solutions of the Picard-Fuchs equation come from the local system (\ref{local-system}) and in Conjecture 1.4 and Proposition 1.5 the Hodge numbers for the solution space are given. This brings us to the next section where we will investigate this in detail.

\section{Statements on the cohomology of the solution space}\label{cohom}

We want to relate already known statements to the work we have done so far in the thesis. First we continue the last section. We will relate our results to work of Corti and Golyshev \cite{CG}. In their paper there is a result that calculates the Hodge numbers of the solution space of the Picard-Fuchs equation. We will state their result in form, which is compatible with our setting.

\begin{prop}\emph{(\cite{CG} Conjecture 1.4 and Proposition 1.5)}\label{hodge-numbers-goly}
Consider the sets $A:=\bigsqcup_{i=1}^n \left(\left(Q_i\setminus \{\widehat{d}\}\right)\cup\{0\}\right)$ and $D_0=(D\setminus \{\widehat{d}\})\cup\{0\}$ (cf. Definition \ref{vundu}). Set $\{\alpha_1,\dots,\alpha_u\}:=A\setminus(A\cap D)$ with $\alpha_i\leq\alpha_{i+1}$ for all $i$ and $\{\beta_1,\dots,\beta_u\}:=D\setminus(A\cap D)$ with $\beta_i<\beta_{i+1}$ for all $i$.
Now consider the differential equation (\ref{pifu}), which is with the above notation given by

\begin{align*}
s^{\widehat{d}}\prod_{i=1}^n\widehat{q}_i^{\widehat{q}_i}\prod_{i=1}^u(\delta+\alpha_i)-(-\widehat{d})^{\widehat{d}}\prod_{i=1}^u(\delta-\beta_i)=0.
\end{align*}

Now define the following function

\begin{align*}
p(k):=|\{j|\,\alpha_j<\beta_k\}|-(k-1) \quad\text{ for } \quad k=1,\dots,u
\end{align*}

and let $p_+:=\max \{p(k)\}$ and $p_-:=\min \{p(k)\}$.

\pagebreak
Then the local system of solutions of the ordinary differential equation above supports a real polarised variation of Hodge structure of weight $p_+-p_-$ and Hodge numbers
\begin{align*}
h^{j-p_-,p_+-j}=|p^{-1}(j)|.
\end{align*}

\end{prop}

\begin{kor} It follows easily from the calculations in Chapter \ref{haupt} that the following numbers coincide:
\begin{itemize}
\item $p_+=p(1)=n-1$
\item $p_-=p(u)=1$
\item $\sum_{j=1}^{n-1} h^{j-1,n-1-j}=\sum_{j=1}^{n-1} |p^{-1}(j)|=u$.
\end{itemize}
\end{kor}

We are able to make the relation between $u$ and the above Hodge numbers even more precise.

\begin{prop}
Let $u=u_0+\dots+u_{n-2}$, where $u_i$ denotes the number of degree $i\cdot d$ basis elements of the $u$ basis elements one needs to write the Picard-Fuchs equation as calculated in the proof of Theorem \ref{thm}. Then
\begin{align*}
u_{i-1}=h^{i-1,n-i-1}=h^{i-p_-,p_+-i}.
\end{align*}
\end{prop}

\begin{rem}
Notice that $u_i\geq 1$ for all $i$, because we have at least one basis element in every degree, and $u_0=u_{n-2}=1$, because in degree $0$ and $n-2$ we have exactly the basis elements $\frac{s\Omega_0}{f}$ and $\frac{s^{n-1}\left(\prod x_i\right)^{n-2}\Omega_0}{f^{n-1}}$ respectively.
\end{rem}

\begin{proof}
We will relate the function $p(k)$ to the $u$ basis elements. In particular we show that for every $1\leq k \leq u$ with $p(k)=i$ we need one basis element in degree $(n-i-1)d$. Notice that $p(k)$ can be written recursively as follows:
\begin{align*}
p(k+1)&=\left|\{j|\,\alpha_j<\beta_{k+1}\}\right|-k=\left|\{j|\,\alpha_j<\beta_k\}\right|-(k-1)+\left|\{j|\,\beta_{k}<\alpha_j<\beta_{k+1}\}\right|-1\\
&=p(k)+\left|\{j|\,\beta_{k}<\alpha_j<\beta_{k+1}\}\right|-1.
\end{align*}
Now we will show the statement via induction. If $p(k)=i$ corresponds to a basis element in degree $(n-i-1)d$, then $p(k+1)$ corresponds to a basis element in degree $(n-i-1)d+(\left|\{j|\,\beta_{k-1}<\alpha_j<\beta_k\}\right|+1)d$. For the correspondence between the function $p$ and the basis elements, we just view the $\alpha_i$ and $\beta_i$ as potential positions on the Jacobi path, where the $\alpha_i$ correspond to positions which have multiple possibilities or which are occupied by a partial derivative that creates an extra vertex and the $\beta_i$ correspond to free positions before shifting. Notice that in contrast to the proof of Theorem \ref{thm} the $p(k)$ count the end of a connected part on the Jacobi path and not the beginning.

First we investigate $k=1$. We know that $\alpha_1=\dots=\alpha_{n-1}=0$ and $\beta_1=1$, so $p(1)=n-1$. This makes sense, because at position $1$ we have the vertex $(n-1,\dots,n-1)$ and when we have used the Griffiths formula $(n-1)$-times, we reach $(0,0,0,0)$ and there is definitely a connected part of the path ending here and we need one basis element in degree $0=(n-(n-1)-1)d$.\\
Now assume that $p(k)=i$ and we already know that this corresponds to the fact that after using the Griffiths formula the appropriate number of times, there is a connected path ending at position $\beta_k$ with a vertex of degree $(n-i-1)d$. Consider the next number $k+1$, then
\begin{align*}
p(k+1)=p(k)+\left|\{j|\,\beta_{k}<\alpha_j<\beta_{k+1}\}\right|-1=i+\left|\{j|\,\beta_{k}<\alpha_j<\beta_{k+1}\}\right|-1
\end{align*}
and we want to show that this leads to a basis element in degree $(n-i-\left|\{j|\,\beta_{k}<\alpha_j<\beta_{k+1}\}\right|)d$. To prove this, we have 3 distinct cases:
\begin{enumerate}[label=(\roman*)]
\item $p(k+1)=p(k)-1$
\item $p(k+1)=p(k)$
\item $p(k+1)>p(k)$
\end{enumerate}
In the first case, there are no $\alpha_i$ between $\beta_k$ and $\beta_{k+1}$, so every position in between is covered by exactly one partial derivative, therefore the smallest number drops by one at position $\beta_k+1$, i.e. $\kappa(\beta_k)-1=\kappa(\beta_k+1)$, and stays the same until $\beta_{k+1}$ is reached. This marks the end of the connected part of the path. Since the smallest number is one less in this case we have to use the Griffiths formula one time less before we reach the basis element and therefore the degree of the basis element for this part of the path is $d$ times bigger than before, so the degree of this basis element is $(n-i)d=(n-(i-1)-1)d$ which agrees with the fact that $p(k+1)=i-1$.

In case (ii) there is exactly one $\alpha_i$ between $\beta_k$ and $\beta_{k+1}$. If $\alpha_i\in\Z$, then this means there is one position between $\beta_k$ and $\beta_{k+1}$ that is occupied by two partial derivatives. So again $\kappa(\beta_k)-1=\kappa(\beta_k+1)$, but before this part of the path ends, the smallest number increases by one due to the double occupation. So with the argument from before, we get a basis element of the same degree $(n-i-1)d=(n-p(k+1)-1)d$. If $\alpha_i\in\Q\setminus\Z$, then either the partial derivative corresponding to $\alpha_i$ is at position $\beta_{k}+1$ and the smallest number did not drop, or it is somewhere between $\beta_k$ and $\beta_{k+1}$. In this case it is in the same place as another partial derivative, because every number not in the set $\{\beta_i\}$ is occupied and we are back to the first consideration. So either way we have a basis element with the same smallest number as before, which therefore also has degree $(n-i-1)d=(n-p(k+1)-1)d$.

Finally, the third case is just an expansion of the previous case. Let us define the number of $\alpha_i$ between $\beta_k$ and $\beta_{k+1}$ as $a_k:=\left|\{j|\,\beta_{k}<\alpha_j<\beta_{k+1}\}\right|$, then with the same argumentation as before, we can see that the smallest number increases by $a_k-1$ on the path between $\beta_k$ and $\beta_{k+1}$, i.e. $\kappa(\beta_k)+a_k-1=\kappa(\beta_{k+1})$. It follows that the degree of the basis element of this part of the path is $(a_k-1)d$ times smaller than the previous basis element. So this basis element has degree $(n-i-1)d-(a_k-1)d=(n-(i+a_k-1)-1)d=(n-i-a_k)d$, which agrees with the above formula and ends the proof.
\end{proof}

\begin{rem}
We have $p(k)+p(u-k+1)=n$. This is due to the fact that the $\alpha_i\neq 0$ and the $\beta_i$ are evenly spread between $1$ and $\widehat{d}-1$. This implies 
$\left|\{j|\,\beta_{k}<\alpha_j<\beta_{k+1}\}\right|=\left|\{j|\,\beta_{u-k}<\alpha_j<\beta_{u-k+1}\}\right|$ for $1<k<u/2$ and therefore
\begin{align*}
p(k)+p(u-k+1)&=\left|\{j|\,\alpha_j<\beta_1\}\right|+\sum_{i=1}^{k-1} \left|\{j|\,\beta_{i}<\alpha_j<\beta_{i+1}\}\right|-(k-1)\\
&\quad+\left|\{j|\,\alpha_j<\beta_1\}\right|+\sum_{i=1}^{u-k} \left|\{j|\,\beta_{i}<\alpha_j<\beta_{i+1}\}\right|-(u-k)\\
&=2(n-1)+\sum_{i=u-k+1}^{u-1}\left|\{j|\,\beta_{i}<\alpha_j<\beta_{i+1}\}\right|\\
&\quad+\sum_{i=1}^{u-k} \left|\{j|\,\beta_{i}<\alpha_j<\beta_{i+1}\}\right|-(u-1)\\
&=2(n-1)+u-(n-1)-(u-1)=n.
\end{align*}
This also implies
\begin{align*}
h^{i-1,n-i-1}=h^{n-i-1,i-1}.
\end{align*}
\end{rem}

In \cite{CG} one can also find a more detailed description of the Hodge numbers that appear here. This relies mainly on the work of Danilov \cite{MR873655} on Deligne-Hodge numbers and Newton polyhedra.

\begin{rem}
The Hodge numbers $h^{i-1,n-i-1}=u_i$ that appear in our work as well as in \cite{CG} are the Deligne-Hodge numbers of the cohomology with compact support of a hypersurface defined by a Laurent polynomial with Newton polyhedron $\Delta$, where $\Delta=\left\langle\left(\frac{\widehat{q}_1}{\widehat{d}},\dots,\frac{\widehat{q}_n}{\widehat{d}}\right),(1,0\dots,0),\dots,(0,\dots,0,1)\right\rangle$. In particular from this viewpoint the $u_i$ are Deligne-Hodge numbers of a toric variety with polytope $\Delta$ in the lattice $\Z\left(\frac{\widehat{q}_1}{\widehat{d}},\dots,\frac{\widehat{q}_n}{\widehat{d}}\right)+\Z^n$.
\end{rem}

\section{The case of Arnold's strange duality}\label{ASD}

In this section we will show all the results and some more details for the 14 exceptional unimodal hypersurfaces singularities. This is a class of examples first studied by Arnol'd in \cite{MR0420689} where he among other things discovered that there is a duality among this 14 exceptional unimodal hypersurfaces singularities which is now known as Arnold's strange duality. One can define Gabrielov and Dolgachev numbers for every one of these hypersurface singularities and he showed that for every one of the 14 exceptional unimodal hypersurface singularities there is another singularity in this list with interchanged Dolgachev and Gabrielov numbers. The consequences of this duality between the $14$ exceptional unimodal hypersurface singularities have been studied by a number of people. An overview on a lot of aspects of this duality can be found in a paper by Ebeling \cite{MR1709346}. These examples were also the starting point for the analysis of the Picard-Fuchs equations in this thesis. We want to concentrate in this section on the duality between the invertible polynomials of Arnold's strange duality and the consequences we get from the results achieved so far.

\pagebreak
In the following table we list important data of the 14 exceptional unimodal hypersurface singularities we need. In particular we have a look at the compactification that comes from compactifying the Milnor fibres in the weighted projective space with one additional dimension which has weight one. The table consists of the polynomial defining the compactification of the hypersurface singularity, the degree of this polynomial, the weights and the dual singularity due to Arnol'd. As in Table \ref{ASD-normal} there is still a the duality between the invertible polynomials.

\begin{table}[H]
\begin{center}
\begin{tabular}{lccccl}
\toprule
Name & $g(w,x,y,z)$ & Deg & Weights & Dual\\
\midrule
E$_{12}$ & $w^{42}+x^7+y^3+z^2$ & 42 & (1,6,14,21) & E$_{12}$\\
&&&&\\
E$_{13}$ & $w^{30}+x^5y+y^3+z^2$ & 30 & (1,4,10,15) & Z$_{11}$\\
Z$_{11}$ & $w^{30}+x^5+xy^3+z^2$ & 30 & (1,6,8,15) & E$_{13}$\\
&&&&\\
E$_{14}$ & $w^{24}+x^4z+y^3+z^2$ & 24 & (1,3,8,12) & Q$_{10}$\\
Q$_{10}$ & $w^{24}+x^4+y^3+xz^2$ & 24 & (1,6,8,9) & E$_{14}$\\
&&&&\\
Z$_{12}$ & $w^{22}+x^4y+xy^3+z^2$ & 22 & (1,4,6,11) & Z$_{12}$\\
&&&&\\
W$_{12}$ & $w^{20}+x^5+y^2z+z^2$ & 20 & (1,4,5,10) & W$_{12}$\\
&&&&\\
Z$_{13}$ & $w^{18}+x^3z+xy^3+z^2$ & 18 & (1,3,5,9) & Q$_{11}$\\
Q$_{11}$ & $w^{18}+x^3y+y^3+xz^2$ & 18 & (1,4,6,7) & Z$_{13}$\\
&&&&\\
W$_{13}$ & $w^{16}+x^4y+y^2z+z^2$ & 16 & (1,3,4,8) & S$_{11}$\\
S$_{11}$ & $w^{16}+x^4+y^2z+xz^2$ & 16 & (1,4,5,6) & W$_{13}$\\
&&&&\\
Q$_{12}$ & $w^{15}+x^3z+y^3+xz^2$& 15 & (1,3,5,6) & Q$_{12}$\\
&&&&\\
S$_{12}$ & $w^{13}+x^3y+y^2z+xz^2$ & 13 & (1,3,4,5) & S$_{12}$\\
&&&&\\
U$_{12}$ & $w^{12}+x^4+y^2z+yz^2$ & 12 & (1,3,4,4) & U$_{12}$\\
\bottomrule
\end{tabular}
\end{center}
\caption{The Compactification of Arnold's strange duality}
\label{ASD-compact}
\end{table}

We will now list the sets $\{\alpha_1,\dots,\alpha_u\}=A\setminus(A\cap D_0)$ and $\{\beta_1,\dots,\beta_u\}=D_0\setminus(A\cap D_0)$, where $A=\bigsqcup_{i=1}^n\left(\left(Q_i\setminus \{\widehat{d}\}\right)\cup\{0\}\right)$ and $D_0=(D\setminus \{\widehat{d}\})\cup\{0\}$ as in Proposition \ref{hodge-numbers-goly}, and the resulting order of the Picard-Fuchs equation $u$. Remember from Definition \ref{vundu} that the sets $Q_i$ and $D$ are defined via the dual weights. Notice that $\delta+\alpha_i$ and $\delta-\beta_i$ are the linear factors in the two summands of the Picard-Fuchs equation. Together with the dual weights and the dual degree they completely determine the Picard-Fuchs equation. We want to mention that $\{\beta_1,\dots,\beta_u\}$ contains all numbers $1\leq b \leq \widehat{d}$ which are coprime to $\widehat{d}$. The set $\{\alpha_1,\dots,\alpha_u\}\cup \Z$ on the other hand contains only elements which are not coprime to $\widehat{d}$.

\begin{table}[H]
\begin{center}
\begin{tabular}{lccc}
\toprule
Name & $\alpha_1,\dots,\alpha_u$ & $\beta_1,\dots,\beta_u$ & $u$\\
\midrule
E$_{12}$ & $0,0,0,6,12,14,18,21,24,28,30,36$ & $1,5,11,13,17,19,23,25,29,31,37,41$ & $12$ \\
&&&\\
E$_{13}$ & $0,0,0,\frac{15}{4},\frac{15}{2},10,\frac{45}{4},15,\frac{75}{4},20,\frac{45}{2},\frac{105}{4}$ & $1,3,7,9,11,13,17,19,21,23,27,29$ & $12$ \\
Z$_{11}$ & $0,0,0,6,\frac{15}{2},12,15,18,\frac{45}{2},24$ & $1,5,7,11,13,17,19,23,25,29$ & $10$ \\
&&&\\
E$_{14}$ & $0,0,0,\frac{8}{3},\frac{16}{3},8,\frac{32}{3},12,\frac{40}{3},16,\frac{56}{3},\frac{64}{3},$ & $1,2,5,7,10,11,13,14,17,19,22,23$ & $12$ \\
Q$_{10}$ & $0,0,0,6,8,12,16,18$ & $1,5,7,11,13,17,19,23$ & $8$ \\
&&&\\
Z$_{12}$ & $0,0,0,\frac{11}{3},\frac{11}{2},\frac{22}{3},11,\frac{44}{3},\frac{33}{2},\frac{55}{3}$ & $1,3,5,7,9,13,15,17,19,21$ & $10$ \\
&&&\\
W$_{12}$ & $0,0,0,4,8,10,12,16$ & $1,3,7,9,11,13,17,19$ & $8$ \\
&&&\\
Z$_{13}$ & $0,0,0,\frac{18}{7},\frac{9}{2},\frac{36}{7},\frac{54}{7},9,\frac{72}{7},\frac{90}{7},\frac{27}{2},\frac{108}{7}$ & $1,2,4,5,7,8,10,11,13,14,16,17$ & $12$ \\
Q$_{11}$ & $0,0,0,\frac{18}{5},6,\frac{36}{5},\frac{54}{5},12,\frac{72}{5}$ & $1,3,5,7,9,11,13,15,17$ & $9$ \\
&&&\\
W$_{13}$ & $0,0,0,\frac{8}{3},\frac{16}{5},\frac{16}{3},\frac{32}{5},8,\frac{48}{5},\frac{32}{3},\frac{64}{5},\frac{40}{3}$ & $1,2,3,5,6,7,9,10,11,13,14,15$ & $12$ \\
S$_{11}$ & $0,0,0,4,\frac{16}{3},8,\frac{32}{3},12$ & $1,3,5,7,9,11,13,15$ & $8$ \\
&&&\\
Q$_{12}$ & $0,0,0,\frac{5}{2},5,\frac{15}{2},10,\frac{25}{2}$ & $1,2,4,7,8,11,13,14$ & $8$ \\
&&&\\
S$_{12}$ & $0,0,0,\frac{13}{5},\frac{13}{4},\frac{13}{3},\frac{26}{5},\frac{13}{2},\frac{39}{5},\frac{26}{3},\frac{39}{4},\frac{52}{5}$ & $1,2,3,4,5,6,7,8,9,10,11,12$ & $12$ \\
&&&\\
U$_{12}$ & $0,0,0,3,6,9$ & $1,2,5,7,10,11$ & $6$ \\
\bottomrule
\end{tabular}
\end{center}
\caption{The sets defining the linear factors of the Picard-Fuchs equation}
\label{ASD-sets}
\end{table}

With the table above and Proposition \ref{hodge-numbers-goly} we are able to state how many basis elements we need in every degree. Using the methods from Chapter \ref{haupt}, we can give an exact basis of the part of the cohomology that is used in the calculations. This can be done by relating the numbers $\alpha_i$ and $\beta_i$ to the positions on the Jacobi path, where a basis element can be chosen. In the examples we have $h^{2,0}=h^{0,2}=1$ and we can choose the basis elements $\frac{s\Omega_0}{f}$ and $\frac{s^3(wxyz)^2\Omega_0}{f^3}$ respectively. This means we have to find $u-2$ basis elements in degree $d$. We will only list a basis for the part of the middle cohomology we used and not a basis for the whole cohomology, or equivalently Milnor ring. We show how to calculate these basis elements in an example. The rest of the basis elements in Table \ref{ASD-basis-elements} can be calculated the same way.

\pagebreak
\begin{ex}
We concentrate on the singularity S$_{11}$, so $f(\x)=w^{16}+x^4+y^2z+xz^2$ with weights $(1,4,5,6)$ and degree $d=16$. The dual weights are $(1,3,8,4)$ and $\widehat{d}=d=16$. We have
\begin{align*}
Q_1=\{16\},\;Q_2=\{\frac{16}{3},\frac{22}{3},16\},\;Q_3=\{2,4,6,8,10,12,14,16\}\,\text{ and }\,Q_4=\{4,8,12,16\}.
\end{align*}
We can read off from the above table what the $\alpha_i$ are and we can see that we need $6$ basis elements in the middle cohomology. Apart from $0$, the $\alpha_i$ consist of the elements in the disjoint union of the $Q_i$ which are rational or appear twice. So the number $\alpha_i$ tells us the position of the basis element on the Jacobi path. This becomes clear in the calculation:\\
We start with $\alpha_4=4$. To calculate the corresponding basis elements we need to check which arrows have already been used, in other words how many numbers in $Q_i$ are $\leq 4$. In $Q_1$ there is no element, so $\partial_w$ does not appear on the Jacobi path before the vertex of the basis element we are looking for. The set $Q_2$ has also no element $\leq 4$, but $Q_3$ contains $2$ and $4$ and $Q_4$ contains also $4$. This means $\partial_y$ was used twice and $\partial_z$ was used once on the Jacobi path before we arrive at the basis element. So starting at $(1,1,1,1)$ we have to add up $(1,1,-1,0)$ twice and $(1,0,1,-1)$ once and we end up at the vertex

\begin{align*}
(1,1,1,1)+2\cdot\partial_y+\partial_z=(1,1,1,1)+2\cdot(1,1,-1,0)+(1,0,1,-1)=(4,3,0,0)
\end{align*}

which means that the basis element we are looking for is given by $w^4x^3$.\\
For $\alpha_5=\frac{16}{3}$, we check that counting the elements $\leq \frac{16}{3}$ in every $Q_i$, we have to use $\partial_x$ once, $\partial_y$ twice and $\partial_z$ once to get the basis element. We calculate that

\begin{align*}
(1,1,1,1)+\partial_z+2\cdot\partial_y+\partial_z&=(1,1,1,1)+(1,-3,1,1)+2\cdot(1,1,-1,0)+(1,0,1,-1)\\
&=(5,0,1,1)
\end{align*}

which leads to $w^5yz$ as basis element. In the same way one gets that $\alpha_6=8$ corresponds to the basis element $w^8x^2$, $\alpha_7=\frac{32}{3}$ gives $w^{10}y$ and finally $\alpha_8=12$ leads to $w^{12}x$. Now we have $5$ basis elements and in addition we have $wxyz$ which corresponds to one of the zeroes in $A$.
\end{ex}

With this construction we are able to calculate the basis in the middle cohomology for all examples, and this is listed in the following table.
So the table includes the name of the singularity, the number $h^{1,1}=u-2$ from Proposition \ref{hodge-numbers-goly} and a basis for the part of the Milnor ring in degree $d$, which gives also a basis of the part of the middle cohomology we are using.

\begin{table}[H]
\begin{center}
\begin{tabular}{lcl}
\toprule
Name & $u-2$ & basis elements in the Milnor ring in degree $d$ \\
\midrule
E$_{12}$ & $10$ & $wxyz,w^{36}x,w^{30}x^2,w^{24}x^3,w^{18}x^4,w^{12}x^5,w^6x^6,w^{28}y,w^{14}y^2,w^{21}z$ \\
&&\\
E$_{13}$ & $10$ & $wxyz,w^{26}x,w^{22}x^2,w^{18}x^3,w^{20}y,w^{10}y^2,w^{15}z,w^{11}xz,w^7x^2z,w^3x^3z$ \\
Z$_{11}$ & $8$ & $wxyz,w^{24}x,w^{18}x^2,w^{12}x^3,w^6x^4,w^{22}y,w^{15}z,w^7yz$ \\
&&\\
E$_{14}$ & $10$ & $wxyz,w^{21}x,w^{18}x^2,w^{16}y,w^8y^2,w^{12}z,w^{13}xy,w^{10}x^2y,w^5xy^2,w^2x^2y^2$ \\
Q$_{10}$ & $6$ & $wxyz,w^{18}x,w^{12}x^2,w^6x^3,w^{16}y,w^8y^2$ \\
&&\\
Z$_{12}$ & $8$ & $wxyz,w^{18}x,w^{14}x^2,w^{16}y,w^{11}z,w^7xz,w^3x^2z,w^5yz$ \\
&&\\
W$_{12}$ & $6$ & $wxyz,w^{16}x,w^{12}x^2,w^8x^3,w^4x^4,w^{10}y^2$ \\
&&\\
Z$_{13}$ & $10$ & $wxyz,w^{15}x,w^{12}x^2,w^{13}y,w^9z,w^{10}xy,w^7x^2y,w^5xy^2,w^2x^2y^2,w^4yz$ \\
Q$_{11}$ & $7$ & $wxyz,w^{14}x,w^{10}x^2,w^{12}y,w^6y^2,w^7xz,w^3x^2z$ \\
&&\\
W$_{13}$ & $10$ & $wxyz,w^{13}x,w^{10}x^2,w^7x^3,w^{12}y,w^9xy,w^6x^2y,w^3x^3y,w^5xz,w^2x^2z$ \\
S$_{11}$ & $6$ & $wxyz,w^{12}x,w^8x^2,w^4x^3,w^{10}z,w^5yz$ \\
&&\\
Q$_{12}$ & $6$ & $wxyz,w^{12}x,w^{10}y,w^5y^2,w^7xy,w^2xy^2$ \\
&&\\
S$_{12}$ & $10$ & $wxyz,w^{10}x,w^7x^2,w^9y,w^8z,w^6xy,w^3x^2y,w^5xz,w^2x^2z,w^4yz$ \\
&&\\
U$_{12}$ & $4$ & $wxyz,w^9x,w^6x^2,w^3x^3$ \\
\bottomrule
\end{tabular}
\end{center}
\caption{Basis elements for the middle cohomology}
\label{ASD-basis-elements}
\end{table}

Of course from the previous work we can immediately calculate the Picard-Fuchs equation, either with the Griffiths-Dwork method as shown in Section \ref{detailed-example-gd}, with a computer algebra system as shown in Appendix \ref{singular}, or by inserting in the $\alpha_i$ and $\beta_j$ as linear factors as in Theorem \ref{conj}. The output for all singularities we investigated in this section is shown in the next table.

\begin{landscape}
\begin{table}
\centering
\begin{tabular}{ll}
\toprule
Name& Picard-Fuchs equation for $f$\\
\midrule
\multirow{2}{*}{E$_{12}$} & $s^{42}\delta^3(\delta+6)(\delta+12)(\delta+14)(\delta+18)(\delta+21)(\delta+24)(\delta+28)(\delta+30)(\delta+36)$\\
& $-2^{22}3^{15}7^7(\delta-1)(\delta-5)(\delta-11)(\delta-13)(\delta-17)(\delta-19)(\delta-23)(\delta-25)(\delta-29)(\delta-31)(\delta-37)(\delta-41)$\\
\midrule
\multirow{2}{*}{E$_{13}$} & $s^{30}\delta^3(\delta+\frac{15}{4})(\delta+\frac{15}{2})(\delta+10)(\delta+\frac{45}{4})(\delta+15)(\delta+\frac{75}{4})(\delta+20)(\delta+\frac{45}{2})(\delta+\frac{105}{4})$\\
& $-3^{9}5^{15}(\delta-1)(\delta-3)(\delta-7)(\delta-9)(\delta-11)(\delta-13)(\delta-17)(\delta-19)(\delta-21)(\delta-23)(\delta-27)(\delta-29)$\\
\midrule
\multirow{2}{*}{Z$_{11}$} & $s^{30}\delta^3(\delta+6)(\delta+\frac{15}{2})(\delta+12)(\delta+15)(\delta+18)(\delta+\frac{45}{2})(\delta+24)$\\
& $-2^{12}3^{15}5^5(\delta-1)(\delta-5)(\delta-7)(\delta-11)(\delta-13)(\delta-17)(\delta-19)(\delta-23)(\delta-25)(\delta-29)$\\
\midrule
\multirow{2}{*}{E$_{14}$} & $s^{24}\delta^3(\delta+\frac{8}{3})(\delta+\frac{16}{3})(\delta+8)(\delta+\frac{32}{3})(\delta+12)(\delta+\frac{40}{3})(\delta+16)(\delta+\frac{56}{3})(\delta+\frac{64}{3})$\\
& $-2^{42}(\delta-1)(\delta-2)(\delta-5)(\delta-7)(\delta-10)(\delta-11)(\delta-13)(\delta-14)(\delta-17)(\delta-19)(\delta-22)(\delta-23)$\\
\midrule
Q$_{10}$ & $s^{24}\delta^3(\delta+6)(\delta+8)(\delta+12)(\delta+16)(\delta+18)-2^{24}3^{9}(\delta-1)(\delta-5)(\delta-7)(\delta-11)(\delta-13)(\delta-17)(\delta-19)(\delta-23)$\\
\midrule
\multirow{2}{*}{Z$_{12}$} & $3^6s^{22}\delta^3(\delta+\frac{11}{3})(\delta+\frac{11}{2})(\delta+\frac{22}{3})(\delta+11)(\delta+\frac{44}{3})(\delta+\frac{33}{2})(\delta+\frac{55}{3})$\\
& $-2^811^{11}(\delta-1)(\delta-3)(\delta-5)(\delta-7)(\delta-9)(\delta-13)(\delta-15)(\delta-17)(\delta-19)(\delta-21)$\\
\midrule
W$_{12}$ & $s^{20}\delta^3(\delta+4)(\delta+8)(\delta+10)(\delta+12)(\delta+16)-2^{22}5^5(\delta-1)(\delta-3)(\delta-7)(\delta-9)(\delta-11)(\delta-13)(\delta-17)(\delta-19)$\\
\midrule
\multirow{2}{*}{Z$_{13}$} & $7^7s^{18}\delta^3(\delta+\frac{18}{7})(\delta+\frac{9}{2})(\delta+\frac{36}{7})(\delta+\frac{54}{7})(\delta+9)(\delta+\frac{72}{7})(\delta+\frac{90}{7})(\delta+\frac{27}{2})(\delta+\frac{108}{7})$\\
& $-2^43^{30}(\delta-1)(\delta-2)(\delta-4)(\delta-5)(\delta-7)(\delta-8)(\delta-10)(\delta-11)(\delta-13)(\delta-14)(\delta-16)(\delta-17)$\\
\midrule
\multirow{2}{*}{Q$_{11}$} & $5^5s^{18}\delta^3(\delta+\frac{18}{5})(\delta+6)(\delta+\frac{36}{5})(\delta+\frac{54}{5})(\delta+12)(\delta+\frac{72}{5})$\\
& $-2^{18}3^{15}(\delta-1)(\delta-3)(\delta-5)(\delta-7)(\delta-9)(\delta-11)(\delta-13)(\delta-15)(\delta-17)$\\
\bottomrule
\end{tabular}
\end{table}
\end{landscape}

\begin{landscape}
\begin{table}
\centering
\begin{tabular}{ll}
\toprule
Name & Picard-Fuchs equation for $f$\\
\midrule
\multirow{2}{*}{W$_{13}$} & $3^65^5s^{16}\delta^3(\delta+\frac{8}{3})(\delta+\frac{16}{5})(\delta+\frac{16}{3})(\delta+\frac{32}{5})(\delta+8)(\delta+\frac{48}{5})(\delta+\frac{32}{3})(\delta+\frac{64}{5})(\delta+\frac{40}{3})$\\
& $-2^{50}(\delta-1)(\delta-2)(\delta-3)(\delta-5)(\delta-6)(\delta-7)(\delta-9)(\delta-10)(\delta-11)(\delta-13)(\delta-14)(\delta-15)$\\
\midrule
S$_{11}$ & $3^3s^{16}\delta^3(\delta+4)(\delta+\frac{16}{3})(\delta+8)(\delta+\frac{32}{3})(\delta+12)-2^{32}(\delta-1)(\delta-3)(\delta-5)(\delta-7)(\delta-9)(\delta-11)(\delta-13)(\delta-15)$\\
\midrule
Q$_{12}$ & $2^6s^{15}\delta^3(\delta+\frac{5}{2})(\delta+5)(\delta+\frac{15}{2})(\delta+10)(\delta+\frac{25}{2})-3^{6}5^{10}(\delta-1)(\delta-2)(\delta-4)(\delta-7)(\delta-8)(\delta-11)(\delta-13)(\delta-14)$\\
\midrule
\multirow{2}{*}{S$_{12}$} & $3^34^45^5s^{13}\delta^3(\delta+\frac{13}{5})(\delta+\frac{13}{4})(\delta+\frac{13}{3})(\delta+\frac{26}{5})(\delta+\frac{13}{2})(\delta+\frac{39}{5})(\delta+\frac{26}{3})(\delta+\frac{39}{4})(\delta+\frac{52}{5})$\\
& $+13^{13}(\delta-1)(\delta-2)(\delta-3)(\delta-4)(\delta-5)(\delta-6)(\delta-7)(\delta-8)(\delta-9)(\delta-10)(\delta-11)(\delta-12)$\\
\midrule
U$_{12}$ & $s^{12}\delta^3(\delta+3)(\delta+6)(\delta+9)-2^{8}3^{9}(\delta-1)(\delta-2)(\delta-5)(\delta-7)(\delta-10)(\delta-11)$\\
\bottomrule
\end{tabular}
\end{table}
\end{landscape}

We give another viewpoint on the Picard-Fuchs equation. From Theorem \ref{conj} we know that the Picard-Fuchs equation always consists of exactly two summands. They can be separated by setting $s^{\widehat{d}}=0$ or $s^{\widehat{d}}=\infty$. We already know that if we view the Picard-Fuchs equation as a polynomial with variable $\delta$ then the zeroes of these polynomials after setting $s^{\widehat{d}}=0$ are given by $\beta_1,\dots,\beta_u$ and the zeroes for $s^{\widehat{d}}=\infty$ are given by $-\alpha_1,\dots,-\alpha_u$. Now we want to focus on a polynomial that has related zeroes. Namely, we define $\chi_0$ to be the polynomial with zeroes $\exp\left(2\pi\textrm{i}\frac{\beta_i}{\widehat{d}}\right)$ for $i=0,\dots,n$ and $\chi_{\infty}$ the polynomial with roots $\exp\left(2\pi\textrm{i}\frac{\alpha_i}{\widehat{d}}\right)$ for $i=0,\dots,n$ and notice that multiple roots in the Picard-Fuchs equation lead to multiple roots of $\chi_{\infty}$ and $\chi_0$. Equivalently, we can first write the Picard-Fuchs equation for the variable $\lambda=(-s)^{-\widehat{d}}$ and then start with the zeroes of this equation for $\lambda=\infty$ and $\lambda=0$.

\begin{notation}
We will shorten the notation for a rational function with only roots of unity as zeroes and poles. We will write $\nu_1\cdots\nu_{m_1}/\eta_1\cdots\eta_{m_2}$ for the rational function
\begin{align*}
\chi(t)=\frac{(1-t^{\nu_1})\cdot\dots\cdot(1-t^{\nu_{m_1}})}{(1-t^{\eta_1})\cdot\dots\cdot(1-t^{\eta_{m_2}})}
\end{align*}
\end{notation}

With this notation we write down the functions $\chi_0$ and $\chi_{\infty}$ in the following table.

\begin{table}[H]
\begin{center}
\begin{tabular}{lcccc}
\toprule
Name & Deg & Weights & $\chi_0$ & $\chi_{\infty}$\\
\midrule
E$_{12}$ & 42 & (1,6,14,21) & $2\cdot 3\cdot 7\cdot 42/1\cdot 6\cdot 14\cdot 21$ & $2\cdot 3\cdot 7$ \\
&&&&\\
E$_{13}$ & 30 & (1,4,10,15) & $3\cdot 30/6\cdot 15$ & $1\cdot3\cdot8$ \\
Z$_{11}$ & 30 & (1,6,8,15) & $5\cdot 30/10\cdot 15$ & $1\cdot4\cdot5$ \\
&&&&\\
E$_{14}$ & 24 & (1,3,8,12) & $2\cdot 24/6\cdot 8$ & $1\cdot2\cdot9$ \\
Q$_{10}$ & 24 & (1,6,8,9) & $4\cdot 24/8\cdot 12$ & $1\cdot3\cdot4$ \\
&&&&\\
Z$_{12}$ & 22 & (1,4,6,11) & $1\cdot 22/2\cdot 11$ & $1\cdot 1\cdot4\cdot6/2$ \\
&&&&\\
W$_{12}$ & 20 & (1,4,5,10) & $2\cdot20/4\cdot10$ & $1\cdot2\cdot5$ \\
&&&&\\
Z$_{13}$ & 18 & (1,3,5,9) & $18/6$ & $1\cdot4\cdot7$ \\
Q$_{11}$ & 18 & (1,4,6,7) & $18/9$ & $1\cdot3\cdot5$ \\
&&&&\\
W$_{13}$ & 16 & (1,3,4,8) & $16/4$ & $1\cdot5\cdot6$ \\
S$_{11}$ & 16 & (1,4,5,6) & $16/8$ & $1\cdot3\cdot4$ \\
&&&&\\
Q$_{12}$ & 15 & (1,3,5,6) & $1\cdot15/3\cdot5$ & $1\cdot 1\cdot6$ \\
&&&&\\
S$_{12}$ & 13 & (1,3,4,5) & $13/1$ & $3\cdot4\cdot5$ \\
&&&&\\
U$_{12}$ & 12 & (1,3,4,4) & $1\cdot 12/3\cdot 4$ & $1\cdot 1\cdot4$ \\
\bottomrule
\end{tabular}
\end{center}
\end{table}

The functions $\chi_0$ and $\chi_{\infty}$ are in all cases a little bit different, but the interesting thing is that the quotient of the two functions is always the same.

\begin{rem}
The rational functions $\chi_0$ and $\chi_{\infty}$ described in the above table have always the property that
\begin{align*}
\frac{\chi_0(t)}{\chi_{\infty}(t)}=\frac{(1-t^{\widehat{d}})}{(1-t^{\widehat{q}_1})(1-t^{\widehat{q}_2})(1-t^{\widehat{q}_3})(1-t^{\widehat{q}_4})}
\end{align*}
\end{rem}

In the next section we will look at this phenomenon in more generality and we will also see that the roots of $\chi_0$ and $\chi_{\infty}$ are the eigenvalues of the local monodromy around $(-1)^{\widehat{d}}\lambda^{-1}=s^{\widehat{d}}=0$ and $(-1)^{\widehat{d}}\lambda^{-1}=s^{\widehat{d}}=\infty$ respectively.

\section{Relations to the Poincar\'{e} series and monodromy}\label{poincare}
In this section we want to relate the numbers in the Picard-Fuchs equation of $f(\x)$ to the Poincar\'{e} series of $g^t(\x)$ and to the monodromy around $0$ and $\infty$ in the solution space of the Picard-Fuchs equation. The last remark in the previous chapter already showed us the direction. 

\subsection*{Poincar\'{e} series}
First we want to investigate the relation to the Poincar\'{e} series. Therefore we consider the Picard-Fuchs equation in the form of (\ref{pf-in-t}) which is a differential equation with parameter $\lambda=(-s)^{-\widehat{d}}$. If we view this differential equation as a polynomial with variable $\mathcal{D}$, then we can immediately read off the zeroes for $\lambda=0$ and $\lambda=\infty$:

\begin{align*}
\lambda&=0:\quad \frac{\alpha_1}{\widehat{d}},\dots,\frac{\alpha_u}{\widehat{d}}\\
\lambda&=\infty:\quad -\frac{\beta_1}{\widehat{d}},\dots,-\frac{\beta_u}{\widehat{d}}
\end{align*}

\begin{rem}
Because of the symmetry of the $\alpha_j$ and $\beta_j$, the sets $\left\{\exp\left(2\pi \textrm{i} \frac{\alpha_j}{\widehat{d}}\right)\right\}$ and $\left\{\exp\left(2\pi\textrm{i} \frac{\beta_j}{\widehat{d}}\right)\right\}$ are closed under complex conjugation.
\end{rem}

We will now relate these numbers $\alpha_j$ and $\beta_j$ or $\exp\left(2\pi\textrm{i} \frac{\alpha_j}{\widehat{d}}\right)$ and $\exp\left(2\pi\textrm{i} \frac{\beta_j}{\widehat{d}}\right)$ respectively to the Poincar\'{e} series of $g^t(\x)$. Let us recall first how the Poincar\'{e} series is defined.

\begin{defi}
Let $A:=\C[\x]/(g(\x))$ be the coordinate algebra of the hypersurface $\{g(\x)=0\}$. Then $A$ admits naturally a grading $A=\bigoplus_{m=0}^\infty A_m$, where $A_m$ is generated by the monomials in $A$ of weighted degree $m$. The Poincar\'{e} series for this hypersurface is given by
\begin{align*}
p_A(t):=p_g(t):=\sum_{m=0}^{\infty} \dim_\C A_m t^m
\end{align*}
\end{defi}

\begin{rem}\emph{(cf. \cite{MR777682})}
If $g(\x)$ is quasihomogeneous with weights $q_1,\dots,q_n$ and weighted degree $d$, then the Poincar\'{e} series is given by
\begin{align*}
p_g(t)=\frac{(1-t^{d})}{(1-t^{q_1})\cdot\dots\cdot(1-t^{q_n})}
\end{align*}
\end{rem}

A rational function of this form is of course uniquely determined by the set of poles and zeroes. So we study these sets for the Poincar\'{e} series of $g^t(\x)$, because as mentioned before this will be related to the Picard-Fuchs equation of $f(\x)$. So we study the zeroes and poles of the function
\begin{align*}
p_{g^t}(t)=\frac{(1-t^{\widehat{d}})}{(1-t^{\widehat{q}_1})\cdot\dots\cdot(1-t^{\widehat{q}_n})}.
\end{align*}
The zeroes of $(1-t^{\widehat{d}})$ are given by the set $\left\{\exp\left(2\pi\textrm{i}\frac{j}{\widehat{d}}\right)|\,0\leq j\leq \widehat{d}-1\right\}$ and the zeroes of $(1-t^{\widehat{q}_1})\cdot\dots\cdot(1-t^{\widehat{q}_n})$ are given by the set $\bigcup_{k=1}^n \left\{\exp\left(2\pi\textrm{i}\frac{j}{\widehat{q}_k}\right)|\,0\leq j\leq \widehat{q}_k-1\right\}$. So putting this together, the zeroes of the Poincar\'{e} series of $g^t(\x)$ are given by
\begin{align*}
&\left\{\exp\left(2\pi\textrm{i}\frac{j}{\widehat{d}}\right)|\,j\in\Z\right\}\setminus\left(\left\{\exp\left(2\pi\textrm{i}\frac{j}{\widehat{d}}\right)|\,j\in\Z\right\}\cap\bigcup_{k=1}^n \left\{\exp\left(2\pi\textrm{i}\frac{j}{\widehat{q}_k}\right)|\,j\in\Z\right\}\right)\\
&=\left\{\exp\left(2\pi\textrm{i}\frac{b}{\widehat{d}}\right)|\,b\in D\right\}\setminus\left(\left\{\exp\left(2\pi\textrm{i}\frac{\beta}{\widehat{d}}\right)|\,b\in D\right\}\cap \left\{\exp\left(2\pi\textrm{i}\frac{a}{\widehat{d}}\right)|\,a\in A\right\}\right)\\
&=\left\{\exp\left(2\pi\textrm{i}\frac{\beta_j}{\widehat{d}}\right)|\,j=0,\dots,u\right\}
\end{align*}

and the poles are given by the set
\begin{align*}
&\bigsqcup_{k=1}^n \left\{\exp\left(2\pi\textrm{i}\frac{j}{\widehat{q}_k}\right)|\,j\in\Z\right\}\setminus\left(\left\{\exp\left(2\pi\textrm{i}\frac{j}{\widehat{d}}\right)|\,j\in\Z\right\}\cap\bigcup_{k=1}^n \left\{\exp\left(2\pi\textrm{i}\frac{j}{\widehat{q}_k}\right)|\,j\in\Z\right\}\right)\\
&=\left\{\exp\left(2\pi\textrm{i}\frac{a}{\widehat{d}}\right)|\,a\in A\right\}\setminus\left(\left\{\exp\left(2\pi\textrm{i}\frac{b}{\widehat{d}}\right)|\,b\in D\right\}\cap \left\{\exp\left(2\pi\textrm{i}\frac{a}{\widehat{d}}\right)|\,a\in A\right\}\right)\\
&=\left\{\exp\left(2\pi\textrm{i}\frac{\alpha_j}{\widehat{d}}\right)|\, j=0,\dots,u\right\},
\end{align*}

where the disjoint union indicates that poles occur in this set counted with multiplicity. Notice that the notation $\left\{\exp\left(2\pi\textrm{i}\frac{a}{\widehat{d}}\right)|\,a\in A\right\}$, where $A=\bigsqcup_{k=1}^n \left(\left(Q_k\setminus\{\widehat{d}\}\right)\cup\{0\}\right)$, is short for $\bigsqcup_{k=1}^n\left\{\exp\left(2\pi\textrm{i}\frac{a_k}{\widehat{d}}\right)|\,a_k\in\left(Q_k\setminus\{\widehat{d}\}\right)\cup\left\{0\right\}\right\}$.

In the above we can see clearly the relation between the zeroes of the Picard-Fuchs equation of $f(\x)$ for $\lambda=0$ and $\lambda=\infty$ and the Poincar\'{e} series of $g^t(\x)$. We summarize this in the following corollary.

\begin{kor}\label{poincare-zeroes}
The zeroes of the Poincar\'{e} series of $g^t(\x)$ are in $1-1$ correspondence with the zeroes of the Picard-Fuchs equation of $f(\x)$ for $\lambda=\infty$ or $s=0$ and the poles of the Poincar\'{e} series of $g^t(\x)$ are in $1-1$ correspondence to the zeroes of the Picard-Fuchs equation of $f(\x)$ for $\lambda=0$ or $s=\infty$.\\
Equivalently the same holds for the Picard-Fuchs equation of $f^t(\x)=g^t(\x)+s\prod x_i$ and the Poincar\'{e} series of $g(\x)$.
\end{kor}

\subsection*{Monodromy}
Now we want to explain why the roots of the Picard-Fuchs equation for $\lambda=(-s)^{-\widehat{d}}=0$ and $\lambda=(-s)^{-\widehat{d}}=\infty$ are in 1-1 correspondence with the eigenvalues of the local monodromy around $0$ and $\infty$ in the solution space of the Picard-Fuchs equation, i.e. the space of the period integrals. More precisely, the eigenvalues of the monodromy around $0$ and $\infty$ are equal to the poles and zeroes of the Poincar\'{e} series respectively. First we recall monodromy in the context of Picard-Fuchs equations in as much generality as we need. References for the relation between monodromy and the Picard-Fuchs equation are \cite{MR1677117}, \cite{MR1191426} and \cite {De}. 

In this subsection we will always regard the Picard-Fuchs equation in $\mathcal{D}=\lambda\frac{\partial}{\partial \lambda}$, so we are working with the differential equation (\ref{pf-in-t})
\begin{align*}
0=\prod_{i=1}^n \widehat{q}_i^{\widehat{q}_i}\prod_{i=1}^{n}\prod_{j=0}^{\widehat{q}_i-1}(\mathcal{D}-\frac{j}{\widehat{q}_i})\prod_{\ell\in I}(\mathcal{D}-\frac{\ell}{\widehat{d}})^{-1}\Phi-\widehat{d}^{\widehat{d}}\lambda\prod_{j=0}^{\widehat{d}-1}(\mathcal{D}+\frac{j}{\widehat{d}})\prod_{\ell\in I}(\mathcal{D}+\frac{\ell}{\widehat{d}})^{-1}\Phi.
\end{align*}
Due to \cite{De} this Picard-Fuchs equation has only regular singular points. This can for example be seen by the fact that in the Picard-Fuchs equation, written as

\begin{align}\label{pf-ausmultipliziert}
\mathcal{D}^u\Phi+\sum_{i=0}^{u-1} h_i(\lambda)\mathcal{D}^i \Phi=0,
\end{align}

all coefficients $h_i(\lambda)$ are holomorphic functions of $\lambda$. Now we can define the residue matrix for $\lambda$.
\begin{defi}
Let $\omega_1,\dots,\omega_u$ be a basis of the solution space of the Picard-Fuchs equation and define the connection matrix $(\Gamma)_{ij}$ via $\mathcal{D}\omega_i=\sum_j \Gamma_{ij} \omega_j$. Then the residue matrix is given by $\text{Res}=\text{Res}_{\lambda=0}\left((\Gamma)_{ij}\right)$.
\end{defi}

\begin{rem}
In the cases we consider $(\Gamma)_{ij}$ has no poles at $\lambda=0$, so the residue matrix is just given by $\text{Res}=\left((\Gamma)_{ij}\right)_{\lambda=0}$.
\end{rem}

\begin{thm}\emph{(\cite{De})}\label{Deligne}
The following relations between the residue matrix and the monodromy around $\lambda=0$ in the solution space of the Picard-Fuchs equation hold.
\begin{enumerate}[label=(\roman*)]
\item $\eta$ is an eigenvalue of \emph{Res} $\Leftrightarrow$ $\exp(2\pi \emph{i}\eta)$ is an eigenvalue of the monodromy.
\item $\exp(-2\pi\emph{i}\emph{Res})$ is conjugate to the monodromy.
\item The monodromy is unipotent $\Leftrightarrow$ \emph{Res} is nilpotent.
\end{enumerate}
\end{thm}

We cannot be sure that $\omega,\mathcal{D}\omega,\dots,\mathcal{D}^{u-1}\omega$, with $\omega$ a solution of the Picard-Fuchs equation, is a basis for the solution space, but we can easily write down the connection matrix for this basis:

\begin{align*}
\Gamma=\begin{pmatrix}0&1&0&\cdots&~&0\\0&0&1&0&\cdots&0\\\vdots&~&\ddots&\ddots&~&\vdots \\0&\cdots&~&0&1&0\\0&0&\cdots&~&0&1 \\-h_1(\lambda)&-h_2(\lambda)&-h_3(\lambda)&\cdots&~&-h_{u-1}(\lambda)\end{pmatrix}
\end{align*}

A theorem by Morisson gives a condition for these elements $\omega,\mathcal{D}\omega,\dots,\mathcal{D}^{u-1}\omega$ to be a basis of the solution space. The condition depends on the eigenvalues of the matrix $\Gamma$.

\begin{thm}\emph{(\cite{MR1191426})}\label{Morrison}
Let $\mathcal{D}\vec{\omega}(\lambda)=\Gamma\vec{\omega}(\lambda)$ be a system of ordinary differential equations with a regular singular point at $\lambda=0$. If distinct eigenvalues of $\Gamma_{\lambda=0}$ do not differ by integers, then $\omega_1,\dots,\omega_u$ with $\vec{\omega}=(\omega_1,\dots,\omega_u)$ is a basis for the solution space of the system of ordinary differential equations.
\end{thm}

So we calculate the eigenvalues of $\Gamma_{\lambda=0}$. For this purpose we only have to remember that the equation (\ref{pf-ausmultipliziert}) or equally equation (\ref{pf-in-t}) has the following solutions for $\lambda=0$:

\begin{align*}
\bigsqcup_{i=1}^n \left\{0,\frac{1}{\widehat{q}_i},\dots,\frac{\widehat{q}_i-1}{\widehat{q}_i}\right\}\setminus\left(\left\{0,\frac{1}{\widehat{d}},\dots,\frac{\widehat{d}-1}{\widehat{d}}\right\}\cap\bigcup_{i=1}^n \left\{0,\frac{1}{\widehat{q}_i},\dots,\frac{\widehat{q}_i-1}{\widehat{q}_i}\right\}\right).
\end{align*}

This means that no distinct eigenvalues differ by an integer and therefore $\Gamma_{\lambda=0}=\text{Res}$ is a residue matrix by Theorem \ref{Morrison}. In addition it follows from Theorem \ref{Deligne} that for every eigenvalue $\eta$ of $\Gamma$ we get an eigenvalue $\exp(2\pi\textrm{i}\eta)$ of the monodromy. So together with Corollary \ref{poincare-zeroes}, we get the following statement.

\begin{kor}
The poles of the Poincar\'{e} series of $g^t(\x)$ are the eigenvalues of the monodromy around $\lambda=(-s)^{-\widehat{d}}=0$ in the solution space of the Picard-Fuchs equation of $f(\x)=g(\x)+s\prod_i x_i$ and the zeroes of the Poincar\'{e} series of $g^t(\x)$ are the eigenvalues of the monodromy around $\lambda=(-s)^{-\widehat{d}}=\infty$.
\end{kor}

The second part of this statement is proved analogously to the first part, substituting only $\lambda$ by $\lambda^{-1}$.

\begin{rem}
For the calculations in the last section this means that the eigenvalues of the monodromy around $\lambda=0$ are given by the roots of $\chi_{\infty}$ and the eigenvalues of the monodromy around $\lambda=\infty$ are given by the roots of $\chi_0$.
\end{rem}

\begin{rem}
Notice that the monodromy around $0$ and $\infty$ is not unipotent, but it is quasi-unipotent, i.e. a power of the monodromy is unipotent. This agrees with Theorem 2.3 in \cite{De}.
\end{rem}

We want to mention that the points $0$ and $\infty$ are not the only points with monodromy. At $\lambda=\prod \widehat{q}_i^{\widehat{q}_i}/\widehat{d}^{\widehat{d}}$ the Picard-Fuchs equation degenerates and therefore we can consider monodromy around this point in the solution space as well. But the monodromy around this point is just a combination of the monodromy around the other two points. This can be seen from the fact that the parameter can be considered on a projective line (cf. \cite{MR1852194}).

Also we want to mention that the critical points of $\lambda$ in the solution space of the Picard-Fuchs equation apart from $\lambda=\infty$ are in 1-1 correspondence with the critical values of $f(\x)$ in $s$. Namely $\lambda=(-s)^{-\widehat{d}}=0$ and $\lambda=(-s)^{-\widehat{d}}=\frac{\prod \widehat{q}_i^{\widehat{q}_i}}{\widehat{d}^{\widehat{d}}}$ are the critical values of $f(\x)$ in $s$.

\newpage
\thispagestyle{empty}

\appendix

\chapter{Examples: Simple K3 singularities}\label{k3}
In this first part of the appendix we will state some famous examples. These were additional examples leading to Theorem \ref{conj}.\\
The following examples were all calculated with the Griffiths-Dwork method using Singular. The code for these calculations can be found in the second part of the appendix. The polynomials given below are those from the list of 95 polynomials in \cite{MR1066667} that can be described as invertible polynomials by considering an involution on the corresponding hypersurface. The polynomials that can be achieved by the involution are due to personal communication with Noriko Yui and will be published soon in a joint paper with Yasuhiro Goto and Ron Livn\'{e} \cite{GLY}.\\
In all the examples below Theorem \ref{conj} can be checked. This can easily be done by computing the dual weights and dual degree and comparing them to the numbers appearing in the Picard-Fuchs equation we calculated.\\
The number in the first column of the table is the index given in the article of Yonemura \cite{MR1066667}. The second column contains the invertible polynomial that was stated as $g(x_1,x_2,x_3,x_4)$ before. The third and fourth column contain the weights of the invertible polynomial and the order of the Picard-Fuchs equation respectively. The order can also be calculated with the result of Theorem \ref{thm}, which is a fast computation once the dual weights and the dual degree are known. Finally in the last column the result of the calculations, namely the Picard-Fuchs equation of $f(x_1,x_2,x_3,x_4)=g(x_1,x_2,x_3,x_4)+sx_1x_2x_3x_4$, is given and this is exactly the formula given in Theorem \ref{conj}.

\begin{landscape}
\begin{table}
\centering
\begin{tabular}{llccl}
\toprule
Nr. & invertible polynomial & weights & order PF & Picard-Fuchs equation\\
\midrule
$1$ & $x_1^4+x_2^4+x_3^4+x_4^4$ & $(1,1,1,1)$ & $3$ & $s^4\delta^3-2^8(\delta-1)(\delta-2)(\delta-3)$\\
\midrule
$2$ & $x_1^6+x_2^4+x_3^4+x_4^3$ & $(2,3,3,4)$ & $6$ & $s^{12}\delta^3(\delta+4)(\delta+6)(\delta+8)-2^{14}3^6(\delta-1)(\delta-2)(\delta-5)(\delta-7)(\delta-10)(\delta-11)$\\
\midrule
$3$ & $x_1^6+x_2^6+x_3^3+x_4^3$ & $(1,1,2,2)$ & $4$ & $s^6\delta^3(\delta+3)-2^23^6(\delta-1)(\delta-2)(\delta-4)(\delta-5)$\\
\midrule
$4$ & $x_1^{12}+x_2^4+x_3^3+x_4^3$ & $(1,3,4,4)$ & $6$ & $s^{12}\delta^3(\delta+3)(\delta+6)(\delta+9)-2^{10}3^9(\delta-1)(\delta-2)(\delta-5)(\delta-7)(\delta-10)(\delta-11)$\\
\midrule
$5$ & $x_1^6+x_2^6+x_3^6+x_4^2$ & $(1,1,1,3)$ & $3$ & $s^6\delta^3-2^63^3(\delta-1)(\delta-3)(\delta-5)$\\
\midrule
$6$ & $x_1^{10}+x_2^5+x_3^5+x_4^2$ & $(1,2,2,5)$ & $4$ & $s^{10}\delta^3(\delta+5)-2^65^5(\delta-1)(\delta-3)(\delta-7)(\delta-9)$\\
\midrule
$7$ & $x_1^8+x_2^8+x_3^4+x_4^2$ & $(1,1,2,4)$ & $4$ & $s^8\delta^3(\delta+4)-2^{14}(\delta-1)(\delta-3)(\delta-5)(\delta-7)$\\
\midrule
$8$ & $x_1^{12}+x_2^6+x_3^4+x_4^2$ & $(1,2,3,6)$ & $6$ & $s^{12}\delta^3(\delta+4)(\delta+6)(\delta+8)-2^{16}3^3(\delta-1)(\delta-3)(\delta-5)(\delta-7)(\delta-9)(\delta-11)$\\
\midrule
\multirow{2}{*}{$9$} & \multirow{2}{*}{$x_1^{20}+x_2^5+x_3^4+x_4^2$} & \multirow{2}{*}{$(1,4,5,10)$} & \multirow{2}{*}{$8$} & $s^{20}\delta^3(\delta+4)(\delta+6)(\delta+8)(\delta+10)(\delta+12)$\\
&&&&$-2^{22}5^5(\delta-1)(\delta-3)(\delta-7)(\delta-9)(\delta-11)(\delta-13)(\delta-17)(\delta-19)$\\
\midrule
$10$ & $x_1^{12}+x_2^{12}+x_3^3+x_4^2$ & $(1,1,4,6)$ & $4$ & $s^{12}\delta^3(\delta+6)-2^{10}3^6(\delta-1)(\delta-5)(\delta-7)(\delta-11)$\\
\midrule
\multirow{3}{*}{$11$} & \multirow{3}{*}{$x_1^{15}+x_2^{10}+x_3^3+x_4^2$} & \multirow{3}{*}{$(2,3,10,15)$} & \multirow{3}{*}{$10$} & $s^{30}\delta^3(\delta+6)(\delta+10)(\delta+12)(\delta+15)(\delta+18)(\delta+20)(\delta+24)$\\
&&&&$-2^{18}3^{12}5^5(\delta-1)(\delta-5)(\delta-7)(\delta-11)(\delta-13)(\delta-17)(\delta-19)(\delta-23)$\\
&&&&$(\delta-25)(\delta-29)$\\
\midrule
$12$ & $x_1^{18}+x_2^9+x_3^3+x_4^2$ & $(1,2,6,9)$ & $6$ & $s^{18}\delta^3(\delta+6)(\delta+9)(\delta+12)-2^{10}3^{12}(\delta-1)(\delta-5)(\delta-7)(\delta-11)(\delta-13)(\delta-17)$\\
\midrule
\multirow{2}{*}{$13$} & \multirow{2}{*}{$x_1^{24}+x_2^8+x_3^3+x_4^2$} & \multirow{2}{*}{$(1,3,8,12)$} & \multirow{2}{*}{$8$} & $s^{24}\delta^3(\delta+6)(\delta+8)(\delta+12)(\delta+16)(\delta+18)$\\
&&&&$-2^{24}3^9(\delta-1)(\delta-5)(\delta-7)(\delta-11)(\delta-13)(\delta-17)(\delta-19)(\delta-23)$\\
\midrule
\multirow{3}{*}{$14$} & \multirow{3}{*}{$x_1^{42}+x_2^7+x_3^3+x_4^2$} & \multirow{3}{*}{$(1,6,14,21)$} & \multirow{3}{*}{$12$} & $s^{42}\delta^3(\delta+6)(\delta+12)(\delta+14)(\delta+18)(\delta+21)(\delta+24)(\delta+28)(\delta+30)(\delta+36)$\\
&&&&$-2^{22}3^{15}7^7(\delta-1)(\delta-5)(\delta-11)(\delta-13)(\delta-17)(\delta-19)(\delta-23)(\delta-25)$\\
&&&&$(\delta-29)(\delta-31)(\delta-37)(\delta-41)$\\%
\bottomrule
\end{tabular}
\end{table}
\end{landscape}

\begin{landscape}
\begin{table}
\centering
\begin{tabular}{llccl}
\toprule
Nr. & invertible polynomial & weights & order PF & Picard-Fuchs equation\\
\midrule
\multirow{2}{*}{$15$} & \multirow{2}{*}{$x_1^5+x_2^5+x_2x_3^3+x_4^3$} & \multirow{2}{*}{$(3,3,4,5)$} & \multirow{2}{*}{$8$} & $2s^{15}\delta^3(\delta+3)(\delta+6)(\delta+\frac{15}{2})(\delta+9)(\delta+12)$\\
&&&&$+3^{12}5^5(\delta-1)(\delta-2)(\delta-4)(\delta-7)(\delta-8)(\delta-11)(\delta-13)(\delta-14)$\\
\midrule
$16$ & $x_1^8+x_2^4+x_1x_3^3+x_4^3$ & $(3,6,7,8)$ & $6$ & $s^{12}\delta^3(\delta+3)(\delta+6)(\delta+9)-2^83^9(\delta-1)(\delta-2)(\delta-5)(\delta-7)(\delta-10)(\delta-11)$ \\
\midrule
$19$ & $x_1^8+x_2^4+x_3^4+x_2x_4^2$ & $(1,2,2,3)$ & $4$ & $s^8\delta^3(\delta+4)-2^{14}(\delta-1)(\delta-3)(\delta-5)(\delta-7)$ \\
\midrule
\multirow{2}{*}{$20$} & \multirow{2}{*}{$x_1^{24}+x_2^4+x_3^3+x_2x_4^2$} & \multirow{2}{*}{$(1,6,8,9)$} & \multirow{2}{*}{$8$} & $s^{24}\delta^3(\delta+6)(\delta+8)(\delta+12)(\delta+16)(\delta+18)$\\
&&&&$-2^{24}3^9(\delta-1)(\delta-5)(\delta-7)(\delta-11)(\delta-13)(\delta-17)(\delta-19)(\delta-23)$ \\
\midrule
\multirow{3}{*}{$22$} & \multirow{3}{*}{$x_1^{15}+x_2^5+x_3^3+x_2x_4^2$} & \multirow{3}{*}{$(1,3,5,6)$} & \multirow{3}{*}{$10$} & $s^{30}\delta^3(\delta+6)(\delta+10)(\delta+12)(\delta+15)(\delta+18)(\delta+20)(\delta+24)$\\
&&&&$-2^{18}3^{12}5^5(\delta-1)(\delta-5)(\delta-7)(\delta-11)(\delta-13)(\delta-17)(\delta-19)$\\
&&&&$(\delta-23)(\delta-25)(\delta-29)$ \\
\midrule
$24$ & $x_1^{12}+x_2^6+x_3^3+x_2x_4^2$ & $(1,2,4,5)$ & $4$ & $s^{12}\delta^3(\delta+6)-2^{10}3^6(\delta-1)(\delta-5)(\delta-7)(\delta-11)$ \\
\midrule
\multirow{2}{*}{$26$} & \multirow{2}{*}{$x_1^{10}+x_2^5+x_3^4+x_1x_4^2$} & \multirow{2}{*}{$(2,4,5,9)$} & \multirow{2}{*}{$8$} & $s^{20}\delta^3(\delta+4)(\delta+8)(\delta+10)(\delta+12)(\delta+16)$\\
&&&&$-2^{22}5^5(\delta-1)(\delta-3)(\delta-7)(\delta-9)(\delta-11)(\delta-13)(\delta-15)(\delta-17)(\delta-19)$ \\
\midrule
\multirow{2}{*}{$27$} & \multirow{2}{*}{$x_1^{12}+x_2^8+x_3^3+x_1x_4^2$} & \multirow{2}{*}{$(2,3,8,11)$} & \multirow{2}{*}{$8$} & $s^{24}\delta^3(\delta+6)(\delta+8)(\delta+12)(\delta+16)(\delta+18)$\\
&&&&$-2^{24}3^9(\delta-1)(\delta-5)(\delta-7)(\delta-11)(\delta-13)(\delta-17)(\delta-19)(\delta-23)$ \\
\midrule
\multirow{3}{*}{$28$} & \multirow{3}{*}{$x_1^{21}+x_2^7+x_3^3+x_1x_4^2$} & \multirow{3}{*}{$(1,3,7,10)$} & \multirow{3}{*}{$12$} & $s^{42}\delta^3(\delta+6)(\delta+12)(\delta+14)(\delta+18)(\delta+21)(\delta+24)(\delta+28)(\delta+30)(\delta+36)$\\
&&&&$-2^{22}3^{15}7^7(\delta-1)(\delta-5)(\delta-11)(\delta-13)(\delta-17)(\delta-19)(\delta-23)(\delta-25)(\delta-29)$\\
&&&&$(\delta-31)(\delta-37)(\delta-41)$ \\
\midrule
$29$ & $x_1^6x_3+x_2^6+x_3^5+x_4^2$ & $(4,5,6,15)$ & $3$ & $s^6\delta^3-2^63^3(\delta-1)(\delta-3)(\delta-5)$ \\
\midrule
$30$ & $x_1^8+x_1x_2^5+x_3^5+x_4^2$ & $(5,7,8,20)$ & $4$ & $s^{10}\delta^3(\delta+5)-2^65^5(\delta-1)(\delta-3)(\delta-7)(\delta-9)$ \\
\midrule
$31$ & $x_1^8+x_2^6+x_2x_3^4+x_4^2$ & $(3,4,5,12)$ & $4$ & $s^8\delta^3(\delta+4)-2^{14}(\delta-1)(\delta-3)(\delta-5)(\delta-7)$ \\
\bottomrule
\end{tabular}
\end{table}
\end{landscape}

\begin{landscape}
\begin{table}
\centering
\begin{tabular}{llccl}
\toprule
Nr. & invertible polynomial & weights & order PF & Picard-Fuchs equation\\
\midrule
\multirow{3}{*}{$32$} & \multirow{3}{*}{$x_1^7+x_2^7+x_2x_3^4+x_4^2$} & \multirow{3}{*}{$(2,2,3,7)$} & \multirow{3}{*}{$12$} & $3^3s^{28}\delta^3(\delta+4)(\delta+8)(\delta+\frac{28}{3})(\delta+12)(\delta+14)(\delta+16)(\delta+\frac{56}{3})(\delta+20)(\delta+24)$\\
&&&&$-2^{34}7^7(\delta-1)(\delta-3)(\delta-5)(\delta-9)(\delta-11)(\delta-13)(\delta-15)(\delta-17)(\delta-19)$\\
&&&&$(\delta-23)(\delta-25)(\delta-27)$ \\
\midrule
$33$ & $x_1^9+x_2^6+x_1x_3^4+x_4^2$ & $(2,3,4,9)$ & $6$ & $s^{12}\delta^3(\delta+4)(\delta+6)(\delta+8)-2^{16}3^3(\delta-1)(\delta-3)(\delta-5)(\delta-7)(\delta-9)(\delta-11)$ \\
\midrule
\multirow{2}{*}{$34$} & \multirow{2}{*}{$x_1^{15}+x_2^5+x_1x_3^4+x_4^2$} & \multirow{2}{*}{$(2,6,7,15)$} & \multirow{2}{*}{$8$} & $s^{20}\delta^3(\delta+4)(\delta+8)(\delta+10)(\delta+12)(\delta+16)$\\
&&&&$-2^{22}5^5(\delta-1)(\delta-3)(\delta-7)(\delta-9)(\delta-11)(\delta-13)(\delta-17)(\delta-19)$ \\
\midrule
\multirow{2}{*}{$35$} & $x_1^8x_2+x_2^7+x_3^4+x_4^2$ & \multirow{2}{*}{$(3,4,7,14)$} & $4$ & $s^8\delta^3(\delta+4)-2^{14}(\delta-1)(\delta-3)(\delta-5)(\delta-7)$ \\
& $x_1^7x_3+x_2^7+x_3^4+x_4^2$ && $6$ & $3^3s^{14}\delta^3(\delta+\frac{14}{3})(\delta+7)(\delta+\frac{28}{3})-2^{10}7^7(\delta-1)(\delta-3)(\delta-5)(\delta-9)(\delta-11)(\delta13)$ \\
\midrule
\multirow{2}{*}{$36$} & $x_1^{10}+x_1x_2^6+x_3^4+x_4^2$ & \multirow{2}{*}{$(2,3,5,10)$} & $6$ & $s^{12}\delta^3(\delta+4)(\delta+6)(\delta+8)-2^{16}3^3(\delta-1)(\delta-3)(\delta-5)(\delta-7)(\delta-9)(\delta-11)$ \\
& $x_1^{10}+x_2^5x_3+x_3^4+x_4^2$ && $4$ & $s^{10}\delta^3(\delta+5)-2^65^5(\delta-1)(\delta-3)(\delta-7)(\delta-9)$ \\
\midrule
\multirow{4}{*}{$37$} & \multirow{2}{*}{$x_1^{16}+x_1x_2^5+x_3^4+x_4^2$} & \multirow{4}{*}{$(1,3,4,8)$} & \multirow{2}{*}{$8$} & $s^{20}\delta^3(\delta+4)(\delta+8)(\delta+10)(\delta+12)(\delta+16)$\\
&&&&$-2^{22}5^5(\delta-1)(\delta-3)(\delta-7)(\delta-9)(\delta-11)(\delta-13)(\delta-17)(\delta-19)$ \\
& \multirow{2}{*}{$x_1^{16}+x_3x_2^4+x_3^4+x_4^2$} && \multirow{2}{*}{$8$} & $3^3s^{16}\delta^3(\delta+4)(\delta+\frac{16}{3})(\delta+8)(\delta+\frac{32}{3})(\delta+12)$\\
&&&&$-2^{32}(\delta-1)(\delta-3)(\delta-5)(\delta-7)(\delta-9)(\delta-11)(\delta-13)(\delta-15)$ \\
\midrule
\multirow{3}{*}{$38$} & \multirow{3}{*}{$x_1^{30}+x_2^5+x_2x_3^3+x_4^2$} & \multirow{3}{*}{$(1,6,8,15)$} & \multirow{3}{*}{$10$} & $2^2s^{30}\delta^3(\delta+6)(\delta+\frac{15}{2})(\delta+12)(\delta+15)(\delta+18)(\delta+\frac{45}{2})(\delta+24)$\\
&&&&$-2^{14}3^{15}5^5(\delta-1)(\delta-5)(\delta-7)(\delta-11)(\delta-13)(\delta-17)(\delta-19)$\\
&&&&$(\delta-23)(\delta-25)(\delta-29)$ \\
\midrule
$39$ & $x_1^{18}+x_2^6+x_2x_3^3+x_4^2$ & $(1,3,5,9)$ & $6$ & $s^{18}\delta^3(\delta+6)(\delta+9)(\delta+12)-2^{10}3^{12}(\delta-1)(\delta-5)(\delta-7)(\delta-11)(\delta-13)(\delta-17)$ \\
\midrule
\multirow{4}{*}{$40$} & \multirow{4}{*}{$x_1^{14}+x_2^7+x_2x_3^3+x_4^2$} & \multirow{4}{*}{$(1,2,4,7)$} & \multirow{4}{*}{$14$} & $2^2s^{42}\delta^3(\delta+6)(\delta+\frac{21}{2})(\delta+12)(\delta+14)(\delta+18)(\delta+21)(\delta+24)(\delta+28)$\\
&&&&$(\delta+30)(\delta+\frac{63}{2})(\delta+36)$\\
&&&&$-2^{22}3^{18}7^7(\delta-1)(\delta-5)(\delta-7)(\delta-11)(\delta-13)(\delta-17)(\delta-19)$\\
&&&&$(\delta-23)(\delta-25)(\delta-29)(\delta-31)(\delta-35)(\delta-37)(\delta-41)$ \\
\bottomrule
\end{tabular}
\end{table}
\end{landscape}

\begin{landscape}
\begin{table}
\centering
\begin{tabular}{llccl}
\toprule
Nr. & invertible polynomial &weights & order PF & Picard-Fuchs equation\\
\midrule
$41$ & $x_1^{12}+x_2^8+x_2x_3^3+x_4^2$ & $(2,3,7,12)$ & $4$ & $s^{12}\delta^3(\delta+6)-2^{10}3^6(\delta-1)(\delta-5)(\delta-7)(\delta-11)$ \\
\midrule
\multirow{3}{*}{$42$} & \multirow{3}{*}{$x_1^{10}+x_2^{10}+x_2x_3^3+x_4^2$} & \multirow{3}{*}{$(1,1,3,5)$} & \multirow{3}{*}{$10$} & $s^{30}\delta^3(\delta+6)(\delta+10)(\delta+12)(\delta+15)(\delta+18)(\delta+20)(\delta+24)$\\
&&&&$-2^{18}3^{12}5^5(\delta-1)(\delta-5)(\delta-7)(\delta-11)(\delta-13)(\delta-17)(\delta-19)$\\
&&&&$(\delta-23)(\delta-25)(\delta-29)$ \\
\midrule
$43$ & $x_1^{12}+x_2^9+x_1x_3^3+x_4^2$ & $(3,4,11,18)$ & $6$ & $s^{18}\delta^3(\delta+6)(\delta+9)(\delta+12)-2^{10}3^{12}(\delta-1)(\delta-5)(\delta-7)(\delta-11)(\delta-13)(\delta-17)$ \\
\midrule
\multirow{2}{*}{$44$} & \multirow{2}{*}{$x_1^{16}+x_2^8+x_1x_3^3+x_4^2$} & \multirow{2}{*}{$(1,2,5,8)$} & \multirow{2}{*}{$8$} & $s^{24}\delta^3(\delta+6)(\delta+8)(\delta+12)(\delta+16)(\delta+18)$\\
&&&&$-2^{24}3^9(\delta-1)(\delta-5)(\delta-7)(\delta-11)(\delta-13)(\delta-17)(\delta-19)(\delta-23)$ \\
\midrule
\multirow{3}{*}{$45$} & \multirow{3}{*}{$x_1^{28}+x_2^7+x_1x_3^3+x_4^2$} & \multirow{3}{*}{$(1,4,9,14)$} & \multirow{3}{*}{$12$} & $s^{42}\delta^3(\delta+6)(\delta+12)(\delta+14)(\delta+18)(\delta+21)(\delta+24)(\delta+28)(\delta+30)(\delta+36)$\\
&&&&$-2^{22}3^{15}7^7(\delta-1)(\delta-5)(\delta-11)(\delta-13)(\delta-17)(\delta-19)(\delta-23)(\delta-25)(\delta-29)$\\
&&&&$(\delta-31)(\delta-37)(\delta-41)$ \\
\midrule
$46$ & $x_1^{12}x_2+x_2^{11}+x_3^3+x_4^2$ & $(5,6,22,33)$ & $4$ & $s^{12}\delta^3(\delta+6)-2^{10}3^6(\delta-1)(\delta-5)(\delta-7)(\delta-11)$ \\
\midrule
$47$ & $x_1^{14}+x_2^7x_3+x_3^3+x_4^2$ & $(3,4,14,21)$ & $6$ & $2^2s^{14}\delta^3(\delta+\frac{7}{2})(\delta+7)(\delta+\frac{21}{2})-2^67^7(\delta-1)(\delta-3)(\delta-5)(\delta-9)(\delta-11)(\delta-13)$ \\
\midrule
$48$ & $x_1^{16}+x_1x_2^9+x_3^3+x_4^2$ & $(3,5,16,24)$ & $6$ & $s^{18}\delta^3(\delta+6)(\delta+9)(\delta+12)-2^{10}3^{12}(\delta-1)(\delta-5)(\delta-7)(\delta-11)(\delta-13)(\delta-17)$ \\
\midrule
\multirow{2}{*}{$49$} & \multirow{2}{*}{$x_1^{21}+x_1x_2^8+x_3^3+x_4^2$} & \multirow{2}{*}{$(2,5,14,21)$} & \multirow{2}{*}{$8$} & $s^{24}\delta^3(\delta+6)(\delta+8)(\delta+12)(\delta+16)(\delta+18)$\\
&&&&$-2^{24}3^9(\delta-1)(\delta-5)(\delta-7)(\delta-11)(\delta-13)(\delta-17)(\delta-19)(\delta-23)$ \\
\midrule
\multirow{4}{*}{$50$} & \multirow{4}{*}{$x_1^{30}+x_2^5x_3+x_3^3+x_4^2$} & \multirow{4}{*}{$(1,4,10,15)$} & \multirow{4}{*}{$12$} & $2^{10}s^{30}\delta^3(\delta+\frac{15}{4})(\delta+\frac{15}{2})(\delta+10)(\delta+\frac{45}{4})(\delta+15)(\delta+\frac{75}{4})(\delta+15)$\\
&&&&$(\delta+\frac{75}{2})(\delta+\frac{105}{4})$\\
&&&&$-2^{10}3^95^{15}(\delta-1)(\delta-3)(\delta-7)(\delta-9)(\delta-11)(\delta-13)(\delta-17)(\delta-19)(\delta-21)$\\
&&&&$(\delta-23)(\delta-27)(\delta-29)$ \\
\midrule
\multirow{3}{*}{$51$} & \multirow{3}{*}{$x_1^{36}+x_1x_2^7+x_3^3+x_4^2$} & \multirow{3}{*}{$(1,5,12,18)$} & \multirow{3}{*}{$12$} & $s^{42}\delta^3(\delta+6)(\delta+12)(\delta+14)(\delta+18)(\delta+21)(\delta+24)(\delta+28)(\delta+30)(\delta+36)$\\
&&&&$-2^{22}3^{15}7^7(\delta-1)(\delta-5)(\delta-11)(\delta-13)(\delta-17)(\delta-19)(\delta-23)(\delta-25)(\delta-29)$\\
&&&&$(\delta-31)(\delta-37)(\delta-41)$ \\
\bottomrule
\end{tabular}
\end{table}
\end{landscape}

\begin{landscape}
\begin{table}
\centering
\begin{tabular}{llccl}
\toprule
Nr. & invertible polynomial & weights & order PF & Picard-Fuchs equation\\
\midrule
$52$ & $x_1^4x_2+x_2^3x_4+x_3^4+x_4^3$ & $(7,8,9,12)$ & $3$ & $s^4\delta^3-2^8(\delta-1)(\delta-2)(\delta-3)$ \\
\midrule
$55$ & $x_1^{10}+x_2^4+x_1x_3^3+x_3x_4^2$ & $(2,5,6,7)$ & $6$ & $s^{12}\delta^3(\delta+4)(\delta+6)(\delta+8)-2^{16}3^3(\delta-1)(\delta-3)(\delta-5)(\delta-7)(\delta-9)(\delta-11)$ \\
\midrule
$56$ & $x_1^6+x_2^5+x_2x_3^3+x_3x_4^2$ & $(5,6,8,11)$ & $3$ & $s^6\delta^3-2^63^3(\delta-1)(\delta-3)(\delta-5)$ \\
\midrule
\multirow{2}{*}{$59$} & \multirow{2}{*}{$x_1^{21}+x_1x_2^4+x_3^3+x_2x_4^2$} & \multirow{2}{*}{$(1,5,7,8)$} & \multirow{2}{*}{$8$} & $s^{24}\delta^3(\delta+6)(\delta+8)(\delta+12)(\delta+16)(\delta+18)$\\
&&&&$-2^{24}3^9(\delta-1)(\delta-5)(\delta-7)(\delta-11)(\delta-13)(\delta-17)(\delta-19)(\delta-23)$ \\
\midrule
\multirow{2}{*}{$60$} & \multirow{2}{*}{$x_1^{18}+x_2^3x_3+x_3^3+x_2x_4^2$} & \multirow{2}{*}{$(1,4,6,7)$} & \multirow{2}{*}{$9$} & $s^{18}\delta^3(\delta+\frac{18}{5})(\delta+6)(\delta+\frac{36}{5})(\delta+\frac{54}{5})(\delta+12)(\delta+\frac{72}{5})$\\
&&&&$-2^{18}3^{15}(\delta-1)(\delta-3)(\delta-5)(\delta-7)(\delta-9)(\delta-11)(\delta-13)(\delta-15)(\delta-17)$ \\
\midrule
$61$ & $x_1^7+x_1x_2^4+x_3^4+x_2x_4^2$ & $(4,6,7,11)$ & $4$ & $s^8\delta^3(\delta+4)-2^{14}(\delta-1)(\delta-3)(\delta-5)(\delta-7)$ \\
\midrule
$65$ & $x_1^{11}+x_1x_2^6+x_3^3+x_2x_4^2$ & $(3,5,11,14)$ & $4$ & $s^{12}\delta^3(\delta+6)-2^{10}3^6(\delta-1)(\delta-5)(\delta-7)(\delta-11)$ \\
\midrule
$68$ & $x_1^{10}+x_2^5x_3+x_3^3+x_2x_4^2$ & $(3,4,10,13)$ & $5$ & $3^3s^{10}\delta^3(\delta+\frac{10}{3})(\delta+\frac{20}{3})-2^{10}5^5(\delta-1)(\delta-3)(\delta-5)(\delta-7)(\delta-9)$ \\
\midrule
$73$ & $x_1^6x_2+x_2^5x_3+x_3^5+x_4^2$ & $(7,8,10,25)$ & $3$ & $s^6\delta^3-2^63^3(\delta-1)(\delta-3)(\delta-5)$ \\
\midrule
$74$ & $x_1^8+x_2^5x_3+x_1x_3^4+x_4^2$ & $(4,5,7,16)$ & $4$ & $s^{10}\delta^3(\delta+5)-2^65^5(\delta-1)(\delta-3)(\delta-7)(\delta-9)$ \\
\midrule
\multirow{2}{*}{$76$} & \multirow{2}{*}{$x_1^{13}+x_2^4x_3+x_1x_3^4+x_4^2$} & \multirow{2}{*}{$(2,5,6,13)$} & \multirow{2}{*}{$8$} & $3^3s^{16}\delta^3(\delta+4)(\delta+\frac{16}{3})(\delta+8)(\delta+\frac{32}{3})(\delta+16)$\\
&&&&$-2^{32}(\delta-1)(\delta-3)(\delta-5)(\delta-7)(\delta-9)(\delta-11)(\delta-13)(\delta-15)$ \\
\midrule
\multirow{3}{*}{$77$} & \multirow{3}{*}{$x_1^{26}+x_1x_2^5+x_2x_3^3+x_4^2$} & \multirow{3}{*}{$(1,5,7,13)$} & \multirow{3}{*}{$10$} & $2^2s^{30}\delta^3(\delta+6)(\delta+\frac{15}{2})(\delta+12)(\delta+15)(\delta+18)(\delta+\frac{45}{2})(\delta+24)$\\
&&&&$-2^{14}3^{15}5^5(\delta-1)(\delta-5)(\delta-7)(\delta-11)(\delta-13)(\delta-17)(\delta-19)$\\
&&&&$(\delta-23)(\delta-25)(\delta-29)$ \\
\midrule
\multirow{3}{*}{$78$} & \multirow{3}{*}{$x_1^{22}+x_2^4x_3+x_2x_3^3+x_4^2$} & \multirow{3}{*}{$(1,4,6,11)$} & \multirow{3}{*}{$10$} & $2^23^6s^{22}\delta^3(\delta+\frac{11}{3})(\delta+\frac{11}{2})(\delta+\frac{22}{3})(\delta+11)(\delta+\frac{44}{3})(\delta+\frac{33}{2})(\delta+\frac{55}{3})$\\
&&&&$-2^{10}11^{11}(\delta-1)(\delta-3)(\delta-5)(\delta-7)(\delta-9)(\delta-13)(\delta-15)(\delta-17)$\\
&&&&$(\delta-19)(\delta-21)$ \\
\bottomrule
\end{tabular}
\end{table}
\end{landscape}

\begin{landscape}
\begin{table}
\centering
\begin{tabular}{llccl}
\toprule
Nr. & invertible polynomial & weights & order PF & Picard-Fuchs equation\\
\midrule
\multirow{2}{*}{$79$} & \multirow{2}{*}{$x_1^{16}+x_1x_2^6+x_2x_3^3+x_4^2$} & \multirow{2}{*}{$(2,5,9,16)$} & \multirow{2}{*}{$6$} & $s^{18}\delta^3(\delta+6)(\delta+9)(\delta+12)$\\
&&&&$-2^{10}3^{12}(\delta-1)(\delta-5)(\delta-7)(\delta-11)(\delta-13)(\delta-17)$ \\
\midrule
$80$ & $x_1^{11}+x_1x_2^8+x_2x_3^3+x_4^2$ & $(4,5,13,22)$ & $4$ & $s^{12}\delta^3(\delta+6)-2^{10}3^6(\delta-1)(\delta-5)(\delta-7)(\delta-11)$ \\
\midrule
\multirow{2}{*}{$81$} & \multirow{2}{*}{$x_1^{13}+x_2^6x_3+x_1x_3^3+x_4^2$} & \multirow{2}{*}{$(2,3,8,13)$} & \multirow{2}{*}{$9$} & $5^5s^{18}\delta^3(\delta+\frac{18}{5})(\delta+6)(\delta+\frac{36}{5})(\delta+\frac{54}{5})(\delta+12)(\delta+\frac{72}{5})$\\
&&&&$-2^{18}3^{15}(\delta-1)(\delta-3)(\delta-5)(\delta-7)(\delta-9)(\delta-11)(\delta-13)(\delta-15)(\delta-17)$ \\
\midrule
\multirow{4}{*}{$82$} & \multirow{4}{*}{$x_1^{22}+x_2^5x_3+x_1x_3^3+x_4^2$} & \multirow{4}{*}{$(1,3,7,11)$} & \multirow{4}{*}{$12$} & $2^{10}s^{30}\delta^3(\delta+\frac{15}{4})(\delta+\frac{15}{2})(\delta+10)(\delta+\frac{45}{4})(\delta+15)(\delta+\frac{75}{4})(\delta+20)$\\
&&&&$(\delta+\frac{45}{2})(\delta+\frac{105}{4})$\\
&&&&$-2^{10}3^95^{15}(\delta-1)(\delta-3)(\delta-7)(\delta-9)(\delta-11)(\delta-13)(\delta-17)(\delta-19)$\\
&&&&$(\delta-21)(\delta-23)(\delta-27)(\delta-29)$ \\
\midrule
$83$ & $x_1^9x_3+x_1x_2^{10}+x_3^3+x_4^2$ & $(4,5,18,27)$ & $5$ & $3^3s^{10}\delta^3(\delta+\frac{10}{3})(\delta+\frac{20}{3})-2^{10}5^5(\delta-1)(\delta-3)(\delta-5)(\delta-7)(\delta-9)$ \\
\midrule
$84$ & $x_1^4x_3+x_2^3x_4+x_2x_3^3+x_4^3$ & $(5,6,7,9)$ & $3$ & $s^4\delta^3-2^8(\delta-1)(\delta-2)(\delta-3)$ \\
\midrule
$85$ & $x_1^7+x_2^3x_4+x_1x_3^3+x_3x_4^2$ & $(2,3,4,5)$ & $6$ & $2^2s^9\delta^3(\delta+3)(\delta+\frac{9}{2})(\delta+6)+3^{12}(\delta-1)(\delta-2)(\delta-4)(\delta-5)(\delta-7)(\delta-8)$ \\
\midrule
\multirow{4}{*}{$87$} & \multirow{4}{*}{$x_1^{13}+x_2^3x_3+x_3^2x_4+x_2x_4^2$} & \multirow{4}{*}{$(1,3,4,5)$} & \multirow{4}{*}{$12$} & $2^83^35^5s^{13}\delta^3(\delta+\frac{13}{5})(\delta+\frac{13}{4})(\delta+\frac{13}{3})(\delta+\frac{26}{5})(\delta+\frac{13}{2})(\delta+\frac{39}{5})(\delta+\frac{26}{3})$\\
&&&&$(\delta+\frac{39}{4})(\delta+\frac{52}{5})$\\
&&&&$+13^{13}(\delta-1)(\delta-2)(\delta-3)(\delta-4)(\delta-5)(\delta-6)(\delta-7)(\delta-8)(\delta-9)$\\
&&&&$(\delta-10)(\delta-11)(\delta-12)$ \\
\midrule
$92$ & $x_1^9x_3+x_1x_2^7+x_2x_3^3+x_4^2$ & $(3,5,11,19)$ & $5$ & $3^3s^{10}\delta^3(\delta+\frac{10}{3})(\delta+\frac{20}{3})-2^{10}5^5(\delta-1)(\delta-3)(\delta-5)(\delta-7)(\delta-9)$ \\
\midrule
$94$ & $x_1^4x_4+x_1x_2^4+x_2x_3^3+x_3x_4^2$ & $(3,4,5,7)$ & $4$ & $2^2s^5\delta^3(\delta+\frac{5}{2})+5^5(\delta-1)(\delta-2)(\delta-3)(\delta-4)$ \\
\bottomrule
\end{tabular}
\end{table}
\end{landscape}

\newpage
\thispagestyle{empty}

\chapter{The Griffiths-Dwork method in Singular}\label{singular}

Below we show the algorithm for the Griffiths-Dwork method in Singular, which can be used for calculating the Picard-Fuchs-equation of a special one-parameter family associated to an arbitrary polynomial in $\C\left[w,x,y,z\right]$. Of course it can easily be adjusted to a different number of variables, but because most of the examples in this thesis are K3-surfaces this is not necessary here. This method of calculation was used for all computations in this thesis unless the calculations are given explicitly. First we fix the notation:\\
Let $g(x_1,x_2,x_3,x_4)$ be any polynomial defining a hypersurface in $\mathbb{P}(q_1,q_2,q_3,q_4)$. Then in this case the algorithm calculates the Picard-Fuchs equation of the one-parameter family $f(x_1,x_2,x_3,x_4)=g(x_1,x_2,x_3,x_4)+sx_1x_2x_3x_4$ using the Griffiths-Dwork method. But one can do it the same way for any other one-parameter family.\\
The output is the polynomial $pf(x_1)$, which is the Picard-Fuchs equation if one replaces the variable $x_1$ by the differential operator $\delta=s\frac{\partial}{\partial s}$. The splitting of the summands of the Picard-Fuchs equation into linear factors has to be done by hand afterwards.

\begin{verbatim}
> ring r=(0,s),(x1,x2,x3,x4),wp(q1,q2,q3,q4); // the ring r is
     	      // C(s)[x1,x2,x3,x4] with a weighted order
> LIB "general.lib";
> intvec q=(q1,q2,q3,q4); // the weight vector is defined
> poly f=g(x1,x2,x3,x4)+s*x1*x2*x3*x4; //defining the one-parameter family f
> ideal j=jacob(f); // defining the Jacobian ideal of the family
> ideal sj=std(j); // calculates a Gröbner basis of the Jacobian ideal
> int d=q[1]+q[2]+q[3]+q[4]; // this number is the degree of f
> ideal kb0=weightKB(sj,0,q); // here the basis of the Milnor ring
> ideal kb1=weightKB(sj,d,q); // in degree 0,d and 2d is calculated
> ideal kb2=weightKB(sj,2*d,q);
> list kl=kb0,kb1,kb2;
> int kn=ncols(kb0)+ncols(kb1)+ncols(kb2);
> matrix m[kn+1][kn]; // the matrix m stores the important 
> m[1,1]=s;						// information derived below
> for(int k=1;k<=kn;k++)
. {
. poly p=factorial(k)*(-1)^k*s^(k+1)*(x1*x2*x3*x4)^k; // These are the
. while(deg(p)>0)      // polynomials that need to be written in the basis
. {								     // of the Milnor ring
. while(reduce(p,sj)==0) // If p is in the Jacobian,
. {								// use the Griffiths formula to reduce the degree
. poly h=0;
. ideal l=lift(j,p);
. for(int jj=1;jj<=4;jj++)
. {
. h=h+diff(l[jj],var(jj)); // this is the Griffiths formula
. }
. p=h*1/(deg(h)/d+1);
. if(deg(p)==0){break;}
. }
. if(deg(p)==0){break;}
. int u=deg(p)/d+1; // If p is not in the Jacobian, calculate the degree to
. ideal kb=kl[u]; // use the correct degree of the Milnor ring and
. ideal li=lift(kb,reduce(p,sj)); //reduce with respect to the basis
. if(u==1) // In the following part we store in the matrix m the coefficients
. {        // of the basis elements we have used so far
. m[k+1,1]=m[k+1,1]+li[1];
. }
. if(u==3)
. {
. m[k+1,kn]=m[k+1,kn]+li[1];
. }
. if(u==2)
. {
. for(int jl=1;jl<=ncols(kb);jl++)
. {
. m[k+1,jl+1]=m[k+1,jl+1]+li[jl];
. }
. }		
. p=p-reduce(p,sj);
. }
. m[k+1,1]=m[k+1,1]+p;
. }
> matrix w[kn+1][kn+1]; // The matrix w consists of the coefficients of the
> for(int kk=1;kk<=kn+1;kk++) // partial derivatives of the holomorphic form 
. {	                        	// omega we began with 
. w[1,kk]=1;
. w[kk,kk]=1;
. if(kk>=3)
. {
. w[2,kk]=w[2,kk-1]+2^(kk-2);
. }
. for(int ll=3;ll<=kk-1;ll++)
. {
. w[ll,kk]=ll*w[ll,kk-1]+w[ll-1,kk-1];
. }
. }
> matrix en=transpose(w)*m; //the matrix en now contains the coefficients of
         // the partial derivatives of omega in the basis of the Milnor ring
. module end=transpose(en);
. module ende=syz(end); // The module ende gives all linear relations between
			              	// the partial derivatives
> poly pf=0;
> for(int lk=1;lk<=nrows(ende);lk++) // The coefficients of the Picard-Fuchs
. {                    // equation are put in a polynomial with variable x1
. pf=pf+ende[lk,1]*x1^(lk-1);
. }
> pf;
\end{verbatim}
\newpage
\thispagestyle{empty}

%--------------------------   Ende des Hauptteils, Literaturverzeichnis und Lebenslauf  -------------------------------

\bibliographystyle{alpha}
\clearpage
\addcontentsline{toc}{chapter}{Bibliography}
\bibliography{Bibi}

\end{document}